\PassOptionsToPackage{dvipsnames}{xcolor}
\documentclass[a4paper,11pt]{article}

\usepackage[
    a4paper,
    top=1.05in,
    bottom=1.05in,
    left=1.00in,
    right=1.00in
]{geometry}
\usepackage[T1]{fontenc}

\usepackage{amsmath,amsthm,mathtools}
\usepackage{amssymb}
\usepackage{newtxtext,newtxmath}
\usepackage{microtype}
\microtypesetup{protrusion=true,expansion=true}

\usepackage{mathrsfs}
\DeclareFontFamily{U}{rsfs}{\skewchar\font127}
\DeclareFontShape{U}{rsfs}{m}{n}{<5-7>rsfs5<7-10>rsfs7<10->rsfs10}{}
\usepackage{slashed}
\usepackage{upgreek}
\allowdisplaybreaks[4]
\usepackage{enumitem}
\setlist{itemsep=0.25em, topsep=0.45em, parsep=0pt}
\setlist[enumerate]{label=(\roman*)}
\usepackage{scalerel}
\usepackage[usestackEOL]{stackengine}
\usepackage[thinc]{esdiff}

\usepackage[all]{xy}
\usepackage{graphicx,psfrag}
\usepackage{pstool}
\usepackage{tikz-cd}
\usetikzlibrary{arrows.meta,positioning}
\usepackage{float}
\usepackage[font=small,labelfont=bf,labelsep=period]{caption}
\usepackage{titlesec}
\titleformat{\section}{\normalfont\Large\bfseries}{\thesection}{0.75em}{}
\titleformat{\subsection}{\normalfont\large\bfseries}{\thesubsection}{0.75em}{}
\titlespacing*{\section}{0pt}{2.6ex plus 0.8ex minus 0.2ex}{1.2ex plus 0.2ex}
\titlespacing*{\subsection}{0pt}{2.0ex plus 0.6ex minus 0.2ex}{0.9ex plus 0.2ex}

\usepackage{fancyhdr}
\fancypagestyle{plain}{%
  \fancyhf{}%
  \fancyfoot[C]{\thepage}%
}

\usepackage[dvipsnames]{xcolor}
\definecolor{linkblue}{RGB}{20,67,128}
\definecolor{citegreen}{RGB}{25,105,80}
\definecolor{urlpurple}{RGB}{100,55,125}
\usepackage[colorlinks=true,
            linkcolor=linkblue,
            citecolor=citegreen,
            urlcolor=urlpurple,
            plainpages=false]{hyperref}
\usepackage{bookmark}

 \usepackage{setspace}

\usepackage[titletoc,title]{appendix}
\usepackage{titletoc}
\usepackage{makeidx}

\usepackage[norefs,nocites,ignoreunlbld]{refcheck}
\usepackage[capitalize,nameinlink,noabbrev]{cleveref}

\newcommand{\R}{\ensuremath{\mathbb{R}}}

\def\xx{\mathbf{x}}
\def\yy{\mathbf{y}}

\newcommand{\ba}{\begin{align*}}
	\newcommand{\ea}{\end{align*}}

\def\Vol{\mathrm{Vol}}

\def\x{\mathbf{x}}

\newcommand{\tr}{\text{tr}}

\newcommand{\Rm}{\ensuremath{\mathrm{Rm}}}
\newcommand{\Ric}{\ensuremath{\mathrm{Ric}}}

\makeatletter
\def\ExtendSymbol#1#2#3#4#5{\ext@arrow 0099{\arrowfill@#1#2#3}{#4}{#5}}

\makeatother

\makeatletter
\def\ExtendSymbol#1#2#3#4#5{\ext@arrow 0099{\arrowfill@#1#2#3}{#4}{#5}}

\makeatother

\DeclarePairedDelimiter\abs{\lvert}{\rvert}%
\makeatletter
\let\oldabs\abs
\def\abs{\@ifstar{\oldabs}{\oldabs*}}
\makeatother

\def\aint{\,\ThisStyle{\ensurestackMath{%
			\stackinset{c}{.2\LMpt}{c}{.5\LMpt}{\SavedStyle-}{\SavedStyle\phantom{\int}}}%
		\setbox0=\hbox{$\SavedStyle\int\,$}\kern-\wd0}\int}

\numberwithin{equation}{section}

\theoremstyle{plain}
\newtheorem{thm}{Theorem}[section]
\newtheorem{cor}[thm]{Corollary}
\newtheorem{prop}[thm]{Proposition}
\newtheorem{lem}[thm]{Lemma}
\newtheorem{conj}[thm]{Conjecture}

\theoremstyle{definition}
\newtheorem{defn}[thm]{Definition}

\newtheorem{claim}[thm]{Claim}

\theoremstyle{remark}
\newtheorem{rem}[thm]{Remark}

\titlecontents{part}[0pt]  
{\addvspace{10pt}\bfseries} 
{\makebox[3em][l]{PART-\thecontentslabel}}      
{}                          
{}                          

\makeindex

\newcommand{\Hess}{\nabla^2}
\DeclareMathOperator{\Area}{Area}

\newcommand{\Lie}{\mathcal L}

\newcommand{\sphere}{\mathbb S^2}
\newcommand{\Disk}{\mathbb D}
\newcommand{\Pheat}{\mathcal P}
\newcommand{\krot}{\varkappa}
\hypersetup{pdftitle={O(3)-symmetry of four-dimensional steady Ricci solitons},
            pdfauthor={Hanbing Fang}}

\title{$O(3)$-symmetry of four-dimensional steady Ricci solitons}
\author{Hanbing Fang}
\date{\today}

\begin{document}
\maketitle
\begin{abstract}
We prove that a complete connected four-dimensional gradient steady
Ricci soliton with positive curvature operator that is $\kappa$-noncollapsed
at all scales admits an effective isometric $O(3)$ action preserving the
potential function and having principal orbits diffeomorphic to $S^2$. 
\end{abstract}
\tableofcontents

\section{Introduction}

\subsection{Background and main theorem}
Steady Ricci solitons are self-similar solutions to the
Ricci flow. Their evolution is generated by diffeomorphisms, so they
exist for all time and provide models for singularity formation.
Understanding their geometry is therefore a natural problem
in the study of singularity formation of Ricci flows. Specifically, a complete Riemannian manifold $(M,g)$ is called a \textbf{gradient steady Ricci soliton} if there exists a smooth potential function $f$ satisfying
\begin{equation}
                    \Ric_g=\nabla_g^2f.
\end{equation}
Here $\Ric_g$ is the Ricci tensor and $\nabla_g^2f$ is the Hessian
of the potential $f$; we write $\mathrm{R}$ for scalar curvature and $\Rm$
for the full curvature tensor. Thus the canonical Ricci flow is $g(t):=\Phi_t^*g$, where
$\frac{\mathrm{d}}{\mathrm{d}t}\Phi_t=-\nabla f\circ\Phi_t$.
Our question is whether positive curvature and noncollapsing force
a four-dimensional steady soliton to carry an isometric $O(3)$ action
with two-dimensional principal orbits. 

The lower-dimensional theory explains both the expected symmetry and
the role of noncollapsing. Hamilton's cigar \cite{HamiltonSurfaces}
is the basic two-dimensional steady example. In Euclidean coordinates,
one normalization is
\[
 g_{\mathrm{cig}}=\frac{\mathrm{d}x^2+\mathrm{d}y^2}{1+x^2+y^2},
 \qquad f_{\mathrm{cig}}=\log(1+x^2+y^2),
 \qquad \mathrm{R}_{\mathrm{cig}}=\frac4{1+x^2+y^2}.
\]
Its cylindrical end is collapsed at large curvature-controlled scales.
Bryant \cite{Bryant} constructed the rotationally symmetric
three-dimensional steady soliton. Brendle \cite{Brendle3} proved
that every complete nonflat three-dimensional noncollapsed gradient steady soliton is Bryant,
up to scaling and isometry. This proved the steady case of
Perelman's conjecture that every complete noncompact
three-dimensional ancient $\kappa$-solution with bounded positive
sectional curvature is the Bryant soliton
\cite[Remark~11.9]{Perelman}. Perelman also announced uniqueness
within the gradient steady class in that remark, without providing
a proof there. Brendle \cite[Theorem~1.3]{BrendleAncient}
subsequently proved the full noncompact classification: every such
solution is, up to parabolic scaling, time translation, and isometry,
the Bryant soliton or a
family of shrinking round cylinders $S^2\times\mathbb R$, possibly
with a quotient. Here an ancient $\kappa$-solution is complete,
nonflat, and $\kappa$-noncollapsed at all scales, with nonnegative
curvature operator and bounded curvature on compact time intervals.
In the compact case, Brendle--Daskalopoulos--\v Se\v sum
\cite[Theorem~1.5]{BDS} proved that an ancient $\kappa$-solution
on $S^3$ is, up to diffeomorphism, parabolic scaling, and time
translation, either a family of shrinking round spheres or
Perelman's ancient oval.

The noncollapsing assumption is essential for this classification.
Without it, a three-dimensional steady soliton need not be Bryant. Lai constructed the collapsed flying wings
\cite{LaiWings} and proved $O(2)$ symmetry for complete nonflat
three-dimensional gradient steady solitons \cite{LaiSymmetry}.
The latter work provides an important precedent for recovering
symmetry from different geometric regions at infinity, rather than
from a single rotational model.

In higher dimensions, Brendle
\cite{BrendleHigher} proved rotational symmetry under positive
sectional curvature and asymptotic cylindricity; Deng--Zhu \cite{DZ}
obtained rotational symmetry under nonnegative curvature operator,
noncollapsing, and linear scalar curvature decay. In the nonflat
case these hypotheses force Bryant geometry.
The unrestricted problem has more possibilities: Lai \cite{LaiWings}
constructed noncollapsed positively curved examples with
$\mathbb Z_2\times O(n-1)$ symmetry for $n\ge4$, and
Chan--Lai--Lee \cite{CLL} constructed further higher-dimensional
families, including noncollapsed subfamilies. Thus $O(3)$ symmetry
is a natural target in dimension four even when full $O(4)$ symmetry
fails. Here ``$O(3)$ symmetry'' always means an effective action
with two-sphere principal orbits.

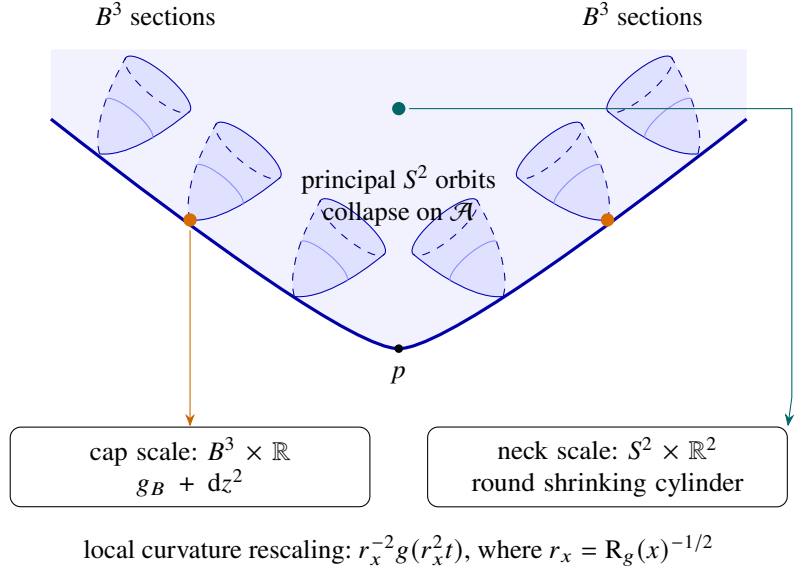
\begin{figure}[htbp]
\centering
\begin{tikzpicture}[font=\small,>=Stealth,x=.92cm,y=.92cm]
 \fill[blue!5] (-5,3.4)
  .. controls (-2.7,1.6) and (-.6,.1) .. (0,.1)
  .. controls (.6,.1) and (2.7,1.6) .. (5,3.4)
  --(5,4.4)--(-5,4.4)--cycle;
 \draw[very thick,blue!65!black] (-5,3.4)
  .. controls (-2.7,1.6) and (-.6,.1) .. (0,.1)
  .. controls (.6,.1) and (2.7,1.6) .. (5,3.4);
 \foreach \side in {-1,1} {
  \foreach \xx/\yy in {1.5/.85,3.0/1.95,4.3/2.9} {
   \begin{scope}[shift={(\side*\xx,\yy)},xscale=-\side,rotate=50]
    \fill[blue!12] (0,0)
     .. controls (.2,-.3) and (.75,-.5) .. (1.3,-.6)
     arc[start angle=-90,end angle=90,x radius=.18,y radius=.6]
     .. controls (.75,.5) and (.2,.3) .. (0,0);
    \draw[blue!70!black] (0,0)
     .. controls (.2,-.3) and (.75,-.5) .. (1.3,-.6)
     arc[start angle=-90,end angle=90,x radius=.18,y radius=.6];
    \draw[dashed,blue!60!black] (0,0)
     .. controls (.2,.3) and (.75,.5) .. (1.3,.6)
     arc[start angle=90,end angle=270,x radius=.18,y radius=.6];
    \draw[blue!45] (.65,-.44) arc[start angle=-90,end angle=90,x radius=.13,y radius=.44];
   \end{scope}
  }
 }
 \fill (0,.1) circle (1.6pt) node[below=3pt] {$p$};
 \fill[orange!85!black] (-3,1.95) circle (2.5pt);
 \fill[orange!85!black] (3,1.95) circle (2.5pt);
 \fill[teal!80!black] (0,3.55) circle (2.5pt);
 \node[align=center] at (0,2.3) {principal $S^2$ orbits\\collapse on $\mathcal A$};
 \node[above] at (-3.5,4.6) {$B^3$ sections};
 \node[above] at (3.5,4.6) {$B^3$ sections};
 \node[draw,rounded corners,align=center,text width=4.5cm] (cap) at (-3,-1.6)
  {cap scale: $B^3\times\mathbb R$\\$g_B+\mathrm dz^2$};
 \node[draw,rounded corners,align=center,text width=4.5cm] (neck) at (3,-1.6)
  {neck scale: $S^2\times\mathbb R^2$\\round shrinking cylinder};
 \draw[->,orange!80!black] (-3,1.8)--(cap.north);
 \draw[->,teal!80!black] (.15,3.55)--(5.65,3.55)--(5.65,-.6)--(neck.north east);
 \node[align=center] at (0,-2.8)
  {local curvature rescaling: $r_x^{-2}g(r_x^2t)$, where $r_x=\mathrm R_g(x)^{-1/2}$};
\end{tikzpicture}
\caption{The two cap regions and the neck geometry of a four-dimensional
non-Bryant steady soliton. The wing layout follows the three-dimensional
sketch in \cite[Figure~1]{LaiWings}, with transverse three-dimensional Bryant soliton
sections replacing cigars and two-sphere orbits replacing circles.
The diagram depicts the geometry after the angular action has been
constructed. The blue boundary represents its fixed axis $\mathcal A$, opened out to
show its two escaping ends. Orange points indicate cap base points;
the green point indicates a neck base point away from the cap cores.
The boxes describe local curvature scale flow models.}
\label{fig:intro-four-dimensional-wing}
\end{figure}

Dimension reduction gives the local geometric setting.
Chow--Deng--Ma \cite{CDM} studied four-dimensional steady solitons
under dimension reduction hypotheses. Chan--Ma--Zhang \cite{CMZ}
proved that the tangent flow at infinity is
$S^3\times\mathbb R$ or $S^2\times\mathbb R^2$, with the first
case forcing the four-dimensional Bryant soliton. In the non-Bryant case, the curvature scale
models are shrinking $S^2\times\mathbb R^2$ cylinders and
$B^3\times\mathbb R$, where $B^3$ is the three-dimensional Bryant
soliton. Hebbar--\v Se\v sum \cite{HS} described the two edges and
large potential levels under hypotheses including a critical point
of the potential. We use the local models of \cite{CMZ} for the collar estimates in
Section~\ref{sec:collars} and construct the intrinsic exterior geometry
in Section~\ref{sec:exterior} without assuming a critical point.
The global edge and potential level results of \cite{HS} require a
critical point; we discuss their relation to our metric blow-down
criterion in Remark~\ref{rem:ray-criterion-generality}, after
constructing the angular action. The asymptotic geometry is illustrated in Figure \ref{fig:intro-four-dimensional-wing}.

\paragraph{Main result.}
We say that $(M^4,g)$ is \textbf{$\kappa$-noncollapsed} at all scales if
\begin{equation}
 \sup_{B_g(x,r)}|\Rm|\le r^{-2}
 \quad\Longrightarrow\quad
 \Vol_g B_g(x,r)\ge\kappa r^4
 \qquad(x\in M,\ r>0).
\end{equation}
Under the curvature assumptions below, this condition is equivalent
to parabolic noncollapsing of the canonical flow at all scales,
with a possible dimensional change in the constant $\kappa$.
We now state the main theorem.

\begin{thm}\label{thm:main}
Let $(M^4,g,f)$ be a complete connected gradient steady Ricci soliton
with positive curvature operator. Suppose that it is
$\kappa$-noncollapsed at all scales. Then $(M,g)$ admits an effective
isometric $O(3)$ action preserving $f$. The action has
principal orbits diffeomorphic to $S^2$. The potential $f$ has a unique critical point, which is its global minimum, and is proper.
\end{thm}

Theorem~\ref{thm:main} recovers a common angular rotation group. Neither
uniform scalar decay nor a critical point of $f$ is included among
the hypotheses. We establish a unique critical point and properness
of the potential as part of the proof, after constructing the angular
action.

The steady soliton part of Haslhofer's classification conjecture
\cite[Conjecture~1.3]{Haslhofer4D} predicts the four-dimensional
Bryant soliton, the product $B^3\times\mathbb R$, and Lai's
one-parameter family as the noncollapsed steady models, up to scaling
and finite quotients. Lai's Oberwolfach report
\cite[Conjecture~2, p.~1817]{LaiOberwolfach} likewise lists the
four-dimensional Bryant soliton and Lai's constructed family, and
additionally conjectures a ray blow-down for those examples. The
report's list does not explicitly include the product
$B^3\times\mathbb R$. For the class with positive curvature
operator considered here, the classification prediction reads as
follows.
\begin{conj}
Every complete four-dimensional gradient steady Ricci soliton with
positive curvature operator that is $\kappa$-noncollapsed at all scales
is, up to scaling and isometry, the four-dimensional Bryant soliton
or a member of Lai's $\mathbb Z_2\times O(3)$-symmetric family.
\end{conj}
The product $B^3\times\mathbb R$ belongs to the nonnegative-curvature
formulation; its mixed sectional curvatures vanish, so it is excluded
by the positive curvature operator assumption above.
Theorem~\ref{thm:main} addresses the symmetry aspect of that problem;
identification with a particular parameterized family needs additional
uniqueness information. In a forthcoming work \cite{Fang26IV}, based on the angular action, we will study the quantitative geometry of the compact potential levels, prove the $\mathbb Z_2\times O(3)$-symmetry and discuss the uniqueness problem. The nonidentity element of the $\mathbb Z_2$ factor reverses the fixed axis of the $O(3)$ factor.

The angular action also allows us to determine the large-scale metric
limit. The following result is proved in
Theorem~\ref{lem:metric-blowdown}. Here, for any $\alpha\in (0,\pi)$, $\mathcal C_\alpha$ is the closed Euclidean sector
of opening angle $\alpha$, with its intrinsic metric and vertex $o$; $\mathcal C_0=[0,\infty)$ is the ray.

\begin{thm}\label{thm:intro-metric-blowdown}
Under the hypotheses of Theorem~\ref{thm:main}, there is an angle
$\alpha\in[0,\pi)$ such that, for every fixed $p\in M$,
\[
 (M,\lambda g,p)\longrightarrow(\mathcal C_\alpha,o)
 \qquad\text{as }\lambda\downarrow0
\]
in the pointed Gromov--Hausdorff sense along the full family of scales.
Moreover, we have
\[
 \alpha=0\quad\Longleftrightarrow\quad
 \mathrm{R}(x)\longrightarrow0\quad\text{uniformly as }x\to\infty.
\]
In the non-Bryant case, let $S:=M^{O(2)}$ be the meridian surface for
the subgroup $O(2)\subset O(3)$ that fixes a chosen line pointwise in the
three-dimensional rotation space and acts orthogonally on its
perpendicular plane. The surface $S$ consists of the points fixed
by every element of this subgroup; it is a complete totally geodesic
two-dimensional surface, rather than the three-dimensional quotient
by $O(2)$. Its Gaussian
curvature $K_S$ determines the angle by
\[
                 \alpha=\pi-\frac12\int_S K_S\,\mathrm{d}\mu_S.
\]
\end{thm}

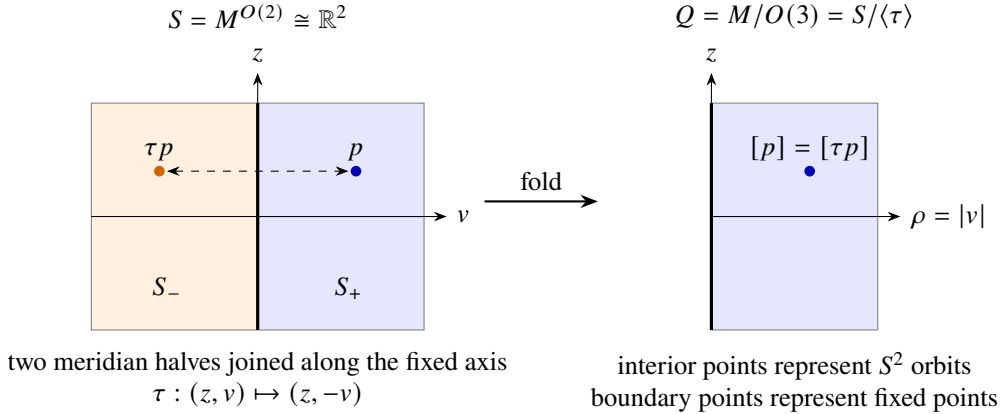
\begin{figure}[!htbp]
\centering
\begin{tikzpicture}[x=1cm,y=1cm,>=Stealth,font=\small]
 \fill[orange!12] (-2.2,-1.5) rectangle (0,1.5);
 \fill[blue!10] (0,-1.5) rectangle (2.2,1.5);
 \draw[gray] (-2.2,-1.5) rectangle (2.2,1.5);
 \draw[very thick] (0,-1.5)--(0,1.5);
 \draw[->] (0,1.5)--(0,1.9) node[above] {$z$};
 \draw[->] (-2.2,0)--(2.5,0) node[right] {$v$};
 \node at (-1.2,-.95) {$S_-$};
 \node at (1.2,-.95) {$S_+$};
 \fill[blue!70!black] (1.3,.6) circle (2pt);
 \fill[orange!80!black] (-1.3,.6) circle (2pt);
 \node[above] at (1.3,.6) {$p$};
 \node[above] at (-1.3,.6) {$\tau p$};
 \draw[<->,dashed] (-1.2,.6)--(1.2,.6);
 \node[above,align=center] at (0,2.35)
  {$S=M^{O(2)}\cong\mathbb R^2$};
 \node[below,align=center] at (0,-1.65)
  {two meridian halves joined along the fixed axis\\
   $\tau:(z,v)\mapsto(z,-v)$};
 \draw[->,thick] (3,.2)--(4.5,.2)
  node[midway,above] {fold};
 \begin{scope}[xshift=6cm]
  \fill[blue!10] (0,-1.5) rectangle (2.2,1.5);
  \draw[gray] (0,-1.5)--(2.2,-1.5)--(2.2,1.5)--(0,1.5);
  \draw[very thick] (0,-1.5)--(0,1.5);
  \draw[->] (0,1.5)--(0,1.9) node[above] {$z$};
  \draw[->] (0,0)--(2.5,0) node[right] {$\rho=|v|$};
  \fill[blue!70!black] (1.3,.6) circle (2pt);
  \node[above] at (1.3,.6) {$[p]=[\tau p]$};
  \node[above,align=center] at (1.1,2.35)
   {$Q=M/O(3)=S/\langle\tau\rangle$};
  \node[below,align=center] at (1.1,-1.65)
   {interior points represent $S^2$ orbits\\
    boundary points represent fixed points};
 \end{scope}
\end{tikzpicture}
\caption{The meridian surface $S$ in the non-Bryant case of
Theorem~\ref{thm:intro-metric-blowdown}. In the equivariant model
$\mathbb R\oplus\mathbb R^3$, the subgroup $O(2)$ fixes
$S=\{(z,(v,0,0))\}$. Each principal two-sphere orbit meets $S$
in the two points $p$ and $\tau p$; the fixed axis is $v=0$.
Here $v$ is a signed coordinate in the equivariant model. The surface is complete
and totally geodesic, with $K_S>0$.
Folding its two halves by the angular involution gives the orbit
quotient $Q$, with the quotient length metric. The involution $\tau$ fixes the
axis pointwise.}
\label{fig:intro-meridian}
\end{figure}

Theorem~\ref{thm:intro-metric-blowdown} describes the geometry seen
at distance scales much larger than the local curvature scale.
The two-sphere orbits disappear in this limit, leaving a sector
that records the geometry transverse to the rotations. Its angle
is independent of the base point and the sequence of scales.
Thus uniform scalar decay has a metric characterization: it holds
exactly when this sector degenerates to a ray. If uniform decay
fails, the limit retains a strictly positive opening angle.
The curvature formula gives a second description of this distinction:
in the non-Bryant case the ray corresponds to total meridian
curvature $2\pi$, whereas
a positive angle corresponds to a strictly smaller total curvature.

\paragraph{The mean curvature flow analogy.}
For mean curvature flow, a translating soliton is the counterpart of
a steady Ricci soliton. Haslhofer \cite{Haslhofer} proved
uniqueness of the bowl translator under noncollapsing and uniform
two-convexity. Brendle--Choi \cite{BC} proved that a
noncompact, strictly convex, noncollapsed ancient solution in
$\mathbb R^3$ is a translating bowl. Their neck improvement method,
like the corresponding Ricci flow arguments, improves approximate
rotations on a large parabolic region and then promotes the improved
fields to exact symmetries.

Higher-dimensional translating solitons display a broader range of
geometry. Hoffman--Ilmanen--Mart\'in--White \cite{HIMW} constructed
families of graphical translators in higher dimensions. Zhu proved
$SO(2)$ symmetry for translators in $\mathbb R^4$ arising as blow-up
limits of embedded mean-convex flows \cite{ZhuFour}, and
$SO(n-1)$ symmetry for uniformly three-convex translating blow-up
limits in higher dimensions \cite{ZhuHigher}. Choi--Haslhofer--Hershkovits \cite{CHH} classified noncollapsed
translators in $\mathbb R^4$, including the oval bowls. These results support
the expectation that cylindrical regions and two distinguished cap
regions should determine a common lower-dimensional rotation group.

\subsection{Strategy and organization}
In this subsection, we explain the stages of the symmetry proof.
The main difficulty is that a single model does not describe all
regions at infinity. After excluding the four-dimensional Bryant
soliton, the local models are shrinking $S^2\times\mathbb R^2$
neck regions and $B^3\times\mathbb R$ cap regions, where $B^3$ is
the three-dimensional Bryant soliton. This division into an edge
model and an intervening neck parallels Yi Lai's symmetry argument
\cite{LaiSymmetry}, where cigar-product regions and a flat cylindrical
region must be treated together. Here the noncollapsed
three-dimensional factor is Bryant and the angular algebra has three
generators. We must improve the local symmetry errors, make the
choices of rotation fields compatible, and finally turn them into
exact Killing fields. The local improvement is given by Theorem~\ref{cylinder:improvement}
and Corollary~\ref{cor:alltiltheat}; Proposition~\ref{prop:globalrotations}
combines the resulting triples into global approximate rotations.

We normalize by
\begin{equation}\label{eq:normalization}
 \Ric=\Hess f,\qquad X:=\nabla f,\qquad \mathrm{R}+|X|^2=1,
 \qquad g_t=\Phi_t^*g.
\end{equation}
The bound $|X|\le1$ makes the canonical flow complete.
Chan--Ma--Zhang \cite[Theorem~1.2]{CMZ} show that, for every
sequence $x_i\to\infty$, a subsequence of the curvature-rescaled
canonical flows converges smoothly to a complete ancient product
$N^3\times\mathbb R$. Their tangent flow alternative
\cite[Corollary~5.4]{CMZ} states that the tangent flow at infinity
is either $S^3\times\mathbb R$ or $S^2\times\mathbb R^2$.
In the first case, \cite[Theorem~5.2]{CMZ} identifies $(M,g)$ with
the four-dimensional Bryant soliton, which already has the action
in Theorem~\ref{thm:main}. In the second case,
\cite[Theorem~1.3]{CMZ} identifies $N^3$ as the three-dimensional
Bryant soliton or the shrinking round cylinder $S^2\times\mathbb R$.
Thus the four-dimensional local models are $B^3\times\mathbb R$
and shrinking $S^2\times\mathbb R^2$ cylinders. The construction
below concerns this second branch.

We divide the argument into five stages. An approximate rotation
is a vector field whose flow changes the metric by a small amount;
this change is measured by its metric Lie derivative $\Lie_Ug$. The local triples
also have small components normal to the relevant two-spheres.
The aim is to improve these errors, align the triples on overlaps,
and finally obtain fields with zero metric Lie derivative. All constructions before
the angular action use local geometry or an auxiliary exhaustion;
they do not require the potential to be proper.

\paragraph{(A) Uniform improvement on necks and tilted caps.}
We begin with Brendle's neck improvement method
\cite[Theorem~8.6]{BrendleAncient}. First, we align nearby approximate
rotations and glue them on a larger region. We then evolve the glued
fields by the vector heat equation. Their metric Lie derivatives
satisfy the linearized Ricci flow equation, and the spectrum of the
shrinking cylinder improves the components that break sphere symmetry.
The slowly decaying terms are handled separately, and Lie brackets
produce a triple with a smaller symmetry error. The precise conclusion
is Theorem~\ref{cylinder:improvement}; its alignment and gluing estimates
are Lemmas~\ref{cylinder:lem:alignment} and~\ref{cylinder:lem:gluing}.

In our four-dimensional cylinder, the spheres have two normal
directions. One comes from the potential; the other lies within a
potential level and is normal to its constant mean curvature (CMC)
spheres. We estimate the potential pairing alongside the metric
Lie derivative. A first bracket controls this pairing; after freezing the
fields in time, the CMC variation equation and a second bracket
control the remaining normal component. We also obtain one constant
change of basis that normalizes the triple on every sphere
and at every time.

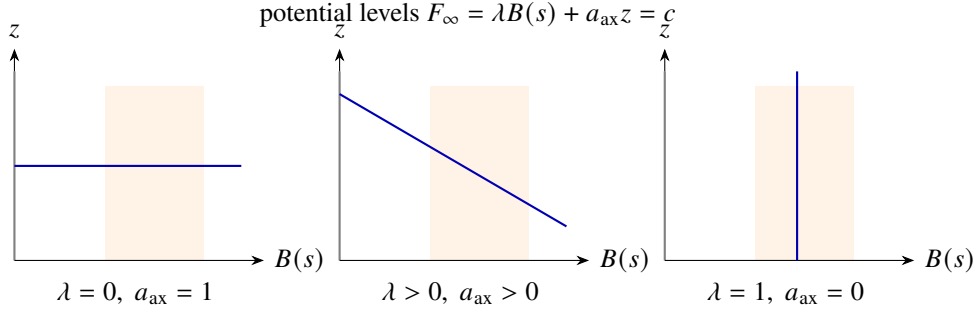
\begin{figure}[!htbp]
\centering
\begin{tikzpicture}[x=1cm,y=1cm,>=Stealth,font=\small]
\foreach \shift/\title in {0/{$\lambda=0,\ a_{\rm ax}=1$},4.3/{$\lambda>0,\ a_{\rm ax}>0$},8.6/{$\lambda=1,\ a_{\rm ax}=0$}} {
 \begin{scope}[xshift=\shift cm]
 \fill[orange!10] (1.2,0) rectangle (2.5,2.3);
 \draw[thick,gray] (0,0)--(0,2.5);
 \draw[->] (0,0)--(3.3,0) node[right] {$B(s)$};
 \draw[->] (0,2.5)--(0,2.8) node[above] {$z$};
 \node[below] at (1.6,-.15) {\title};
 \end{scope}
}
\draw[blue!70!black,thick] (0,1.25)--(3,1.25);
\draw[blue!70!black,thick] (4.3,2.2)--(7.3,.45);
\draw[blue!70!black,thick] (10.35,0)--(10.35,2.5);
\node at (6,3.25) {potential levels $F_\infty=\lambda B(s)+a_{\rm ax}z=c$};
\end{tikzpicture}
\caption{Potential levels in the local Bryant-product quotient.
Each point with $B(s)>0$ represents a two-sphere; $B(s)=0$ is the product
of the Bryant tip with the line. The axes are product coordinates. The shading marks a fixed collar away
from the tip. When $a_{\rm ax}=0$ and $c>0$, the level is $B(s)=c$ and is
parametrized by $z$, so the graph $z=(c-\lambda B(s))/a_{\rm ax}$ is unavailable.
The collar estimates of Proposition~\ref{prop:alltiltcollar} and the
product domains of Proposition~\ref{prop:productdomains} remain
uniform across these cases.}
\label{fig:introtilt}
\end{figure}

Near a cap the local model is $B^3\times\mathbb R$. Before applying
the heat improvement, we choose domains and measurements that remain
uniform when the potential tilts away from the product line.
The potential direction in a Bryant-product chart need
not approach the additional line. Figure~\ref{fig:introtilt}
illustrates the resulting potential levels. In a tip-normalized
chart, the time-zero slice of the limiting normalized potential has the form
\begin{equation}\label{eq:introtilt}
 F_\infty=\lambda B(s)+a_{\rm ax}z,\qquad
 \lambda^2+a_{\rm ax}^2=1,\qquad 0\le\lambda,a_{\rm ax}\le1,
\end{equation}
Here $B^3$ denotes the three-dimensional Bryant soliton, with tip $o$
and scalar curvature $\mathrm R_B(o)=1$. The variable $s=d_{g_B}(o,\cdot)$
is geodesic distance from the tip, and $B(s)$ is its soliton potential,
characterized by $\Hess_{g_B}B=\Ric_{g_B}$ and $B(0)=0$.
The coordinate $z$ is arclength on the additional $\mathbb R$ factor,
oriented so that $a_{\rm ax}\ge0$. We write
$g_B:=\mathrm{d}s^2+\omega(s)^2g_{S^2}$, where $g_{S^2}$ is the unit
round metric and $\omega(s)$ is the radius of its two-sphere orbits. The gradient
of $F_\infty$ at the tip has length $a_{\rm ax}$, which may vanish.
For $a_{\rm ax}>0$ the induced metric on a potential level is
\begin{equation}
 \left(1+\frac{\lambda^2}{a_{\rm ax}^2}\frac{\mathrm d B}{\mathrm ds}(s)^2\right)\mathrm{d}s^2
                   +\omega(s)^2g_{S^2}.
\end{equation}
This expression loses uniform control as $a_{\rm ax}\downarrow0$ because
solving $F_\infty=c$ for $z$ introduces the factor $a_{\rm ax}^{-1}$.
Instead, we fix $s_+>s_-\ge s_0>0$ and $Z_*>0$ and work on the annulus
$\mathcal A_{\rm col}=\{s_-\le s\le s_+,\ |z|\le Z_*\}$ in the product metric
$g_0=g_B+\mathrm dz^2$. There, we have
\begin{equation}
 q^2:=|\nabla F_\infty|_{g_0}^2
 =a_{\rm ax}^2+\lambda^2\frac{\mathrm d B}{\mathrm ds}(s)^2
 \ge\min\{1,\inf_{[s_-,s_+]}\frac{\mathrm d B}{\mathrm ds}(s)^2\}>0,
\end{equation}
uniformly for $a_{\rm ax}\in[0,1]$. Thus the two normal directions
\[
 N_0=q^{-1}(\lambda \frac{\mathrm d B}{\mathrm ds}\partial_s+a_{\rm ax}\partial_z),\qquad
 \nu_0=q^{-1}(a_{\rm ax}\partial_s-\lambda \frac{\mathrm d B}{\mathrm ds}\partial_z)
\]
remain controlled on the collar even when $a_{\rm ax}=0$.
Proposition~\ref{prop:alltiltcollar} proves the corresponding uniform
CMC foliation and Jacobi inverse estimates for nearby metrics and
potentials, with fixed interior margins. Here $s_0$ is the threshold
in Proposition~\ref{prop:jacobi}. At the tip, however, $\frac{\mathrm d B}{\mathrm ds}(0)=0$ and $q=a_{\rm ax}$ can vanish.
The core estimates therefore use fixed product domains containing
$s=0$ and measure the potential pairing by
\[
 UF=\langle U,\nabla F\rangle,
 \qquad\text{rather than}\qquad
 \langle U,n_f\rangle=\frac{UF}{|\nabla F|}.
\]
Here $F$ is the rescaled heat potential, and
$n_f:=\nabla F/|\nabla F|$ is defined only on its regular levels.
The pairing $UF$ remains defined at a critical point. The heat estimates
on these domains and their restriction to regular collars are given
in Corollary~\ref{cor:alltiltheat} and
Proposition~\ref{prop:productdomains}; their constants require no
positive lower bound for $a_{\rm ax}$.

On these product domains, a tensor comparison estimate replaces the
cylindrical spectral argument and improves the metric Lie derivative
on a fixed smaller cap. The extra line requires control of mixed
and line components as well as tensors tangent to the Bryant factor;
Proposition~\ref{cap:prop:product-quotient} supplies this estimate.
For collars bounding controlled domains inside regular potential
levels, Proposition~\ref{prop:normal-component} controls the CMC
normal component. Proposition~\ref{prop:heatcap} gives the heat
improvement when the axial coefficient is bounded below.
Corollary~\ref{cor:alltiltheat} treats all tilts, including
$a_{\rm ax}=0$, using the annular bracket estimate of
Proposition~\ref{prop:potentialbracket} without assuming bounded
potential level interiors.
Proposition~\ref{prop:commongram} supplies one constant normalization
on all cap collar leaves and times.

\paragraph{(B) Intrinsic geometry and exterior iteration.}
To repeat the local improvements, the triple constructed in one
application must satisfy the hypotheses of the next. The potential is not yet known
to be proper, so its levels cannot serve as an exhaustion. We instead
locate the cap tips using the scalar curvature and the Ricci tensor.
The use of Ricci eigenvalues relative to scalar curvature to distinguish
caps from necks follows Brendle's approach
\cite[Proposition~9.1 and Definition~9.2]{BrendleAncient}.
Here we also use the scalar curvature gradient transverse to the
additional product line to select the tip curve.
In a region where the lowest Ricci eigenvalue is separated from the
others, let $L$ denote its \textbf{lowest Ricci eigenline} and set $E:=L^\perp$,
with $P_E$ the orthogonal projection onto $E$. The intrinsic tip set
$\Gamma:=\{P_E\nabla\mathrm R=0\}$, restricted to the eigenvalue-gap
region of Proposition~\ref{prop:tipcurves}, is a union of smooth
curves outside a compact set: the model scalar curvature has a nondegenerate transverse Hessian there.
Propositions~\ref{prop:tipcurves} and~\ref{prop:tipnormal} establish
these curves and their normal tubes.

We define the cap domains intrinsically and reserve spatial and time
margins, including the overlap where a cap triple supplies a neck
triple. Common actual CMC leaves make the two constructions
compatible. Lemma~\ref{lem:compact-exclusion} and
Propositions~\ref{prop:exteriorleaves} and~\ref{prop:intrinsicwitnesses}
provide the compact exclusion sets, the leaves, and the compatible
local triples.

We then measure a weighted exterior symmetry error $\mathscr E(x)$. If
$\mathcal I(x)$ is the set of centers whose local triples are needed
for one improvement at $x$, and $Q$ is the assigned curvature scale,
the neck and cap recurrences can be written schematically as
\[
 \mathscr E(x)\le\tfrac12\sup_{z\in\mathcal I(x)}\mathscr E(z),\qquad
 Q(z)\le(1+\sigma)Q(x).
\]
The entire finite cap construction is included in one application
of this recurrence. For sufficiently small $\sigma$, the decrease
of the symmetry error dominates the possible increase of the scale.
More precisely, Proposition~\ref{prop:typedclosure} gives, for any
$p>0$,
\[
 e_N(x,0)\le C_p(1+d_g(o,x))^{-p},\qquad
 e_C(q)\le C_p(1+d_g(o,q))^{-p},
\]
for neck centers outside the specified cap core and for intrinsic
cap centers $q\in\Gamma$, respectively, after excluding a compact
set. Here $o\in M$ is a fixed reference point, not an assumed
critical point. Thus the local rotations become
increasingly accurate as their centers escape the compact core.

\paragraph{(C) From local triples to global fields.}
The improved triples must be aligned before they can be combined into global fields.
A change of basis between two triples is an orthogonal matrix.
Choosing such matrices independently on different overlaps can lead
to an inconsistency around a loop. We need a consistent choice of
basis, together with an exhaustion along which to propagate it. Greene--Wu \cite{GW}
supply a smooth proper strictly convex function $\psi$ with bounded
gradient. This auxiliary function controls the exterior geometry.

To control these choices, we record each approximate rotation together
with the skew-adjoint part of its covariant derivative. These pairs
determine a nearby rank-three subbundle of
$TM\oplus\mathfrak{so}(TM)$. A frame of this bundle represents a
choice of three rotations, and its connection records how that choice
varies between charts. Proposition~\ref{prop:rotationconnection} constructs
a metric connection with rapidly decaying curvature. Outward parallel
transport along $\psi$ then gives an $L^1$ limiting connection
(Proposition~\ref{prop:radialconnectionlimit}).
Lemma~\ref{lem:Loneparallel} and
Proposition~\ref{prop:smoothexteriorgauge} provide a weak parallel
frame and its quantitative smoothing. Taking the vector field component of this smooth frame gives three
global approximate rotations $\mathcal U_1,\mathcal U_2,\mathcal U_3$. Their size and derivatives
grow at most polynomially. Their metric Lie derivatives and potential
pairings have any prescribed polynomial decay, while their cyclic
bracket defects tend uniformly to zero. In Proposition~\ref{prop:globalrotations},
we will obtain the following quantitative bounds. For every $p>0$, outside a
compact set depending on $p$, and with $r(x):=d_g(o,x)$,
\[
 |\mathcal U_i|+|\nabla \mathcal U_i|\le C_p(1+r),\qquad
 \sum_{j=0}^2\bigl(|\nabla^j\Lie_{\mathcal U_i}g|
                  +|\nabla^j(\mathcal U_if)|\bigr)
 \le C_p(1+r)^{-p},\qquad i=1,2,3.
\]
For cyclic $(i,j,k)$, the bracket defects satisfy
$|[\mathcal U_i,\mathcal U_j]-\mathcal U_k|\to0$ uniformly as $r\to\infty$.
The arbitrary polynomial decay absorbs the finite losses in the
alignment and smoothing constructions.

\paragraph{(D) Exact rotations and the angular action.}
Recall that $X:=\nabla f$ and set $\mathcal D:=\Delta+\nabla_X$,
where $\Delta$ is the rough Laplacian on vector fields. Choose
$p>9$ in Proposition~\ref{prop:globalrotations}, and extend each
exterior field $\mathcal U_a$ smoothly across the compact core
using a cutoff in $\psi$; denote the global extension by $U_a$.
The cutoff changes the error only on a compact set. We subtract a
decaying solution $V_a$ of the drift equation
\[
                      \mathcal D V_a=\mathcal D U_a
\]
and set $W_a:=U_a-V_a$. Proposition~\ref{prop:globalcorrection}
gives $|V_a|+|\nabla V_a|+|\nabla^2V_a|\le C(1+r)^{-(p-5)}$.
The metric Lie derivative $\Lie_{W_a}g$ now solves the
stationary drift--Lichnerowicz equation. Our tensor rigidity result
shows that an $L^1(\mathrm d\mu_g)$ solution of this equation
vanishes, using
positive curvature and integration in the ordinary volume measure.
Thus the $W_a$ are Killing fields; the precise statements are
Proposition~\ref{prop:distanceLone} and
Corollary~\ref{cor:globalKilling}.

This correction follows the approach in Brendle's proofs of rotational
symmetry for steady solitons \cite{Brendle3,BrendleHigher}.
In the asymptotically cylindrical setting, his tensor rigidity theorem
identifies solutions with the prescribed pointwise decay as constant
multiples of $\Ric$ \cite[Theorem~4.4]{BrendleHigher}. Subtracting the
corresponding multiple of $X$ from the corrected vector field removes
this remaining mode and gives a Killing field. Here
Proposition~\ref{prop:distanceLone} gives vanishing directly from
unweighted integrability. Its distance-cutoff argument requires
neither a proper potential nor an asymptotically cylindrical end.
The choice $p>9$ makes the corrected metric Lie derivatives decay
faster than $r^{-4}$, supplying the needed unweighted integrability.

The exact bracket relations integrate the three fields to a compact
group action. Averaging the auxiliary strictly convex exhaustion
supplies a fixed point. The existence of escaping points where a unit linear combination
of the generators tends to zero, as in
\eqref{eq:almostisotropy}, rules out the four-dimensional quaternionic
representation of $SU(2)$. The remaining representation is the usual three-dimensional
rotation representation plus a fixed line. It gives an effective
$SO(3)$ action with two-sphere orbits, and the invariant metric extends
this action to $O(3)$. Proposition~\ref{prop:orbitcompletion}
carries out this group action step.

After constructing the angular action, we prove
Theorem~\ref{lem:metric-blowdown} in three steps. First, the curvature scale size
bound \eqref{eq:exact-generator-size} for the exact Killing fields, together with
Proposition~\ref{prop:curvaturelength}, shows that the rotation
orbits have diameter $o(r)$ in balls of radius $r$. They therefore
collapse under distance rescaling by $r^{-1}$, and the soliton has
the same metric limit as its two-dimensional quotient.

Second, the standard action from
Proposition~\ref{prop:orbitcompletion} allows the proof of
Theorem~\ref{lem:metric-blowdown} to identify this quotient with
one half of a complete totally geodesic meridian
surface $S$, with its intrinsic induced metric. Positive sectional curvature gives $K_S>0$, and the
surface ideal-boundary description used in
Theorem~\ref{lem:metric-blowdown} determines a unique limiting cone
of total angle $2\pi-\int_S K_S\,\mathrm{d}\mu_S$. The angular
involution folds this cone onto the sector of half that angle;
Figure~\ref{fig:intro-meridian} illustrates the folding.

Finally, Lemma~\ref{lem:decay-ray-blowdown} shows that uniform
scalar decay collapses the distance spheres to a ray. Conversely,
Lemma~\ref{lem:ray-direct-converse} uses their sublinear diameter,
convexity of the potential, and the soliton identity to recover
uniform scalar decay from a ray limit. These two implications do
not require angular symmetry or an assumed critical point.

\paragraph{(E) Critical point geometry.}
Suppose for contradiction that
$\nabla f$ never vanishes. Lemma~\ref{reflcrit:foliation} then identifies
every potential level with a complete positively curved rotational
$\mathbb R^3$ having one tip on the fixed axis. The normal-potential flow
identifies these levels, and the restriction of $f$ to the axis is
strictly increasing. We study this hypothetical geometry through the
orbit radius $u$ and the quotient metric $h$.

The first step is to rule out a positive limiting scalar curvature at
the outgoing end of the axis. The decorated sector limit records both
the potential and $u^2$ on the metric blow-down. A positive end limit
would produce a ridge with two distinct slopes for the limiting squared
radius. The heat equation comparison in Lemma~\ref{reflcrit:heat-ridge}
excludes that ridge, giving Corollary~\ref{reflcrit:Lzero}.
The radius envelope and the ambient singular-profile barrier then rule
out a finite infimum of $f$, as shown in
Corollary~\ref{reflcrit:infinite-infimum}.

In the remaining case, $\inf_M f=-\infty$, high-potential levels have
uniform scalar decay and contain a fixed normal trajectory with
$u^2\sim2c$, where $c:=f$. The one-cap profile analysis in
Lemma~\ref{reflcrit:onecap} sharpens this to
$F(0,c)=\sqrt{2c}+O(1)$ and ${\partial_{z} F}(0,c)=O(c^{-1/2})$.
These estimates allow the ambient barrier of
Lemma~\ref{reflcrit:barrier} to be applied on the tip-to-reference
segments. At a fixed radius, this barrier bounds the level slope
strictly below one, whereas scalar decay and the meridian Jacobi
equation force that slope to tend to one. This contradiction proves
existence of a critical point in Proposition~\ref{reflcrit:critical}.
Finally, $\nabla^2f=\Ric>0$ gives uniqueness and a linear lower bound
for $f$ in distance from that point; completeness then gives properness
by Lemma~\ref{lem:critical-properness-II}.

\paragraph{Organization of the paper.}
The paper is organized as follows.
Section~\ref{sec:collars} develops the model geometry and collar
estimates. Sections~\ref{sec:cylinder} and~\ref{sec:capimprovement}
prove the cylinder and cap improvements. Section~\ref{sec:exterior}
constructs the intrinsic exterior geometry and closes the recurrence.
Section~\ref{sec:globalization} aligns the local rotations on overlaps.
Section~\ref{sec:exact} proves the global analytic correction and
constructs the angular action, determines the metric blow-down,
and establishes the critical point geometry. 

Theorem~\ref{thm:intro-metric-blowdown} is proved as
Theorem~\ref{lem:metric-blowdown} in
Subsection~\ref{sec:exact-action}. The critical point assertion of
Theorem~\ref{thm:main} is proved in
Proposition~\ref{reflcrit:critical} of that subsection. The proof of
Theorem~\ref{thm:main} is completed at the end of
Subsection~\ref{sec:exact-action}.

\paragraph{Acknowledgements.} The author thanks Yu Li for helpful discussions on this problem. The author is supported by Hong Kong RGC grants No. 14300623
and No. 14304225.

\paragraph{AI disclosure.}
The author prepared an initial note containing the main ideas, proof
strategy, intermediate results, and proof sketches for the
$O(3)$-symmetry of steady solitons with uniform scalar decay. He
subsequently used ChatGPT-6 Astra and the Danus system to help expand
the technical arguments and improve the presentation. This assistance
also included work on removing the scalar decay hypothesis. The author
reviewed and verified all AI-assisted material and takes full
responsibility for the mathematical claims, proofs, and source
attributions in this manuscript.

\section{Model geometry and uniform collar estimates}\label{sec:collars}
In this section, we describe the local models supplied by dimension
reduction and obtain collar estimates uniform in the coefficients of
the limiting potential, including a vanishing axial coefficient.
Lemma~\ref{lem:boundedlocalgeometry} gives the local coefficient bounds,
and Proposition~\ref{prop:alltiltcollar} constructs the uniform collars.
The transverse estimates in Propositions~\ref{prop:potentialbracket}
and~\ref{prop:normal-component} will control the two normal components
of approximate rotations in the cap improvement.
The main results are Propositions~\ref{prop:alltiltcollar},
\ref{prop:potentialbracket}, and~\ref{prop:normal-component}.

\subsection{Soliton identities, model geometry, and notation}
We fix the operators and normalization used throughout the paper. The identities \eqref{eq:basicidentities} and~\eqref{eq:scaledpotentialidentities}, Lemma~\ref{lem:boundedlocalgeometry}, and the Bryant estimates \eqref{eq:BryantODE}--\eqref{eq:Bryantasymptotics} provide the model and coefficient bounds used in the collar estimates.

Throughout, $g_{S^2}$ denotes the round metric of sectional curvature one.
The \textbf{scalar-normalized round cylinder} is
\begin{equation}
 g_{\rm cyl}:=2g_{S^2}+\mathrm dz_1^2+\mathrm dz_2^2
 \quad\text{on }S^2\times\mathbb R^2,\qquad
 \mathrm R_{g_{\rm cyl}}=1.
\end{equation}
Thus its spherical factor has radius $\sqrt2$. Every scalar-normalized
round-product limit below uses this metric; it is the time $-1$ slice
of the shrinking cylinder in \eqref{cylinder:eq:cylinder}.
For scalar functions $v$, we use $\nabla_gv$ for the gradient and
$|\nabla_gv|_g$ for its norm, with the metric subscript omitted when
unambiguous. On the orbit quotient the gradient is $\nabla_hv$.
The differential $\mathrm dv$ is retained when a covector is required.

The \textbf{Laplacian} is $\Delta_g:=\nabla^i\nabla_i$, with
nonpositive spectrum on closed manifolds. On vector fields and tensors
it is the rough connection Laplacian. We suppress the metric subscript
only when the current metric is unambiguous. The
\textbf{Lichnerowicz Laplacian} is defined by
\begin{equation}\label{eq:lichnerowicz-convention}
 (\Delta_Lh)_{ij}:=\Delta h_{ij}+2\Rm(h)_{ij}
                        -\Ric_i{}^k h_{kj}-\Ric_j{}^k h_{ik}.
\end{equation}

For a symmetric two-tensor $h$, the \textbf{curvature action} is
\begin{equation}\label{eq:curvature-action}
 \Rm(h)_{ij}:=R_{ikjl}h^{kl},\qquad
 \Rm(h,h):=\langle\Rm(h),h\rangle.
\end{equation}
Indices are raised with the current metric. Our four-tensor convention
satisfies $\Rm(v,w,v,w)=\sec(v\wedge w)$ for orthonormal $v,w$.
For a covariant two-tensor $T$, $T^2$ means
$(T^2)_{ij}:=T_i{}^kT_{kj}$. The notation $\Rm*U$ in the Killing transport equations denotes a schematic contraction with a vector field $U$. We write $\Ric^\sharp$ for its \textbf{Ricci endomorphism}.
For a vector field $U$, set $UF:=\mathrm dF(U)$. The
\textbf{metric Lie derivative} has components
\[
 (\Lie_Ug)_{ij}=\nabla_iU_j+\nabla_jU_i.
\]
Thus $\Lie_Ug=2\operatorname{sym}\nabla U$ after lowering the vector
index with $g$; $\operatorname{sym}$ and $\operatorname{skew}$ denote
the symmetric and skew parts under this metric identification. The symbols $\mathrm d\mu_g$ and $\mathrm d\mu_S$
denote Riemannian volume and induced hypersurface area.
The notation $p\asymp q$ for positive quantities means
$C^{-1}q\le p\le Cq$ with fixed positive constants.
All distances defining escape or polynomial weights use the original
metric $g$. A fixed base point may be changed when the
auxiliary exhaustion is chosen: for any fixed $o,o'$,
\[
 (1+d(o,o'))^{-1}(1+d(o,x))\le1+d(o',x)
       \le(1+d(o,o'))(1+d(o,x)).
\]
Changing the base point therefore alters only the constants in the
polynomial estimates. Unless a scale is displayed, local norms are
measured in the stated normalized chart. A parabolic derivative counts one time
derivative as two spatial derivatives. The \textbf{error notation} $O_{C^k}(e)$ means
the estimate
\[
 \|A\|_{C^k(\Omega)}
 :=\sum_{j=0}^k\sup_\Omega|\nabla^j A|\le Ce
\]
on the stated interior domain $\Omega$. For a spacetime domain $Q$, a plain $C^k(Q)$ norm uses spatial
derivatives, uniformly in time. When time derivatives are required,
we write
\[
 \|A\|_{\mathrm{par},k;Q}
 :=\sum_{j+2\ell\le k}\sup_Q|\nabla^j\partial_t^\ell A|.
\]
For separate spatial and time orders we use the \textbf{mixed derivative norm}
\[
 \|A\|_{C_t^\ell C_x^k(Q)}
 :=\sum_{j=0}^k\sum_{b=0}^\ell
          \sup_Q|\nabla_x^j\partial_t^b A|.
\]
Thus $C^k(Q)$ and $O_{C^k}$ contain no time derivatives, whereas
$C_t^\ell C_x^k$ includes exactly the displayed time orders; the
parabolic norm above instead uses the total order $j+2\ell\le k$.
These conventions also apply when the evolution parameter is $c$, $T$,
or $\tau$. In particular, $\partial_x^j\partial_T H$ denotes one $T$ derivative
and $j$ spatial derivatives. The connection and tensor norm are those
of the specified actual or model metric. Time derivatives are taken
at fixed spatial coordinates in the indicated chart. Derivatives are
written explicitly as $\partial_t A$ or $\frac{\mathrm d}{\mathrm dt}A$;
a subscript on a quantity denotes an index, component, or parameter
dependence, never differentiation.
The domain and its margins are fixed before the comparison tolerance
is chosen.
The norms $C^{k,\alpha}$ are the usual H\"older norms in the fixed
comparison charts, or in the explicitly specified leaf charts.
Unless otherwise stated, $C>1$ denotes a constant depending on the fixed
geometric bounds, derivative order, and domain margins; its value may change
from line to line. Small positive comparison constants are written as $C^{-1}$;
$c$ is reserved for the level parameter. An error parameter such as $e$, $\epsilon$, or
$\eta$ is nonnegative; $\eta$ measures geometric comparison (see \eqref{cylinder:eq:eta-cylindrical}), whereas
$e$ or $\epsilon$ measures the error of the vector fields.
Smallness thresholds are chosen after the domains and margins are fixed.
A \textbf{controlled geometry} on a fixed normalized domain means that,
in the indicated charts, the metric and its inverse have uniform
$C^{m,\alpha}$ bounds through the orders required in the statement,
with a uniform positive lower bound for the metric. For a
\textbf{controlled family} of collars or leaves, the parametrizations,
their inverse coordinate maps, and their derivatives in the leaf
parameters have the same uniform bounds on the prescribed enlargements.
Regular potential coordinates also require a fixed lower bound for
$|\nabla F|$, and a Jacobi inverse bound is required whenever stated.
These bounds, together with the fixed domain margins, are the background
data in Propositions~\ref{prop:alltiltcollar} and~\ref{prop:potentialbracket}
and Lemma~\ref{lem:rotation-gram}; constants are uniform only while
these data remain fixed.

We use the constants $\kappa_{\mathrm{pin}}$ for the
cap quotient absorption constant, $c_{\mathrm{pin}}$ for the distinct
perturbation threshold in Proposition~\ref{cap:prop:product-quotient}, and
$\krot=\sqrt{6\pi}$ for the
\textbf{bracket constant} of area-normalized round rotations. The \textbf{alternating symbol} $\epsilon_{ijk}$ has
$\epsilon_{123}=1$, changes sign under an exchange of two indices,
and vanishes whenever two indices coincide; repeated algebra indices are summed unless a
sum is written explicitly. For example,
\[
 [J_i,J_j]=\epsilon_{ijk}J_k
 :=\sum_{k=1}^3\epsilon_{ijk}J_k,
 \qquad [J_1,J_2]=J_3.
\]
Here the middle expression abbreviates the explicit sum; the bracket
identity fixes the \textbf{model rotation convention}. We write $J_i$
for model rotations with cyclic bracket constant one and
$K_i:=\krot J_i$ for those with normalized Gram matrix $I$.
Throughout, a \textbf{normalized Gram matrix} means
\[
 (Q_S(U))_{ij}:=\Area(S)^{-2}
       \int_S\langle U_i,U_j\rangle\,\mathrm d\mu_S.
\]
For a real matrix $A=(A_{ij})$, the Frobenius norm is
\[
 |A|_{\mathrm F}:=\left(\sum_{i,j}|A_{ij}|^2\right)^{1/2}
 =\bigl(\operatorname{tr}(A^{\mathsf T}A)\bigr)^{1/2}.
\]
For \textbf{spherical projections} we use
$\Pi_0$ for averaging and $P:=I-\Pi_0$ for mean-zero projection;
the relevant measure is specified locally. The \textbf{canonical potential} is $f_t:=f\circ\Phi_t$.
Where the potential gradient is nonzero, we write
\[n_f:=\frac{\nabla f_t}{|\nabla f_t|}.\]
This is the \textbf{potential level unit normal}.
Equivalently, $n_f=\nabla F/|\nabla F|$ in the heat gauge. We use
$\nu$ for the unit normal to a CMC sphere within that level. The term CMC means constant mean curvature with respect to the
metric induced on the indicated three-dimensional level. The ambient
normal tube coordinate $\beta$ is defined intrinsically from the tip curve
in Proposition~\ref{prop:tipnormal}; the model radial weight is $B_t$.

We distinguish the \textbf{geometric comparison error} from the
\textbf{rotation error}: $\eta$ controls the metric and normalized
potential on an enlarged model chart, while $\epsilon$ or $e$ controls
$\Lie_Ug$, $UF$, transverse pairings, and the normalized Gram matrix
on the stated measurement sets. Constants denoted $C_E$ may depend on
the comparison size $E$; the fixed domains are fixed before $E$ and
then $\eta$ are chosen. On a chart of curvature length $r$, all such
norms use $\widehat g=r^{-2}g$ unless the original metric is specified.
In particular, model radial position $B(s)$ and the actual ambient
tube coordinate $\beta$ are compared only on the specified charts.

We next record the soliton identities and the rescaled Bryant-product
model used in the local improvement arguments.
The contracted Bianchi identity and $\Ric=\Hess f$ give
\begin{align}
 \Delta f&=\mathrm{R},& \mathrm d\mathrm{R}&=-2\Ric(X,\cdot),&
 (\Delta+X)\mathrm{R}&=-2|\Ric|^2.\label{eq:basicidentities}
\end{align}
It follows that $\mathrm{R}+|X|^2$ is constant. Since the soliton
has positive scalar curvature, this constant is positive. Scaling the metric gives
\eqref{eq:normalization}, i.e. 
\[\mathrm{R}+|X|^2=1.\]
 If $f_t=f\circ\Phi_t$, then we obtain
\[
 \partial_tf_t=-|\nabla f_t|_{g_t}^2=\mathrm{R}_{g_t}-1,
 \qquad (\partial_t-\Delta_{g_t})(f_t+t)=0.
\]
At scale $r$, we set $\widehat g(t):=r^{-2}g(r^2t)$ and define
\[
  F(t):=r^{-1}f_{r^2t}+rt,\qquad \lambda:=r^{-1}.
\]
Additive constants or time translations may be included when centering
a chart. Direct scaling yields
\begin{equation}\label{eq:scaledpotentialidentities}
 \Hess_{\widehat g}F=\lambda\Ric_{\widehat g},\qquad
 |\nabla F|_{\widehat g}^2=1-\lambda^2\mathrm{R}_{\widehat g},\qquad
 (\partial_t-\Delta_{\widehat g})F=0.
\end{equation}
These formulas specify the potential normalization throughout the local
heat arguments.

\begin{lem}\label{lem:boundedlocalgeometry}
Use the normalization \eqref{eq:normalization}, and set $X:=\nabla f$.
Under the hypotheses of Theorem~\ref{thm:main}, for every fixed
integer $m\ge0$ there are constants $C_m,i_0,v_0>0$ such that
\[
 |\nabla^m\Rm|+|\nabla^{m+1}X|\le C_m,\qquad
 \operatorname{inj}_g(M)\ge i_0,\qquad
 \operatorname{Vol}_gB_g(x,1)\ge v_0.
\]
\end{lem}
\begin{proof}
Positive curvature operator and $\mathrm{R}\le1$ give $|\Rm|\le C$ on $M$.
Since each $g_t$ is isometric to $g$, the same bound holds on the
entire eternal flow. We apply Shi's derivative estimates \cite{Shi}
on a unit backward interval ending at any chosen time to obtain
\[
 \sup_M|\nabla^m\Rm|_{g_t}\le C_m.
\]
The identities $|X|\le1$ and $\nabla X=\Ric^\sharp$ control the
drift and its derivatives. We choose $r_0\in(0,1)$ so that
$C\le r_0^{-2}$. Noncollapsing applies to every radius-$r_0$
parabolic ball and gives
\[
 \operatorname{Vol}_gB_g(x,1)
 \ge\operatorname{Vol}_gB_g(x,r_0)\ge\kappa r_0^4.
\]
The curvature and volume bounds give a uniform injectivity radius
by Cheeger--Gromov--Taylor \cite[Section 4]{CGT}.
This finishes the proof.
\end{proof}

We write the three-dimensional Bryant soliton, normalized by $\mathrm{R}_B(o)=1$,
as
\[
 g_B=\mathrm{d}s^2+\omega(s)^2\gamma,\qquad
 \gamma:=g_{S^2},\qquad B(0)=0.
\]
Here $o$ is the Bryant tip, $s:=d_{g_B}(o,\cdot)$ is radial
geodesic distance, and $\gamma$ is the unit round two-sphere metric.
We use the metric-specific abbreviations
$\Ric_B:=\Ric_{g_B}$ and $\mathrm{R}_B:=\mathrm{R}_{g_B}$.
The function $\omega(s)$ is the orbit radius; $B(s)$ is the radial
soliton potential, satisfying $\Hess_{g_B}B=\Ric_{g_B}$ and $B(0)=0$.
All derivatives of $B$ and $\omega$ below are taken with respect to the radial variable $s$.
Let 
\begin{align*}
	k:=-\frac{\mathrm d^2 \omega}{\mathrm ds^2}/\omega\qquad\text{and}\qquad l:=(1-\bigl(\frac{\mathrm d}{\mathrm d s}\omega\bigr)^2)/\omega^2
\end{align*}
denote
the radial and tangential sectional curvatures. The soliton equation
and normalization become
\begin{equation}\label{eq:BryantODE}
 \frac{\mathrm d^2 B}{\mathrm ds^2}=2k=-\frac{2\frac{\mathrm d^2 \omega}{\mathrm ds^2}}{\omega},\qquad
 \frac{\mathrm d B}{\mathrm ds}\frac{\frac{\mathrm d \omega}{\mathrm ds}}{\omega}=k+l
     =\frac{1-\bigl(\frac{\mathrm d}{\mathrm d s}\omega\bigr)^2-\omega\frac{\mathrm d^2 \omega}{\mathrm ds^2}}{\omega^2},\qquad
 (\frac{\mathrm d B}{\mathrm ds})^2+4k+2l=1.
\end{equation}
Smoothness at the tip gives $\omega(0)=0$, $\frac{\mathrm d \omega}{\mathrm ds}(0)=1$ and
$\frac{\mathrm d B}{\mathrm ds}(0)=0$. In particular, we have that for small $s$, 
\begin{equation}
 \Ric_B(o)=\tfrac13g_B(o),\quad
 B(s)=\tfrac16s^2+O(s^4),\quad
 \omega(s)=s-\tfrac1{36}s^3+O(s^5),\quad
 \Hess \mathrm{R}_B(o)=-\tfrac29g_B(o).
\end{equation}
To see the last identity, we can differentiate $\nabla \mathrm{R}_B=-2\Ric_B^\sharp(\nabla B(s))$
at the critical point and use $\Hess B(s)=\Ric_B$.

We use the standard Bryant asymptotics (see \cite{Bryant}): for large
$s$, with constants depending only on this normalization,
\begin{equation}\label{eq:Bryantasymptotics}
 \begin{gathered}
 \omega(s)\asymp s^{1/2},\quad
 0<\frac{\mathrm d \omega}{\mathrm ds}(s)\le Cs^{-1/2},\quad
 |\frac{\mathrm d^2 \omega}{\mathrm ds^2}(s)|\le Cs^{-3/2},\\
 \frac{\mathrm d B}{\mathrm ds}(s)\longrightarrow1,\quad 0<\frac{\mathrm d^2 B}{\mathrm ds^2}(s)\le Cs^{-2},\quad
 C^{-1}(1+s)^{-1}\le \mathrm{R}_B(s)\le C(1+s)^{-1}.
 \end{gathered}
\end{equation}
The leading coefficients used in the cylindrical limit also follow
from these estimates and \eqref{eq:BryantODE}. In fact, for large $s$, we have
\[
 \frac{\mathrm d B}{\mathrm ds}\omega\frac{\mathrm d \omega}{\mathrm ds}=1-\bigl(\frac{\mathrm d}{\mathrm d s}\omega\bigr)^2-\omega\frac{\mathrm d^2 \omega}{\mathrm ds^2}=1+O(s^{-1}),
 \qquad 1-\frac{\mathrm d B}{\mathrm ds}=\int_s^\infty \frac{\mathrm d^2 B}{\mathrm ds^2}(v)\,\mathrm{d}v=O(s^{-1}).
\]
It follows by integration and substitution in the curvature formulas
that
\[
 \omega(s)^2=2s+O(\log s),\qquad
 \mathrm{R}_B(s)=s^{-1}+O(s^{-2}\log s),\qquad
 \frac{\mathrm d B}{\mathrm ds}(s)=1-\frac1{2s}+O(s^{-2}\log s).
\]
In particular, $B(s)\asymp s$ outside a fixed ball. For a value $v\ge0$, define
$\mathrm{R}_{B,\mathrm{pot}}(v):=\mathrm{R}_B(B^{-1}(v))$.
Here $B^{-1}$ is the inverse of the strictly increasing radial potential
for $s>0$, extended continuously at zero. For $v>0$ this gives
\[
 \frac{\mathrm{d}\mathrm{R}_{B,\mathrm{pot}}}{\mathrm{d}v}=-4k(B^{-1}(v)),
 \qquad \mathrm{R}_{B,\mathrm{pot}}(v)\asymp(1+v)^{-1},\qquad
 \left|\frac{\mathrm{d}}{\mathrm{d}v}\log\mathrm{R}_{B,\mathrm{pot}}\right|\le\frac{C}{1+v}.
\]
On $B^3\times\mathbb R$, any limiting normalized potential satisfying
$\Hess F=\lambda\Ric$ differs from $\lambda B(s)$ by a function with
parallel gradient. The Bryant factor has no nonzero parallel vector
field, since its Ricci curvature is positive. Hence
$F=\lambda B(s)+a_{\rm ax}z+b_0$. The second identity in
\eqref{eq:scaledpotentialidentities}, evaluated at the tip, gives
$\lambda^2+a_{\rm ax}^2=1$. This proves \eqref{eq:introtilt} and explains
precisely which family of tilted potentials the local estimates must
cover.

\subsection{Uniform collar estimates}

In this subsection, we construct collars on which the two transverse
components of approximate rotations can be controlled uniformly.
Proposition~\ref{prop:alltiltcollar} constructs regular potential
collars uniformly down to a vanishing axial coefficient, provided the
collar stays away from the Bryant tip.
The normal integral identity of Proposition~\ref{prop:bracket-normal-integral}
and the Jacobi inverse give the estimates in
Propositions~\ref{prop:potentialbracket} and~\ref{prop:normal-component}.
Proposition~\ref{prop:tipnormal} later compares the ambient tube
coordinate with the model potential $B(s)$.

\paragraph{Collar notation.}
For the collar results in this subsection we use the following
\textbf{leaf notation}. Here $F$ is a smooth function on the ambient
domain, with $|\nabla F|\ge q_0>0$ on the collar, and with the metric
and foliation bounds specified in each statement. In applications it
is the normalized soliton potential, possibly with the time-dependent
additive constant of the heat gauge; the abstract collar estimates
only require the stated smoothness and regularity assumptions.
A regular level is $\Sigma_c:=\{F=c\}$;
the induced metric is $h_c:=g|_{T\Sigma_c}$, and a CMC sphere is
$S_\xi\subset\Sigma_c$, where $\xi$ is parametrized by $(c,\sigma)$.
Its induced area measure is $\mathrm d\mu_\xi$, its unit normal in
$\Sigma_c$ is $\nu_\xi$, and its second fundamental form is
$\mathrm{II}_\xi$. We set
\begin{align}\label{equ:defJacOpe}
J_\xi:=\Delta_{S_\xi}+|\mathrm{II}_\xi|^2
                  +\Ric_{\Sigma_c}(\nu_\xi,\nu_\xi),\qquad
 \Pi_\xi v:=\frac{1}{\Area(S_\xi)}\int_{S_\xi}v\,\mathrm d\mu_\xi,
 \qquad P_\xi:=I-\Pi_\xi.	
\end{align}
When the leaf parameter is fixed we suppress $\xi$. A bounded cap
$\Omega\subset\Sigma_c$ carries volume measure $\mathrm d\mu_{\Sigma_c}$.
A \textbf{collar chart} is a fixed parametrization
$\Theta:S^2\times I_c^+\times I_\sigma^+\to\mathcal U$, where $I_c^+,I_\sigma^+$
are bounded intervals. For closed interior intervals
$I_c\Subset I_c^+$ and $I_\sigma\Subset I_\sigma^+$ we write
\begin{equation}
 \mathcal C^+:=\Theta(S^2\times I_c^+\times I_\sigma^+),\qquad
 \mathcal C:=\Theta(S^2\times I_c\times I_\sigma).
\end{equation}
Thus an interior collar includes fixed positive margins in both leaf
parameters. Unless a different domain is displayed, collar norms use
this convention, with
\[
 \|v\|_{C^k(\mathcal C)}
 :=\sum_{|\alpha|+|\beta|\le k}
    \sup_{S^2\times I_c\times I_\sigma}
    |\nabla_\theta^\alpha\partial_{(c,\sigma)}^\beta(v\circ\Theta)|.
\]
The tensor norm is defined in the same fixed charts. For a fixed
sphere $S\subset(\Sigma_c,h_c)$, we write $H(S,h_c)$ for its mean
curvature and use the \textbf{metric variation of mean curvature}
\[
 (D_hH)_S(k):=\left.\frac{\mathrm d}{\mathrm dt}\right|_{t=0}
                   H(S,h_c+tk).
\]
The notation $DH$ in a leaf calculation abbreviates this operator
at its current induced level metric; it is distinct from the Jacobi
operator governing normal displacement of $S$. 

\begin{prop}
\label{prop:bracket-normal-integral}
Let $(L^3,h)$ be an oriented smooth Riemannian three-manifold, with
volume form $\mu$, and let $S\subset L$
be a closed oriented embedded surface. For smooth vector fields $V,W$
defined near $S$, let $\iota_V$ denote contraction of a differential
form with $V$, let $\nu$ be the unit normal determined by the
orientations, and let $\mathrm d\mu_S$ be its induced area form. Then we obtain
\[
 \int_S\iota_{[V,W]}\mu
 =\int_S \bigl((\operatorname{div}_h W)\iota_V\mu
              -(\operatorname{div}_h V)\iota_W\mu\bigr),
\]
or equivalently,
\[
 \int_S\langle[V,W],\nu\rangle_h\,\mathrm d\mu_S
 =\int_S\bigl((\operatorname{div}_hW)\langle V,\nu\rangle_h
              -(\operatorname{div}_hV)\langle W,\nu\rangle_h\bigr)
       \,\mathrm d\mu_S.
\]
Consequently, if $|V|+|W|\le B$ and
$|\Lie_Vh|+|\Lie_Wh|\le e$ near $S$, then we obtain
\[
 \left|\frac1{\Area(S)}
          \int_S\langle[V,W],\nu\rangle\,\mathrm{d}\mu_S\right|
 \le C B e.
\]
\end{prop}
\begin{proof}
Cartan's formula and $\mathrm{d}\iota_W\mu=(\operatorname{div}_hW)\mu$ give
\[
 \iota_{[V,W]}\mu
 =\mathrm{d}(\iota_V\iota_W\mu)
   +(\operatorname{div}_hW)\iota_V\mu
   -(\operatorname{div}_hV)\iota_W\mu.
\]
The integral of the exact term over the closed surface is zero.
To conclude the inequality, we use
$\operatorname{div}_hV=\tfrac12\operatorname{tr}_h\Lie_Vh$
and the corresponding identity for $W$.
The assertion follows.
\end{proof}

We next estimate the mean-zero Jacobi operator on the radial spheres
of a tilted potential level. The estimate remains uniform as the
axial coefficient tends to zero, even though the graph metric itself
loses uniform control.

We write the tip-normalized Bryant metric and potential as
\[
 g_B=\mathrm{d}s^2+\omega(s)^2g_{S^2},\qquad \Hess_{g_B}B=\Ric_{g_B},
\]
and let $\lambda^2+a_{\rm ax}^2=1$ with $a_{\rm ax}>0$. The level of
$F:=\lambda B(s)+a_{\rm ax}z$ carries the graph metric
\[
 g_{\lambda}=E(s)^2\mathrm{d}s^2+\omega(s)^2g_{S^2},\qquad \mathrm{where}\qquad
 E(s)=\sqrt{1+(\lambda/a_{\rm ax})^2\frac{\mathrm d B}{\mathrm ds}(s)^2}.
\]
The graph metric degenerates as $a_{\rm ax}\downarrow0$. The spectral estimate
on each fixed radial sphere nevertheless has a constant independent of $a_{\rm ax}$.

\begin{prop}\label{prop:jacobi}
Let $(B^3,g_B,B)$ be the Bryant soliton with $\mathrm R_B(o)=1$
and $B(0)=0$. For $0<a_{\rm ax}\le1$, set $\lambda:=\sqrt{1-a_{\rm ax}^2}$ and let
$g_\lambda$ be the metric induced on a level of
$\lambda B(s)+a_{\rm ax}z$ in $(B^3\times\mathbb R,g_B+\mathrm dz^2)$,
as displayed above. On $S_s=\{s=\mathrm{constant}\}$, define
$J:=\Delta_{S_s}+|\mathrm{II}|^2+\Ric_{g_\lambda}(\nu,\nu)$,
where $\mathrm{II}$ and $\nu$ are its second fundamental form and
unit normal in that level. All integrals below use its induced area measure
and $u$ is a smooth real function on $S_s$.
There are $s_0,C>0$, independent of $0<a_{\rm ax}\le1$ and
$\lambda=\sqrt{1-a_{\rm ax}^2}$, such that the mean-zero Jacobi operator $J$
of the radial sphere at $s\ge s_0$ in $g_{\lambda}$ is negative definite,
with
\[
 -\int_{S_s}uJu\,\mathrm{d}\mu_{S_s}\ge \frac{1}{Cs}\int_{S_s}u^2\,\mathrm{d}\mu_{S_s}, \qquad\mathrm{where}\ u \ \mathrm{satisfies}
 \qquad \int_{S_s}u\,\mathrm{d}\mu_{S_s}=0.
\]
\end{prop}

\begin{proof}
We set $\mathrm{d}\sigma:=E\,\mathrm{d}s$ and $w(\sigma):=\omega(s)$. The scalar Jacobi
operator is
\[
 J=\Delta_{S^2(\omega)}+V,\qquad
 V=\frac{2\bigl(\frac{\mathrm d}{\mathrm d s}\omega\bigr)^2}{E^2\omega^2}
   -\frac{2\frac{\mathrm d^2 \omega}{\mathrm ds^2}}{E^2\omega}
   +\frac{2\frac{\mathrm d \omega}{\mathrm ds}\frac{\mathrm d E}{\mathrm ds}}{E^3\omega}.
\]
All derivatives in this expression are with respect to $s$.
For large $s$, the standard Bryant asymptotics give
\[
 \omega\asymp\sqrt{s},\quad |\frac{\mathrm d \omega}{\mathrm ds}|\le Cs^{-1/2},\quad
 |\frac{\mathrm d^2 \omega}{\mathrm ds^2}|\le Cs^{-3/2},\quad \frac{\mathrm d B}{\mathrm ds}\ge c_0>0,\quad
 0<\frac{\mathrm d^2 B}{\mathrm ds^2}\le Cs^{-2}.
\]
Writing $q=\lambda/a_{\rm ax}$, one has the parameter-independent bound
\[
 0\le\frac{\frac{\mathrm d E}{\mathrm ds}}E
 =\frac{q^2\frac{\mathrm d B}{\mathrm ds}\frac{\mathrm d^2 B}{\mathrm ds^2}}{1+q^2\bigl(\frac{\mathrm d}{\mathrm d s}B\bigr)^2}\le\frac{\frac{\mathrm d^2 B}{\mathrm ds^2}}{\frac{\mathrm d B}{\mathrm ds}}\le Cs^{-2}.
\]
Since $E\ge1$, the three terms in $V$ are bounded in absolute value
by $Cs^{-2}$, $Cs^{-2}$, and $Cs^{-3}$, respectively, uniformly in
$a_{\rm ax}$. The first nonzero eigenvalue of $-\Delta_{S^2(\omega)}$ is
$2/\omega^2\ge c_1/s$. For sufficiently large $s_0$, the zeroth-order
term is at most $c_1/(2s)$, proving the assertion.
This finishes the proof.
\end{proof}

We next make the collar estimate uniform over the allowed potential
tilts. This uniformity is needed when different caps have different
orientations of the potential direction.

\begin{prop}
\label{prop:alltiltcollar}
Let $(B^3,g_B,B)$ be the Bryant soliton and put
$g_0:=g_B+\mathrm dz^2$. Fix $s_+>s_-\ge s_0$, $Z_*>0$,
and $0<\delta<\min\{s_-,(s_+-s_-)/2,Z_*\}$, where $s_0$
is as in Proposition~\ref{prop:jacobi}. Define the enlarged annulus
$\mathcal A^+:=\{s_--\delta<s<s_++\delta,\ |z|<Z_*+\delta\}$
and the interior
$\mathcal A^-=\{s_-+\delta<s<s_+-\delta,\ |z|<Z_*-\delta\}$.
Fix an integer $k\ge0$ and $0<\alpha<1$. There are
$\eta_0,C>0$, depending on these fixed data but independent of $a_{\rm ax}$,
with the following property.

 Let $0\le a_{\rm ax}\le1$,
$\lambda:=\sqrt{1-a_{\rm ax}^2}$, and $F_0:=\lambda B(s)+a_{\rm ax}z$.
Suppose $g$ is a smooth Riemannian metric and $F$ a smooth function
on $\mathcal A^+$ such that
\begin{align}\label{equ:comparison}
\|g-g_0\|_{C^{k+8,\alpha}(\mathcal A^+)}
 +\|F-F_0\|_{C^{k+8,\alpha}(\mathcal A^+)}\le\eta_0.	
\end{align}
\begin{enumerate}[label=\textup{(\roman*)},ref=\textup{(\roman*)}]
\item The interior $\mathcal A^-$ admits a controlled
local two-parameter foliation by CMC two-spheres in actual levels of
$F$. The lower bound $|\nabla F|_g\ge q_*$ displayed below, the local coordinate bounds, the
mean-zero Jacobi inverse, and the allowed perturbation size can all be
chosen uniformly for $a_{\rm ax}\in[0,1]$ on this fixed annulus.
Here controlled means that the normal graphs and their derivatives
in both leaf parameters have bounds in $C^{k+2,\alpha}$ depending
only on the fixed data. The graphs are $O(\eta)$ in these norms
when the displayed comparison norm \eqref{equ:comparison} is at most $\eta\le\eta_0$.
After decreasing $\eta_0$, we may take
\[
 q_*:=\tfrac12\min\{1,\inf_{[s_--\delta,s_++\delta]}\frac{\mathrm d B}{\mathrm ds}(s)\}>0,
 \qquad |\nabla F|_g\ge q_*\quad\hbox{on }\mathcal A^-.
\]
This constant depends on the fixed annulus and the Bryant normalization,
but not on $a_{\rm ax}$. The
inverse of $PJP$ on mean-zero functions is bounded from
$C^{k,\alpha}$ to $C^{k+2,\alpha}$ by $C$, where $P,\,J$ are defined in \eqref{equ:defJacOpe}.
\item At the same potential value, the CMC families agree on overlaps
after matching their mean-height parameters in a common graph chart.
Precisely, fix one reference model sphere $S_{\rm ref}$, its induced
model metric $h_{\rm ref}$, and the identification of its neighborhood
within the potential level used for the graph construction.
Use the same unit normal and graph chart for both spheres. If
$u:S_{\rm ref}\to\mathbb R$ is the signed normal graph height in this
chart, its reference-area mean height is
\[
 \overline u_{\rm ref}:=
 \frac{1}{\operatorname{Area}_{h_{\rm ref}}(S_{\rm ref})}
       \int_{S_{\rm ref}}u\,\mathrm d\mu_{h_{\rm ref}}.
\]
The measure in this definition is fixed on the reference sphere. Two sufficiently
small CMC graphs in the same potential level coincide when their
values of $\overline u_{\rm ref}$ agree.
\end{enumerate}
\end{prop}
\begin{proof}
We first establish the uniform estimates and construct the foliation in (i).
We set $q:=(a_{\rm ax}^2+\lambda^2 \frac{\mathrm d B}{\mathrm ds}(s)^2)^{1/2}$. On the fixed annulus, we have
\[
 q\ge\min\{1,\inf \frac{\mathrm d B}{\mathrm ds}\}>0.
\]
The model unit normal to the potential level and the unit normal to a
sphere within that level can therefore be written without division by
$a_{\rm ax}$:
\[
 N_0=q^{-1}(\lambda \frac{\mathrm d B}{\mathrm ds}\partial_s+a_{\rm ax}\partial_z),\qquad
 \nu_0=q^{-1}(a_{\rm ax}\partial_s-\lambda \frac{\mathrm d B}{\mathrm ds}\partial_z).
\]
Their coefficients and all required derivatives are uniformly bounded.
Flow boxes for these fields give local coordinates consisting of the
potential value, a transverse parameter within its level, and the
sphere variable. Continuous dependence on the compact parameter
quarter-circle $a_{\rm ax}^2+\lambda^2=1$, $a_{\rm ax},\lambda\ge0$, and uniform spatial
derivative bounds on the fixed annulus give uniform chart sizes after shrinking
margins. One can equivalently use the flow of
$\nabla F_0/|\nabla F_0|^2$ for the first coordinate and $\nu_0$ for
the second.

For $a_{\rm ax}>0$, Proposition~\ref{prop:jacobi} gives a uniform negative
gap for the mean-zero Jacobi operator. The coefficients of this operator
in the displayed nonsingular frame extend smoothly to $a_{\rm ax}=0$.
At $a_{\rm ax}=0$ the level is $S^2(\omega(s))\times\R$, its sphere is totally
geodesic, and its Jacobi operator is the round Laplacian. The gap and
uniform Schauder bounds thus hold on the full parameter interval.

We use the regular-level implicit function theorem to identify nearby
actual $F$ levels with the model levels. In each such level, solve the
projected CMC equation for a mean-zero normal graph over a model
sphere; its linearization is the preceding invertible Jacobi operator
on mean-zero functions. We now establish uniform estimates for the parameter derivatives. In the
flow boxes above, let $\xi:=(c,\sigma)$ be the two leaf parameters,
let $P v:=v-\overline v_{\rm ref}$, with the fixed reference-area
mean defined in (ii), and write
\[
 \mathcal G(u;a_{\rm ax},\xi,g,F)=P H(u;a_{\rm ax},\xi,g,F).
\]
On mean-zero $C^{k+2,\alpha}$ graphs its derivative at the model
zero graph is $L_{a_{\rm ax},\xi}=-PJP$, with
\[
 \sup_{a_{\rm ax},\xi}\|L_{a_{\rm ax},\xi}^{-1}\|_{C^{k,\alpha}\to C^{k+2,\alpha}}
       \le C,\qquad
 \|\mathcal G(0)\|_{C^{k,\alpha}}\le C\eta.
\]
Uniform smooth coefficient bounds on the compact parameter family
give, for small $u,v$,
\[
 \|\mathcal N(u)-\mathcal N(v)\|_{C^{k,\alpha}}
 \le C(\eta+\|u\|_{C^{k+2,\alpha}}+\|v\|_{C^{k+2,\alpha}})
                      \|u-v\|_{C^{k+2,\alpha}},
\]
where $\mathcal N(u)=\mathcal G(u)-\mathcal G(0)-L_{a_{\rm ax},\xi}u$.
Thus $u=-L_{a_{\rm ax},\xi}^{-1}(\mathcal G(0)+\mathcal N(u))$ is a
contraction on a ball of radius $C_1\eta$, with constants uniform
in $a_{\rm ax}$. Differentiating $\mathcal G(u;a_{\rm ax},\xi,g,F)=0$ gives
\[
 D_u\mathcal G\,\partial_{\xi_j}u
       =-\partial_{\xi_j}\mathcal G.
\]
The same inverse estimate applies. Repetition gives every required
parameter derivative within the stated regularity reserve. The lapse
is $1+O(\eta)$ in the transverse flow coordinate and is positive, so
these spheres foliate the interior collar. This proves (i).

 To prove (ii), note that prescribing the mean height determines the graph
uniquely near a model leaf. For a nonzero prescribed mean, subtract
that constant from the graph function and apply the same mean-zero
implicit function theorem. Thus equality of the potential and mean-height
labels gives equality of the leaves, proving (ii).
All these arguments take place on an annulus where $q$ is bounded below.
They make no assertion about regular potential coordinates at the tip.
This proves the stated conclusion.
\end{proof}

The metric Lie derivative estimate alone does not control motion across potential
levels. We therefore use the bracket relations to improve the
fields' potential tangency.

\begin{prop}
\label{prop:potentialbracket}
Let $(\mathcal U^4,g)$ be a smooth Riemannian domain and let
$F:\mathcal U\to\mathbb R$ be a smooth function. Suppose a collar
$\mathcal C\Subset\mathcal U$ has a smooth two-parameter foliation
$S_\xi$, $\xi=(c,\sigma)$, by closed CMC two-spheres in the
regular levels $\{F=c\}$. Assume $|\nabla F|\ge q_0>0$ and
a uniform inverse bound for $P_\xi J_\xi P_\xi$ on mean-zero
functions. For an integer $k\ge0$ and $\alpha\in(0,1)$, assume
controlled $C^{k+8,\alpha}$ background and foliation geometry
on a fixed enlargement of $\mathcal C$, including bounds for the
metric inverse and leaf parametrizations. Suppose $Y,Z$ are smooth
ambient vector fields with $C^{k+5,\alpha}$ norms bounded by a fixed $B_0$ there and
\[
 \|\Lie_Yg\|_{C^{k+4,\alpha}}
 +\|\Lie_Zg\|_{C^{k+4,\alpha}}
 +\|YF\|_{C^{k+4,\alpha}}+\|ZF\|_{C^{k+4,\alpha}}\le e.
\]
Set ${W}:=[Y,Z]$ and let
$\Theta(\theta,c,\sigma)$ be fixed foliation charts parametrizing
$\mathcal C$. The normals are $n_f=\nabla F/|\nabla F|$ and
the actual CMC normal $\nu$ within the potential level.
Then the metric Lie derivative and the two normal components satisfy
\[
 \|\Lie_{W}g\|_{C^k(\mathcal C)}
 +\|\langle {W},n_f\rangle\|_{C^k(\mathcal C)} +\|\langle {W},\nu\rangle\|_{C^k(\mathcal C)} \le Ce.
\]
The constant depends only on $k$,
$\alpha$, $q_0$, $B_0$, the Jacobi inverse bound, the background
and foliation bounds, and the fixed margins. For the collars of
Proposition~\ref{prop:alltiltcollar}, the constant is uniform in $a_{\rm ax}$.
\end{prop}
\begin{proof}
The identities
\[
 W F=Y(ZF)-Z(YF),\qquad
 \Lie_{W}g=\Lie_Y\Lie_Zg-
                  \Lie_Z\Lie_Yg
\]
give the potential pairing and metric Lie derivative bounds, proving
\begin{align*}
	\|\Lie_{W}g\|_{C^k(\mathcal C)}
 +\|\langle {W},n_f\rangle\|_{C^k(\mathcal C)} \le Ce.
\end{align*}

Next, for
$T\in\{Y,Z,{W}\}$ define its projection onto the $F$ levels by
$\bar T:=T-(TF)\nabla F/|\nabla F|^2$. This changes each field
by $O(e)$ with the required derivatives. In the variation equation
below, $DH$ is the metric derivative of mean curvature on the fixed leaf. Their
projected metric Lie derivatives are $O(e)$, and the projection of ${W}$ differs from
the bracket of the projected $Y,Z$ by $O(e)$.
Write $W_\Sigma:=[\bar Y,\bar Z]$. On each regular
level $\Sigma=\{F=c\}$, both $\bar Y$ and $\bar Z$ are tangent,
and their induced metric Lie derivatives are $O(e)$ on the enlarged
collar. Moreover,
\[
 \bar{W}-W_\Sigma=O_{C^{k+2,\alpha}}(e).
\]
Apply Proposition~\ref{prop:bracket-normal-integral} to
$\bar Y,\bar Z$ on the closed sphere in $(\Sigma,g_\Sigma)$.
It gives
\[
 \left|\Pi_\xi\langle W_\Sigma,\nu\rangle\right|\le Ce,
 \qquad
 \left|\Pi_\xi\langle W,\nu\rangle\right|\le Ce.
\]
The second bound uses the displayed difference and
$\langle\bar{W},\nu\rangle=\langle W,\nu\rangle$.
Set $w:=\langle W,\nu\rangle$. CMC variation gives
$Jw=-DH(\Lie_{\bar {W}}g_{\{F=c\}})$, and the mean-zero Jacobi
inverse bounds its remaining part. We give the details of the mean-zero inversion, including its dependence
on the leaves. Pull back to the reference sphere and write
$m(\xi)=\Pi_\xi w$, $v=w-m(\xi)$, with $P_\xi=I-\Pi_\xi$ the
actual-area mean-zero projection. We set
$q_\xi:=|\mathrm{II}|^2+\Ric_{\{F=c\}}(\nu,\nu)$ and
$H_\xi=-DH(\Lie_{\bar {W}}g_{\{F=c\}})$. Since $J1=q_\xi$, we have
\[
 P_\xi Jv=P_\xi H_\xi-m(\xi)P_\xi q_\xi,
 \qquad \Pi_\xi v=0.
\]
The normal integral formula gives $|m|\le Ce$, so the uniform Jacobi inverse
and Schauder estimate give $\|v\|_{C^{k+2,\alpha}}\le Ce$
whenever the right side is controlled in $C^{k,\alpha}$.
For parameter derivatives, differentiation of an average is explicit:
\[
 \partial_{\xi_j}\Pi_\xi u
   =\Pi_\xi(\partial_{\xi_j}u)
    +\Pi_\xi\bigl((u-\Pi_\xi u)\,
                      \partial_{\xi_j}\log \mathrm{d}\mu_\xi\bigr).
\]
We apply this also to the exact normal integral identity to obtain
$|\partial_\xi^I m|\le C_Ie$. Differentiating the projected
equation puts derivatives of $J$ and $P_\xi$ only against lower
parameter derivatives of $w$. After transferring to fixed mean-zero
spaces, the equation at order $I$ has the same invertible principal
operator and a right side bounded by
\[
 C_Ie+C_I\sum_{|J|<|I|}\|\partial_\xi^Jw\|_{C^{k-|I|+2,\alpha}}.
\]
Here the estimate at parameter order $|I|\le k$ is taken in
$C^{k-|I|+2,\alpha}$ in the angular variables; thus derivatives are
counted by total order. The previously established lower parameter
orders have at least this angular regularity. Induction proves the
full collar bounds within the stated derivative reserve.
The projections cost at most one derivative in their metric Lie derivatives, the
bracket costs one, and $DH$ costs one more; the displayed derivative
margins suffice for the $C^k$ conclusion. No bounding domain is used.
This finishes the proof.
\end{proof}

We now estimate the component normal to a CMC sphere within its
potential level. We first project the field onto that level and then
apply the normal integral identity and the mean-zero Jacobi estimate.
The second transverse pairing on a cap can be estimated without a second
commutator, provided the improved data hold throughout the bounded cap
inside the leaf. For an individual field, collar data alone do not
control the integral of the normal component over its inner boundary. The bracket cancellation
in Proposition~\ref{prop:bracket-normal-integral} provides a different alternative
when an additional bracket is available.

\begin{prop}\label{prop:normal-component}
Let $F$ be a smooth function on a Riemannian manifold and suppose
$|\nabla F|\ge q_0>0$ on a neighborhood of a compact domain
$\Omega\subset\Sigma=\{F=c\}$. Suppose $S=\partial\Omega$ is a smooth
closed CMC hypersurface of $(\Sigma,g_\Sigma)$, with outward unit normal
$\nu$. We fix $m\ge0$ and $0<\alpha<1$. Assume controlled smooth geometry
on this neighborhood, a bound for $\operatorname{Vol}(\Omega)/
\Area(S)$, and an inverse estimate
\begin{equation}\label{eq:projectedjacobi}
 \|v\|_{C^{m+2,\alpha}(S)}
 \le C_J\|PJPv\|_{C^{m,\alpha}(S)},\qquad\mathrm{where}\ v\ \mathrm{satisfies} \int_S v\,\mathrm d\mu_S=0.
\end{equation}
If a smooth ambient vector field $Y$ satisfies
\[
 \|\Lie_Yg\|_{C^{m+1,\alpha}}
 +\|YF\|_{C^{m+2,\alpha}}\le e
\]
on this neighborhood, then we obtain
\begin{equation}\label{eq:normal-component-conclusion}
 \|\langle Y,\nu\rangle\|_{C^{m+2,\alpha}(S)}\le C e.
\end{equation}
The constant depends only on the stated geometry, $q_0$, the bounds for
$F$, $\operatorname{Vol}(\Omega)/
\Area(S)$ and $C_J$.
\end{prop}
\begin{proof}
We set $u:=YF$ and project the field exactly onto the potential level:
\[
 V=Y-u\frac{\nabla F}{|\nabla F|^2}.
\]
The projection is defined on the whole neighborhood, $VF=0$, and
\[
 \|\Lie_Vg_\Sigma\|_{C^{m+1,\alpha}}\le C e.
\]
Indeed, writing $Z=\nabla F/|\nabla F|^2$,
$\Lie_{uZ}g=u\Lie_Zg+\mathrm{d}u\otimes Z^\flat+
Z^\flat\otimes \mathrm{d}u$, with all coefficients bounded by the hypotheses.
Let $w:=\langle V,\nu\rangle=\langle Y,\nu\rangle$.
The divergence theorem inside $\Sigma$ gives
\begin{equation}\label{eq:normal-component-flux}
 \left|\int_S w\,\mathrm d\mu_S\right|
 =\left|\int_\Omega\operatorname{div}_\Sigma V\,\mathrm d\mu_\Sigma\right|
 =\frac12\left|\int_\Omega
        \operatorname{tr}_{g_\Sigma}\Lie_Vg_\Sigma\,\mathrm d\mu_\Sigma\right|
 \le Ce\operatorname{Vol}(\Omega).
\end{equation}
Thus its average $\bar w$ is bounded by $Ce$.

For completeness, let $\psi_s$ be the flow of $V$ on $\Sigma$.
Naturality of mean curvature identifies its metric variation on the
fixed hypersurface $S$, in direction $\Lie_Vg_\Sigma$, with the
variation under the deformation $\psi_s(S)$. Since $H$ is constant on
$S$, the tangential derivative of $H$ vanishes and the hypersurface
variation formula gives
\[
 (D_hH)_S(\Lie_Vg_\Sigma)=-Jw.
\]
The metric variation operator is linear of order one in its argument;
its coefficients are controlled by the hypersurface geometry. Therefore
$\|Jw\|_{C^{m,\alpha}}\le Ce$.
We set $w_0:=w-\bar w$. Since $Pw_0=w_0$, we have
\begin{equation}\label{eq:normal-component-projected}
 PJPw_0=P(Jw)-\bar w\,P(J1).
\end{equation}
Both terms on the right have $C^{m,\alpha}$ norm at most $Ce$.
We apply \eqref{eq:projectedjacobi} and then restore the already controlled
average. This proves \eqref{eq:normal-component-conclusion}.
The assertion follows.
\end{proof}

\begin{lem}
\label{lem:rotation-gram}
Consider a controlled family of near-round two-spheres $(S,\gamma)$,
with area bounded above and below and uniform geometry through eight
derivatives. Let $T_1,T_2,T_3$ be uniformly bounded tangent fields,
close to a nondegenerate oriented round rotational triple, such that
\[
 \|\Lie_{T_a}\gamma\|_{C^6}\le e,\qquad
 [T_a,T_b]=C_{ab}^{\ c}T_c+O_{C^2}(e).
\]
Suppose the coefficients $C_{ab}^{\ c}$ are the same on every
sphere and are close to $c_{\rm rot}\epsilon_{abc}$ for one fixed
$c_{\rm rot}>0$. Define
\[
 M:=\operatorname{sym}\left(
       \frac12\sum_{a,b}\epsilon_{dab}C_{ab}^{\ c}\right)_{dc},
 \qquad Q_*:=\frac{\det M}{6\pi}M^{-1}.
\]
For sufficiently small comparison error, $M,Q_*$ are positive
definite with uniformly controlled inverses, and the normalized Gram matrix satisfies
$Q_S(T)=Q_*+O(e)$ on every sphere. The conclusion also holds for
ambient vector fields differing from $T_a$ along the spheres by $O_{C^2}(e)$.
The error constant and the comparison threshold depend only on the
area bounds, the uniform geometry through eight derivatives, the
field bounds, and the upper and positive lower bounds for the
singular values of the comparison rotational triple and for
$c_{\rm rot}$. They are independent of the sphere in this family.
\end{lem}
\begin{proof}
Naturality gives $T_aK=(DK)_\gamma(\Lie_{T_a}\gamma)=O_{C^4}(e)$.
The fields span each tangent plane with controlled right inverse,
so $\mathrm dK=O_{C^3}(e)$. Integrating on the sphere and using
Gauss--Bonnet yields $K-4\pi/\operatorname{Area}(S)=O_{C^3}(e)$.
After scaling to area $4\pi$, apply \cite[Theorem~2.2]{Le} with
$n=3,p=4$. Its $W^{5,4}$ conformal-factor estimate gives a round
metric $\gamma_0$, of the original area, with
$\gamma_0-\gamma=O_{C^4}(e)$. This choice is made separately on
each sphere; no parameter derivative of the uniformization is used.

For the Killing operator $\mathcal K Z=\Lie_Z\gamma_0$, the symbol
$\xi\otimes Z^\flat+Z^\flat\otimes\xi$ is injective for $\xi\ne0$.
Thus $\mathcal K^*\mathcal K$ is elliptic and has a positive spectral
gap on the orthogonal complement of the round Killing fields.
Its elliptic estimate and the $C^{1,\alpha}$ bound for
$\mathcal K T_a$ give, for fixed $0<\alpha<1$,
\[
 T_a=P_{ai}J_i+O_{C^2}(e),\qquad G:=PP^{\mathsf T},
\]
where $J_i$ have cyclic bracket coefficient one and $P,P^{-1}$
are bounded. Choose the orientation so that $\det P>0$.
Whole-sphere $L^2$ nondegeneracy, rather than pointwise rank,
then compares coefficients in the bracket identity:
\[
 C_{ab}^{\ c}=\det(P)\epsilon_{abd}(G^{-1})_{dc}+O(e),
 \qquad M=\det(P)G^{-1}+O(e).
\]
Since $\det G=(\det P)^2$, it follows that
$G=(\det M)M^{-1}+O(e)$. Finally, on a round sphere of any radius,
\[
 \operatorname{Area}(S)^{-2}
 \int_S\langle J_i,J_j\rangle\,\mathrm d\mu
       =(6\pi)^{-1}\delta_{ij}.
\]
Indeed the unit-sphere integral is $(8\pi/3)\delta_{ij}$ and both
it and squared area scale by the fourth power of radius.
Restoring $\gamma$ and the original fields costs $O(e)$ and proves
the assertion. All constants are uniform on the controlled family.
\end{proof}

The collar estimates will be used in the CMC reconstruction in the
proof of Theorem~\ref{cylinder:improvement} and in the cap reconstruction
of Corollary~\ref{cor:alltiltheat}. In both places,
Propositions~\ref{prop:potentialbracket} and~\ref{prop:normal-component}
control the two transverse components of the corrected rotations,
while Lemma~\ref{lem:rotation-gram} provides a common Gram normalization.

\section{Neck improvement on cylindrical regions}\label{sec:cylinder}
In this section, we prove neck improvement on evolving
$S^2\times\mathbb R^2$ neighborhoods. Theorem~\ref{cylinder:improvement}
reduces the symmetry error $\epsilon$ to $\epsilon/2$ on a smaller
neighborhood, as in Brendle's neck improvement theorem
\cite[Theorem~8.6]{BrendleAncient}. We first define local symmetry and
control the normalized potential; Lemmas~\ref{cylinder:lem:alignment}
and~\ref{cylinder:lem:gluing} then provide compatible fields for the
heat argument. The proof of the improvement theorem is divided into
nine steps. Cap improvement follows in Section~\ref{sec:capimprovement}.

The analytic quantities are the metric Lie derivative
$h:=\Lie_Vg$ and potential pairing $u:=VF$. In the heat gauge of
Lemma~\ref{cylinder:lem:heat-gauge}, the corrected fields satisfy
\[
 \partial_tV=\Delta V+\Ric^\sharp V,\qquad
 (\partial_t-\Delta_L)h=0,\qquad
 (\partial_t-\Delta)u=-\lambda\langle\Ric,h\rangle.
\]
The last identity is Lemma~\ref{cylinder:lem:pairing}; the operator
$\Delta_L$ always has the convention \eqref{eq:lichnerowicz-convention}.
Local alignment, the cylinder spectrum, and two cyclic brackets then
produce time-independent constructed fields. Here the second flat direction introduces a potential-normal
component in addition to the CMC normal component. We estimate it
by a scalar heat equation coupled to the metric Lie derivative. We also keep explicit
polynomial bounds for gluing in two flat variables. For the final
normalization, Brendle's Steps 11--14 compare the Gram matrices by
geometric variation along the neck and in time. We instead recover
a common Gram matrix from constant bracket coefficients, which
controls both transverse parameters at once.


We work in the steady sign convention $\Ric=\Hess f$ and normalize
$\mathrm{R}+|\nabla f|^2=1$ before performing the curvature scale rescaling.
The canonical Ricci flow is generated by $-\nabla f$.
After translating and rescaling the distinguished spacetime point,
we set $(x_0,t_0):=(x_0,-1)$ and $\mathrm{R}(x_0,-1)=1$.
The comparison cylinder is
\begin{equation}\label{cylinder:eq:cylinder}
 \bar g(t):=(-2t)g_{\sphere}+\mathrm{d}z_1^2+\mathrm{d}z_2^2,\qquad t<0,
\end{equation}
where $g_{\sphere}$ has sectional curvature one. We write $z:=(z_1,z_2)$ for the flat coordinates and
$D_a:=\sphere\times\Disk_a$, where $\Disk_a\subset\mathbb R^2$
is the Euclidean disk of radius $a$.

We first specify the evolving comparison domain and the local symmetry
conditions used in the heat correction.

\begin{defn}[Evolving cylindrical neighborhood]\label{cylinder:def:neighborhood}
Let $(M^4,g(s))_{s\in I}$ be a Ricci flow, let $(x_0,t_0)$ be a
spacetime point with $Q:=\mathrm{R}_{g(t_0)}(x_0)>0$, and set
$r:=Q^{-1/2}$. The \textbf{curvature-normalized flow} at this point is
\begin{equation}\label{cylinder:eq:pointed-rescaling}
 \widehat g(t):=Q\,g\bigl(t_0+Q^{-1}(t+1)\bigr),\qquad
 \mathrm{R}_{\widehat g(-1)}(x_0)=1.
\end{equation}
Fix $0<\eta<1/2$, let
$\bar x:=(p_\circ,0)\in S^2\times\mathbb R^2$ be a fixed model point.
We call $(x_0,t_0)$ \textbf{$\eta$-cylindrical} if the normalized flow
is defined on $[-\eta^{-1},0]$ and there is a time-independent
spatial diffeomorphism onto its image
\[
 \Psi:\overline B_{\bar g(-1)}(\bar x,\eta^{-1})\longrightarrow M,
 \qquad \Psi(\bar x)=x_0,
\]
defined on an open neighborhood of this closed model ball, such that
\begin{equation}\label{cylinder:eq:eta-cylindrical}
 \sup_{\substack{y\in\overline B_{\bar g(-1)}(\bar x,\eta^{-1})\\
                 -\eta^{-1}\le t\le-\eta}}
 \sum_{j+2k\le \lfloor \eta^{-1}\rfloor }
 \bigl|\bar\nabla^j\partial_t^k
       (\Psi^*\widehat g(t)-\bar g(t))\bigr|_{\bar g(t)}\le\eta.
\end{equation}
Here $\bar g(t)$ is \eqref{cylinder:eq:cylinder} and $\bar\nabla$ is
its spatial covariant derivative.
We also call the image spacetime region an \textbf{evolving
$\eta$-cylindrical neighborhood centered at $(x_0,t_0)$}.
\end{defn}

The definition is invariant under parabolic rescaling.
For prescribed $R_*>0$, $T_*>1$, and derivative order $m$, if we choose
\begin{equation}\label{cylinder:eq:fixed-domain-threshold}
 \eta\le\min\left\{\frac12,
 \frac{1}{2(R_*+\pi\sqrt2+1)},
 \frac{1}{2(T_*+1)},\frac{1}{m+1}\right\},
\end{equation}
then $D_{R_*}\times[-T_*,-1]$ is contained in the model comparison
region, and
\eqref{cylinder:eq:eta-cylindrical} controls its parabolic derivatives
through order $m$.

If the regular potential levels in a time-$t$ curvature scale
neighborhood admit a smooth foliation by CMC two-spheres, let $\nu$
denote the unit normal to these spheres inside their potential level.
We use this notation only on neighborhoods equipped with such a
foliation; its existence is part of the local symmetry hypothesis below.
The second transverse direction is the unit potential normal $n_f$.
In rescaled coordinates we estimate its pairing through $\nabla F$,
using the exact heat-gauge identity
\begin{equation}
 \langle U,n_f\rangle=\frac{UF}{|\nabla F|},\qquad
 |\nabla F|^2=1-\lambda^2\mathrm{R}.
\end{equation}
Here $\lambda$ is the positive rescaling constant in
Lemma~\ref{cylinder:lem:heat-gauge}. As in
\cite[Definition~8.2]{BrendleAncient}, the following definition is
made on a neighborhood equipped with a CMC foliation.
\begin{defn}[Local symmetry]\label{cylinder:def:symmetry}
Let $\mathrm{R}(x_0,t_0)=r^{-2}$ and let
$\mathcal P:=B_{g(t_0)}(x_0,100r)\times[t_0-100r^2,t_0]$.
Assume that $\nabla f_t\ne0$ on an open neighborhood of
$\overline{\mathcal P}$ and that, at every time, its potential
levels admit a smoothly varying foliation by embedded CMC two-spheres.
In the cylindrical setting we use the small-normal-graph family
described in Remark~\ref{cylinder:rem:cmc-existence}, with whole-leaf
margins in the larger neighborhood. Relative to this foliation, the point $(x_0,t_0)$ is
\textbf{$\epsilon$-symmetric} if there are three smooth time-independent
fields $U_1,U_2,U_3$ on an open neighborhood of the closed spatial
domain such that
\begin{align*}
 \sup_{\mathcal P}\sum_{a=1}^3\sum_{j=0}^2
       r^{2j}|\nabla^j\Lie_{U_a}g(t)|^2&\le\epsilon^2,\\
 \sup_{\mathcal P}r^{-2}\sum_{a=1}^3
       \bigl(|\langle U_a,\nu\rangle|^2
             +|\langle U_a,n_f\rangle|^2\bigr)&\le\epsilon^2,\\
 \sup_S\sum_{a,b=1}^3|Q_{ab}(S)-\delta_{ab}|^2&\le\epsilon^2,
\end{align*}
where all norms and pairings are measured in $g(t)$, and
\[
 Q_{ab}(S):=(Q_S(U))_{ab}=\Area_{g(t)}(S)^{-2}
               \int_S\langle U_a,U_b\rangle_{g(t)}\,\mathrm{d}\mu_{S,g(t)}.
\]
Here $\mathrm d\mu_{S,g(t)}$ is the induced area measure on $S$.
The last supremum is over every whole CMC leaf contained in the
spatial neighborhood at every time; we require this
family to be nonempty at each time.
\end{defn}

\begin{rem}\label{cylinder:rem:cmc-existence}
The CMC foliation in Definition~\ref{cylinder:def:symmetry} is part
of the local geometric data. Its construction on cylindrical regions
is given in Step~6 of the proof of Theorem~\ref{cylinder:improvement}
(see \eqref{cylinder:eq:lapse});
Proposition~\ref{prop:exteriorleaves} constructs compatible families
on the exterior neck neighborhoods and proves their agreement on
overlaps and with the cap collars of
Proposition~\ref{prop:alltiltcollar}.

This additional assumption is natural in the steady soliton setting.
Lemma~\ref{cylinder:lem:heat-gauge} derives potential closeness from
metric closeness and the soliton identities: after normalization and
a rotation of the flat coordinates,
\[
 g\approx(-2t)g_{\sphere}+\mathrm dz_1^2+\mathrm dz_2^2,
 \qquad F\approx z_1.
\]
The regular potential levels are therefore close to
$S^2\times\mathbb R$, whose slices $S^2\times\{z_2\}$ have constant
mean curvature. Within an actual potential level, write a nearby
sphere as a normal graph $u=\sigma+v$, where $\sigma$ is its prescribed
mean height and $v$ has zero average. If $P$ subtracts the
reference-sphere average, the CMC equation is
\[
 P H(\sigma+v)=0,\qquad\mathrm{where}\qquad \int_{S^2}v\,\mathrm d \mu_S=0.
\]
Its linearization in $v$ at a product slice is, up to sign, the
sphere Laplacian on zero-average functions, which is invertible.
The parameter-dependent implicit function theorem gives a smooth
family with positive transverse lapse. Thus the assumption records
the geometric structure used to define $\nu$ and the whole-leaf Gram
normalization; the later construction verifies it on the neighborhoods
used in the proof.
\end{rem}

The local alignment and model estimates now allow us to improve an
approximate rotational triple on a smaller cylindrical neighborhood.
The next theorem gives the quantitative gain needed for iteration.

\begin{thm}[Cylindrical improvement]\label{cylinder:improvement}
Let $(M^4,g,f)$ satisfy the hypotheses of Theorem~\ref{thm:main},
with \eqref{eq:normalization}, and let $g(t):=\Phi_t^*g$.
There are $L<\infty$ and $\epsilon_1>0$ with the following property.
Let $(x_0,t_0)$ be the center of an evolving $\epsilon_1$-cylindrical
neighborhood in the canonical flow of a steady soliton as above,
with $\mathrm{R}(x_0,t_0)=r^{-2}$. Suppose every point of
\(
 B_{g(t_0)}(x_0,Lr)\times[t_0-Lr^2,t_0]
\)
is $\epsilon$-symmetric, where $0\le\epsilon\le\epsilon_1$.
Then $(x_0,t_0)$ is $\epsilon/2$-symmetric.
\end{thm}

We use only the local soliton identities and the displayed smooth
cylinder closeness. The assertion depends only on these local hypotheses.
The constants and domains are determined in the following way.
We choose a fixed large terminal domain
$D_H\times[-\tau_*, -1]$, with, for example,
$H=\tau_*=10^4$. Only fixed smaller subdomains of this domain are
used in the final geometric arguments. Subsequently choose a large
number $T\gg H^2+\tau_*$, and work initially on
$D_{8\sqrt T}\times[-2T,-1]$.
Finally, we choose the cylindrical threshold sufficiently small depending
on $T$ and use $\eta:=\epsilon_1$ as an upper bound for the metric
comparison error. In particular $\epsilon\le\eta$ throughout the proof. We use $C$ for constants independent of $T,\eta,\epsilon$,
$C_q$ for constants depending on a fixed polynomial exponent, and
$C(T)$ for finite constants which may grow with $T$.
Constants depending on the fixed $H,\tau_*$ are included in $C$.
They are chosen before $T$. After choosing $T$, take $L:=1000T$. We then decrease the single
cylindrical threshold $\epsilon_1$ so that
\begin{equation}\label{cylinder:eq:improvement-domain-threshold}
 \epsilon_1\le\min\{1/121,(10^6T)^{-1}\},\qquad
 C(T)\epsilon_1<1/4.
\end{equation}
We impose \eqref{cylinder:eq:fixed-domain-threshold} as well for every
finite spatial radius and backward duration used in the proof.
These are finitely many lower bounds on $\epsilon_1^{-1}$, imposed
after $T$ is fixed. They place the initial heat region, all enlarged
comparison neighborhoods, and the fixed terminal regions inside the
$\epsilon_1$-cylindrical chart. A further fixed enlargement of $L$
only changes these finite threshold requirements. The rest of this section is devoted to proving Theorem~\ref{cylinder:improvement}.

The next lemma derives potential closeness from the metric comparison
and the steady soliton identities.

\begin{lem}\label{cylinder:lem:heat-gauge}
In the setting of
Theorem~\ref{cylinder:improvement}, we fix a finite comparison domain $D_{8\sqrt{T}}\times [-2T,-1]$,
a derivative order $m$, and let $\eta$
be an upper bound for its smooth comparison error as in
\eqref{cylinder:eq:eta-cylindrical}. In the curvature-normalized coordinates there is a positive constant
$\lambda$ and a potential $F(x,t)$ such that
\begin{equation}\label{cylinder:eq:F-identities}
 \Hess F=\lambda\Ric,\qquad
 |\nabla F|^2=1-\lambda^2\mathrm{R},\qquad
 \partial_tF=\Delta F=\lambda \mathrm{R}.
\end{equation}
On any fixed large cylindrical domain, after rotating the flat
coordinates and adding a constant, we have
\begin{equation}\label{cylinder:eq:F-close}
 \lambda\le C\eta,\qquad
 \|F-z_1\|_{C^m}\le C(T)\eta
\end{equation}
for every derivative order needed below, with parabolic derivatives
included when they occur.
\end{lem}
\begin{proof}
Let $Q:=\mathrm{R}_{g(t_0)}(x_0)$ and $r:=Q^{-1/2}$, as in
\eqref{cylinder:eq:pointed-rescaling}, and let
$f_s:=f\circ\Phi_s$ be the original pulled-back potential. On the
normalized flow put
\[
 \lambda:=r^{-1},\qquad
 F(x,t):=r^{-1}\bigl(f_{t_0+r^2(t+1)}(x)-f_{t_0}(x_0)\bigr)+r(t+1).
\]
The same time-independent comparison map $\Psi$ pulls back both
$\widehat g$ and $F$ to the model domain; we suppress it below. Constant metric rescaling leaves the connection
and the covariant Ricci tensor unchanged. Therefore
$\Hess F=r^{-1}\Ric$. The normalized Hamilton identity gives
$|\nabla F|^2=1-r^{-2}\mathrm{R}$. Finally, we get
\[
 \partial_tF=-r|\nabla F|^2+r=r^{-1}\mathrm{R}=\Delta F.
\]
The added function of time changes neither the levels nor the
gradient direction at any fixed time.

To see why $\lambda$ is small, take a whole coordinate sphere
$S=\sphere\times\{z\}$ at time $-1$ near the distinguished point.
Its ambient mean curvature vector is $O(\eta)$ and
$\tr_{TS}\Ric\ge C^{-1}>0$. The restriction formula and integration give
\[
 0=\int_S\Delta_S F\,\mathrm d\mu_S
 =\lambda\int_S\tr_{TS}\Ric\,\mathrm d\mu_S+
   \int_S\langle\mathbf H_S,\nabla F\rangle\,\mathrm d\mu_S.
\]
Since $|\nabla F|\le1$, this implies $\lambda\le C\eta$.
On this sphere the restricted Hessian of $F$ is $O(\eta)$;
the mean-zero sphere estimate gives $|\nabla_S F|\le C\eta$.
The gradient identity then gives a flat component of length
$1+O(\eta)$. Rotate that component to $\partial_{z_1}$.
Integrating $\Hess F=\lambda\Ric$ on the fixed convergence
domain proves the spatial estimates in \eqref{cylinder:eq:F-close}.
The equation $\partial_tF=\lambda \mathrm{R}$ gives the temporal ones.
All constants here may depend on the already fixed domain $T$.
This proves the stated conclusion.
\end{proof}

Let
\begin{equation}
 \Pheat_g V:=\Delta_g V+\Ric_g^\sharp V.
\end{equation}
For every vector field $V$, we have
\begin{equation}\label{cylinder:eq:Bianchi}
 \big(\operatorname{div}\Lie_Vg-
 \tfrac12\mathrm{d}\operatorname{tr}_g\Lie_Vg\big)^\sharp=\Pheat_gV.
\end{equation}
If $\partial_tV=\Pheat_gV$, then for $h=\Lie_Vg$, we obtain
\begin{equation}\label{cylinder:eq:h-equation}
 \partial_th=\Delta_{L,g}h.
\end{equation}
These identities follow respectively by commuting two derivatives
of $V$ and by using
$2(D\Ric)_g(h)=-\Delta_Lh+
\Lie_{(\operatorname{div}h-\frac12\mathrm{d}\operatorname{tr}h)^\sharp}g$.

\begin{lem}\label{cylinder:lem:pairing}
Let $g(t)$ be a smooth Ricci flow on a four-dimensional domain,
let $\lambda>0$ be constant, and let $F$ and $V$ be a smooth
function and a smooth time-dependent vector field.
If $F$ satisfies \eqref{cylinder:eq:F-identities} and
$\partial_tV=\Pheat_gV$, then $u=VF=\langle V,\nabla F\rangle$
satisfies
\begin{equation}\label{cylinder:eq:u-equation}
 (\partial_t-\Delta_g)u=-\lambda\langle\Ric,h\rangle,
 \qquad \mathrm{where}\qquad h=\Lie_Vg.
\end{equation}
\end{lem}
\begin{proof}
In invariant notation the scalar product rule is
\[
 \Delta(VF)=(\Delta V)F+V\Delta F+
 \Ric(V,\nabla F)+2\nabla_iV^j\nabla_i\nabla_jF.
\]
On the other hand,
$\partial_t(VF)=(\Delta V+\Ric^\sharp V)F+V\Delta F$.
The Ricci terms cancel. The symmetry of $\Hess F$ and
$\Hess F=\lambda\Ric$ give
$-2\nabla_iV^j\nabla_i\nabla_jF
=-\lambda\langle\Ric,\Lie_Vg\rangle$.
This finishes the proof.
\end{proof}

Thus, the scalar source is proportional to the symmetry error, with
an additional factor of geometric smallness.

We set $\krot:=\sqrt{6\pi}$. We fix sphere rotation fields $J_a$ with
$[J_a,J_b]=\epsilon_{abc}J_c$, where $\epsilon_{abc}$ is the alternating symbol with $\epsilon_{123}=1$ and repeated indices are summed. We set $K_a:=\krot J_a$.
On every round sphere, we have
\begin{equation}
 \Area(S)^{-2}\int_S\langle K_a,K_b\rangle\,\mathrm{d}\mu_S=\delta_{ab},
 \qquad [K_a,K_b]=\krot\epsilon_{abc}K_c.
\end{equation}

The heat correction requires one triple on the entire comparison
region. We first choose constant matrices on overlaps and then glue
the aligned triples.

\begin{lem}\label{cylinder:lem:alignment}
Let $s>0$ and consider the product neighborhoods
\[
 \mathcal P_s:=\bigl(S^2\times\Disk_{10\sqrt{s}}\bigr)\times[-2s,-s],
 \qquad
 \mathcal P'_s:=\bigl(S^2\times\Disk_{5\sqrt{s}}\bigr)
                       \times[-3s/2,-s].
\]
Suppose the metric, regular potential levels, and a given smooth
foliation of these levels by CMC two-spheres are sufficiently close
to the cylindrical model in the rescaled $C^6$ topology after
the parabolic rescaling
$(z,t,g)\mapsto(s^{-1/2}z,s^{-1}t,s^{-1}g)$.
The potential is rescaled by $F\mapsto s^{-1/2}F$, up to an
additive constant. Let $\eta$ bound the rescaled geometric
perturbation, and assume $0\le\epsilon\le\eta$ with $\eta$
sufficiently small.
Suppose two smooth time-independent triples $U_a,\widetilde U_a$,
$a=1,2,3$, satisfy the metric Lie derivative, two-normal, and
whole-leaf Gram bounds of Definition~\ref{cylinder:def:symmetry}
on $\mathcal P_s$, with $r=\sqrt{s}$ and error $\epsilon$.
The CMC leaves through $\mathcal P'_s$ are assumed to be whole
leaves contained in the larger spatial product.
\begin{enumerate}[label=\textup{(\roman*)},ref=\textup{(\roman*)}]
\item There is a constant matrix $O\in O(3)$ such that
\begin{equation}\label{cylinder:eq:alignment}
 \sup_{\mathcal P'_s}\max_{1\le a\le3}
 \sum_{j=0}^3 s^{(j-1)/2}
       |\nabla^j(U_a-O_{ab}\widetilde U_b)|\le C\epsilon.
\end{equation}
\item For three triples on overlapping such neighborhoods, the alignment
matrices satisfy the cocycle estimate
\eqref{cylinder:eq:alignment-cocycle} whenever the common smaller
overlap contains a whole CMC leaf with fixed interior margins in
these rescaled coordinates.
\end{enumerate}
\end{lem}
\begin{proof}
We first prove (i). For simplicity, we take $s=1$. The general case follows by the
parabolic rescaling in the statement: the norm of the $j$th
covariant derivative of a vector field acquires the factor
$s^{(j-1)/2}$.
We first control the values and first derivatives of the local fields.
We then choose an orthogonal alignment and estimate the difference
between the aligned triples.

We give the compactness argument with the required kernel condition.
For an approximate Killing field with metric Lie derivative $h=\Lie_Ug$, the
Killing transport identity has the form
\begin{equation}\label{cylinder:eq:Ktransport}
 \nabla\nabla U=\Rm*U+
 \tfrac12\big(\nabla h+\nabla h-\nabla h\big),
\end{equation}
with the last three terms carrying the usual permuted indices.
It bounds the value and first derivative along fixed-length paths
by their initial values and $\|h\|_{C^1}$. The initial value and first covariant derivative
can be bounded directly, without invoking a separate Korn estimate.
Let $B(x,a)$ be an interior ball, let $M:=\sup_{B(x,a)}|U|$, and
suppose $|\Rm|\le K$ there. Along a unit-speed geodesic
from $x$, use a parallel frame and integrate
\eqref{cylinder:eq:Ktransport} twice. Taylor's formula gives, for every
unit initial velocity $v$,
\[
 a|\nabla_vU(x)|
 \le 2M+Ca^2\bigl(KM+\|\nabla h\|_{C^0(B(x,a))}\bigr).
\]
Thus a fixed interior ball bounds $|U|+|\nabla U|$ at its center
by $M+C\|h\|_{C^1}$. Along any path of length at most $D$,
the coupled transport system for $(U,\nabla U)$ then gives
\[
 \sup_{\rm path}(|U|+|\nabla U|)
 \le C_D\bigl(|U(x)|+|\nabla U(x)|+\|h\|_{C^1}\bigr)
\]
by Gronwall's inequality, with curvature bounds included in $C_D$.
Finally \eqref{cylinder:eq:Ktransport} bounds $\nabla^2U$, and its first
derivative bounds $\nabla^3U$ by the controlled value and first covariant derivative,
$|\nabla\Rm|$, and $\|h\|_{C^2}$. This supplies both
the compactness and the $C^3$ estimates used below. All paths
stay in the larger domain; all their lengths and
interior ball radii are fixed before the small parameters.

The Gram condition and the full-neighborhood transverse bounds
give a uniform bound for $U$ on a smaller region. Otherwise divide
by an unbounded norm, use the preceding transport and elliptic
bounds, and pass to a nonzero model Killing field. It vanishes
on a central sphere by the normalized Gram bound and has no flat
component everywhere by the transverse bounds. The only model
Killing fields with these properties are zero. Here the full
neighborhood matters: a planar rotation vanishes on the central
sphere but not on neighboring spheres.

We now make the nonlinear alignment estimate explicit. We fix nested
connected full-sphere product regions $\Omega'\Subset\Omega$ and
one whole central leaf $S\subset\Omega'$. For two vector fields $V,W$, we write
\[
 \langle V,W\rangle_S :=\Area(S)^{-2}\int_S\langle V,W\rangle\,\mathrm{d}\mu.
\]
For any bounded rotation triple $W=(W_1,W_2,W_3)$, we have a uniform linear
coercivity estimate
\begin{equation}\label{cylinder:eq:alignment-coercivity}
 \|E\|_{C^0(\Omega')}\le C\left(
 \|\Lie_Eg\|_{C^2(\Omega)}+
 \|EF\|_{C^0(\Omega)}+
 \|\langle E,\nu\rangle\|_{C^0(\Omega)}+
 \sum_b|\langle E,W_b\rangle_S|\right).
\end{equation}
Here $E$ is a single arbitrary smooth field; it need not be a
rotation triple. To prove the uniform estimate, suppose it fails for
perturbations tending to zero, divide by the left-hand side,
and use the explicit transport estimates above. These give
compactness on $\overline{\Omega'}$, and extend the bounds to
slightly larger regions inside $\Omega$. The limit is a nonzero
model Killing field. Its two transverse components vanish on
an open product neighborhood, so it is a sphere rotation.
The limiting $W_b$ form a basis of those rotations by their
Gram normalization. All three last pairings vanish, so this
sphere rotation is zero, a contradiction. This proves a
uniform constant for all sufficiently small perturbations.
The same transport estimate then controls the $C^3$ norm on a
slightly smaller region by the right-hand side of
\eqref{cylinder:eq:alignment-coercivity} and the metric Lie derivative norm.

We choose $O\in O(3)$ minimizing
$\sum_a\langle U_a-(O\widetilde U)_a,
U_a-(O\widetilde U)_a\rangle_S$; a minimizer exists by
compactness of $O(3)$. We set $W:=O\widetilde U$, $E:=U-W$,
$D_{ab}:=\langle E_a,W_b\rangle_S$, and
$Q_{\mathrm{err},ab}:=\langle E_a,E_b\rangle_S$.
Differentiation along every skew-symmetric infinitesimal change
of $O$ gives $D=D^{\mathsf T}$. The exact Gram identity is
\[
 Q_U-Q_W=D+D^{\mathsf T}+Q_{\mathrm{err}}=2D+Q_{\mathrm{err}}.
\]
Both rotation triple Gram matrices are $I+O(\epsilon)$, so, with
$e=\max_a\|E_a\|_{C^0(\Omega')}$, we have
\[
                    |D|\le C\epsilon+Ce^2.
\]
We apply \eqref{cylinder:eq:alignment-coercivity} to each $E_a$. The
metric Lie derivatives and both transverse pairings of $E_a$ are $O(\epsilon)$,
since the basis change is constant and orthogonal. Consequently, it follows that
\[
                        e\le C\epsilon+Ce^2.
\]
It remains to justify absorption uniformly. As the perturbation
and rotation triple error tend to zero, the bounded triples converge,
after subsequences, to Gram-orthonormal bases of the same model
sphere Killing algebra. The minimizing $O$ has a convergent
subsequence and makes their limiting whole-leaf difference zero.
The model difference therefore vanishes on the entire product.
Thus $e=o(1)$ uniformly as $\eta\to0$, for $\epsilon\le\eta$.
We choose a fixed geometric threshold making $Ce\le1/2$.
The last inequality then gives $e\le2C\epsilon$, and the transport
estimate yields the claimed $C^3$ bound.

All these choices can be made at one common time. The aligned
difference and $O$ are time-independent. Uniform metric and
connection comparisons on the fixed backward rotation triple interval
therefore give the same bound throughout that interval.

This proves (i). We now prove (ii). On a triple overlap, write $U_i=(U_{i,1},U_{i,2},U_{i,3})$ for the $i$th triple and
choose $O^{ij}$ so that $U_i-O^{ij}U_j$ satisfies
\eqref{cylinder:eq:alignment}. Put
${\mathfrak C}:=O^{ij}O^{jk}-O^{ik}$. On a whole common CMC leaf $S$, the identity
\[
 {\mathfrak C} U_k=
 O^{ij}(O^{jk}U_k-U_j)
 +(O^{ij}U_j-U_i)
 +(U_i-O^{ik}U_k)
\]
gives $|\mathfrak C U_k|\le C\epsilon$. Its Gram matrix
$Q^{(k)}(S)$ satisfies $Q^{(k)}(S)\ge\tfrac12 I$ for small
$\epsilon$. Hence
\[
 \tfrac12|{\mathfrak C}|_{\mathrm F}^2
 \le\operatorname{tr}\bigl(\mathfrak C Q^{(k)}(S){\mathfrak C}^{\mathsf T}\bigr)
 =\Area(S)^{-2}\int_S\sum_a|(\mathfrak C U_k)_a|^2\,\mathrm d\mu_S
 \le C\epsilon^2.
\]
Here $|\cdot|_{\mathrm F}$ is the Frobenius matrix norm, and the
area of $S$ is bounded above and below in the normalized product.
Consequently, we have
\begin{equation}\label{cylinder:eq:alignment-cocycle}
 |O^{ij}O^{jk}-O^{ik}|_{\mathrm F}\le C\epsilon.
\end{equation}
The Gram matrix and this matrix inequality are unchanged by
parabolic rescaling, so the same estimate holds for every $s>0$.
This proves (ii) and completes the proof.
\end{proof}

Before applying the improvement argument, we must put the local
rotational triples into a common frame and glue them on the
comparison region. The following lemma supplies this construction.

\begin{lem}\label{cylinder:lem:gluing}
In the setting of Theorem~\ref{cylinder:improvement}, we translate
and rescale so that the distinguished point is $(x_0,-1)$ and its
scalar curvature is one. Fix $T\ge2$,
choose its comparison domain to contain $D_{8\sqrt T}\times[-2T,-1]$
with the margins used below, and choose its comparison error $\eta$
after $T$.
Assume all points of the large neighborhood in the theorem are
$\epsilon$-symmetric. Then there is a smooth time-dependent triple $(U_1,U_2,U_3)$ on
$D_{8\sqrt T}\times[-2T,-1]$ such that, with $q=10$, $s=-t$, and $W=2T/s\ge1$, we have
\begin{align}
 |\Lie_{U_i}g|+\sqrt s\,|\nabla\Lie_{U_i}g|
       +s|\nabla^2\Lie_{U_i}g|&\le C W^q\epsilon,
 \label{cylinder:eq:U-metric-variation}\\
 |\partial_tU_i|+|\Pheat_gU_i|&\le C s^{-1/2}W^q\epsilon,
 \label{cylinder:eq:U-source}\\
 |U_i|&\le C\sqrt s,\qquad |U_iF|\le C\sqrt s\,W^q\epsilon.
 \label{cylinder:eq:U-size}
\end{align}
After one fixed orthogonal choice of basis, the spatial derivatives
of orders $0,1,2,3$ satisfy
\begin{equation}\label{cylinder:eq:U-model-closeness}
 \sup_{(x,t)\in D_{8\sqrt T}\times[-2T,-1]}
 \max_{1\le i\le3}\sum_{j=0}^3(-t)^{(j-1)/2}
 \bigl|\bar\nabla^j(U_i-K_i)\bigr|_{\bar g(t)}
 \le C(T)\eta.
\end{equation}
Here $\bar\nabla$ is the spatial connection of the comparison cylinder and $K_i$ are the normalized round rotations. The inequality \eqref{cylinder:eq:U-model-closeness} is understood after pulling back by $\Psi$ as in Definition \ref{cylinder:def:neighborhood}. 
\end{lem}

\begin{proof}
In this proof $a,b,c$ label mesh boxes and $i=1,2,3$ labels
components of a triple. The construction has four stages: choose rotation triples on a fixed
parabolic mesh, align their bases along a tree, bound the error on
non-tree overlaps, and interpolate. The quantitative chain is
\[
 \text{triangle error }O(\epsilon)
 \ \Longrightarrow\ \text{loop error }O(W^4\epsilon)
 \ \Longrightarrow\ \text{aligned overlap error }O(W^4\epsilon)
 \ \Longrightarrow\ \eqref{cylinder:eq:U-metric-variation}\text{--}\eqref{cylinder:eq:U-size}.
\]
The mesh refinement below is used only to estimate loop errors;
the original tree fixes the bases used in the interpolation.

\medskip\noindent\textbf{Step 1. We choose the mesh and rotation triples.}
The potential comparison required by
Lemma~\ref{cylinder:lem:alignment} follows from
Lemma~\ref{cylinder:lem:heat-gauge}, after the local parabolic
rescaling and the corresponding rescaling of $F$.
The CMC families are the small-graph families of
Remark~\ref{cylinder:rem:cmc-existence}; their parameter estimates
and uniqueness follow from the same projected equation used in
Step~6 below. Thus these geometric hypotheses hold on the enlarged
comparison boxes once $\eta$ is sufficiently small for fixed $T$.
We use a dyadic cover in backward time, with boxes of spatial size
$C^{-1}\sqrt{s}$ and time size $C^{-1}s$, where $s=-t$ and $C>1$ is fixed
and large. Each expanded box is covered by a local rotation triple based
near its later time; the comparison domain contains them because
$100\sqrt{2T}+8\sqrt T\ll L$ and $202T<L$. We reserve whole-sphere
overlaps and fixed margins. Lemma~\ref{cylinder:lem:alignment}
gives constant orthogonal transitions whose triangle errors are
$C\epsilon$. The following construction controls their accumulation.

\medskip\noindent\textbf{Step 2. We fix compatible tree matrices.}
It suffices to select the large heat parameter along the sequence
$2T=4^J$. We set $s:=-t$ and $s_j:=2T\,4^{-j}$. We cover the time
slabs $s_{j+1}\le s\le s_j$ by rectangular boxes of flat
spatial side $h\sqrt{s_j}$ and time side $hs_j$, where $h$
is a small fixed dyadic number. We use a fixed flat square of
side $C\sqrt T$ containing the disk needed for the heat
construction. Its corners and the slightly enlarged square
are inside the original improvement chart after increasing
its fixed domain constant $L$.

At an interface the spatial meshes differ by a factor two.
We subdivide the coarse interface faces to agree with the finer
mesh. Triangulate each affected cell by coning its subdivided
faces to its center, using a common triangulation on shared
faces. Only a bounded number of faces and vertices is introduced
per cell, because neighboring scales differ by fixed factors.
The centers and subdivision points have dyadic coordinates.
There are finitely many shapes after parabolic rescaling, so
their local coordinate and partition bounds are uniform.
We choose all interface subdivisions in advance; adding a later
epoch will not alter an earlier cell or its data.

Every vertex receives an $\epsilon$-symmetry rotation triple whose
expanded domain contains its star and a fixed number of
neighboring stars. Base it near the later time of that expanded
star. A sufficiently small fixed $h$ makes this possible at
the curvature scale $r\asymp\sqrt s$. Neighboring clusters
have common whole sphere leaves and comparable scales; the
full-neighborhood alignment lemma supplies constant orthogonal
matrices $P^{ab}$ on their edges. Their triangle holonomies
are $O(\epsilon)$ by the whole-sphere Gram bound. We add
transitions for any bounded-distance overlaps using the same
expanded domains.

Let $K_j$ be the prefix through epoch $j$, including its finer
boundary subdivision. Its edge graph is connected. We construct
a spanning tree recursively: retain the entire old tree and
connect each new vertex to the old connected graph by a tree
inside the new slab. No edge connecting two previously reached
vertices is added. We fix the root gauge at the earliest epoch
and propagate it along these tree edges. We denote the resulting orthogonal matrix at a vertex $a$ by $G^a$.
Thus all gauges are fixed once assigned and do not change at later refinements.

\medskip\noindent\textbf{Step 3. Bound the non-tree overlap error.}
Let $s_*$ be the smallest time scale in $K_j$ and put
$W_*:=2T/s_*$. Refine the whole prefix to a uniform rectangular
grid with flat step $C^{-1}h\sqrt{s_*}$ and time step $C^{-1}h s_*$,
where a fixed dyadic $C^{-1}$ includes the cell centers and interface
vertices. There are $O(\sqrt{W_*})$ intervals in each flat
direction and $O(W_*)$ in time. Consequently, the refined grid
has at most $CW_*^2$ vertices and cells. The original complex
has no more than this number of vertices, up to a fixed factor.
Its tree paths therefore have length at most $CW_*^2$.

We label each refinement vertex by an original local rotation triple whose
expanded box covers its small star. At original vertices retain
their already assigned labels. Such a label can be selected
from the original boxes near that location, so all neighboring
labels have comparable curvature scales. Each small triangle
has a common whole-sphere comparison region and hence an
$O(\epsilon)$ cocycle error. Repeated labels have identity
transition.

We replace an original edge by a coordinate path in its incident
cell on the refined grid. This path has length at most $CW_*$.
All its labels and the original endpoint labels lie in one
bounded cluster of comparable-scale original boxes. Their
expanded rotation triple domains have a common whole sphere and a
common time interval. Triangle fans based at one endpoint
compare the original transition with the refined path product
at a total cost of $CW_*\epsilon$. This comparison does
not redefine the original matrix or its gauge.

We now take the original tree loop for an edge in $K_j$. It has
at most $CW_*^2$ edges. The preceding replacements give a
refined loop of length at most $CW_*^3$, at a total holonomy
error at most $CW_*^3\epsilon$.

On the uniform rectangular refinement, choose the coordinate
tree consisting of one time chain, all rows in the first flat
direction on one side, and all spokes in the second flat
direction. Let $n_1,n_2,n_t$ be the numbers of grid intervals in the two
flat directions and in time. Filling the coordinate-tree edge loops
by rectangles shows that any refined edge's gauge discrepancy is
bounded by
\[
 C(n_1+n_2+n_t)\epsilon\le CW_*\epsilon.
\]
Equivalently, its coordinate-tree loop fills with at most
$CW_*$ elementary squares. We multiply these edge discrepancies
around the refined loop: the coordinate-tree matrices cancel,
and orthogonal conjugations preserve the norm. Its holonomy
is at most $CW_*^4\epsilon$. Including the replacement
error gives the same bound for the original tree loop.
Therefore the original, compatible gauges satisfy
\[
             \|G^aP^{ab}(G^b)^{-1}-I\|
                       \le CW_*^4\epsilon.
\]
Extra local overlap edges are reduced to a bounded number of
original edges and local triangles, at the same order.

This reasoning uses the uniform refinement only to estimate
holonomy. It never replaces the previously chosen tree or its
gauges, which is the required compatibility between epochs.
For an overlap at time $s$, use the prefix ending at its epoch
or the next neighboring epoch. Its smallest scale is comparable
to $s$, so $W_*\le C(2T/s)$.

\medskip\noindent\textbf{Step 4. We keep the alignment fixed in time and interpolate.}
The passage from a single-time alignment to a parabolic overlap
uses time independence of the fields and controlled connections.
More precisely, the following estimate holds. Suppose a fixed spatial set
$\Omega$ lies in both expanded rotation triple domains throughout an
interval $I$ of length at most $C_0r^2$. Assume on a neighborhood
of $\Omega\times I$ that
\[
 |\nabla^j\Rm_{g(t)}|\le C_jr^{-2-j},
                         \qquad 0\le j\le3.
\]
For any time-independent field $E$ and $\tau,t\in I$, we have
\[
 \sum_{j=0}^3 r^{j-1}|\nabla_{g(t)}^jE|_{g(t)}
 \le C\sum_{j=0}^3r^{j-1}|\nabla_{g(\tau)}^jE|_{g(\tau)}
\]
pointwise on $\Omega$, with $C$ depending only on the displayed
bounds and $C_0$. We rescale to $r=1$ to prove this. The metric
evolution and bounded Ricci tensor first give uniform metric
comparability. The connection evolution is
\[
 \partial_t\Gamma^k_{ij}
 =-g^{k\ell}(\nabla_i\Ric_{j\ell}
             +\nabla_j\Ric_{i\ell}
             -\nabla_\ell\Ric_{ij}).
\]
We integrate it, and its first two covariant spatial derivatives,
on the bounded time interval. Induction and Gronwall's inequality
give bounded connection difference and its derivatives through
order two, measured in the fixed connection at $\tau$.
Expanding $\nabla_{g(t)}=\nabla_{g(\tau)}+\mathsf a(t)$ on tensors
then gives the asserted estimate through order three.

We apply this to $E=U_a-P^{ab}U_b$, choosing $P^{ab}$ only once,
at a common time $\tau$. The local alignment estimate at $\tau$
therefore holds throughout the overlap, with the same constant
matrix. No error proportional to the sizes of $U_a,U_b$ is added:
the comparison is applied to their already small difference.

The spatial set $\Omega$ here must contain the union of the small
boxes over the overlap interval. This is why expanded fixed
spatial rotation triple domains were reserved. On cylindrical grid boxes
the flat parameters are fixed and the sphere factor is whole;
neighboring epoch scales differ by fixed factors. On cap boxes
the evolving radial parameter moves at bounded speed over the
short time interval, and a fixed additional spatial margin
contains that union. These are applications of this estimate
on one fixed spatial set, not a new choice of alignment at each
time or an identification of different moving domains.

Let $\chi_a$ be the parabolic partition of unity and write
$\widetilde U_a:=G^aU_a$. We define
\[
 U:=\sum_a\chi_a\widetilde U_a,\qquad \sum_a\chi_a=1.
\]
On a box with reference rotation triple $\widetilde U_{a_0}$, we have
\begin{align*}
 \Lie_Ug
 &=\sum_a\chi_a\Lie_{\widetilde U_a}g
   +2\sum_a\mathrm d\chi_a\mathbin{\odot}
          (\widetilde U_a-\widetilde U_{a_0})^\flat,\\
 \partial_tU
 &=\sum_a(\partial_t\chi_a)
          (\widetilde U_a-\widetilde U_{a_0}),
\end{align*}
where $\odot$ denotes the symmetric tensor product.
These identities show that the terms involving the reference rotation triple cancel
in every cutoff term.
A spatial derivative of a cutoff costs $s^{-1/2}$ and a time
derivative costs $s^{-1}$. Derivatives of the interpolated metric Lie derivative
involve only the metric Lie derivative of the rotation triples and up to three
derivatives of their differences on overlaps. The bounds just
proved give \eqref{cylinder:eq:U-metric-variation}. The overlap bounds cost at most $W^4$; the normalized cutoff derivatives contribute only the displayed powers of $s$. Thus $q=10$ suffices in all these estimates.
Since the rotation triples themselves are time independent, subtracting
one of them in the time-cutoff terms gives
$|\partial_tU_i|\le Cs^{-1/2}W^q\epsilon$.
Equation \eqref{cylinder:eq:Bianchi} gives the other term in
\eqref{cylinder:eq:U-source}. Convex interpolation and the local norm
and transverse bounds prove \eqref{cylinder:eq:U-size}.

Finally, we view each bounded local rotation triple in the exact product
metric. Its metric Lie derivative through order two, transverse errors, and Gram
error are $O(\eta)$ there: replacing the metric and normals costs
$O(\eta)$ times the already bounded $C^3$ field norm, and
$\epsilon\le\eta$. We apply Lemma~\ref{cylinder:lem:alignment} in the exact
product to this triple and the standard $K_a$. It yields an
$O(\eta)$ comparison through order three after a constant
orthogonal change of basis. Over the finitely many boxes for
fixed $T$, the rooted choices change this bound by at most
$C(T)\eta$. This uses the proved alignment estimate, rather
than an additional high-order estimate for the local rotation triples. The assertion follows.
\end{proof}

The exponent $q$ will never determine a derivative order or a number
of brackets. Its only later effect is a finite constant $C_q$ in
a heat estimate. We fix $q_U=10$ for the glued data and
$q_h=11$ for the corrected heat data; neither exponent depends on $T$.

\begin{proof}[Proof of Theorem~\ref{cylinder:improvement}]
\textbf{Step 1. Align and glue the local rotations.}
We apply Lemma~\ref{cylinder:lem:gluing} with the fixed exponent
$q_U=10$. It produces the time-dependent fields $U_a$ on
$D_{8\sqrt T}\times[-2T,-1]$ with metric Lie derivatives, potential pairings,
and vector heat defects proportional to $\epsilon$. The construction
above replaces the one-flat-variable gluing in Brendle's Step 1;
its polynomial bounds are needed before the long heat evolution.

\medskip\noindent\textbf{Step 2. Correct the fields by vector heat flow.}
As in Brendle's Steps 2--3, solve the vector heat equation with
the glued fields as parabolic boundary data. The metric Lie derivative then
solves a homogeneous Lichnerowicz equation. In addition, the
potential pairing satisfies Lemma~\ref{cylinder:lem:pairing}.

\begin{claim}\label{cylinder:lem:heat-correction}
Put $Q_{\mathrm{heat}}:=D_{7\sqrt T}\times[-2T,-1]$.
There are fields $V_a$ solving the initial boundary problem
\[
 \begin{cases}
 (\partial_t-\Pheat_g)V_a=0&\text{on }Q_{\mathrm{heat}},\\
 V_a=U_a&\text{on }\partial_{\mathrm{par}}Q_{\mathrm{heat}},
 \end{cases}
\]
where
\[
 \partial_{\mathrm{par}}Q_{\mathrm{heat}}
 :=(\overline D_{7\sqrt T}\times\{-2T\})
 \cup(\partial D_{7\sqrt T}\times[-2T,-1]).
\]
Write $s:=-t$ and fix $q=q_h=11$. On $D_{6\sqrt T}\times[-3T/2,-1]$, we have
\begin{equation}\label{cylinder:eq:h-u-bounds}
 |h_a|\le C(2T/s)^q\epsilon,\qquad
 |u_a|\le C\sqrt s(2T/s)^q\epsilon,
 \quad h_a=\Lie_{V_a}g,\quad u_a=V_aF.
\end{equation}
For fixed $R,H>1$, chosen before $T$, put
\[
 Q_{R,H}^{\mathrm{term}}:=D_R\times[-H,-1],\qquad
 Q_{R,H}^{\mathrm{term},+}:=D_{R+1}\times[-H-1,-1].
\]
Choose $T$ so that $R+1<5\sqrt T$ and $H+1<T$.
On these smaller terminal domains all required derivatives of
$h_a,u_a$ are bounded by $C(T)\epsilon$, and $V_a-K_a$ is bounded by
$C(T)\eta$ in the required high-order norms.
\end{claim}

In the estimates for $w_a$ in this proof, the source exponent is
$q_U=10$ and the resulting exponent is $q_h=q_U+1=11$.
The linear parabolic Dirichlet problem on the compact product domain
has a solution; one can smooth incompatible corner data and take
an interior limit. Let $w_a:=V_a-U_a$. The norm evolution under
Ricci flow cancels the Ricci term in the vector equation, giving
the Kato inequality
\[
 (\partial_t-\Delta)|w_a|
 \le|\partial_tU_a-\Pheat_gU_a|
 \le Cs^{-1/2}(2T/s)^{q_U}\epsilon.
\]
The zero parabolic boundary values and a time barrier of the form
$C\epsilon\sqrt s(2T/s)^{q_h}$ imply
\begin{equation}\label{cylinder:eq:w-bound}
 |w_a|\le C\epsilon\sqrt s(2T/s)^{q_h}.
\end{equation}
Indeed, for $s=-t$ and $b(s)=C\epsilon\sqrt{s}(2T/s)^{q_h}$,
\[
 \partial_tb=C(q_h-\tfrac12)\epsilon\,s^{-1/2}(2T/s)^{q_h}
       \ge C'\epsilon\,s^{-1/2}(2T/s)^{q_U},
\]
because $1\le s\le2T$ and $q_h=q_U+1$. The barrier is
nonnegative on the parabolic boundary, where $w_a=0$.
Interior gradient estimates on cylinders of radius $C^{-1}\sqrt s$
and backward duration $C^{-1}s$ give
$|\nabla w_a|\le C(2T/s)^{q_h}\epsilon$ on the stated smaller domain.
There is an initial time margin at least $T/2$ and a flat
spatial margin $\sqrt T$. We choose the cylinder radius constant
large enough that $C^{-1}s<T/2$ and $C^{-1}\sqrt s<\sqrt T$ whenever
$s\le3T/2$. Thus, the cylinders are available even at the
new lower endpoint $t=-3T/2$.
Combining this with \eqref{cylinder:eq:U-metric-variation} proves the bound for $h_a$.
The bound for $u_a$ follows from \eqref{cylinder:eq:U-size},
\eqref{cylinder:eq:w-bound} and $|\nabla F|\le1$.

The tensor $h_a$ solves the homogeneous equation
\eqref{cylinder:eq:h-equation}. The bounds just proved on
$D_{6\sqrt T}\times[-3T/2,-1]$ give, by interior parabolic
regularity, its high spatial and time derivative bounds on
$D_{(11/2)\sqrt T}\times[-T,-1]$. In particular the entire
initial slice and lateral boundary of the subsequent
$D_{5\sqrt T}\times[-T,-1]$ model comparison have an earlier
time and spatial margin. The high derivative inhomogeneous term estimate
therefore does not use regularity at the initial time of a
newly posed Dirichlet problem. For $u_a$, apply the same interior estimates to the triangular
coupled system for $(h_a,u_a)$, whose scalar equation is
\eqref{cylinder:eq:u-equation} with zeroth-order coupling to $h_a$.
Equivalently, first estimate $h_a$ on a domain between the two
just displayed, and then estimate $u_a$ on the smaller one.
This retains an earlier time margin for the scalar source as well.
Finally, $K_a$ solves
the vector equation for the exact cylinder. The equation for
$V_a-K_a$ has source $\Pheat_gK_a=O(C(T)\eta)$ and boundary
values $O(C(T)\eta)$. The norm estimate and interior regularity
give the claimed high-order closeness. This finishes the proof of the claim.

\medskip\noindent\textbf{Step 3. We compare the metric Lie derivative and pairing with the exact model.}
The comparison is made with common initial and lateral values;
this ensures that its errors are proportional to $\eta\epsilon$.
On $Q_{\mathrm{cmp}}:=D_{5\sqrt T}\times[-T,-1]$ define
\[
 \begin{cases}
 (\partial_t-\Delta_{L,\bar g})\bar h_a=0,\quad
 (\partial_t-\Delta_{\bar g})\bar u_a=0&\text{on }Q_{\mathrm{cmp}},\\
 (\bar h_a,\bar u_a)=(h_a,u_a)&\text{on }\partial_{\mathrm{par}}Q_{\mathrm{cmp}}.
 \end{cases}
\] On smaller terminal domains, we have
\begin{equation}\label{cylinder:eq:comparison}
 \|h_a-\bar h_a\|_{C^m}+\|u_a-\bar u_a\|_{C^m}
 \le C(T)\eta\epsilon
\end{equation}
for the fixed finite derivative orders used below.
For the tensor difference, the inhomogeneous term is
$(\Delta_{L,g}-\Delta_{L,\bar g})h_a$.
For the scalar difference it is
$(\Delta_g-\Delta_{\bar g})u_a-
\lambda\langle\Ric,h_a\rangle$.
Lemma~\ref{cylinder:lem:heat-gauge} and Claim~\ref{cylinder:lem:heat-correction}
bound both by $C(T)\eta\epsilon$ with all needed derivatives.
We next estimate the norm of this comparison directly, making
its dependence on the symmetry error explicit. We write
$k:=h_a-\bar h_a$ and $v:=u_a-\bar u_a$. They have zero initial
and lateral data. In the model metric the tensor norm obeys
\[
 (\partial_t-\Delta_{\bar g})|k|
       \le C(T)|k|+C(T)\eta\epsilon,
\]
by Kato's inequality and bounded curvature on this fixed domain;
the time variation of the tensor norm contributes only to the
zeroth-order coefficient. The scalar difference obeys the same
inequality without the first term. Comparison with the scalar
ODE $\frac{\mathrm d y}{\mathrm dt}=C(T)y+C(T)\eta\epsilon$, $y(-T)=0$, gives
\[
           \sup(|k|+|v|)\le C(T)\eta\epsilon.
\]
This estimate requires no derivatives of the boundary data.
On the smaller terminal domains, interior parabolic estimates
using the already controlled derivatives of the inhomogeneous term give
\eqref{cylinder:eq:comparison}. In particular no term of size $\eta$
alone can enter from the common boundary or initial data.
These constants need not be polynomial in $T$ because $\eta$
will be chosen last. The earlier polynomial bounds for $h_a,u_a$
are used separately as data bounds for the exact-model equations;
they are not replaced by this potentially much larger $C(T)$.

The boundary and initial values of $\bar h_a,\bar u_a$ have
the polynomial bounds \eqref{cylinder:eq:h-u-bounds}. The following
spectral argument uses those boundary bounds directly; it
does not assume the small comparison error is independent of $T$.

\medskip\noindent\textbf{Step 4. Separate decaying and retained cylinder modes.}
This is the counterpart of the cylinder analysis used in Brendle's
Steps 4--5. We include every tensor block involving either flat
direction, as well as the scalar potential pairing.
The mode calculation below separates decay of angular eigencomponents from interior
regularity. Its polynomial boundary data constants are independent
of the large heat parameter. The retained derivative orders and the
nested terminal domains are recorded together in
Step~9 below.

All norms in this step are measured using the product metric at the
indicated time. Fixed terminal domains can be enlarged by a fixed factor
before applying interior estimates. In particular none of their constants
depend on $T$.

\begin{claim}\label{cylinder:lem:angular-decay}
An angular block means an invariant spherical spectral component of
a coefficient in the decomposition
\[
 h=\omega g_{\sphere}+\chi+\sum B_{ij}\,\mathrm dz_i\otimes\mathrm dz_j
       +\sum\mathrm dz_i\odot\sigma_i,
\]
or of the homogeneous model comparison $\bar u_a$ for the scalar
pairing $u_a=V_aF$, at each fixed flat coordinate $z$.
Here $\chi$ is sphere trace-free and $\sigma_i$ are sphere one-forms;
the relevant spherical projections commute with the flat heat operator.
Suppose such a block of a product heat equation has solution operator with norm
bounded by $C(u/v)^\alpha$ from time $-v$ to time $-u$, where
$1\le u\le v$ and $\alpha>0$. On $D_{5\sqrt T}\times[-T,-1]$,
suppose its initial norm is at most $C\epsilon$ and its lateral boundary
norm at time $-v$ is at most $C(T/v)^q\epsilon$.
The initial norm is the supremum on $D_{5\sqrt T}\times\{-T\}$;
the lateral norm at time $-v$ is the supremum on
$\partial D_{5\sqrt T}\times\{-v\}$, with the fiber norms of the
product metric at that time. On each $Q_{R,H}^{\mathrm{term}}$ defined
in Claim~\ref{cylinder:lem:heat-correction}, with fixed $R,H$ and
$T$ sufficiently large, its norm is at most $C_qT^{-\alpha}\epsilon$.
And for a scalar with initial size $C\sqrt T\epsilon$ and lateral size
$C\sqrt v(T/v)^q\epsilon$, the nonconstant spherical part is at most
$C_qT^{-1/2}\epsilon$ on that domain.
\end{claim}

The flat heat equation and the time-dependent angular operator commute.
Thus, the solution is the integral of angular solution operators against the
initial heat measure and lateral caloric measure of the flat disk.
For points a distance at least $C^{-1}\sqrt T$ from its boundary, the measure
of lateral caloric measure at backward durations at most $a$ is bounded
by $C\exp(-T/(Ca))$. For completeness, this estimate follows by applying
the exponential supersolutions for the one-dimensional heat equation
to each coordinate of the disk: exit from the disk requires one
coordinate to move by at least a fixed multiple of $\sqrt T$.
Optimization of the exponential parameter gives the stated bound.

At a lateral time $-v$, with $v=u+a$, the weighted bound is
\[
 C\epsilon (u/T)^\alpha(T/v)^{q+\alpha}.
\]
We decompose the durations into
$[2^{-j-1}T,2^{-j}T]$. The integral of the factor remaining after
$\epsilon(u/T)^\alpha$ is bounded by
\[
 C\sum_{j\ge0}2^{(j+1)(q+\alpha)}e^{-2^j/C}<\infty.
\]
This is a number depending on $q,\alpha$ and the fixed terminal
domain, not on $T$. To make the constant explicit, let
\[
 \Gamma_u:=\partial\Disk_{5\sqrt T}\times[-T,-u]
\]
be the lateral boundary available to a terminal point $(z,-u)$,
and let $\omega_{z,-u}$ be its lateral caloric measure for the flat
heat equation. This measure has total measure at most one. For $x\ge2$,
the exit estimate gives
\[
 \omega_{z,-u}\bigl(\{(y,-v)\in\Gamma_u:T/v\ge x\}\bigr)
       \le C e^{-x/C},
\]
since $T/v\ge x$ implies that the backward duration $v-u$ is
at most $T/x$. Integrating this tail bound, for each fixed $p>0$,
we obtain
\[
 \int_{\Gamma_u}\left(\frac{T}{v}\right)^p
                  \,\mathrm{d}\omega_{z,-u}(y,-v)
 \le 2^p+Cp\int_2^\infty x^{p-1}e^{-x/C}\,\mathrm{d}x
 \le 2^p+C^{p+1}\Gamma(p+1).
\]
Taking $p=q+\alpha$ (or $q+1/2$ for
the scalar bound) gives an explicit finite choice of the moment
constant. Thus the angular decay constant can be fixed before
choosing $T$, however large it is as a function of the fixed $q$.
The initial contribution has the same bound. The first nonzero scalar
sphere eigenvalue is $2$; consequently its angular solution operator is
$C u/v$. The corresponding lateral bound is
\[
 C\epsilon\frac{u}{\sqrt T}(T/v)^{q+1/2}.
\]
The same convergent series proves the scalar assertion. The finitely
many fixed terminal values of $u$ are absorbed in the constant. This finishes the proof of the claim.

\begin{claim}
For the model solutions constructed above, on
$Q_{R,H}^{\mathrm{term}}=D_R\times[-H,-1]$, with fixed $R,H$ and
interior estimates taken from $Q_{R,H}^{\mathrm{term},+}$, we have
\begin{align}
 \bar h_a&=Q_a(z,t)+(-t)\psi_a(\theta)g_{\sphere}+E_a,
 \label{cylinder:eq:decomp}
 \end{align}
 and
 \begin{align}
 \|(I-\Pi_0)\bar u_a\|_{C^m}+\|E_a\|_{C^m}
 &\le C_qT^{-1/2}\epsilon.\label{cylinder:eq:decay}
\end{align}
Here $\Pi_0v:=(4\pi)^{-1}\int_{S^2}v\,\mathrm d\mu_{g_{\sphere}}$
is spherical averaging at fixed $z,t$,
$\Pi_1$ is the $L^2(g_{\sphere})$ projection onto the first spherical
harmonics, and
$\Pi_{\ge2}=I-\Pi_0-\Pi_1$; $\psi_a$ is a time-independent first
spherical harmonic with norm at most $C_q\epsilon$, and
\[
 Q_a=\alpha_a(z,t)g_{\sphere}+
       \sum_{i,j=1}^2B_{a,ij}(z,t)\,\mathrm{d}z_i\otimes \mathrm{d}z_j
\]
is invariant under every $K_b$. Its fixed terminal $C^m$ norm is at
most $C(T)\epsilon$. The estimates hold for any of the fixed derivative
orders used later.
\end{claim}

We write a symmetric tensor as
\[
 h=\omega g_{\sphere}+\chi+
   \sum B_{ij}\,\mathrm{d}z_i\otimes \mathrm{d}z_j+
   \sum \mathrm{d}z_i\odot\sigma_i,
\]
where $\omega$ and $B_{ij}$ are scalar coefficients, $\chi$ is a
symmetric trace-free tensor on the sphere, and $\sigma_i$ are sphere
one-forms. Here $\Delta_z=\partial_{z_1}^2+\partial_{z_2}^2$ and
$\Delta^{\rm rough}_{\sphere}$ is the rough Laplacian of the unit sphere.
With $u=-t$, substitution in $\partial_th=\Delta_{L,\bar g}h$ gives
\begin{align*}
 \partial_t\omega&=\Delta_z\omega+(2u)^{-1}\Delta_{\sphere}\omega,\quad
 \partial_t B_{ij}=\Delta_zB_{ij}+(2u)^{-1}\Delta_{\sphere}B_{ij},\\
 \partial_t\sigma_i&=\Delta_z\sigma_i+
       (2u)^{-1}(\Delta^{\rm rough}_{\sphere}\sigma_i-\sigma_i),\\
 \partial_t\chi&=\Delta_z\chi+
       (2u)^{-1}(\Delta^{\rm rough}_{\sphere}\chi-4\chi).
\end{align*}
Here the signs and the time-dependent norm factors can be checked
without a tensor-spectrum convention. We set $a:=2u$, so that the
sphere has curvature $a^{-1}$ and $\Ric|_{TS^2}=g_{\sphere}$.
In the convention \eqref{eq:lichnerowicz-convention},
$\Delta_Lh:=\Delta h+2\Rm(h)-\Ric\circ h-h\circ\Ric$,
the curvature terms cancel on $\omega g_{\sphere}$, equal
$-a^{-1}h$ on a mixed tensor, and equal $-4a^{-1}\chi$ on a
sphere trace-free tensor. Moreover, the norms of the sphere trace,
mixed, and sphere trace-free components include respectively the
factors $a^{-1},a^{-1/2},a^{-1}$. Thus a coefficient decay
$(u/v)^\beta$ becomes norm decay $(u/v)^{\beta-1}$ for two
sphere covariant indices, and $(u/v)^{\beta-1/2}$ for one.
This accounts for all the exponents below, including the zero-eigenvalue
first-harmonic sphere trace. It also explains why a tensor
coefficient estimate alone would give an incorrect gain.

Angular time is $\frac12\log(v/u)$. Scalar eigenvalues are
$\ell(\ell+1)$. On one-forms the first eigenvalue of
$-\Delta^{\rm rough}$ is $1$, by the Hodge decomposition and the
Weitzenb\"ock formula. On trace-free tensors the rough operator is
nonpositive, which already suffices here. After also accounting for
the metric factors in tensor norms, the coefficient and norm exponents are related by
\[
 \alpha=\beta-\frac{k}{2},
 \qquad k=\text{number of sphere covariant indices}.
\]
For a scalar harmonic of degree $\ell$, the coefficient exponent
is $\beta=\ell(\ell+1)/2$. Thus, the nonconstant flat block and
the higher-harmonic sphere trace satisfy, respectively,
\begin{align*}
 B_{ij}\,\mathrm dz_i\otimes\mathrm dz_j:\quad
 &\alpha=\frac{\ell(\ell+1)}2\ge1
       &&(\ell\ge1),\\
 \omega g_{\sphere}:\quad
 &\alpha=\frac{\ell(\ell+1)}2-1\ge2
       &&(\ell\ge2).
\end{align*}
The one-form spectral bound and the curvature term in the
trace-free equation give
\[
 \beta_{\mathrm{mixed}}\ge1,
 \qquad \alpha_{\mathrm{mixed}}\ge1-\tfrac12=\tfrac12,
 \qquad
 \beta_{\mathrm{tf}}\ge2,
 \qquad \alpha_{\mathrm{tf}}\ge2-1=1.
\]
These are sufficient lower bounds; sharp exponents are not needed.
The first-harmonic sphere trace instead has
\[
 \beta_{\mathrm{trace},1}=1,\qquad
 \alpha_{\mathrm{trace},1}=1-1=0,
\]
so it must be retained.
The spherical averages are retained
separately in $Q_a$. In particular the first-harmonic trace is zero-eigenvalue
in tensor norm, despite decay of its coefficient. The dyadic
caloric-measure estimate in Claim~\ref{cylinder:lem:angular-decay} handles
the polynomial lateral data before any terminal derivative estimate
is applied.
The asserted uniform angular solution operator bounds follow either from
the sphere spectral expansion, or from its $L^2$ spectral gap followed
by smoothing for one fixed angular time. Short angular times are
bounded directly by the heat semigroup. Claim~\ref{cylinder:lem:angular-decay}
therefore controls all the displayed blocks.

We retain the spherical averages. For the remaining first-harmonic
trace write $\Pi_1\omega=u a(\theta,z,t)$. Its equation is exactly
$\partial_ta=\Delta_z a$. The initial and lateral bounds for $a$
are $C\epsilon$ and $C(T/v)^q\epsilon$. The caloric-measure argument in
Claim~\ref{cylinder:lem:angular-decay}, now with exponent zero, bounds $a$ by
$C_q\epsilon$ on an inner disk of radius comparable to $\sqrt T$.
Interior flat heat estimates on such disks and time intervals give
\[
 |\nabla_z a|\le C_qT^{-1/2}\epsilon,\qquad
 |\partial_ta|\le C_qT^{-1}\epsilon
\]
on smaller inner disks and for $-T/4\le t\le-1$.
The same estimates apply to higher derivatives. We set
$\psi(\theta):=a(\theta,0,-1)$. Integration on the fixed terminal
domain shows that $u(a-\psi)g_{\sphere}$ has $C^m$ norm at most
$C_qT^{-1/2}\epsilon$ there. Initial and boundary bounds for the retained
averages give their stated $C(T)\epsilon$ bounds.
The scalar conclusion is precisely the scalar part of
Claim~\ref{cylinder:lem:angular-decay}. Interior parabolic estimates on slightly
larger fixed terminal domains give all the derivative estimates. This finishes the proof of the claim.

\medskip\noindent\textbf{Step 5. Remove the zero-eigenvalue component, take a bracket, and freeze.}
As in Brendle's Steps 6--8, subtract the first-harmonic correction
and use a cyclic bracket to remove rotation-invariant errors.
Here the same bracket also improves the potential pairing.

In the rest of the proof, we write
\begin{equation}\label{cylinder:eq:delta}
 d:=\big(C_qT^{-1/2}+C(T)\eta\big)\epsilon.
\end{equation}
Constants may be increased a fixed number of times; the final choices
of $T,\eta$ are made after all these increases. By model comparison,
\eqref{cylinder:eq:decomp}--\eqref{cylinder:eq:decay} hold for $h_a,u_a$ with errors
$O_{C^m}(d)$. We set
\[
 \zeta_a:=-\tfrac14\nabla_{\sphere}\psi_a,\qquad W_a:=V_a-\zeta_a.
\]
Since $\Hess_{\sphere}\psi_a=-\psi_a g_{\sphere}$, we know
\[
 \Lie_{\zeta_a}\bar g=(-t)\psi_a g_{\sphere},\qquad
 \Pheat_{\bar g}\zeta_a=0.
\]
It follows that
\begin{equation}
 \Lie_{W_a}g=Q_a+O_{C^m}(d),\qquad
 W_a-K_a=O_{C^m}(C(T)\eta).
\end{equation}
We use $\epsilon\le\eta$ in the second estimate, but never require
$\eta$ to be smaller than $\epsilon$.

The Bianchi identity and the vector heat equation give
\[
 \partial_tW_a=\Pheat_gW_a+\Pheat_g\zeta_a
              =b_a(z,t)+O_{C^{m-1}}(d),
\]
where $b_a$ is a flat, sphere-independent vector field. Indeed it is
the model Bianchi operator applied to $Q_a$. Its derivatives are
$O(C(T)\epsilon)$. In particular, $[K_b,b_a]=0$.

We define cyclically
\[
 X_1:=\krot^{-1}[W_2,W_3],\quad
 X_2:=\krot^{-1}[W_3,W_1],\quad
 X_3:=\krot^{-1}[W_1,W_2].
\]
The identity
$\Lie_{[V,W]}g=\Lie_V\Lie_Wg-\Lie_W\Lie_Vg$ and
$\Lie_{K_b}Q_a=0$ show that
\begin{equation}\label{cylinder:eq:X-est}
 \Lie_{X_a}g=O_{C^{m-1}}(d),\qquad
 \partial_tX_a=O_{C^{m-2}}(d).
\end{equation}
Every remaining product contains either an $O(d)$ error or a factor
$W-K=O(C(T)\eta)$ multiplied by a factor $O(C(T)\epsilon)$.

To control the transverse component, we use the following estimate. Since
$F-z_1=O_{C^m}(C(T)\eta)$, we have
$\zeta_aF=O_{C^{m-1}}(C(T)\eta\epsilon)$. Consequently, with $w_a=W_aF$, we get
\[
 (I-\Pi_0)w_a=O_{C^{m-1}}(d),\qquad w_a=O_{C^{m-1}}(C(T)\epsilon).
\]
For example, we can write
\[
 X_1F=\krot^{-1}(W_2w_3-W_3w_2)
     =\krot^{-1}\big(K_2(I-\Pi_0)w_3-K_3(I-\Pi_0)w_2\big)
       +O_{C^{m-2}}(C(T)\eta\epsilon).
\]
Thus $X_aF=O_{C^{m-2}}(d)$. We set
\begin{equation}
 Y_a(x):=X_a(x,-1).
\end{equation}
Integrating \eqref{cylinder:eq:X-est} over the fixed terminal time interval
proves, on smaller fixed terminal domains,
\begin{equation}\label{cylinder:eq:Y-est}
 \Lie_{Y_a}g=O_{C^{m-3}}(d),\quad
 Y_aF=O_{C^{m-3}}(d),\quad
 Y_a-K_a=O_{C^{m-3}}(C(T)\eta).
\end{equation}
These fields are time-independent. This first bracket suffices to
improve the potential pairing.

\medskip\noindent\textbf{Step 6. Control the CMC normal component by a second bracket.}
We adapt Brendle's Steps 9--10 within each potential level.
The extra term involving $Y_aF$ records the second ambient normal
direction and is controlled by Step 5.

We give the remaining geometric argument on fixed terminal domains.
Thus the constants in the following inverse estimates are fixed before
$T$ is selected. We decrease $\eta(T)$ so that the geometry, $F$, and
the fields are sufficiently close to the product there.

Inside each level of $F$, the CMC spheres are a smooth foliation.
This follows by writing them as zero-average graphs over the product
spheres and applying the implicit function theorem to the mean-zero
part of mean curvature. Its linearization is the sphere Laplacian
on zero-average functions, which is invertible. The same construction
with parameters gives bounds for all required derivatives in the
level parameter, the leaf parameter, and time. It also gives the
usual local uniqueness of the CMC foliation.

Let $\sigma$ be the leaf parameter inside each $F$ level, chosen so
that its positive normal lapse $\vartheta$ satisfies
\begin{equation}\label{cylinder:eq:lapse}
 \partial_\sigma=\vartheta\nu,\qquad \int_S\vartheta\,\mathrm{d}\mu=1.
\end{equation}

The first identity denotes the normal part of the parametrization;
tangential reparametrization has no effect below. Such a parameter
is obtained by integrating the volume between successive leaves.
We extend $\sigma$ smoothly between $F$ levels. We set
$D_F:=\nabla F/|\nabla F|^2$ and $\chi_F:=D_F\sigma$. For every ambient field $Y$, we have
\begin{equation}\label{cylinder:eq:sigma-identity}
 Y\sigma=\vartheta^{-1}\langle Y,\nu\rangle+\chi_FYF.
\end{equation}
In particular, $\vartheta^{-1}\nu$ is the gradient of $\sigma$ \emph{inside}
the level. Equation~\eqref{cylinder:eq:sigma-identity} accounts for the other
ambient direction and avoids identifying it with the ambient gradient.

We project $Y_a$ onto the level by
$Y_a^\top=Y_a-(Y_aF)D_F$. Its induced metric Lie derivative and all
needed derivatives are $O(d)$, by \eqref{cylinder:eq:Y-est}. Let $\mathrm{II}$ be the second fundamental form of a leaf within its
potential level, and let $\Delta_S$ be its intrinsic Laplacian. We write
$a_a:=\langle Y_a,\nu\rangle$ and
\[
 J:=\Delta_S+|\mathrm{II}|^2+\Ric_{\{F=\mathrm{const}\}}(\nu,\nu).
\]
The mean curvature variation formula gives
\begin{equation}
 Ja_a=H_a,\qquad H_a=O_{C^k}(d).
\end{equation}
Indeed, the metric variation is the small tensor
$\Lie_{Y_a^\top}g|_{T\{F=\mathrm{const}\}}$; its contribution to
the mean curvature variation uses at most its first derivative.
The tangential derivative of the leaf's constant mean curvature
vanishes. Applying the same formula to the leaf variation gives
$J\vartheta=b$, a function constant on each sphere (the sign of $b$ depends
on the convention for mean curvature and is immaterial).

We set $G_a(F,\sigma,t):=\int_Sa_a\,\mathrm{d}\mu$. The function
$e_a:=a_a-G_a\vartheta$ has zero integral. Since $J\vartheta=b$ and $\int_S\vartheta\,\mathrm d\mu=1$,
self-adjointness gives
\[
 bG_a=\int_Sa_aJ\vartheta\,\mathrm d\mu
     =\int_S\vartheta Ja_a\,\mathrm d\mu
     =\int_S\vartheta H_a\,\mathrm d\mu.
\]
Consequently,
\[
 Je_a=H_a-\int_SH_a\vartheta\,\mathrm{d}\mu.
\]
The restriction of $J$ to zero-average functions, with its range
projected to zero-average functions, is uniformly invertible near
the product sphere. Projecting the last equation to that range and
using elliptic estimates gives $e_a=O_{C^k}(d)$. Differentiating the
projected equation and the zero-average condition in its parameters
gives the same bound for all required parameter derivatives.

The divergence theorem on the slab between two leaves gives
\begin{equation}\label{cylinder:eq:normal-integral-variation}
 \partial_\sigma G_a=\int_S \vartheta\operatorname{div}_{\{F=\mathrm{const}\}}
                   Y_a^\top\,\mathrm{d}\mu=O_{C^{k-1}}(d).
\end{equation}
The remaining derivatives of $G_a$, including its $F$ derivative,
are uniformly bounded. Equations~\eqref{cylinder:eq:sigma-identity}--\eqref{cylinder:eq:normal-integral-variation}
therefore imply
\begin{equation}\label{cylinder:eq:Ysigma}
 Y_a\sigma=G_a(F,\sigma,t)+O_{C^{k-1}}(d),\qquad
 \partial_\sigma G_a=O_{C^{k-1}}(d).
\end{equation}

We take one more cyclic bracket and get
\begin{equation}
 Z_1:=\krot^{-1}[Y_2,Y_3],\quad
 Z_2:=\krot^{-1}[Y_3,Y_1],\quad
 Z_3:=\krot^{-1}[Y_1,Y_2].
\end{equation}

These fields are time-independent. Their metric Lie derivatives
and their pairings with $\mathrm{d}F$ are $O(d)$. To verify the remaining
pairing, for any $b,c$, we expand
\begin{align*}
 [Y_b,Y_c]\sigma=
 (\partial_FG_c)Y_bF-(\partial_FG_b)Y_cF +(\partial_\sigma G_c)Y_b\sigma
    -(\partial_\sigma G_b)Y_c\sigma+O_{C^{k-2}}(d).
\end{align*}
Each term is $O(d)$, by \eqref{cylinder:eq:Y-est} and \eqref{cylinder:eq:Ysigma}.
Using \eqref{cylinder:eq:sigma-identity} once again proves
\begin{equation}\label{cylinder:eq:Z-est}
 \Lie_{Z_a}g=O_{C^{k-3}}(d),\quad Z_aF=O_{C^{k-3}}(d),\quad
 \langle Z_a,\nu\rangle=O_{C^{k-3}}(d),\quad
 Z_a-K_a=O_{C^{k-3}}(C(T)\eta).
\end{equation}
There are exactly two brackets in the construction, on nested fixed
terminal domains. Their number does not increase with $T$ or $q$.

\medskip\noindent\textbf{Step 7. Recover one Gram matrix from constant brackets.}
The normalization must be the same on all leaves and times.
The following lemmas replace the variation argument in the final
steps of Brendle's proof by matching values and first derivatives, Killing transport,
and the algebra of the sphere rotations.

The following argument uses the time independence of $Z_a$ in an
essential way. Leafwise choices of an orthogonal frame or a Gram
normalization would not suffice for Definition~\ref{cylinder:def:symmetry}.

\begin{claim}
There are real constants $C_{ab}^{\ \ c}$, independent of space and
time, such that on a smaller fixed terminal domain
\begin{equation}\label{cylinder:eq:structure}
 [Z_a,Z_b]=C_{ab}^{\ \ c}Z_c+O_{C^k}(d),\qquad
 C_{ab}^{\ \ c}=\krot\epsilon_{abc}+O(C(T)\eta).
\end{equation}
Here and below we may lower $k$ by a fixed number from its value in
\eqref{cylinder:eq:Z-est}.
\end{claim}

We choose one point of one leaf at time $-1$ and a tangent orthonormal
pair $e_1,e_2$ there. On the span of $Z_1,Z_2,Z_3$, the map consisting
of the two tangential values and the one tangential skew derivative
\[
 Z\longmapsto
 \left(\langle Z,e_1\rangle,\langle Z,e_2\rangle,
 \tfrac12(\langle\nabla_{e_1}Z,e_2\rangle-
               \langle\nabla_{e_2}Z,e_1\rangle)\right)
\]
is invertible. This is true for the three round rotations and
persists under a small perturbation. We define $C_{ab}^{\ \ c}$ by
matching these three components of $[Z_a,Z_b]$.
For $E=[Z_a,Z_b]-C_{ab}^{\ \ c}Z_c$, its metric Lie derivative, $EF$, and
$E\sigma$ are $O_{C^k}(d)$. At the chosen point the two tangential
values and tangential skew derivative vanish. The two transverse
values are $O(d)$. Differentiating $EF,E\sigma$ bounds the normal
components of $\nabla E$ by $O(d)$ there, since $E=O(d)$.
The symmetric part of $\nabla E$ is controlled by its metric Lie derivative.
Thus $|E|+|\nabla E|=O(d)$ at that point.

The Killing transport identity expresses $\nabla^2E$ as
$\Rm*E$ plus fixed permutations of $\nabla\Lie_Eg$.
Along paths of bounded length, the resulting linear differential
inequality gives $|E|+|\nabla E|=O(d)$ throughout the connected
fixed spatial domain at time $-1$. Differentiation of the identity
gives the higher derivative bounds. Since $E$ and the coefficients
$C_{ab}^{\ \ c}$ are time-independent, these coordinate bounds hold
at every time in the fixed terminal interval. Uniform metric and
connection bounds convert them to the asserted geometric estimates. This finishes the proof of the claim.

\begin{claim}\label{cylinder:lem:Gram}
There is a single positive definite matrix $Q_*$, close to $I$, such
that on every final CMC leaf at every final time, we have
\begin{equation}
 \Area(S)^{-2}\int_S\langle Z_a,Z_b\rangle\,\mathrm{d}\mu
       =(Q_*)_{ab}+O(d).
\end{equation}
\end{claim}

Project $Z_a$ to each final CMC leaf. Equation
\eqref{cylinder:eq:Z-est} gives an $O_{C^8}(d)$ change on the
retained domains, and the projected fields have small intrinsic
metric Lie derivatives and the constant bracket coefficients in
\eqref{cylinder:eq:structure}. The final near-round family has
controlled geometry and area range. Lemma~\ref{lem:rotation-gram}
therefore gives the same matrix
$Q_*=(6\pi)^{-1}(\det M)M^{-1}$ on every leaf and time.
Since $M=\krot I+O(C(T)\eta)$, $Q_*$ is close to $I$.
This proves the claim.

\begin{claim}\label{cylinder:lem:normalized-brackets}
We set $\widehat Z:=Q_*^{-1/2}Z$. On the fixed terminal domain, we have
\[
 [\widehat Z_a,\widehat Z_b]
       =\krot\epsilon_{abc}\widehat Z_c+O_{C^k}(d).
\]
Consequently $\mathcal U_a=\widehat Z_a/\krot$ obey
$[\mathcal U_a,\mathcal U_b]=\epsilon_{abc}\mathcal U_c+O_{C^k}(d)$.
The coefficients used in this normalization are constant in space
and time. The same assertion holds throughout a fixed cap whenever
its constant structure-coefficient estimate holds there and its
common Gram estimate holds on the collar.
\end{claim}

We set $\mathsf M:=Q_*^{-1/2}$ and select one collar or neck leaf at one time.
In its round parametrization the projected fields have coefficient
matrix $P$ supplied by the Killing projection in
Lemma~\ref{lem:rotation-gram}. Since
$PP^{\mathsf T}=6\pi Q_*+O(d)$, we get
\[
 (\mathsf M P)(\mathsf M P)^{\mathsf T}=6\pi I+O(d).
\]
All matrices and their inverses are uniformly bounded, and
$\det(\mathsf M P)>0$. The polar decomposition therefore yields
$\mathsf M P=\krot O+O(d)$ with $O\in SO(3)$. An $SO(3)$ change of
basis preserves the standard cyclic brackets of the $J_i$.
Thus on this leaf the projected $\widehat Z_a$ satisfy the
claimed bracket relation up to $O(d)$.

On the other hand the original constant coefficient estimate,
under this constant change of basis, gives on the whole domain
\[
 [\widehat Z_a,\widehat Z_b]
 =\widehat C_{ab}^{\ \ c}\widehat Z_c+O_{C^k}(d),\qquad
 \widehat C_{ab}^{\ \ c}
 =\mathsf M_{ai}\mathsf M_{bj}C_{ij}^{\ \ m}(\mathsf M^{-1})_{mc}.
\]
Comparison on the chosen whole leaf and its positive definite
$L^2$ Gram matrix imply
$\widehat C_{ab}^{\ \ c}=\krot\epsilon_{abc}+O(d)$.
These are estimates for constant numbers. Substituting them back
into the full-domain relation, with the bounded field derivatives,
proves the conclusion on every leaf and time, and also in the
cap interior. Division by $\krot$ gives the second assertion.
The argument uses a polar factor on one leaf, which determines the
constant error. Thus no smooth choice of polar factors across the
family of leaves is required. This finishes the proof of the claim.

\medskip\noindent\textbf{Step 8. We normalize the constructed triple and choose the constants.}

We set $\widehat Z_a:=(Q_*^{-1/2})_{ab}Z_b$. This is one constant
change of basis, so the fields remain time-independent. By
\eqref{cylinder:eq:Z-est} and Claim~\ref{cylinder:lem:Gram}, on the entire required
final neighborhood their metric Lie derivative through order two, their two
transverse pairings, and their area-squared normalized Gram error are all
bounded by $Cd$. The unit potential normal in
Definition~\ref{cylinder:def:symmetry} is
$n_f=\nabla F/|\nabla F|$. Thus it follows that
\[
 |\langle\widehat Z_a,n_f\rangle|
   =\frac{|\widehat Z_aF|}{\sqrt{1-\lambda^2\mathrm{R}}}
   \le (1+C\eta)\,|\widehat Z_aF|\le Cd.
\]
Adding the time function in the heat gauge changes neither the
potential levels nor this normal.

The fixed terminal domains were chosen large enough to contain the
whole ball and backward interval required in Definition~\ref{cylinder:def:symmetry}, with
room for each fixed interior restriction above. The metric is
uniformly comparable with the model there. Thus, after including the
fixed scale factors and the finite square sums in that definition,
the new symmetry parameter is bounded by
\[
 \epsilon_{\rm new}\le\big(C_qT^{-1/2}+C(T)\eta\big)\epsilon.
\]
We use the fixed heat data exponent $q=q_h=11$. We choose $T$
large enough that $C_qT^{-1/2}<1/4$. Set $L:=1000T$ and choose $\epsilon_1$ after $T$ according to
\eqref{cylinder:eq:improvement-domain-threshold} and all remaining
finite smallness conditions. Restrict the
$\epsilon_1$-cylindrical chart using
\eqref{cylinder:eq:fixed-domain-threshold}; its comparison error is
at most $\eta:=\epsilon_1$, and $\lfloor\epsilon_1^{-1}\rfloor\ge120$ supplies the
entire derivative reserve. In particular $C(T)\eta<1/4$.
The normalized upper time $-\epsilon_1$ is later than $-1$, so every
terminal domain used above has a positive future margin as well.
These choices are independent of the rotation error
$\epsilon\in(0,\epsilon_1]$.
It follows that $\epsilon_{\rm new}\le\epsilon/2$, as required. If
$\epsilon=0$, the hypothesis at the distinguished point already gives
the conclusion. The remaining regularity and domain estimates verify that each operation
has the stated regularity and fits within the prescribed neighborhood.

\medskip\noindent\textbf{Step 9. Verify the derivative orders and domain margins.}

For a spacetime tensor $\mathsf T$, write
\[
 \|\mathsf T\|_{\mathrm{par},m;Q}
 :=\sum_{j+2k\le m}\sup_Q|\nabla^j\partial_t^k\mathsf T|.
\]
One time derivative thus costs two spatial orders. The comparison and heat equation estimates are
\begin{gather*}
 \|g-\bar g\|_{\mathrm{par},120}\le C(T)\eta,\qquad
 \sum_i\|\Lie_{U_{a,i}}g\|_{C^2}\le C\epsilon
 \ \Longrightarrow\
 \sum_i\|U_{a,i}-(P^{ab}U_b)_i\|_{C^3}\le C\epsilon,\\
 q_U=10,\qquad q_h=q_U+1=11,\qquad s=-t,\\
 |V-U|\le C\epsilon\sqrt{s}(2T/s)^{11},\qquad
 |h|+s^{-1/2}|u|\le C\epsilon(2T/s)^{11}.
\end{gather*}
Only the $C^3$ overlap estimate is needed to glue the given fields;
no higher derivative of a local rotation triple is assumed. We fix $q_h=11$
before choosing $T$.

On successively smaller interior domains, parabolic regularity
and the exact-model comparison give
\begin{align*}
 \|h_a\|_{\mathrm{par},64}+\|u_a\|_{\mathrm{par},64}
     &\le C(T)\epsilon,\\
 \|V_a-K_a\|_{\mathrm{par},64}
     &\le C(T)\eta,\\
 \|\text{model comparison inhomogeneous term}\|_{\mathrm{par},60}
     &\le C(T)\eta\epsilon,\\
 \|h_a-\bar h_a\|_{\mathrm{par},60}
       +\|u_a-\bar u_a\|_{\mathrm{par},60}
     &\le C(T)\eta\epsilon,\\
 \|E_a\|_{\mathrm{par},60}
       +\|(I-\Pi_0)\bar u_a\|_{\mathrm{par},60}
     &\le C_qT^{-1/2}\epsilon.
\end{align*}
Here $E_a$ is the decaying tensor remainder in
\eqref{cylinder:eq:decomp}. We set
$d:=(C_qT^{-1/2}+C(T)\eta)\epsilon$ as in
\eqref{cylinder:eq:delta}. After the first bracket and freezing,
\[
 \|\Lie_{Y_a}g\|_{C^{56}}+\|Y_aF\|_{C^{56}}\le Cd.
\]
The projected CMC equation retains spatial order $52$ for the
normal component error $e_a=a_a-G_a\vartheta$, including the required
parameter derivatives. The second bracket then gives
\[
 \|\Lie_{Z_a}g\|_{C^{50}}
   +\|Z_aF\|_{C^{50}}+\|Z_a\sigma\|_{C^{50}}\le Cd.
\]
Finally the structure coefficients and common Gram normalization
retain order $40$ for the metric Lie derivative and both transverse errors of
$\widehat Z_a$, and order $39$ for their bracket error, all with
bound $Cd$ on smaller fixed domains.
The inhomogeneous term uses at most $60+2=62<64$ orders of $h,u$.
Angular eigencomponent decay gives the zeroth-order gain first; interior
regularity for the homogeneous model equations gives the higher
orders on smaller terminal domains with the same gain. The two
brackets, freezing, and CMC variation fit within the displayed
derivative requirements. In particular, we have
\[
 (2,0,0)_{\text{symmetry functional}},\qquad
 (2,2,0)_{\text{global gluing}}
 \quad<\quad (40,40,39)_{\text{retained}}.
\]
The Gram estimate uses the leafwise round metric itself and does
not differentiate its choice.

We choose the terminal domains as follows. We keep
$H=\tau_*=10^4$. Obtain the initial terminal bounds on the larger
fixed domain $D_{2H}\times[-2\tau_*,-1]$, which lies in the heat
comparison domain after increasing $T$. For $j=0,\ldots,10$ put
\[
 Q_j:=D_{H(1-j/20)}\times[-\tau_*(1-j/20),-1].
\]
The margins and final containment are explicit:
\[
 R_j-R_{j+1}=H/20=500,\qquad
 B_j-B_{j+1}=\tau_*/20=500,\qquad
 Q_{10}=D_{5000}\times[-5000,-1],
\]
where $R_j=H(1-j/20)$ and $B_j=\tau_*(1-j/20)$.
Thus estimates on $Q_j$ control $Q_{j+1}$ with fixed spatial
and earlier-time margins. Assign these
restrictions to spectral regularity, freezing, CMC construction,
CMC parameter estimates, second brackets, matching values and first derivatives,
Killing transport, uniformization, Killing projection, and final
normalization, in that order. The last domain still has flat
radius $5000$ and backward endpoint $-5000$, containing the
radius $100$ and backward endpoint $-101$ required by the rotation triple
with ample metric comparison margins. Whole sphere factors remain
intact. All constants for these terminal operations are fixed
before $T$; the geometric tolerance is chosen last so that the
geometry belongs to their common fixed bounded neighborhood.
This proves the stated conclusion.
\end{proof}

\begin{cor}\label{cor:intrinsicneckcontraction}
Let $g(t)=\Phi_t^*g$ be the normalized canonical flow of a soliton
satisfying the hypotheses of Theorem~\ref{thm:main}. Fix $\rho\ge100$, $\ell_{\rm neck}\ge10\sqrt{\rho+1}$ and
$0<\theta<1$.
For a spacetime center $(y,t)$ put $r_{y,t}:=\mathrm R(y,t)^{-1/2}$
and define
\begin{align*}
 \mathcal P_\rho(y,t)
 &:=B_{g(t)}(y,\rho r_{y,t})
                  \times[t-\rho r_{y,t}^2,t],\\
 \mathcal P_\rho^+(y,t)
 &:=B_{g(t)}(y,(\rho+\ell_{\rm neck})r_{y,t})
                  \times[t-(\rho+\ell_{\rm neck})r_{y,t}^2,t].
\end{align*}
Define $e_N(y,t)$ as the infimum of the numbers $e\ge0$ for which
a smooth time-independent triple satisfies all three inequalities
of Definition~\ref{cylinder:def:symmetry}, with $\epsilon=e$,
$r=r_{y,t}$, and $\mathcal P$ replaced by
$\mathcal P_\rho(y,t)$. The Gram inequality is measured on every
whole CMC leaf contained in this spatial ball at each time.
The fields are defined on an open neighborhood of its closure. There exist $L<\infty$, a cylindrical threshold
$\eta_{\rm neck}>0$, and a rotation-error threshold
$\epsilon_{\rm neck}>0$, depending only on $\rho,\ell_{\rm neck},\theta$
and the fixed constants of the cylindrical argument, with the
following property.

 For $(x,s)$, we put $r:=r_{x,s}$ and set
\begin{equation}\label{cylinder:eq:neck-input-set}
 \mathcal I(x,s):=
 \overline B_{g(s)}(x,Lr)\times[s-Lr^2,s].
\end{equation}
Assume $(x,s)$ is $\eta_{\rm neck}$-cylindrical. On a neighborhood
of the union of $\mathcal P_\rho^+(y,t)$ over
$(y,t)\in\mathcal I(x,s)$, we use one smooth family of CMC two-spheres
in regular potential levels as in Remark~\ref{cylinder:rem:cmc-existence}. We also require that each leaf
meeting a domain extends to a whole sphere in its enlarged
domain. If two domains overlap, their foliations are
restrictions of this same family: the leaves are the same subsets,
although their parameter labels may differ. All normals and Gram
matrices in the definition of $e_N$ refer to this family.
If
\[
 e_{\rm in}:=\sup_{(y,t)\in\mathcal I(x,s)}e_N(y,t)
                     \le\epsilon_{\rm neck},
\]
then
\begin{equation}\label{cylinder:eq:intrinsic-neck-contraction}
 e_N(x,s)\le\theta e_{\rm in}.
\end{equation}
\end{cor}

The comparison sizes and thresholds may be chosen as follows. Start with
$H=\tau_*:=10^4(1+\rho+\ell_{\rm neck})^2$ in the preceding proof and retain
$q=q_h=11$. Choose $T$, with $2T=4^J$, sufficiently large that
\[
 T\ge100(H^2+\tau_*),\qquad
 C_{\rm eq}C_qT^{-1/2}\le\theta/4,
\]
where $C_{\rm eq}$ bounds the conversion from the terminal
$C^2$ estimates for $\Lie_Ug$, $UF$, and the CMC-normal component
to the three errors in Definition~\ref{cylinder:def:symmetry} on
$\mathcal P_\rho$. It includes the fixed factors $r_{y,t}^j$,
$r_{y,t}^{-1}$, and $\operatorname{Area}(S)^{-2}$, the finite sums
over the three fields, the conversion $UF/|\nabla F|$, and the
comparison of actual and model metric norms on the fixed enlarged
domain. Thus $C_{\rm eq}$ depends on $\rho,\ell_{\rm neck}$ and
these fixed margins, but not on $T$ or the rotation error. Set $L:=1000T$ and
\begin{align*}
 R_*&:=4\bigl[L+(\rho+\ell_{\rm neck}+1)\sqrt{2(L+1)}
                         +\pi\sqrt{2(L+1)}\bigr],\\
 T_*&:=4(L+1)(\rho+\ell_{\rm neck}+1).
\end{align*}
One may take $\epsilon_{\rm neck}:=\eta_{\rm neck}/2$, where
\begin{equation}\label{cylinder:eq:neck-threshold}
 0<\eta_{\rm neck}\le\min\left\{
 \eta_{\rm geom},\frac1{121},\frac1{10^6T},
 \frac1{2(R_*+\pi\sqrt2+1)},\frac1{2(T_*+1)},
 \frac{\theta}{4C_{\rm eq}C(T)}\right\}.
\end{equation}
Here $\eta_{\rm geom}>0$ is a fixed smallness threshold for the
potential, CMC, alignment, and metric comparison estimates on these
specified domains. All constants in this choice are fixed before
$\eta_{\rm neck}$ and independently of $e_{\rm in}$.

\begin{proof}
For simplicity, we translate and rescale the base point to $(x,-1)$ with scalar curvature
one. We repeat the proof of Theorem~\ref{cylinder:improvement},
keeping the measurement size $\rho$ and margin $\ell_{\rm neck}$ fixed.
The larger terminal values $H,\tau_*$ ensure that the final domain
contains $\mathcal P_\rho^+(x,-1)$ with all ten interior restrictions
used in Step~9. The required centers in the finite parabolic mesh of
Lemma~\ref{cylinder:lem:gluing} lie in
\eqref{cylinder:eq:neck-input-set}. Their domains contain
the expanded mesh boxes because $\rho\ge100$; the dyadic-box
constant is chosen sufficiently small for the fixed margins.

We verify that the larger geometric chart contains the domains
appearing in the statement. Decrease $\eta_{\rm geom}$ so that
on the normalized comparison region
\[
 \tfrac12\le(-t)\mathrm R(y,t)\le2.
\]
For $-L-1\le t\le-1$ this implies
$r_{y,t}\le\sqrt{2(L+1)}$.
The enlarged comparison balls therefore have spatial extent at most
$L+(\rho+\ell_{\rm neck})\sqrt{2(L+1)}$ in the final-time metric, and their
earlier endpoints are no earlier than
$-(L+1)-2(\rho+\ell_{\rm neck})(L+1)>-T_*$.
Here $g(t)\ge g(-1)$ for $t\le-1$, since $\Ric\ge0$.
The fixed metric comparison and first-exit margins place all these
balls inside $D_{R_*}$. The bounds
\eqref{cylinder:eq:neck-threshold} and the restriction rule
\eqref{cylinder:eq:fixed-domain-threshold} supply this entire
geometric region, including all required derivatives through order
$120$. No rotation triple is required at an additional center outside
$\mathcal I(x,-1)$ merely to provide these margins.

Lemma~\ref{cylinder:lem:heat-gauge} derives the potential comparison
from metric closeness and the soliton identities. The prescribed
CMC family has the small-graph and parameter estimates from the
projected CMC equation. On overlaps the same whole spheres are used
for the two-normal bounds and the Gram integrals, so
Lemma~\ref{cylinder:lem:alignment} applies to any two selected local triples. No derivative of a normal component of a local field is assumed.
Only the metric Lie derivative through order two is needed to
obtain the $C^3$ overlap comparison and glue the fields. The heat
argument, two brackets, and constant Gram normalization then give,
for local triples of error at most $e'\le\eta_{\rm neck}$,
\[
 e_{\rm out}\le
 C_{\rm eq}\bigl(C_qT^{-1/2}+C(T)\eta_{\rm neck}\bigr)e'
 \le\tfrac\theta2 e'.
\]
The constants may depend on the fixed measurement sizes and
margins. Their choices above do not depend on the symmetry error.
Every comparison error is multiplied by $e'$; there is no additive
geometric error.

Finally let $e=e_{\rm in}\le\eta_{\rm neck}/2$. For each
$0<\zeta<\eta_{\rm neck}/2$, choose at every center of the finite
mesh a triple with error less than $e_N(y,t)+\zeta\le e+\zeta$.
Such triples exist by the definition of the infimum. The construction
produces an admissible time-independent constructed triple on the
prescribed domain, with error at most
$\tfrac\theta2(e+\zeta)$. Hence
\[
 e_N(x,-1)\le\tfrac\theta2(e+\zeta).
\]
Letting $\zeta\downarrow0$ proves
\eqref{cylinder:eq:intrinsic-neck-contraction}, including $e=0$;
no attainment of an infimum is needed. Undoing the parabolic
normalization gives the assertion at $(x,s)$.
\end{proof}

\section{Symmetry improvement on Bryant-product caps}\label{sec:capimprovement}
In this section, we improve symmetry on Bryant-product caps.
Proposition~\ref{prop:productgrid} aligns the initial fields,
Proposition~\ref{cap:prop:product-quotient} gives the tensor estimate,
and Proposition~\ref{prop:heatcap} proves heat improvement.
The main result, Corollary~\ref{cor:alltiltheat},
expresses the resulting estimates uniformly for $a_{\rm ax}\in[0,1]$.

We first recall the Bryant normalization and flow from
Section~\ref{sec:collars}. On the Bryant soliton
$(B^3,g_B)$ with tip $o$, put
\[
 s(x):=d_{g_B}(o,x),\qquad
 g_B=\mathrm ds^2+\omega(s)^2g_{S^2},\qquad
 B_0(x):=B(s(x)),\qquad B(0)=0.
\]
The Bryant diffeomorphisms and canonical metric are
\begin{equation}
 \frac{\mathrm d}{\mathrm dt}\Psi_t(x)
   =-\nabla^{g_B}B_0(\Psi_t(x)),\qquad
 \Psi_0=\operatorname{Id},\qquad g_B(t)=\Psi_t^*g_B.
\end{equation}
Extend $\Psi_t$ by the identity on the product line. In its fixed
spatial coordinates define
\begin{equation}\label{eq:cap-current-radial}
 s_t:=s\circ\Psi_t=d_{g_B(t)}(o,\cdot),\qquad
 B_t:=B_0\circ\Psi_t=B(s_t),\qquad
 \bar g(t):=g_B(t)+\mathrm dz^2.
\end{equation}
Thus $s$ and $B_0$ are time-zero functions, while $s_t$ and $B_t$
are the time-$t$ distance and potential. The time derivatives at fixed
spatial points satisfy
\begin{equation}
 \partial_tB_t
 =-|\nabla^{g_B}B_0|_{g_B}^2\circ\Psi_t
 =-|\nabla^{g_B(t)}B_t|_{g_B(t)}^2\in[-1,0],\qquad
 0\le B_t-B_{t_*}\le t_*-t\quad(t\le t_*).
\end{equation}
In particular $B_{t_*}$ is a frozen spatial function when the final
time $t_*$ is fixed. We use $B_t$ for all evolving radial weights
and domains below. The line coordinate is unchanged by this flow:
$z\circ\Psi_t=z$. Where the heat argument uses $Z$, it denotes the
same line arclength with the specified centering, $Z=z-z_*$.

For constants $\lambda,a_{\rm ax}\ge0$ with $\lambda^2+a_{\rm ax}^2=1$, the model
heat-normalized potential is
\begin{equation}
 F_{\rm mod}(x,z,t):=\lambda B_t(x)+a_{\rm ax}z+\lambda t.
\end{equation}
Its exact identities are
\[
 \Hess_{\bar g(t)}F_{\rm mod}=\lambda\Ric_{\bar g(t)},\qquad
 |\nabla F_{\rm mod}|_{\bar g(t)}^2=1-\lambda^2\mathrm R_{\bar g(t)},\qquad
 \partial_tF_{\rm mod}=\Delta_{\bar g(t)}F_{\rm mod}
                      =\lambda\mathrm R_{\bar g(t)}.
\]
For a comparison with an actual flow, identify the chart with the
model by a time-independent spatial comparison map and denote the
pulled-back actual data by $g(t),F$. The function $F$ is in the heat
gauge \eqref{cylinder:eq:F-identities}; its additive constant is chosen
for comparison with $F_{\rm mod}$. Model coordinates such as $B_t,Z$
on this chart retain their definitions in
\eqref{eq:cap-current-radial}; they are not actual distance functions.

For an enlarged comparison region $\mathcal Q^+$, define
\[
 \|T\|_{\mathrm{par},m;\mathcal Q^+}
 :=\sum_{j+2k\le m}\sup_{\mathcal Q^+}
       |\bar\nabla^j\partial_t^kT|_{\bar g(t)}.
\]
Here $\bar\nabla$ is the spatial model connection. A comparison
error $\eta$ means, at the derivative order specified in the result,
\begin{equation}\label{eq:cap-comparison-error}
 \|g-\bar g\|_{\mathrm{par},m;\mathcal Q^+}
 +\|F-F_{\rm mod}\|_{\mathrm{par},m;\mathcal Q^+}\le\eta;
\end{equation}
the heat estimates use $m=60$.
For the cap heat size parameter $E\ge2\Lambda$, define
\begin{equation}\label{eq:cap-heat-enlargement}
 \begin{aligned}
 Q_E&:=\{B_t<E,\ |Z|<E^2,\ -E^2<t\le0\},\\
 Q_E^+&:=\{B_t<2E,\ |Z|<2E^2,\ -2E^2<t\le0\},\\
 \widehat Q_E&:=\{B_t<2E+1,\ |Z|<2E^2+1,\ -2E^2-1<t\le0\}.
 \end{aligned}
\end{equation}
 All are sets in the fixed
product chart, with moving radial coordinate $B_t$; the tip $B_t=0$
is interior. All final-time boundaries
are relative to $t\le0$.
A spatial $C^k$ norm means
$\sup_t\sum_{j=0}^k\sup_{\mathcal Q^+(t)}|\nabla^jT|_{g(t)}$.
The parameters $\epsilon,e$ measure vector field errors, such as
$\Lie_Ug$, $UF$, normal pairings, and the Gram error
$|Q_S(U)-I|$, with $Q_S$ defined in
Definition~\ref{cylinder:def:symmetry}. On regular potential levels, we have
\[
 n_f=\frac{\nabla F}{|\nabla F|},\qquad
 \langle U,n_f\rangle=\frac{UF}{|\nabla F|};
\]
$\nu$ denotes the CMC-sphere normal inside that level. Constants denoted by $C_E$ may depend on the enlarged domain size
$E$; constants without that subscript are independent of $E$ once
the fixed domains, derivative orders, and margins are fixed.

We first establish the overlap alignment estimate needed to glue
local triples in the product chart. The weighted gluing construction
will be given in Proposition~\ref{prop:productgrid} below.

\begin{prop}\label{prop:loworderalignment}
Fix a length $r>0$, connected spatial domains
$\Omega\Subset\Omega^+$, a compact time interval $I$, and positive
constants $d_0,r_0,D_0,T_0,K$. Assume, for every $t\in I$,
\[
 \operatorname{dist}_{g(t)}(\Omega,\partial\Omega^+)\ge d_0r,
 \qquad |I|\le T_0r^2.
\]
In fixed spatial charts of radius comparable to $r$, assume that the
coefficients of $r^{-2}g(t)$ and their inverses have $C^6$ bounds
at most $K$ when the coordinate variables are divided by $r$.
These charts and the metrics at different times are uniformly
comparable. Fix $t_*\in I$ and a whole collar sphere
$S\subset\Omega$. For every $t\in I$ and $x\in\Omega$, assume
there is a path $\gamma\subset\Omega^+$ from $S$ to $x$ of
$g(t)$-length at most $D_0r$. More precisely, $\gamma$ is covered by
at most $K$ balls $B_{g(t)}(x_\ell,r_0r)$, successive balls meet,
and $B_{g(t)}(x_\ell,2r_0r)\Subset\Omega^+$ for every center.

Let $\mathcal C\Subset\Omega^+$ be an open collar containing $S$,
parametrized by whole potential level CMC spheres and two transverse parameters on
fixed intervals after division of lengths by $r$. Its parametrizations
and inverse charts have rescaled $C^4$ bounds at most $K$.
Let $(n_f,\nu)$ be their orthonormal normal frame.
Fix a compact family of connected model neck or Bryant-product cap
regions with the above domains and margins. Their regular collars
consist of whole potential level spheres.
After the parabolic rescaling
$\widehat g(\tau)=r^{-2}g(t_*+r^2\tau)$, assume that the geometry
on a fixed enlargement of $\Omega^+\times I$ is $\eta$-close to
one of the above models: the metric comparison is spatial $C^8$,
uniform in rescaled time, and the collar parametrizations and
unit normal frames are $C^4$-close in the same fixed charts.
The model family may include the
zero-axial-coefficient Bryant product.

Let $U=(U_1,U_2,U_3)$ and $V=(V_1,V_2,V_3)$ be smooth
time-independent triples on a neighborhood of
$\overline{\Omega^+}$. For a vector field $X$, a covariant
$2$-tensor $h$, and a scalar function $u$, write
\[
 \begin{aligned}
 \|X\|_{k,r;D}&:=\sup_{t\in I}\sum_{j=0}^k
       r^{j-1}\sup_D|\nabla^jX|_{g(t)},\\
 \|h\|_{k,r;D}&:=\sup_{t\in I}\sum_{j=0}^k
       r^j\sup_D|\nabla^jh|_{g(t)},\\
 \|u\|_{k,r;D}&:=\sup_{t\in I}\sum_{j=0}^k
       r^j\sup_D|\nabla^ju|_{g(t)}.
 \end{aligned}
\]
All derivatives here are spatial. Suppose that, for some $e\ge0$, the following estimates hold:
\[
 \sum_{i=1}^3\bigl(\|\Lie_{U_i}g\|_{2,r;\Omega^+}
                 +\|\Lie_{V_i}g\|_{2,r;\Omega^+}\bigr)\le e,
\]
\[
 r^{-1}\sum_{i=1}^3\sum_{n\in\{n_f,\nu\}}
 \bigl(\|\langle U_i,n\rangle\|_{0,r;\mathcal C}
       +\|\langle V_i,n\rangle\|_{0,r;\mathcal C}\bigr)\le e,
\]
\[
 |Q_S(U)-I|_{\mathrm F}+|Q_S(V)-I|_{\mathrm F}\le e,
 \qquad\mathrm{where}\qquad
 Q_S(U)_{ij}:=\Area_{g(t_*)}(S)^{-2}
       \int_S\langle U_i,U_j\rangle_{g(t_*)}\,\mathrm d\mu_{g(t_*)}.
\]
\begin{enumerate}[label=\textup{(\roman*)},ref=\textup{(\roman*)}]
\item There are constants $\eta_0,e_0,C>0$, depending only on the fixed
model family, domain data, and geometric bounds, such that
$0\le\eta\le\eta_0$ and $0\le e\le e_0$ imply the existence of a
constant matrix $P\in O(3)$ satisfying
\[
 \sum_{i=1}^3\left\|U_i-\sum_{j=1}^3P_{ij}V_j
                   \right\|_{3,r;\Omega}\le Ce.
\]
\item For three triples $U_d=(U_{d,1},U_{d,2},U_{d,3})$,
$d\in\{a,b,c\}$, satisfying these hypotheses on common domains
and a common whole sphere, choose the pairwise comparison matrices
$P^{ab}$, acting by $(P^{ab}U_b)_i=\sum_jP^{ab}_{ij}U_{b,j}$,
with $P^{ba}=(P^{ab})^{-1}$. Then
\[
 |P^{ab}P^{bc}P^{ca}-I|_{\mathrm F}\le Ce.
\]
\end{enumerate}
\end{prop}
\begin{proof}
We first prove (i). Set $\widehat g(\tau)=r^{-2}g(t_*+r^2\tau)$ and leave the vector fields unchanged. The normalized Gram matrices are unchanged, and the stated norms become the corresponding norms with $r=1$. We work at $\tau=0$ with this rescaled metric until the final time-comparison step and suppress the hat.
We first verify the model Killing field property used below. On a
regular collar of either a cylindrical model or a Bryant-product
model, write
\[
 \bar g=g_{\mathrm{base}}(y)+\rho(y)^2g_{\mathbb S^2},
 \qquad y=(y^1,y^2).
\]
For the cylinder, the base is the two-dimensional flat factor and
$\rho$ is spatially constant; for the Bryant product, one can take
$y=(s,z)$, where $s$ is the Bryant radial distance. The sphere fibers
are the model collar leaves. Their two unit normals
$(\bar n_f,\bar\nu)$ span the transverse plane, including when the
Bryant product has zero axial coefficient. Consequently a Killing
field $X$ whose two normal pairings vanish throughout this collar
has no transverse component. The sphere--sphere components of
$\Lie_X\bar g=0$ imply
\[
 X(y,\theta)=\sum_{i=1}^3c_i(y)J_i(\theta),
\]
where $J_1,J_2,J_3$ are a fixed basis of the Killing fields of the
unit round sphere. The mixed components, with angular indices
$\mu,\nu$, give
\[
 0=(\Lie_X\bar g)_{\alpha\mu}
   =\rho(y)^2(g_{\mathbb S^2})_{\mu\nu}
        \partial_{y^\alpha}X^\nu,
 \qquad \alpha=1,2.
\]
Thus $\sum_i(\partial_{y^\alpha}c_i)J_i=0$ and
$\partial_{y^\alpha}c_i=0$ for every $i,\alpha$. The coefficients
are constant on the connected collar. The difference between $X$
and this constant linear combination of the model rotations
vanishes on an open set. The Killing field identity
$\nabla^2X=\Rm*X$, applied to that difference along paths,
extends the equality over the connected model region, including
its tip when present. We have therefore proved, for both models,
\[
 \{X:\Lie_X\bar g=0,\ \langle X,\bar n_f\rangle
                  =\langle X,\bar\nu\rangle=0
                  \text{ on the open model collar}\}
       =\operatorname{span}\{J_1,J_2,J_3\}.
\]

We next derive the estimates used to compare the triples.
If $h=\Lie_{\mathcal W}g$, the differentiated metric Lie derivative identity gives
\[
 \nabla^2{\mathcal W}=\Rm*{\mathcal W}+\mathcal A(\nabla h),
\]
where $\mathcal A$ is a fixed linear combination of permuted indices.
An interior ball, Taylor's formula along geodesics, and this identity
bound the value and first covariant derivative at its center by the $C^0$ norm of ${\mathcal W}$ on the
ball and the $C^1$ norm of $h$. Along a unit-speed connecting path $\gamma$, this identity gives
\[
 \frac{\mathrm d}{\mathrm d\ell}
 (|\mathcal W|+|\nabla\mathcal W|)(\gamma(\ell))
 \le C(|\mathcal W|+|\nabla\mathcal W|)(\gamma(\ell))
       +C|\nabla h|(\gamma(\ell))
\]
in the upper-derivative sense. Gronwall's inequality and one further
spatial differentiation therefore give
\[
 \|\mathcal W\|_{C^3(\Omega)}
 \le C\bigl(\|\mathcal W\|_{C^0(B_0)}
                    +\|h\|_{C^2(\Omega^+)}\bigr),
\]
where $B_0$ is a fixed interior ball meeting the common collar.
The same estimate holds on the intermediate domains containing
$\overline\Omega\cup\overline{\mathcal C}$. Constants depend only
on the domain, chart, ball, and path bounds in the statement.

The Gram and transverse conditions first imply a uniform bound for
each rotation triple on a smaller connected domain containing the
common sphere. Otherwise normalize a field with unbounded supremum
there. The transport bounds give a subsequential, nonzero limiting
model Killing field. Its normal components vanish throughout the
open collar, so it is a rotation. Its $L^2$ norm on the common sphere
is zero, since the unnormalized Gram matrices are bounded. This is
impossible for a nonzero rotation. Propagation by the same transport
system gives the bound on the rest of the overlap.

We write $\langle\cdot,\cdot\rangle_S$ for the area-squared normalized
inner product on the common sphere. Let
$W=(W_1,W_2,W_3)$ satisfy the small Lie derivative, transverse,
and Gram bounds of the statement. For every smooth field
$\mathcal W$, the following uniform estimate holds on $\Omega$:
\begin{equation}\label{eq:lowcoercivity}
 \|{\mathcal W}\|_{C^3(\Omega)}\le C\left(
 \|\Lie_{\mathcal W}g\|_{C^2(\Omega^+)}
 +\sum_{\alpha=1}^2\|\langle {\mathcal W},n_\alpha\rangle\|_{C^0(\mathcal C)}
 +\sum_{j=1}^3|\langle {\mathcal W},W_j\rangle_S|\right).
\end{equation}
Here $n_1=n_f$, $n_2=\nu$, and all terms are evaluated at
the selected time. The same estimate holds on fixed intermediate
domains containing $\overline\Omega\cup\overline{\mathcal C}$
with the prescribed interior margins.
To justify the estimate uniformly near the model family, first prove
its $C^0$ version by contradiction, normalizing that norm to one.
The transport bounds give compactness on the smaller domain and
propagate it to the larger interior sets needed for the argument.
A failing sequence approaching the model family has a nonzero
limiting Killing field tangent to all collar spheres. The limiting
$W_j$ form a basis of rotations, since their Gram matrices tend to
$I$. Orthogonality to all three on $S$ forces the limiting field to
vanish. This proves the $C^0$ estimate for sufficiently small common
geometric and rotation triple thresholds. The transport identity then gives
\eqref{eq:lowcoercivity}. A quantitative relation between those two
small parameters was not used.

We choose $P\in O(3)$ minimizing
$\sum_{i=1}^3\|U_i-\sum_jP_{ij}V_j\|_{L^2(S)}^2$, with the same area normalization.
Set $W_i:=\sum_jP_{ij}V_j$, $\mathcal W_i:=U_i-W_i$,
and $D_{ij}:=\langle\mathcal W_i,W_j\rangle_S$.
In the following triple estimates, $W=(W_1,W_2,W_3)$,
$\mathcal W=(\mathcal W_1,\mathcal W_2,\mathcal W_3)$,
and norms mean sums over the three components.
Variation through orthogonal matrices gives $D=D^{\mathsf T}$.
Writing $Q_U=Q_S(U)$, $Q_W=Q_S(W)$, and
$Q_{\mathcal W}=Q_S(\mathcal W)$, the Gram identity gives
\[
 Q_U-Q_W=2D+Q_{\mathcal W},
 \qquad |D|\le Ce+C\|{\mathcal W}\|_{C^0}^2.
\]
Applying \eqref{eq:lowcoercivity} row by row therefore yields
\[
 \|{\mathcal W}\|_{C^3}\le Ce+C\|{\mathcal W}\|_{C^0}^2.
\]
The minimizing $P$ makes $\|{\mathcal W}\|_{C^0}=o(1)$ as the geometric and
rotation triple thresholds tend to zero: any limiting triples are orthonormal
bases of the same rotational algebra, and minimization identifies
them on $S$, hence everywhere. To make the absorption explicit, choose the common threshold so
that $C\|{\mathcal W}\|_{C^0}\le1/2$. Since
$\|{\mathcal W}\|_{C^0}\le\|{\mathcal W}\|_{C^3}$, the preceding inequality gives
\[
 \|{\mathcal W}\|_{C^3}\le Ce+\tfrac12\|{\mathcal W}\|_{C^3},
 \qquad \|U-PV\|_{C^3}\le2Ce.
\]
This absorption requires only separate upper bounds for
$\eta$ and $e$.

Since the fields and
$P$ are time-independent, controlled comparison of metrics and
connections transfers the spatial estimate to the bounded overlap
time interval. Returning to $g(t)$ gives the stated scale-invariant norms and constants independent of $r$. This proves (i).

 We now prove (ii). On a triple overlap, composing the three comparisons gives
\[
 \sum_{i=1}^3\left\|\sum_j
 (P^{ab}P^{bc}P^{ca}-I)_{ij}U_{a,j}\right\|_{C^0(S)}\le Ce.
\]
The normalized Gram matrix $Q_S(U_a)$ has eigenvalues bounded
below by a fixed positive constant. Taking the normalized $L^2(S)$
norm therefore yields $|P^{ab}P^{bc}P^{ca}-I|_{\mathrm F}\le Ce$.  This proves (ii).
\end{proof}

In the next proposition, for an actual annular CMC sphere $S$ at time $t$, define its
radial and axial labels by
\[
 v(S,t):=\frac{1}{\Area_{g(t)}(S)}\int_S B_t\,\mathrm d\mu_{g(t)},\qquad
 \zeta(S,t):=\frac{1}{\Area_{g(t)}(S)}\int_S Z\,\mathrm d\mu_{g(t)}.
\]
Extend $v$ and $\zeta$ as functions constant on each leaf. On the
model these labels are exactly $(B_t,Z)$. The potential label remains
$c=F|_S$; in the model it satisfies
\[
 c=\lambda v+a_{\rm ax}(\zeta+z_*)+\lambda t,
 \qquad Z=z-z_*.
\]
Thus $(v,\zeta)$ are radial and axial leaf labels, distinct from the
potential and mean-height labels $(c,\sigma)$ of
Section~\ref{sec:collars}. 

\begin{prop}\label{prop:productgrid}
Fix $N\ge1$, $c_N\ge1$, a sufficiently large $\Lambda$, and
$J\ge4$, $B>4J\Lambda$, $H>0$. Let
\[
 E\ge\max\{2\Lambda,B+\Lambda/8+1,\sqrt{H+2}\},
 \qquad w(v):=(v+c_N)^{-N}.
\]
Suppose $(g,F)$ is $\eta$-close in parabolic $C^{60}$-norm on
$\widehat Q_E$ to a tilted Bryant product with
$a_{\rm ax}\in[0,1]$. Fix $h\in(0,1)$ and overlapping
intervals $I_k,J_\ell,K_m$ of lengths comparable to $h$. Cover
$\overline{Q_E}$ by the annular and core boxes respectively
\[
 \begin{aligned}
 \mathcal B_{k\ell m}
   &:=\{v\in I_k,\ \zeta\in J_\ell,\ t\in K_m\},\\
 \mathcal B_{0\ell m}
   &:=\{B_t<2\Lambda+h,\ Z\in J_\ell,\ t\in K_m\}.
 \end{aligned}
\]
Use $\alpha,\beta,\gamma$ for box indices and $i,j\in\{1,2,3\}$ for field components. Let $\mathcal B_{\alpha}^+$ be its enlargement
by fixed multiples of $h$ in these parameters, with bounded
multiplicity inside $\widehat Q_E$; at $t=0$ enlarge only toward
earlier times.
For each $\alpha$, let
$U_{\alpha}=(U_{\alpha,1},U_{\alpha,2},U_{\alpha,3})$ be a time-independent triple on
an open spatial domain $\mathcal W_{\alpha}$. If $\mathcal B_\alpha$ is an annular box, then
\[
 \mathcal B_{\alpha}^+\subset\mathcal W_{\alpha}\times T_{\alpha}\subset\widehat Q_E,
 \qquad \mathcal C_{\alpha}^+\subset\mathcal W_{\alpha}\times T_{\alpha},
\]
Here $\mathcal C_{\alpha}^+$ is an open collar formed by whole regular
CMC leaves, containing the leaves meeting each annular box.
If $\mathcal B_\alpha$ is a core box, let $Z_{\alpha}$ be the midpoint of
$J_\ell$ and $t_{\alpha}$ the final time of its expanded time interval.
For these core boxes only, we require
\[
 \begin{aligned}
 \{B_{t_{\alpha}}<B+\Lambda/8,\ |Z-Z_{\alpha}|<H+1\}&\subset\mathcal W_{\alpha},\\
 [t_{\alpha}-H-1,t_{\alpha}]&\subset T_{\alpha}.
 \end{aligned}
\]
On each pair or triple overlap, assume that the common field
domains and collars satisfy the geometric hypotheses of
Proposition~\ref{prop:loworderalignment}, with a common whole
sphere and constants uniform in $E$ and $a_{\rm ax}$, including
those controlling the field domains and connecting paths.

At a box center $(x_{\alpha},t_{\alpha})$, we put
$r_{\alpha}:=\mathrm R(x_{\alpha},t_{\alpha})^{-1/2}$ and
$e_{\alpha}:=C_N\epsilon w(B_{t_{\alpha}}(x_{\alpha}))$.
Assume
\begin{align*}
 \sup_{\mathcal W_{\alpha}\times T_{\alpha}}\sum_{i=1}^3\sum_{j=0}^2
       r_{\alpha}^{2j}|\nabla^j\Lie_{U_{\alpha,i}}g|^2&\le e_{\alpha}^2,\\
 \sup_{\mathcal C_{\alpha}^+}r_{\alpha}^{-2}\sum_{i=1}^3
 \bigl(|\langle U_{\alpha,i},n_f\rangle|^2
      +|\langle U_{\alpha,i},\nu\rangle|^2\bigr)&\le e_{\alpha}^2,\\
 \sup_{t,S}|Q_S(U_{\alpha})-I|_{\mathrm F}&\le e_{\alpha},
\end{align*}
where $S$ ranges over whole leaves of $\mathcal C_{\alpha}^+$.
Assume also, on $\mathcal B_{\alpha}^+$,
\[
 \sum_{i=1}^3\sum_{j=0}^3|\nabla^jU_{\alpha,i}|\le C\sqrt E,
 \qquad
 \sum_{i=1}^3|U_{\alpha,i}F|
       \le C_N\sqrt E\,\epsilon w(B_t).
\]

\begin{enumerate}[label=\textup{(\roman*)},ref=\textup{(\roman*)}]
\item For every fixed $E$, there is $\eta_0(E)>0$, uniform for
$a_{\rm ax}\in[0,1]$, such that, if
$0\le\epsilon\le\eta\le\eta_0(E)$, there exists a smooth
time-dependent triple $(U_1,U_2,U_3)$ on $Q_E$ satisfying
\begin{equation}\label{eq:productgridestimate}
 \sum_{i=1}^3\left(\sum_{j=0}^2|\nabla^j\Lie_{U_i}g|
 +|\partial_tU_i|+|\Delta U_i+\Ric^\sharp U_i|+|U_iF|\right)
 \le C_NE^{5/2}\epsilon w(B_t).
\end{equation}
The constant $C_N$ is independent of $E$ and $a_{\rm ax}$.
\item After one constant orthogonal change of basis, the same triple satisfies
\begin{equation}\label{eq:productgrid-model-closeness}
 \sum_{i=1}^3\|U_i-K_i\|_{C^3(Q_E)}\le C_E\eta,
 \qquad\mathrm{where}\qquad K_i:=\krot J_i.
\end{equation}
Here $C_E$ may
depend on $E$ but is uniform in $a_{\rm ax}$.
\end{enumerate}
\end{prop}
\begin{proof}
We first construct the triple and prove (i).
Proposition~\ref{prop:alltiltcollar} constructs the annular CMC
families uniformly in $a_{\rm ax}$, and uniqueness identifies them
on overlaps. Use the labels and boxes defined in the statement,
with the mesh and enlargement factors independent of $E$ and
$a_{\rm ax}$.
To justify these labels, parametrize the model angular spheres by
$(b,z_0)=(B_t,Z)$ and use the regular-level identifications and
mean-zero graph construction of Proposition~\ref{prop:alltiltcollar}.
Its parameter estimates show that the resulting map
$(b,z_0)\mapsto(v,\zeta)$ is $O(\eta)$-close to the identity,
with three parameter derivatives and one time derivative, on the
fixed enlarged annulus. The estimates are uniform in $a_{\rm ax}$.
After shrinking the comparison threshold, this map is invertible
on the retained domains. Since the averages are defined on the
whole actual sphere, they agree for overlapping descriptions of
the same leaf. They therefore give compatible leaf coordinates.
For each fixed $E$, decrease $\eta_0(E)$ so that $(v,\zeta)$ and
their derivatives through order three and their first time
derivatives differ from $(B_t,Z)$ by a fixed small fraction of $h$.
The radial, axial, and temporal index counts satisfy
\[
 n_{\rm rad}\le CE,\qquad n_Z\le CE^2,\qquad n_t\le CE^2.
\]
On the core layer choose the radial cutoff constant near the tip;
its other cutoffs use only $Z$ and time. The domains $\mathcal W_{\alpha}$ contain the core, collars, and
connecting paths.
No inversion of the potential height or division by $a_{\rm ax}$
is used, and boxes ending at zero require only earlier-time data.

On each pair or triple overlap choose common nested spatial
domains and a time interval as in
Proposition~\ref{prop:loworderalignment}. By hypothesis their
interior margins and connecting paths are controlled after
rescaling by $r_{\alpha}$. The metric Lie derivative bounds hold on
these entire domains, and the normal and Gram bounds hold on the
common collar. Proposition~\ref{prop:loworderalignment} therefore
applies.
Write $b_{\alpha}:=B_{t_{\alpha}}(x_{\alpha})$. Model comparison gives
$r_{\alpha}\asymp\sqrt{b_{\alpha}+1}$. On intersecting boxes, $r_{\beta}\asymp r_{\alpha}$ and
$w(b_{\beta})\asymp w(b_{\alpha})\asymp w(B_t)$. The scale-invariant
conclusion therefore gives
\[
 \sum_{i=1}^3\sum_{j=0}^3r_{\alpha}^{j-1}
 |\nabla^j(U_{\alpha,i}-(P^{\alpha\beta}U_{\beta})_i)|
       \le C_N\epsilon w(B_t).
\]
Since $C^{-1}\le r_{\alpha}\le C\sqrt E$, this implies
\begin{equation}\label{eq:productgrid-alignments}
 \begin{aligned}
 \sum_{i=1}^3\sum_{j=0}^3
 |\nabla^j(U_{\alpha,i}-(P^{\alpha\beta}U_{\beta})_i)|
     &\le C_N\sqrt E\,\epsilon w(B_t),\\
 |P^{\alpha\beta}P^{\beta\gamma}P^{\gamma\alpha}-I|_{\mathrm F}
     &\le C_N\epsilon w(B_t).
 \end{aligned}
\end{equation}
Here the first inequality holds on expanded pair overlaps and the
second on expanded triple overlaps, with
$(P^{\alpha\beta}U_{\beta})_i=\sum_jP^{\alpha\beta}_{ij}U_{\beta,j}$.
Choose a matrix on each
unoriented edge once and define its reverse by inversion. The
whole-sphere Gram bound gives the stated cocycle estimate for
these choices. Changing the radial label by a bounded mesh amount
changes $w$ by a factor depending only on the fixed data and $N$.
Only the $C^2$ metric Lie derivative, $C^0$ collar-normal, and
whole-leaf Gram bounds are used to obtain these comparisons.

We recall the discrete estimate with its weight direction. Root the
grid at its outer radial boundary and one corner of its axial and time
ranges. We use a tree consisting of the time chain at this corner, all
axial rows on the outer radial boundary, and radial spokes directed
inwards from these rows. We define constant orthogonal matrices $G^{\alpha}$ by
transport along this tree. On a radial tree edge,
$G^{\alpha}P^{\alpha\beta}(G^{\beta})^{-1}=I$.
The tree loop of an axial edge fills a radial--axial rectangle of at
most $CE$ elementary faces. For a time edge, first move it radially
outward and then axially to the root time chain; the filling uses at
most $CE^2$ faces. All these faces have radial index at least that
of the original edge, up to a fixed mesh margin.

Subdivide each elementary rectangular face into two triangles of
expanded boxes with common overlap. The cocycle estimate in
\eqref{eq:productgrid-alignments} bounds its matrix product by
$C_N\epsilon w(B_t)$. Orthogonal conjugation preserves the matrix norm, and for a
product of orthogonal matrices,
\[
 \|H_1\cdots H_k-I\|\le\sum_{j=1}^k\|H_j-I\|.
\]
Since $w$ decreases outwards, this proves
\[
 \|G^{\alpha}P^{\alpha\beta}(G^{\beta})^{-1}-I\|
 \le C_NE^2\epsilon w(B_t)
\]
on all grid edges. Additional overlaps have bounded index separation
and reduce to these edges and boundedly many local triangles. Changing
$B_t$ by a fixed mesh amount changes $w(B_t)$ by a factor depending on $N$
but not on $E$. This is the only weight comparison used.

We write $V_{\alpha}:=G^{\alpha}U_{\alpha}$, meaning $V_{\alpha,i}=\sum_j(G^{\alpha})_{ij}U_{\alpha,j}$. The local alignment bound, the matrix discrepancy,
and the $C^3$ field size bound give
\[
 \|V_{\alpha}-V_{\beta}\|_{C^3}
 \le C_NE^{5/2}\epsilon w(B_t).
\]
Choose a partition of unity $\chi_{\alpha}$ subordinate to
$\mathcal B_{\alpha}^+$ with
\[
 \sum_{\alpha}\chi_{\alpha}=1,\qquad
 \sum_{\alpha}\left(\sum_{j=0}^3|\nabla^j\chi_{\alpha}|
                       +|\partial_t\chi_{\alpha}|\right)\le C.
\] Such cutoffs are obtained from fixed bump functions
in the grid coordinates, normalized by their sum. The model derivative
bounds and the small coordinate perturbations established above justify
these uniform bounds. We set $U_i:=\sum_{\alpha}\chi_{\alpha}V_{\alpha,i}$, $i=1,2,3$.

In every differentiated cutoff term subtract one active reference
field $V_{\beta}$, using $\sum_{\alpha} \mathrm{d}\chi_{\alpha}=0$ and its derivatives. Local metric Lie derivatives
and the aligned differences then give the first term of
\eqref{eq:productgridestimate}. Since the local fields and gauges are
time-independent, $\sum_{\alpha}\partial_t\chi_{\alpha}=0$ gives the bound for
$\partial_tU$. The identity
\[
 (\Delta U+\Ric^\sharp U)^\flat
 =\operatorname{div}(\Lie_Ug)
       -\tfrac12\mathrm{d}\operatorname{tr}_g(\Lie_Ug)
\]
gives the vector heat bound. Finally, $U_iF=\sum_{\alpha}\chi_{\alpha}V_{\alpha,i}F$ is bounded
directly by the assumed local pairings; the orthogonal matrices
preserve their triple norm. No derivatives of the partition occur in
this last identity. This proves \eqref{eq:productgridestimate}.
Thus we have proved (i). We next prove (ii). For fixed $E$,
view each local triple in the exact product metric. The metric,
normal, and whole-leaf Gram comparisons change its error by at most
$C_E\eta$, using the field bounds on the enlarged domains from the proof of
Proposition~\ref{prop:loworderalignment} and $\epsilon\le\eta$.
That proposition, applied
to this triple and $K$, gives a constant orthogonal comparison with
$C^3$ error $C_E\eta$ on its interior domain. Fix the model basis at
one initial grid vertex. Propagate it along the same finite tree
used above. The overlap comparisons imply $C_E\eta$ error at every
vertex, since the number of edges is finite for fixed $E$.
Interpolating with the fixed partition preserves this $C^3$ bound,
which proves \eqref{eq:productgrid-model-closeness}. Its constant
may depend on $E$; the constant in the weighted error estimate
\eqref{eq:productgridestimate} remains independent of $E$. This proves (ii).

To apply Proposition~\ref{prop:heatcap} at heat size $E$, perform
this gluing construction at size $2E$. The domains satisfy
\[
 Q_E^+\subset Q_{2E},
\]
so the resulting fields are defined on the full domain required by
that proposition. Their bound there is
$C_N(2E)^{5/2}\epsilon w(B_t)\le C'_NE^4\epsilon w(B_t)$.
The geometric comparison is taken on $\widehat Q_{2E}$, with its
margins fixed before choosing its tolerance. The factor $2$ changes
only the constants and does not change the radial exponent $N$.
\end{proof}

Anderson--Chow prove a drift inequality for $|h|^2/\mathrm{R}^2$
on three-dimensional Ricci flows with positive scalar curvature
\cite{AC}. Their more general denominator $\mathrm{R}+\rho$
allows a zeroth-order term of order $\rho$ when the scalar curvature is
not positive. A related estimate with the downward-shifted denominator
$\mathrm{R}-\rho$ is due to K.~Choi and appears in
\cite[Theorem~B.1]{BrendleAncient}. On a
three-dimensional Ricci flow with nonnegative Ricci curvature, it
controls $e^{2c_{\#}\rho t}|h|^2/(\mathrm{R}-\rho)^2$ by a drift
inequality whenever $\mathrm{R}\ge C_{\#}\rho$, for universal
constants $c_{\#}>0$ and $C_{\#}\ge10$.
Here we establish an estimate with the same denominator on a
four-dimensional Bryant product. The product curvature calculation
gives a strictly positive quadratic form after the downward shift,
and hence an absorption term. The estimate also holds for
sufficiently small curvature perturbations of the product.

\begin{prop}\label{cap:prop:product-quotient}
Let $g(t)$ be a four-dimensional Ricci flow on a spacetime domain
$\mathcal Q$, and let $h$ satisfy $\partial_th=\Delta_Lh$ on
$\mathcal Q$. Fix $\rho>0$. For each $(x,t)\in\mathcal Q$, suppose
there are a scaled Bryant product $(N,g_N)$, a point $q\in N$,
and a linear isometry
\[
 \iota_{x,t}:(T_xM,g(t))\longrightarrow(T_qN,g_N).
\]
Set $\Rm_N:=\iota_{x,t}^*\Rm_{g_N}(q)$ and
$\mathrm R_N:=\mathrm R_{g_N}(q)$, and assume
\[
 \mathrm R_N>\rho,\qquad \mathrm R_g>\rho,\qquad
 |\Rm_g-\Rm_N|\le\delta\mathrm R_N,\qquad
 \delta\le c_{\mathrm{pin}}\rho/\mathrm R_N.
\]
Here $c_{\mathrm{pin}}>0$ is a sufficiently small numerical constant.
Put
\[
 \mathfrak q_h:=\frac{|h|^2}{(\mathrm R_g-\rho)^2},\quad
 S_\rho(h):=|\Ric_g|^2|h|^2-(\mathrm R_g-\rho)\Rm_g(h,h),\quad
 \mathscr D:=\partial_t-\Delta-
 \frac{2\nabla\mathrm R_g}{\mathrm R_g-\rho}\cdot\nabla.
\]
Then on $\mathcal Q$, we have
\[
 S_\rho(h)\ge\tfrac16\rho\mathrm R_N|h|^2,
 \qquad (\mathscr D+\tfrac14\rho)\mathfrak q_h\le0.
\]
\end{prop}
Note that on the exact product metric, $1/6$ can be replaced by $1/3$.
\begin{proof}
Put $\mathcal H:=\partial_t-\Delta$ and $\mathrm R:=\mathrm R_g$.
For $\partial_th=\Delta_Lh$, the norm calculation gives
\begin{equation}\label{cap:eq:quotient}
 \mathscr D \mathfrak q_h=-\frac{2}{(\mathrm{R}-\rho)^4}
       |(\mathrm{R}-\rho)\nabla h-h\otimes\nabla \mathrm{R}|^2
       -\frac{4S_\rho(h)}{(\mathrm{R}-\rho)^3}.
\end{equation}
To check it, use
$\mathcal H |h|^2=-2|\nabla h|^2+4\Rm(h,h)$ and
$\mathcal H \mathrm{R}=2|\Ric|^2$, differentiate the quotient, and complete the
gradient square. In detail, put $\phi:=\mathrm{R}-\rho$ and $H:=|h|^2$.
The product rule gives
\begin{align*}
 \mathcal H(H\phi^{-2})={}&
 \phi^{-2}\bigl(-2|\nabla h|^2+4\Rm(h,h)\bigr)
 -4\phi^{-3}H|\Ric|^2\\
 &-6\phi^{-4}H|\nabla\phi|^2
       +4\phi^{-3}\langle\nabla H,\nabla\phi\rangle,\\
 -2\phi^{-1}\langle\nabla\phi,\nabla(H\phi^{-2})\rangle
 ={}&-2\phi^{-3}\langle\nabla\phi,\nabla H\rangle
       +4\phi^{-4}H|\nabla\phi|^2.
\end{align*}
Adding these equations and using
$\nabla H=2\langle h,\nabla h\rangle$ gives exactly
\eqref{cap:eq:quotient}. No dimension restriction is involved in
this differential identity; the dimension enters the curvature
inequality for $S_\rho$ proved next.

For the Bryant product the required algebra has a direct proof.
We choose an orthonormal frame on the Bryant factor with $e_1$ radial
and $e_2,e_3$ tangent to its round orbit. Its sectional curvatures
are $K_{12}=K_{13}=k\ge0$ and $K_{23}=l\ge0$. At the tip take
$k=l$ by isotropy. We write
\[
 \mathrm{R}_N=4k+2l,\qquad s:=|\Ric_N|^2=6k^2+4kl+2l^2.
\]
We use the curvature action
$\Rm_N(H)_{ij}:=(\Rm_N)_{ikjl}H^{kl}$, as in
\eqref{eq:curvature-action}. We first prove that
$sI-\mathrm{R}_N\Rm_N$ is positive semidefinite
on symmetric two-tensors. On the off-diagonal symmetric tensors of the Bryant factor,
$\Rm_N$ has eigenvalues $-k,-k,-l$, so the assertion holds there.
On diagonal tensors the curvature contraction has matrix
\[
 C=\begin{pmatrix}0&k&k\\k&0&l\\k&l&0\end{pmatrix}.
\]
The vector $(0,1,-1)$ has eigenvalue $-l$. On its orthogonal
complement, with orthonormal coordinates $(x,(y+z)/\sqrt2)$,
$sI-\mathrm{R}_NC$ is
\[
 \begin{pmatrix}s&-\sqrt2\mathrm{R}_N k\\-\sqrt2\mathrm{R}_N k&s-\mathrm{R}_Nl\end{pmatrix}.
\]
Its second diagonal entry is $6k^2$, and its determinant is
\[
                  6k^2s-2\mathrm{R}_N^2k^2=4k^2(k-l)^2\ge0.
\]
Its first diagonal entry is positive when $\mathrm{R}_N>0$. This proves
\[
          s|H|^2-\mathrm{R}_N\Rm_N(H,H)\ge0
\]
for every symmetric $H$ on the three-dimensional factor.
For $0<\rho\le \mathrm{R}_N$, the exact convex decomposition
\begin{align*}
 s|H|^2-(\mathrm{R}_N-\rho)\Rm_N(H,H)
 =(1-\rho/\mathrm{R}_N)\bigl(s|H|^2-\mathrm{R}_N\Rm_N(H,H)\bigr)+(\rho/\mathrm{R}_N)s|H|^2
 \ge \tfrac13\rho \mathrm{R}_N|H|^2
\end{align*}
uses $s\ge \mathrm{R}_N^2/3$. Thus, the shifted estimate is uniform even
at the cylindrical limit $k=0$ and at the round tip $k=l$.

For the four-dimensional product, decompose
$h=H+\mathrm{d}z\otimes\alpha+\alpha\otimes \mathrm{d}z+\chi_0\,\mathrm{d}z^2$. Then
\[
 S_\rho(h)=S_{\rho,N}(H)+s(2|\alpha|^2+\chi_0^2)
               \ge\tfrac13\rho \mathrm{R}_N|h|^2.
\]
The same inequality therefore controls the Bryant, mixed, and
line components of $h$.

If the actual curvature tensor differs from the product tensor
by at most $\delta \mathrm{R}_N$ in an orthonormal identification, expansion
of $S_\rho$ changes it by at most $C\delta \mathrm{R}_N^2|h|^2$.
Thus, the inequality persists with half its constant when
\begin{equation}
 \delta\le C^{-1}\rho/\mathrm{R}_N.
\end{equation}
Choosing $c_{\mathrm{pin}}$ sufficiently small therefore gives
$S_\rho(h)\ge\rho\mathrm{R}_N|h|^2/6$ and
$\mathrm{R}\le2\mathrm{R}_N$. Substituting these bounds into
\eqref{cap:eq:quotient} yields
\[
 \mathscr D\mathfrak q_h
 \le-\frac23\rho\frac{\mathrm{R}_N}{\mathrm{R}-\rho}\mathfrak q_h
 \le-\frac13\rho\mathfrak q_h.
\]
The nonpositive gradient square can be discarded, and $1/4<1/3$
gives the stated absorption inequality.
\end{proof}

We now improve the metric Lie derivative of the glued fields by the vector heat
equation. The proposition is local: it assumes the comparison chart
and weighted fields specified in its statement.
Section~\ref{sec:exterior} verifies these hypotheses on the exterior
of the original soliton before applying the estimate repeatedly.

\begin{prop}\label{prop:heatcap}
Fix $a_0,H,d_{\rm out},\tau_{\rm out}>0$ and connected spatial
cap domains $D_{\rm out}\Subset D_{\rm out}^+$ with
$\operatorname{dist}_{g(0)}(D_{\rm out},\partial D_{\rm out}^+)
\ge d_{\rm out}$. For $t\in[-H-\tau_{\rm out},0]$, let
$C_{\rm out}(t)\Subset C_{\rm out}^+(t)\subset D_{\rm out}^+$
be regular potential level CMC collars.
Set
\[
 \begin{aligned}
 Q_{\rm out}&:=D_{\rm out}\times[-H,0],\\
 Q_{\rm collar}&:=\{(x,t):x\in C_{\rm out}(t),\ -H\le t\le0\},\\
 \mathcal Q_{\rm out}^+&:=D_{\rm out}^+
                              \times[-H-\tau_{\rm out},0].
 \end{aligned}
\]
Assume the geometry, volume-to-area, and mean-zero Jacobi inverse
bounds of Proposition~\ref{prop:normal-component} hold uniformly
on the enlarged collars and their bounded potential level interiors.
These sets and margins are fixed before $E$. There is $E_0$ such that $\mathcal Q_{\rm out}^+\Subset Q_E$
relative to $t\le0$ for every $E\ge E_0$. For each such $E$
there are $\eta_0(E)>0$ and $C_E<\infty$ with the following property.

 Suppose $(g,F)$ has comparison error at most $\eta$
in \eqref{eq:cap-comparison-error}, with $m=60$ on $\widehat Q_E$,
for a Bryant product with $a_{\rm ax}\in[a_0,1]$, and satisfies
\[
 \partial_tg=-2\Ric,\qquad \Hess F=\lambda\Ric,\qquad
 |\nabla F|^2=1-\lambda^2\mathrm R,\qquad
 (\partial_t-\Delta)F=0,\qquad 0\le\lambda\le1.
\]
Fix a weight shift $c_{16}\ge1$, and let $J=(J_1,J_2,J_3)$ be
model rotations with cyclic bracket constant one. Suppose $U=(U_1,U_2,U_3)$ on $Q_E^+$ satisfies
\[
 \sum_{i=1}^3\|U_i-J_i\|_{C^3(Q_E^+)}\le C_E\eta
\]
and, pointwise on $Q_E^+$,
\[
 \sum_{i=1}^3\left(\sum_{j=0}^2|\nabla^j\Lie_{U_i}g|
 +|\partial_tU_i|+|\Delta U_i+\Ric^\sharp U_i|+|U_iF|\right)
 \le CE^4\epsilon(B_t+c_{16})^{-16}.
\]
For $0<\epsilon\le\eta\le\eta_0(E)$, the following hold.
\begin{enumerate}[label=\textup{(\roman*)},ref=\textup{(\roman*)}]
\item There is a time-independent
triple $Y=(Y_1,Y_2,Y_3)$ near $D_{\rm out}$ such that, for $i=1,2,3$,
\begin{align}
 \|\Lie_{Y_i}g\|_{C^2(Q_{\rm out})}
 +\|Y_iF\|_{C^2(Q_{\rm out})}
 +\|\langle Y_i,\nu\rangle\|_{C^2(Q_{\rm collar})}
 &\le(CE^{-9}+C_E\eta)\epsilon.
\end{align}
\item The same triple satisfies for $i=1,2,3$,
\[
 \|Y_i-J_i\|_{C^2(Q_{\rm out})}\le C_E\eta.
\]
All $C^k$ norms are spatial and uniform in time; $\nu$ is the
actual CMC normal inside the potential level. The constants depend
on the fixed domains, field bounds, and $a_0$; $C$ is independent of $E$.
\end{enumerate}
\end{prop}

\begin{proof}
We construct the triple in (i) and keep the model comparison in (ii)
throughout the heat correction and freezing. For glued triples normalized by $Q_S=I$, first divide the fields by
$\sqrt{6\pi}$ to match the bracket-one model rotations. This changes
only the fixed constants in the assumed bounds.
Choose connected domains $D_1,D_2$ with
$D_{\rm out}\Subset D_1\Subset D_2\Subset D_{\rm out}^+$, and put
\[
 Q_1:=D_1\times[-H-\tau_{\rm out}/3,0],\qquad
 Q_2:=D_2\times[-H-2\tau_{\rm out}/3,0].
\]
Thus $Q_{\rm out}\Subset Q_1\Subset Q_2
\Subset\mathcal Q_{\rm out}^+$, relative to $t\le0$.
Choose $D_1$ also to contain fixed enlargements of all output
collars. Compactness of their parameter ranges and the assumed
spatial margins permit this choice. All intermediate estimates below
use these fixed spatial and earlier-time margins.

\medskip\noindent\textbf{Step 1. Solve the vector heat problem with weighted data.}
Put $\mathcal H:=\partial_t-\Delta$ and
$\mathcal V W:=\Delta W+\Ric^\sharp W$ for a vector field $W$.
For fixed $E$ and sufficiently small $\eta$, model comparison gives
on $Q_E^+$
\[
 -2\le\mathcal H B_t\le-\tfrac12,\quad |\nabla B_t|^2\le2,
 \quad \frac{1}{C(B_t+1)}\le \mathrm{R}\le\frac{C}{B_t+1}.
\]
Take $p=14$ and $b_{\rm wt}\ge\max\{c_{16},100(p+1)\}$.
The weights $B_t+c_{16}$ and $B_t+b_{\rm wt}$ differ by a fixed factor,
which is absorbed in $C$.
We solve
\[
 \begin{cases}
 (\partial_t-\mathcal V)V_i=0&\text{on }Q_E^+,\\
 V_i=U_i&\text{on }\partial_{\mathrm{par}}Q_E^+.
 \end{cases}
\]
Here the parabolic boundary is the initial face at $t=-2E^2$ and
the spatial faces $B_t=2E$, $|Z|=2E^2$; the Bryant tip is interior. Smooth exhaustion handles its
corners; alternatively all estimates below can first be proved on
smooth approximating domains with the same interior margins.

Set $w_i:=V_i-U_i$ and
$f_i:=-(\partial_t-\mathcal V)U_i$. Then
$(\partial_t-\mathcal V)w_i=f_i$, $w_i=0$ on
$\partial_{\rm par}Q_E^+$, and the assumed bounds give
$|f_i|\le CE^4\epsilon(B_t+b_{\rm wt})^{-16}$.
Under Ricci flow,
\[
 \mathcal H|w_i|^2=-2|\nabla w_i|^2+2\langle f_i,w_i\rangle.
\]
Kato's inequality therefore gives
\[
 \mathcal H|w_i|\le C E^4\epsilon(B_t+b_{\rm wt})^{-16}.
\]
For $G:=(B_t+b_{\rm wt})^{-p}$,
\begin{align*}
 \mathcal H G
 &=-p(B_t+b_{\rm wt})^{-p-1}\mathcal H B_t
   -p(p+1)(B_t+b_{\rm wt})^{-p-2}|\nabla B_t|^2\\
 &\ge \frac p4(B_t+b_{\rm wt})^{-p-1}.
\end{align*}
Since $w_i$ vanishes on the parabolic boundary, comparison yields
\[
 |w_i|\le C E^4\epsilon(B_t+b_{\rm wt})^{-14}.
\]
On a unit parabolic neighborhood of $\overline{Q_E}$ contained
in $Q_E^+$, the metric coefficients are uniformly controlled and
$B_t+b_{\rm wt}$ varies by a bounded factor. Interior $W^{2,1}_r$
estimates, with $r$ above the parabolic Sobolev threshold, give
\[
 |\nabla w_i|\le CE^4\epsilon(B_t+b_{\rm wt})^{-14}.
\]
These estimates require only the displayed $C^0$ bound for $f_i$. Consequently, on a neighborhood of the closure of
$Q_E$ from \eqref{eq:cap-heat-enlargement},
\begin{equation}\label{eq:weightedh}
 h_i=\Lie_{V_i}g,\qquad
 \partial_th_i=\Delta_Lh_i,\qquad
 |h_i|\le C E^4\epsilon(B_t+b_{\rm wt})^{-14}.
\end{equation}
The exact scalar pairing satisfies
\begin{equation}\label{eq:tiltedscalar}
 \mathcal H(V_iF)=-\lambda\langle\Ric,h_i\rangle.
\end{equation}
Moreover, $|V_iF|\le|U_iF|+|w_i|$, since $|\nabla F|\le1$.
Consequently $|V_iF|\le CE^4\epsilon(B_t+b_{\rm wt})^{-14}$ on $Q_E$.
Interior estimates for \eqref{eq:weightedh} and
\eqref{eq:tiltedscalar} give $C_E\epsilon$ bounds for the spatial
derivatives used below on $Q_2$.
The equation for $V_i-J_i$ has boundary values $U_i-J_i$ and
forcing $\mathcal V_gJ_i$, both bounded by $C_E\eta$.
Comparison and interior regularity yield
\[
 \sum_{i=1}^3\|V_i-J_i\|_{C^{24}(Q_2)}\le C_E\eta.
\]

We next improve $h_i$ on $Q_2$ with a constant independent of $E$.
We set $\rho:=\rho_{\mathrm{pin}}/E$, with
$\rho_{\mathrm{pin}}>0$ sufficiently small, and
$q_i:=|h_i|^2/(\mathrm{R}-\rho)^2$. On $Q_E$, the model scalar curvature
satisfies $\mathrm{R}_B\ge (CE)^{-1}$. We first fix $E$, and then make the
curvature comparison small relative to $\rho/\mathrm{R}_B$.
Proposition~\ref{cap:prop:product-quotient} gives
\[
 \mathcal Pq_i\le0,\qquad
 \mathcal P:=\mathcal H-\frac{2\nabla\mathrm{R}}{\mathrm{R}-\rho}
       \cdot\nabla+\kappa_{\mathrm{pin}}\rho,\qquad
 \kappa_{\mathrm{pin}}:=\tfrac14.
\]
The proposition applies to every symmetric tensor on the product,
including its mixed and line components. Its use here is therefore
independent of the potential tilt.

We use the comparison line coordinate $Z$, not $F$, for the axial barrier.
Let $\mathbf b:=2\nabla \mathrm{R}/(\mathrm{R}-\rho)$ and
$\gamma^2:=\kappa_{\mathrm{pin}}\rho/8$. Model comparison gives
\[
 |\nabla Z|^2\le2,\qquad
 |\mathcal H Z|+|\mathbf b\cdot\nabla Z|\le C_E\eta.
\]
Choosing $\eta$ after $E$ makes
$\gamma C_E\eta\le\kappa_{\mathrm{pin}}\rho/4$. Direct differentiation shows
\[
 \mathcal P\frac{\cosh(\gamma Z)}{\cosh(\gamma E^2)}\ge0.
\]
By \eqref{eq:weightedh}, initial and axial boundary data for $q_i$
are bounded by $CE^8\epsilon^2$, while the radial boundary data are
bounded by $CE^{-18}\epsilon^2$. The nonnegative supersolution
\[
 C\epsilon^2\left[
 E^8e^{-\kappa_{\mathrm{pin}}\rho(t+E^2)}+E^{-18}
       +E^8\frac{\cosh(\gamma Z)}{\cosh(\gamma E^2)}\right]
\]
dominates all these data. This comparison includes the moving face $B_t=E$; the tip is
interior. On the fixed enlargement $\mathcal Q_{\rm out}^+$,
$B_t,|Z|,|t|$ are bounded independently of $E$. The barrier gives
\[
 q_i\le C\epsilon^2
       (E^8e^{-cE}+E^{-18}+E^8e^{-cE^{3/2}}).
\]
Since $E^8e^{-cE}+E^8e^{-cE^{3/2}}\le CE^{-18}$ for large $E$,
interior regularity of $\partial_th_i=\Delta_Lh_i$ gives
\[
 \sum_{i=1}^3\|h_i\|_{C^{24}(Q_2)}\le CE^{-9}\epsilon.
\]
The constants depend on the fixed domains and margins, not on $E$.

\medskip\noindent\textbf{Step 2. Freeze the fields on the fixed region.}
Set $W_i:=V_i(\cdot,0)$ on $D_1$ and keep $W_i$ independent
of time. The identity
\[
 \partial_tV_i=
 (\operatorname{div}h_i-\tfrac12\mathrm{d}\operatorname{tr}h_i)^\sharp
\]
can be integrated from $t$ to zero. With $d:=CE^{-9}\epsilon$,
the resulting estimates on $Q_1$ are
\[
 \|\Lie_{W_i}g\|_{C^{22}(Q_1)}\le d,\qquad
 \|W_iF\|_{C^{22}(Q_1)}\le C_E\epsilon,\qquad
 \|W_i-J_i\|_{C^{22}(Q_1)}\le C_E\eta.
\]
\medskip\noindent\textbf{Step 3. Improve potential tangency by a bracket.}
The core estimate must not divide by $|\nabla F|$. Instead, the soliton
identity and its Ricci variation control the angular derivatives of
$W_aF$; the bracket then removes the remaining invariant component.

For $a=1,2,3$, put $u_a:=W_aF$ and $H_a:=\Lie_{W_a}g$.
The indices $i,j,k$ in the next two identities are tensor indices. The variation of the
Levi--Civita connection is
\[
 (\delta\Gamma[H_a])_{ij}^{\ k}
 =\frac12g^{k\ell}
       \bigl(\nabla_i(H_a)_{j\ell}+\nabla_j(H_a)_{i\ell}
                                      -\nabla_\ell(H_a)_{ij}\bigr).
\]
Differentiating $\Hess F=\lambda\Ric$ under the flow of $W_a$
therefore gives the indexed identity
\begin{equation}
 \nabla_i\nabla_j u_a
 =\lambda(D\Ric)_g[H_a]_{ij}
                      +(\delta\Gamma[H_a])_{ij}^{\ k}\nabla_kF.
\end{equation}
The Ricci variation is a second-order linear operator in $H_a$.
Since $\lambda\le1$ and the terminal geometry is fixed, this implies
$\|\Hess u_a\|_{C^{20}}\le Cd$. Contracting the scalar commutation
identity gives
\begin{equation}
 \Ric_{jk}\nabla^k u_a
     =\nabla^i(\nabla_i\nabla_j u_a)-\nabla_j\Delta u_a,
 \qquad
 \|\Ric^\sharp\nabla u_a\|_{C^{19}}\le Cd.
\end{equation}
Let $P_B$ denote projection onto the three Bryant directions in the
comparison product. On $Q_1$,
$\Ric_B^\sharp|_{TB}\ge c_B I$ for a constant $c_B>0$. Comparison
of the actual and model endomorphisms, together with
$\|\nabla u_a\|_{C^{21}}\le C_E\epsilon$, yields
\begin{align*}
 \|P_B\nabla u_a\|_{C^{19}}
 &\le C\bigl(\|\Ric^\sharp\nabla u_a\|_{C^{19}}
                             +\eta\|\nabla u_a\|_{C^{19}}\bigr)\\
 &\le C d+C_E\eta\epsilon,\qquad
 \|J_j u_a\|_{C^{19}}\le C d+C_E\eta\epsilon.
\end{align*}
Only the Bryant restriction is inverted; no inverse of the vanishing
line eigenvalue or of the axial potential coefficient is used.
For cyclic indices $(a,b,c)$ define $A_c:=[W_a,W_b]$. The two
identities needed for the improved errors are
\begin{align*}
 \Lie_{A_c}g&=\Lie_{W_a}H_b-\Lie_{W_b}H_a,\\
 A_cF&=W_a u_b-W_bu_a\\
     &=J_a u_b-J_bu_a+(W_a-J_a)u_b-(W_b-J_b)u_a.
\end{align*}
Each term in the last line contains either an improved angular
pairing or the product of an $O(C_E\eta)$ coefficient with an
$O(C_E\epsilon)$ derivative. Thus, with $e_E:=d+C_E\eta\epsilon$, on $Q_1$,
\[
 \|\Lie_{A_i}g\|_{C^{16,\alpha}(Q_1)}
 +\|A_iF\|_{C^{16,\alpha}(Q_1)}\le Ce_E,\qquad
 \|A_i-J_i\|_{C^{16,\alpha}(Q_1)}\le C_E\eta.
\] Here and below $0<\alpha<1$ is fixed. The additional four
retained heat derivatives allow these H\"older estimates on the fixed
interior domains. Every term containing $W_a-J_a$ also contains a
derivative of $W_bF$, so its error is $C_E\eta\epsilon$.

\medskip\noindent\textbf{Step 4. Control the within-level normal.}
For the present positive-axial-coefficient statement, set $Y=A$.
The preceding estimates hold on the bounded level interiors as well
as on the collars. Comparison gives $|\nabla F|\ge a_0/2$ there.
Proposition~\ref{prop:normal-component} therefore controls
$\langle Y_i,\nu\rangle$. For its transverse derivatives, pull the
projected Jacobi equation \eqref{eq:normal-component-projected} to
the fixed sphere chart and differentiate it together with its
zero-average constraint. The same projected inverse controls each
successive derivative; differentiating the flux identity
\eqref{eq:normal-component-flux} controls its average. Induction
through order two proves (i). The first bracket preserves comparison
with $J$ and time independence, giving (ii).

For later use, Steps 1--3 require neither the lower bound for
$a_{\rm ax}$ nor bounded level interiors: their constants are uniform
for $a_{\rm ax}\in[0,1]$. They give the displayed
$C^{16,\alpha}$ estimates for $A$ and retain comparison with $J$
through at least seventeen derivatives on fixed interior domains.
\end{proof}

The following proposition normalizes the triple by one constant matrix
on the whole fixed cap and time interval.

\begin{prop}\label{prop:commongram}
Fix connected cap domains $D_{\rm out}\Subset D^+$ and a compact
time interval $I$. Let
$C_{\rm out}(t)\Subset C^+(t)\subset D^+$ be regular potential level
CMC collars of near-round spheres, with the leaf-parameter margins
of Proposition~\ref{prop:heatcap}. Assume
uniform $C^8$ bounds for the metric, its inverse, and the collar
charts, and uniformly bounded paths in $D^+$ from a fixed collar
sphere to every point of $D_{\rm out}$.
Put $\mathcal C^+:=\{(x,t):x\in C^+(t),\ t\in I\}$. All
collar norms include angular and transverse spatial derivatives
in these foliation charts, uniformly in $t$.

Let $Y=(Y_1,Y_2,Y_3)$ be time-independent on $D^+$ and sufficiently
close in spatial $C^9$ to a standard oriented Bryant rotational triple.
Suppose, for $i=1,2,3$,
\[
 \|\Lie_{Y_i}g\|_{C^8(D^+\times I)}
 +\|\langle Y_i,n_f\rangle\|_{C^8(\mathcal C^+)}
 +\|\langle Y_i,\nu\rangle\|_{C^8(\mathcal C^+)}\le e.
\]
\begin{enumerate}[label=\textup{(\roman*)},ref=\textup{(\roman*)}]
\item There is one constant positive definite matrix $Q_*$, with
uniformly bounded inverse, such that every whole leaf
$S\subset C_{\rm out}(t)$, $t\in I$, satisfies
\begin{equation}
 \sum_{i,j=1}^3\left|\Area(S)^{-2}
 \int_S\langle Y_i,Y_j\rangle\,\mathrm d\mu-(Q_*)_{ij}\right|\le Ce.
\end{equation}
\item The triple $\widehat Y_i:=\sum_j(Q_*^{-1/2})_{ij}Y_j$ satisfies
$|Q_S(\widehat Y)-I|_{\mathrm F}\le Ce$ on the same leaves and
satisfies, for $i=1,2,3$,
\begin{align*}
	\|\Lie_{\widehat Y_i}g\|_{C^8(D^+\times I)}
 +\|\langle \widehat Y_i,n_f\rangle\|_{C^8(\mathcal C^+)}
 +\|\langle \widehat Y_i,\nu\rangle\|_{C^8(\mathcal C^+)}\le Ce.
\end{align*}
Constants depend only on the fixed geometry and frame nondegeneracy.
\end{enumerate}
\end{prop}

\begin{proof}
We first prove (i). Fix $t_0\in I$, a collar sphere, and $p$ on that sphere.
Choose an orthonormal basis $e_1,e_2$ of its tangent plane. For an
ambient field $Z$, define its tangential value and derivative data by
\[
 \mathcal J_pZ:=\left(\langle Z,e_1\rangle,\langle Z,e_2\rangle,
 \tfrac12(\langle\nabla_{e_1}Z,e_2\rangle
                    -\langle\nabla_{e_2}Z,e_1\rangle)\right)_p.
\]
The map $v\mapsto\mathcal J_p(\sum_{i=1}^3v_iY_i)$ is invertible,
with uniformly bounded inverse. On a standard round sphere the
first two entries detect rotations moving $p$, and the third
detects the rotation fixing $p$; the $C^9$ comparison preserves
this invertibility. For each $a,b\in\{1,2,3\}$, choose the constant
coefficients $C_{ab}^{\ c}$ so that
\[
 E_{ab}:=[Y_a,Y_b]-\sum_{c=1}^3C_{ab}^{\ c}Y_c,\qquad
 \mathcal J_pE_{ab}=0.
\]
In the following formulas repeated field indices range from one
to three. Cap estimates are on $D_{\rm out}\times I$; leaf estimates
use the induced metric of the specified sphere.

The metric Lie derivative of $E_{ab}$ is $O(e)$ by the Lie derivative commutator
identity. Its two transverse pairings, with their derivatives, are
also $O(e)$. To justify this latter assertion, we use controlled local leaf labels $(\chi,\sigma)$ for the actual
foliation. The hypotheses imply $Y_i\chi,Y_i\sigma=O(e)$ with derivatives,
so $[Y_a,Y_b]\chi,[Y_a,Y_b]\sigma=O(e)$. Their controlled coefficients then convert these estimates into
bounds for the actual orthonormal conormals. At the chosen point,
$E_{ab}=O(e)$ since its tangential value is zero. Differentiating the
two small transverse pairings determines the normal components of
$\nabla E_{ab}$ up to $O(e)$. Small metric Lie derivative then determines the other
mixed components, and the matched tangential skew derivative supplies
the only remaining tangential component. Hence $E_{ab}$ and its full first covariant derivative
are $O(e)$ there.

The differentiated metric Lie derivative identity expresses $\nabla^2E_{ab}$ as a
curvature contraction with $E_{ab}$ plus a linear combination of first
derivatives of $\Lie_{E_{ab}}g$. Integrating the coupled equations
for $(E_{ab},\nabla E_{ab})$ along bounded paths inside the enlarged
fixed cap and applying Gronwall gives $E_{ab}=O(e)$ with the required
fixed derivative bounds at $t_0$. We perform this propagation on a fixed
connected interior neighborhood of
\[
 \overline{D_{\rm out}}\ \cup\
           \bigcup_{t\in I}\overline{C_{\rm out}(t)}.
\]
The fixed collar margins and compact time interval retain this set
inside $D^+$; connectedness and a finite covering supply connecting
paths of bounded length. The constants may depend on these fixed
domains. Thus the structure estimate holds on the output collars
as well as on the cap interior. All fields and
coefficients are time-independent. Uniform comparison of the metrics
and connections over the fixed time interval therefore gives
\begin{equation}\label{eq:constantstructure}
 [Y_a,Y_b]=C_{ab}^{\ c}Y_c+O_{C^2}(e)
\end{equation}
on the full fixed cap at every time. Higher bounds used below
are obtained by the same differentiated metric Lie derivative identity; the eight
assumed derivatives suffice.

On each collar sphere project $Y_a$ to its tangent bundle:
\[
 T_a=Y_a-\langle Y_a,n_f\rangle n_f-\langle Y_a,\nu\rangle\nu.
\]
The hypotheses give $T_a-Y_a=O_{C^8}(e)$, small intrinsic metric
Lie derivatives, and the same bracket relation
\eqref{eq:constantstructure} up to $O_{C^2}(e)$. The projected
triple is close to the oriented round rotations. Apply
Lemma~\ref{lem:rotation-gram} to this controlled family, using the
single collection of constants $C_{ab}^{\ c}$. Its matrix $Q_*$ is
independent of leaf and time and proves (i).
For (ii), conjugate each leaf Gram matrix by $Q_*^{-1/2}$.
Since that matrix is constant and uniformly bounded, the metric Lie
derivatives and both transverse pairings obey the same estimates. This finishes the proof.
\end{proof}

\begin{cor}\label{cor:alltiltheat}
Use the fixed cap $D_{\rm out}$, time interval $[-H,0]$, and
spacetime domains $Q_{\rm out},Q_{\rm collar}$ of
Proposition~\ref{prop:heatcap}. Let $(g,F)$ satisfy its flow and
potential identities on $\widehat Q_E$. Regular potential levels
and CMC normals are required on the annular collars; no bounded
potential level interiors are assumed. For $a_{\rm ax}\in[0,1]$, assume
\[
 \|g-\bar g\|_{\mathrm{par},60;\widehat Q_E}
 +\|F-F_{\rm mod}\|_{\mathrm{par},60;\widehat Q_E}\le\eta.
\]
Let $J=(J_1,J_2,J_3)$ be bracket-one model rotations and
$U=(U_1,U_2,U_3)$ a triple on $Q_E^+$ satisfying
\[
 \sum_i\|U_i-J_i\|_{C^3(Q_E^+)}\le C_{\rm rot}(E)\eta,
\]
and, pointwise on $Q_E^+$,
\[
 \sum_i\left(\sum_{j=0}^2|\nabla^j\Lie_{U_i}g|
 +|\partial_tU_i|+|\Delta U_i+\Ric^\sharp U_i|+|U_iF|\right)
 \le C_{\rm gl}E^4\epsilon(B_t+c_{16})^{-16}.
\]
Here $c_{16},C_{\rm gl}$ are fixed independently of $E,a_{\rm ax}$,
and $C_{\rm rot}(E)$ is independent of $a_{\rm ax}$.
Fix $s_+>s_-\ge s_0$, $Z_*>0$, and
$0<\delta<\min\{s_-,(s_+-s_-)/2,Z_*\}$ as in
Proposition~\ref{prop:alltiltcollar}. For every
$t\in[-H-\tau_{\rm out},0]$, require
\[
 \begin{aligned}
 C_{\rm out}(t)&\subset
 \{s_-+\delta<s_t<s_+-\delta,\ |Z|<Z_*-\delta\},\\
 C_{\rm out}^+(t)&\subset
 \{s_--\delta<s_t<s_++\delta,\ |Z|<Z_*+\delta\}.
 \end{aligned}
\]
Use the actual CMC spheres of Proposition~\ref{prop:alltiltcollar},
with the fixed spatial and leaf-parameter margins of
Proposition~\ref{prop:heatcap}.
For every sufficiently large $E$ there is $\eta_0(E)>0$, uniform
in $a_{\rm ax}$, such that $0<\epsilon\le\eta\le\eta_0(E)$ yields
a time-independent triple $Z=(Z_1,Z_2,Z_3)$ near $D_{\rm out}$ with
\begin{align*}
 &\sum_i\bigl(
 \|\Lie_{Z_i}g\|_{C^8(Q_{\rm out})}
 +\|Z_iF\|_{C^8(Q_{\rm out})}
 +\|\langle Z_i,n_f\rangle\|_{C^8(Q_{\rm collar})}
 +\|\langle Z_i,\nu\rangle\|_{C^8(Q_{\rm collar})}\bigr)\\
 &\qquad+\sup_{t,S}|Q_S(Z)-I|_{\mathrm F}
       \le(CE^{-9}+C_E\eta)\epsilon.
\end{align*}
The supremum is over whole collar leaves for $t\in[-H,0]$. The constant $C$ depends on
the fixed domains, $c_{16},C_{\rm gl}$, but not on $E,a_{\rm ax}$;
$C_E$ may also depend on $E,C_{\rm rot}(E)$, but not on $a_{\rm ax}$.
\end{cor}
\begin{proof}
Apply Steps 1--3 of the proof of Proposition~\ref{prop:heatcap},
which are uniform for $a_{\rm ax}\in[0,1]$: the product metric is
independent of the tilt, the scalar identities use $\lambda\le1$,
the axial barrier uses $Z$, and only the positive Bryant Ricci block
is inverted. They yield a time-independent first bracket triple $A$
with metric Lie derivative and potential pairing bounded by
$e_E:=(CE^{-9}+C_E\eta)\epsilon$ in $C^{16,\alpha}$, and
$C_E\eta$ model comparison through seventeen derivatives.

Set $T_1=[A_2,A_3]$, $T_2=[A_3,A_1]$, and $T_3=[A_1,A_2]$.
The commutator identities give the required $C^8$ metric and
potential bounds on fixed enlargements. Proposition~\ref{prop:alltiltcollar}
with $k=16$ requires comparison order $24<60$ and supplies the
background, foliation, gradient, and Jacobi inverse bounds for
Proposition~\ref{prop:potentialbracket} with $k=8$. Apply the latter
to each pair $A_i$ to obtain both normal bounds in $C^8$.
The second bracket preserves time independence and model comparison
through at least nine derivatives. Proposition~\ref{prop:commongram}
now supplies one constant positive definite matrix $Q_*$.
Set $Z=Q_*^{-1/2}T$ and absorb the fixed matrix factors into the
constants. This gives the whole-leaf Gram estimate and all other
asserted bounds.
\end{proof}

At final time $t_*$ and model line center $z_*$, we use the frozen Bryant
coordinate $B_{t_*}$ and $z-z_*$. A \textbf{product domain cap rotation triple} has
metric Lie derivative and potential pairing bounds on
\[
 D_* =\{B_{t_*}\le B,\ |z-z_*|\le H\},
\]
and CMC-normal and area-squared normalized Gram bounds on the actual collar
leaves contained in
\[
 C_* =\{\Lambda/2\le B_{t_*}\le2J\Lambda,\ |z-z_*|\le H\},
\]
for times in $[t_*-H,t_*]$, in tip units. Here $B,H,J,\Lambda$
are fixed before the heat parameter $E$. We also use the spatial enlargements
\begin{equation}\label{eq:cap-frozen-domain-margins}
 \begin{aligned}
 D_*^+&:=\{B_{t_*}<B+\Lambda/8,\ |z-z_*|<H+1\},\\
 C_*^+&:=\{3\Lambda/8<B_{t_*}<2J\Lambda+\Lambda/8,
                                  \ |z-z_*|<H+1\}.
 \end{aligned}
\end{equation}
\begin{prop}\label{prop:productdomains}
Fix $\rho\ge100$, $\ell_{\rm neck}\ge10\sqrt{\rho+1}$, and an integer
$2\le k\le8$. Use the neck domains of
Corollary~\ref{cor:intrinsicneckcontraction} and the frozen cap
domains $D_*,C_*,D_*^+,C_*^+$ of
\eqref{eq:cap-frozen-domain-margins}, centered at $(z_*,t_*)$.
There are $\Lambda,J,B,H$, chosen before $E$, such that on
sufficiently accurate enlarged Bryant-product charts, uniformly
for $a_{\rm ax}\in[0,1]$, the following holds.
\begin{enumerate}[label=\textup{(\roman*)},ref=\textup{(\roman*)}]
\item If $3\Lambda/4\le B_{t_*}(x)\le2\Lambda$ and
$z_*=z(x)$, then
\[
 \mathcal P_\rho^+(x,t_*)\subset D_*\times[t_*-H,t_*],
\]
and all whole collar leaves through this neck domain lie in
$C_*\times[t_*-H,t_*]$. A cap triple $U=(U_1,U_2,U_3)$ with
error $e$ through order $k$ on these cap domains restricts to a
neck triple on $\mathcal P_\rho(x,t_*)$ with error at most
$C_k\Lambda^{k/2}e$ in neck curvature normalization.
\item Let $\mathcal B_\alpha^+$ be an enlarged core box of
Proposition~\ref{prop:productgrid}, with line center $z_*$ and
final time $t_*$. Then
\[
 \overline{\mathcal B_\alpha^+}\subset D_*\times[t_*-H,t_*].
\]
Its common whole collar spheres lie in $C_*\times[t_*-H,t_*]$;
the connecting paths and their fixed neighborhoods lie in
$D_*\times[t_*-H,t_*]$.
\item After translating $t_*,z_*$ to zero, choose the fixed domains
in Corollary~\ref{cor:alltiltheat} so that $D_*^+\subset D_{\rm out}$
and every whole sphere through $C_*^+$ lies in $C_{\rm out}(t)$
for $-H\le t\le0$. For a fixed margin constant $C_{\rm marg}$, the conditions
\[
 E>B+H+C_{\rm marg},\qquad E^2>H+C_{\rm marg}
\]
put these spacetime domains and their enlargements inside $Q_E$.
For sufficiently small comparison error, the triple of
Corollary \ref{cor:alltiltheat} therefore restricts to $D_*\times[-H,0]$ and
$C_*\times[-H,0]$ with
its stated error bounds.
\end{enumerate}
All restrictions use the same actual CMC spheres. The constants
depend on $\rho,\ell_{\rm neck},k$ and the fixed margins.
\end{prop}
\begin{proof}
We first prove (i) in the exact model and set $t_*:=0,z_*=0$.
The radial potential satisfies
\[
 \partial_t B_t=-|\nabla^{g_B}B_0|_{g_B}^2\circ\Psi_t,\qquad
 0\le B_{t}-B_{0}\le -t\quad(t\le0),
\]
and $B_{t}$ is an increasing function of $B_{0}$. The line coordinate is
independent of time. These facts do not involve $a_{\rm ax}$.

Put $r_p:=\mathrm R(p,t_*)^{-1/2}$ at the cap tip and
$r_x:=\mathrm R(x,t_*)^{-1/2}$ at the neck center, and write
$\vartheta:=r_x/r_p$.
On the overlap band,
$C^{-1}\sqrt\Lambda\le\vartheta\le C\sqrt\Lambda$.
In tip curvature normalization, $\mathcal P_\rho^+(x,t_*)$
has spatial radius at most $L\vartheta$ and backward duration
at most $T\vartheta^2$, where $L,T$ depend only on $\rho,\ell_{\rm neck}$.
The bounds $|\nabla B_{0}|\le1$ and $|\nabla z|=1$ show that its frozen
radial and line coordinates vary by at most $CL\sqrt\Lambda$.
We choose $\Lambda$ so large that the radial variation is less than
$\Lambda/8$ and the collar is beyond the radius in
Proposition~\ref{prop:jacobi}. We choose $J\ge4$, then $B>4J\Lambda$,
and choose $H$ larger, with strict margins, than both
$CL\sqrt\Lambda$ and $T\sup\vartheta^2$ and than the fixed core-box
and alignment-path requirements. The cap's line center is the line
coordinate of $x$, not a point selected by its potential value.
Thus no division by $a_{\rm ax}$ is needed.

Model CMC spheres have constant $B_{t}$ and $z$, hence also constant
$B_{0}$ and $z$. All whole leaves through the prescribed neck domain
therefore lie in the collar, with strict margins.
For the perturbed chart the same containment follows by smooth
comparison and the uniform collar construction. The actual leaves
of cap and neck rotation triples coincide by local small-graph CMC uniqueness,
when the neck rotation triple is specified using this common foliation.

On these collars $|\nabla F|$ has a lower bound independent of $a_{\rm ax}$.
Its inverse and the required derivatives have uniform bounds on the
fixed enlarged annulus and time interval. Under
$g_{\mathrm{neck}}=\vartheta^{-2}g_{\mathrm{cap}}$, an order-$j$ metric Lie derivative
norm acquires $\vartheta^j$, while a unit-normal pairing with $j$
derivatives acquires $\vartheta^{j-1}$. The pairing $UF$ gives the
potential-normal component by division by $|\nabla F|$ only on the
collar, with the same controlled derivative bounds. The normalized
whole-leaf Gram matrix is invariant under this constant metric scaling.
For $j\le k$ these conversions are bounded by $C_k\Lambda^{k/2}$.
All comparisons here are between the same vector fields and leaves;
there is no additive symmetry error from geometric comparison.
This proves (i).

We now prove (ii). A short core box in its own centered frozen
coordinates has $B_{t_*}\le2\Lambda$ and
$|z-z_*|\le Ch$, with time length $Ch$ for a fixed small mesh $h$.
Choose its common collar near $B_{t_*}=\Lambda$, inside
$\{\Lambda/2<B_{t_*}<2J\Lambda\}$. Radial and angular paths joining this collar to the core
stay in a fixed radial enlargement and a fixed line slab. The chosen
$B,H$ contain all these sets. Time translation and line translation
of the model provide the recentered chart; its frozen radial coordinate
is $B_{t_*}$, so the same bounds apply at every such grid center.
On an actual finite comparison chart, center normalization and these
coordinate changes are small perturbations of the model ones.
The strict spatial and time inequalities preserve these containments
for sufficiently small comparison error.

This proves (ii). For (iii), enlarge $D_*^+$ and the whole spheres through $C_*^+$
by fixed spatial and earlier-time margins. This enlarged spacetime
region satisfies, for $-H\le t\le0$,
\[
 B_{t}\le B+H+C,\qquad |z|\le H+C,
\]
with a fixed margin constant $C$. Its collar stays away from
the tip in the current radial coordinate. The second-bracket normal integral
argument needs these collar neighborhoods, but no potential level
cap interiors. We choose $E$ so that $Q_E$ contains this entire fixed
auxiliary spacetime region with the margins required for heat estimates,
freezing, two brackets, and normalization. For example, strict versions
of $E>B+H+C$ and $E^2>H+C$ suffice after enlarging $C$.
We then choose the comparison error to preserve these inclusions.
Corollary~\ref{cor:alltiltheat} gives the required bounds on the
auxiliary region; restricting its fields gives exactly the original
$D_*,C_*,[t_*-H,t_*]$ bounds. The constants $B,H$ therefore remain unchanged.
This proves (iii).
\end{proof}

The proposition proves the local domain containments and restriction
estimates for this definition of cap symmetry. To iterate them, we must
also supply symmetry triples at every required center and verify that
the cap and neck constructions use common actual leaves on their
overlaps. These are the geometric tasks of Section~\ref{sec:exterior}.

\section{Exterior iteration and decay of local symmetry errors}\label{sec:exterior}
In this section, we combine the local improvements into an exterior
iteration and obtain polynomial decay of the local symmetry errors. We first define
intrinsic cap centers and domains in Propositions~\ref{prop:tipcurves}
and~\ref{prop:tipnormal}, then close the cap/neck recurrence in
Proposition~\ref{prop:typedclosure} which is the main result of this section and gives arbitrary polynomial decay of the local symmetry errors.

Throughout this section, let $(M,g,f)$ be a complete steady soliton as in Theorem \ref{thm:main} with
$\Ric=\Hess f$, $\mathrm R>0$, and $\mathrm R+|\nabla f|^2=1$. We put $X=\nabla f$, let $\Phi_t$ be the flow of $-X$, and set
$g_t=\Phi_t^*g$. We assume throughout that $(M,g)$ is not the four-dimensional Bryant soliton.
We first show that the comparison regions in the local improvement
procedures remain outside a prescribed compact set when their base centers are sufficiently remote. We use a lower bound for scalar curvature on that compact set
to estimate how far a local construction can reach towards it.

\begin{lem}
\label{lem:compact-exclusion}
 Fix $L,T,C>0$. For a base point $x$, put
$r_x:=\mathrm{R}(x)^{-1/2}$. Suppose that for some point $y$ and $t$, we have
\[
 -Tr_x^2\le t\le0,\qquad d_{g_t}(x,y)\le Lr_x,\qquad
 \mathrm{R}_{g_t}(y)\le C \mathrm{R}(x).
\]
For every nonempty compact set $K$, put $m_K=\min_K \mathrm{R}>0$. If $\Phi_t(y)\in K$, then
\[
 d_g(x,K)\le D_K:=L\sqrt{C/m_K}+TC/m_K.
\]
\end{lem}
\begin{proof}
If $\Phi_t(y)\in K$, then
\[
 m_K\le\mathrm R(\Phi_t(y))=\mathrm R_{g_t}(y)
       \le C\mathrm R(x),\qquad r_x\le\sqrt{C/m_K}.
\]
Since $|\nabla f|\le1$ and $\Phi_t:(M,g_t)\to(M,g)$ is an isometry,
\[
 \begin{aligned}
 d_g(x,\Phi_t(y))
 &\le d_g(x,\Phi_t(x))+d_g(\Phi_t(x),\Phi_t(y))\\
 &\le |t|+d_{g_t}(x,y)\le Tr_x^2+Lr_x\le D_K.
 \end{aligned}
\]
Thus $d_g(x,K)\le D_K$.
\end{proof}

On any fixed enlarged modeled parabolic domain, the distance, time,
and scalar-comparison bounds in this lemma hold with constants
depending on that domain. The scalar bound must be checked on the
actual comparison domain; an arbitrary distant point is not a permitted comparison center. For a fixed finite collection of local domain types one takes
the largest constants. A recentering through finitely many local
steps can instead be handled by repeated use of the lemma.

\begin{defn}\label{def:recurrence-reachability}
Let $\mathscr R\subset S\times S$ be a relation. A pair
$(u,v)\in\mathscr R$ is a \textbf{recurrence edge} from $u$ to $v$.
We say that $v$ is \textbf{reachable from $u$ in at most $J$ steps}
if there are $0\le m\le J$ and $u_0,\ldots,u_m\in S$ with
\[
 u_0=u,\qquad u_m=v,\qquad
 (u_i,u_{i+1})\in\mathscr R\quad(0\le i<m).
\]
For the canonical soliton flow, we use the following relations.
For $L,T,C>0$, define $\mathscr R_{L,T,C}$ on $M$ by
\[
 \begin{aligned}
 x\,\mathscr R_{L,T,C}\,z\quad\Longleftrightarrow\quad
 &\text{there are }(y,t)\text{ with }z=\Phi_t(y),\\
 &-Tr_x^2\le t\le0,\quad d_{g_t}(x,y)\le Lr_x,\quad
       \mathrm R_{g_t}(y)\le C\mathrm R(x),
 \end{aligned}
\]
where $r_x:=\mathrm R(x)^{-1/2}$.
For a finite family $\mathcal F=\{(L_\ell,T_\ell)\}$ of positive
pairs, define $\mathscr R_{\mathcal F}^{\rm par}$ on spacetime by
\[
 \begin{aligned}
 (x,s)\,\mathscr R_{\mathcal F}^{\rm par}\,(y,t)
 \quad\Longleftrightarrow\quad
 &s-T_\ell r_{x,s}^2\le t\le s,\\
 &d_{g(s)}(x,y)\le L_\ell r_{x,s}
                   \quad\text{for some }\ell,
 \end{aligned}
\]
where $r_{x,s}:=\mathrm R_{g(s)}(x)^{-1/2}$.
\end{defn}

\begin{cor}\label{cor:finite-exclusion}
Fix $o\in M$, $L,T,C>0$, $\bar C>\max\{1,C\}$, and a nonempty
compact set $K$. There are compact sets $K_j$ and $D_*<\infty$ such that
\[
 K\subset K_0\subset K_1\subset\cdots,\qquad
 \bigcup_{j=0}^\infty K_j=M,\qquad
 K_j\subset\overline B_g(o,D_*\bar C^j).
\]
Every compact subset of $M$ lies in $K_j$ for all sufficiently large
$j$. If $x\notin K_j$, then every point reachable from $x$ in at most
$j$ steps of $\mathscr R_{L,T,C}$ lies outside $K$.
The constant $D_*$ depends on $o,K,L,T,\bar C$.
\end{cor}
\begin{proof}
Put $m_0:=\min_K\mathrm R>0$ and choose $D_0\ge1$ with
$K\subset\overline B_g(o,D_0)$. Define
\[
 \begin{aligned}
 m_j&:=m_0\bar C^{-j},\\
 D_{j+1}&:=D_j+Lm_{j+1}^{-1/2}+Tm_{j+1}^{-1},\\
 K_j&:=\{x:d_g(o,x)\le D_j,\ \mathrm R(x)\ge m_j\}.
 \end{aligned}
\]
Completeness makes each $K_j$ compact. Since $D_j$ increases and
$m_j$ decreases, $K\subset K_0\subset K_1\subset\cdots$.
Consider a chain $x_0,\ldots,x_m$ in $\mathscr R_{L,T,C}$, as in
Definition~\ref{def:recurrence-reachability}. If $x_{i+1}\in K_\ell$,
the defining curvature inequality gives
\[
 \mathrm R(x_i)\ge\frac{\mathrm R(x_{i+1})}{C}
                    \ge m_{\ell+1}.
\]
The distance estimate in Lemma~\ref{lem:compact-exclusion} then gives
\[
 d_g(x_i,x_{i+1})\le Lm_{\ell+1}^{-1/2}+Tm_{\ell+1}^{-1},
 \qquad x_i\in K_{\ell+1}.
\]
By contraposition, $x_0\notin K_j$ implies
$x_i\notin K_{j-i}$ for $0\le i\le m$. Since
$K\subset K_0\subset K_{j-i}$, all these points lie outside $K$.

Summing the geometric series yields
\[
 D_j\le D_0+
 \frac{L\sqrt{\bar C}}{\sqrt{m_0}(\sqrt{\bar C}-1)}\bar C^{j/2}
 +\frac{T\bar C}{m_0(\bar C-1)}\bar C^j
 \le D_*\bar C^j.
\]
Finally, $D_j\to\infty$ and $m_j\to0$. Every compact subset has
bounded distance from $o$ and positive minimum scalar curvature,
so it lies in $K_j$ for all sufficiently large $j$.
\end{proof}

With the rotation triple domains compatible, the local improvement estimates
can be iterated along the exterior. The following recurrence
quantifies the resulting decay.

\begin{prop}\label{prop:polynomialrecurrence}
Fix $o\in M$ and $L,T,C>0$. For each $x\in M$, put
$r_x:=\mathrm R(x)^{-1/2}$ and choose a nonempty set
\[
 \mathcal A_x\subset
 \{(y,t):-Tr_x^2\le t\le0,\ d_{g_t}(x,y)\le Lr_x,
                         \ \mathrm R_{g_t}(y)\le C\mathrm R(x)\}.
\]
Define its representatives in the original soliton by
\[
 \mathcal Z(x):=\{\Phi_t(y):(y,t)\in\mathcal A_x\}\subset M.
\]
Let $K$ be nonempty and compact, and let $e:M\to[0,\infty)$ satisfy
$e\le e_0<\infty$ on $M\setminus K$. Suppose, for some $0<\theta<1$,
\[
 e(x)\le\theta\sup_{z\in\mathcal Z(x)}e(z)
 \quad\text{whenever }\{x\}\cup\mathcal Z(x)\subset M\setminus K.
\]
Then, for every $\bar C>\max\{1,C\}$, outside a compact set, we have
\[
 e(x)\le C_1e_0(1+d_g(o,x))^{-p},\qquad\mathrm{where}\qquad
 p:=\frac{-\log\theta}{\log\bar C}.
\]
More generally, fix an integer $J\ge1$ and set
\[
 \mathcal Z_{\le J}(x):=
 \{x_m:1\le m\le J,\ x_0=x,\ x_{i+1}\in\mathcal Z(x_i)
                                      \ (0\le i<m)\}.
\]
If the assumed inequality instead
uses $\sup_{z\in\mathcal Z_{\le J}(x)}e(z)$ whenever
$\{x\}\cup\mathcal Z_{\le J}(x)\subset M\setminus K$, the same
conclusion holds with $p=-\log\theta/(J\log\bar C)$.
\end{prop}
\begin{proof}
Use $K_j,D_*$ from Corollary~\ref{cor:finite-exclusion}, with
$K\subset K_0$. Its proof gives
\[
 x\notin K_{j+1}\quad\Longrightarrow\quad
 \mathcal Z(x)\subset M\setminus K_j.
\]
Since $e\le e_0$ outside $K_0$, induction yields
$e(x)\le\theta^je_0$ on $M\setminus K_j$.
For $d:=d_g(o,x)\ge2D_*$, take
\[
 j:=\left\lfloor\frac{\log(d/(2D_*))}{\log\bar C}\right\rfloor.
\]
Then $K_j\subset\overline B_g(o,d/2)$, so $x\notin K_j$, and
\[
 e(x)\le\theta^je_0\le\theta^{-1}e_0(d/(2D_*))^{-p}.
\]
For $\mathcal Z_{\le J}$, the same corollary gives
$\mathcal Z_{\le J}(x)\subset M\setminus K_{Jj}$ when
$x\notin K_{J(j+1)}$. Induction gives $e\le\theta^je_0$ outside
$K_{Jj}$; use $J\log\bar C$ in the denominator of the displayed
choice of $j$ to obtain the stated exponent.
\end{proof}

\begin{prop}\label{prop:matchedscale}
In the tip-normalized Bryant product with canonical metric
$g_t=g_B(t)+\mathrm dz^2$ and radial function $B_t$, for $v\ge0$, we put
\[
 \mathrm R_{B,\mathrm{pot}}(v):=\mathrm R_B(B^{-1}(v)),\qquad
 b_*:=2\Lambda,\qquad Q_C:=\mathrm R_{B,\mathrm{pot}}(b_*).
\]
Let $\mathcal T:=\{o_B\}\times\mathbb R$ be the tip line.
A cap center $p\in\mathcal T$ has assigned scale $Q_C$;
a neck center $(y,t)$ with $B_t(y)\ge b_*$ has assigned scale
$Q_N(y,t):=\mathrm R_{B,\mathrm{pot}}(B_t(y))$.
For any $(y,t)$ define
\[
 Q_{\mathrm{req}}(y,t):=
 \begin{cases}
 Q_N(y,t),&B_t(y)\ge b_*,\\
 Q_C,&B_t(y)<b_*.
 \end{cases}
\]
In the second case, the associated cap center is
$\pi_{\mathcal T}(y):=(o_B,z(y))$, rather than $y$.
Define $Q_{\mathrm{base}}(x_0,0)$ at the chosen time-zero base center $x_0$ by
\[
 Q_{\mathrm{base}}(x_0,0):=
 \begin{cases}
 Q_C,&\text{for a cap center }p\in\mathcal T,\\
 Q_N(x,0),&\text{for a neck center }x\text{ with }B_0(x)\ge b_*.
 \end{cases}
\]
Fix $L\ge1$. For sufficiently large $\Lambda$, the estimate
\[
 Q_{\mathrm{req}}(y,t)
 \le(1+CL/\sqrt\Lambda)Q_{\mathrm{base}}(x_0,0)
\]
holds with $x_0=p$ for every $t\le0$ and $y$ in the cap center
case, and with $x_0=x$ for
all $(y,t)$ satisfying
\[
 t\le0,\qquad d_{g_t}(x,y)\le LQ_N(x,0)^{-1/2}
\]
in the neck center case. The constant $C$ is independent of the
potential tilt.
\end{prop}
\begin{proof}
The Bryant asymptotics and soliton identity give, for large $v$,
\[
 (Cv)^{-1}\le\mathrm{R}_{B,\mathrm{pot}}(v)\le C/v,\qquad
 0\le-{\frac{\mathrm d}{\mathrm d v}\mathrm{R}_{B,\mathrm{pot}}}(v)\le C/v^2.
\]
Indeed, $v$ is comparable to geodesic radius, and
${\frac{\mathrm d}{\mathrm d v}\mathrm{R}_{B,\mathrm{pot}}}(v)=-2\Ric(\partial_s,\partial_s)$;
the radial Ricci component is $O(s^{-2})$. In particular, $|\partial_v\log\mathrm{R}_{B,\mathrm{pot}}|\le C/v$.
For a cap base center, monotonicity gives
\[
 Q_{\mathrm{req}}(y,t)\le\mathrm R_{B,\mathrm{pot}}(b_*)
       =Q_C=Q_{\mathrm{base}}
\]
for every $y,t$.

For a neck center write $v=B_{0}(x)$. Since $B_{t}(x)\ge v$ at earlier
times and $|\nabla^{g_t} B_{t}|\le1$, we have
\[
 B_{t}(y)\ge v-L/\sqrt{\mathrm{R}_{B,\mathrm{pot}}(v)}\ge v-C L\sqrt v.
\]
We choose $\Lambda$ large relative to $L^2$, so the last expression is
at least $v/2$. If $y$ is a neck comparison center, integration of the logarithmic
derivative bound gives
\[
 \frac{\mathrm{R}_{B,\mathrm{pot}}(B_{t}(y))}{\mathrm{R}_{B,\mathrm{pot}}(v)}
 \le \exp(C L/\sqrt v)\le1+C L/\sqrt\Lambda.
\]
If $y$ lies below the threshold, then
$v-C L\sqrt v\le b_*$, while $v\ge b_*$. Solving this quadratic
inequality in $\sqrt v$ yields
$0\le v-b_*\le C L\sqrt\Lambda$ after enlarging the fixed lower
bound for $\Lambda$. The same logarithmic derivative estimate gives
$\mathrm{R}_{B,\mathrm{pot}}(b_*)/\mathrm{R}_{B,\mathrm{pot}}(v)\le1+C L/\sqrt\Lambda$.
This is precisely the cap replacement bound.
The assertion follows.
\end{proof}

The proofs of Lemma~\ref{lem:compact-exclusion} and
Corollary~\ref{cor:finite-exclusion} also apply to a positive assigned
scale $Q$ if
\[
 -TQ(x)^{-1}\le t\le0,\qquad
 d_{g_t}(x,y)\le LQ(x)^{-1/2},\qquad
 Q(\Phi_t(y))\le CQ(x),
\]
and $Q$ has a positive lower bound on each compact set of centers.
For cap and neck centers, retain the type of each center and use
the original-soliton distance between their spatial representatives.

\begin{prop}
\label{prop:sublinearcurvscale}
Let $(M^4,g,f)$ satisfy the hypotheses of Theorem~\ref{thm:main}
and \eqref{eq:normalization}. Then for every
fixed $o\in M$, we have
\[
 d_g(o,x)\sqrt{\mathrm{R}(x)}\longrightarrow\infty
 \quad\text{as }x\to\infty.
\]
Equivalently, $\mathrm{R}(x)^{-1/2}=o(d_g(o,x))$.
\end{prop}

\begin{proof}
Suppose there are $x_i\to\infty$ and $A<\infty$ with
$d_g(o,x_i)\sqrt{\mathrm R(x_i)}\le A$. Since
$d_g(o,x_i)\to\infty$,
\[
 0<\mathrm R(x_i)\le\frac{A^2}{d_g(o,x_i)^2}\longrightarrow0.
\]
Set $g_i:=\mathrm R(x_i)g$. By the dimension reduction compactness
theorem \cite[Theorem 1.2]{CMZ}, a subsequence of $(M,g_i,x_i)$
converges smoothly to a complete pointed limit. Moreover,
$d_{g_i}(x_i,o)\le A$. Choose a limit ball of radius greater than
$A+1$ and a larger relatively compact comparison chart. Smooth
convergence bounds scalar curvature on the corresponding charts.
A minimizing $g_i$-geodesic of length at most $A$ from $x_i$
cannot leave the comparison image: at its first exit, metric
comparison with the larger limit ball would force length greater
than $A$. Thus $o$ lies in the controlled image, and
$\mathrm R_{g_i}(o)$ is uniformly bounded.

For a scale factor $\lambda_{\rm sc}>0$, scalar curvature satisfies
$\mathrm R_{\lambda_{\rm sc}g}=\lambda_{\rm sc}^{-1}\mathrm R_g$. Hence
\[
 \mathrm R_{g_i}(o)=\frac{\mathrm R_g(o)}{\mathrm R_g(x_i)}
                   \longrightarrow\infty,
\]
since $\mathrm R_g(o)>0$ is fixed. This is a contradiction.
\end{proof}

We next define cap centers intrinsically, so that centers obtained
in overlapping product charts can be compared. In the non-Bryant
alternative, the Ricci eigenspaces and the scalar curvature gradient
determine an exterior set of tip curves. This is analogous to the
intrinsic use of tips in Brendle's analysis of positively curved
ancient solutions; see \cite[Proposition~9.1 and Definition~9.2]{BrendleAncient}.

\begin{prop}
\label{prop:tipcurves}
Let $\mu_1\le\mu_2\le\mu_3\le\mu_4$ be the Ricci
endomorphism eigenvalues, and define
\[
 \mathcal U:=\{\mu_2-\mu_1>\mathrm{R}/8\}.
\]
On $\mathcal U$, let $L$ be the smallest Ricci eigenline, let
$E:=L^\perp$, and let $P_E$ be orthogonal projection onto $E$. We set
\[
 \Gamma:=\{x\in\mathcal U:P_E\nabla \mathrm{R}(x)=0\}.
\]
Outside a compact set $K$, $\Gamma$ is a relatively closed smooth
embedded one-dimensional submanifold, possibly empty. Along every
sequence $p_i\in\Gamma$ escaping to infinity, pointed limits of $\mathrm R(p_i)g$ are $B^3\times\mathbb R$ products based at the tip, and
\[
 \frac{\mu_1(p_i)}{\mathrm{R}(p_i)}\to0,\qquad
 \frac{\mu_j(p_i)}{\mathrm{R}(p_i)}\to\frac13\quad(j=2,3,4),\qquad
 \angle(T_{p_i}\Gamma,L_{p_i})\to0.
\]
In a sufficiently accurate fixed-size chart close to a
Bryant product, the local set $\Gamma$ is the unique small graph over
the product tip line. Its $C^k$ graph norm is $O(\eta)$ when the metric
comparison is $O(\eta)$ in $C^{k+6}$, for fixed $k$ and fixed chart
margins.
\end{prop}
\begin{proof}
The gap makes $L$ and $P_E$ smooth, even when eigenvalues within $E$
coincide. We consider any escaping sequence in $\Gamma$.
By \cite[Theorems 1.2, 1.3, 1.5 and Corollary 5.4]{CMZ}, the
non-Bryant case has $S^2\times\mathbb R^2$ as its tangent flow
at infinity, and its curvature scale reductions are $B^3\times\mathbb R$ products or shrinking $S^2\times\mathbb R^2$ cylinders. We pass to one of these subsequential limits. The cylinder has two
zero Ricci eigenvalues, contradicting the limiting gap at least $1/8$
when basepoint scalar curvature is one. In a Bryant product, $E$ is
the Bryant tangent space and $P_E\nabla \mathrm{R}=\nabla \mathrm{R}_B$. Its only zero
is the Bryant tip. This proves the asserted location of every limit
and the eigenvalue limits.

At the normalized Bryant tip, $\Ric_B=\frac13g_B$ and $\nabla B(s)=0$.
Differentiating $\mathrm{d}\mathrm{R}_B=-2\Ric_B(\nabla B(s),\cdot)$ therefore gives
\[
 \Hess \mathrm{R}_B=-2\Ric_B\circ\Ric_B=-\tfrac29g_B
\]
on the three Bryant directions. We use smooth normal coordinates $y$ on
a fixed small Bryant ball and the product coordinate $z$ on a fixed
line interval. Trivialize $E$ by projection from the three model
coordinate directions. We write the components of the section as
$\mathcal S(y,z):=(P_E\nabla \mathrm{R})(y,z)$. For the model, we have
\[
 \mathcal S_B(0,z)=0,\qquad
 D_y\mathcal S_B(0,z)=-\frac29I_3,\qquad
 \partial_z\mathcal S_B(0,z)=0.
\]
We choose a fixed small transverse radius $\delta>0$. After decreasing
$\delta$ and then the comparison error, the maps
\[
 \mathcal T_z(y)=y+\frac92\mathcal S(y,z)
\]
satisfy $\|D_y\mathcal T_z\|\le1/2$ and
$|\mathcal T_z(0)|\le C\eta$, uniformly on the fixed line interval.
They preserve $\{|y|\le2C\eta\}$ when $2C\eta<\delta$.
The contraction theorem therefore gives a unique zero $y=\psi(z)$
with $|\psi(z)|\le2C\eta$. The derivative formula is
\[
 \frac{\mathrm d \psi}{\mathrm dz}(z)
 =-\bigl(D_y\mathcal S(\psi(z),z)\bigr)^{-1}
                                      \partial_z\mathcal S(\psi(z),z).
\]
The inverse is uniformly bounded. Differentiating this identity
repeatedly gives $\|\psi\|_{C^k}\le C_k\eta$ on smaller fixed
intervals; the assumed $C^{k+6}$ metric control bounds every
coefficient and leaves a derivative margin. In particular, we have
\[
 T\Gamma=\operatorname{span}\{\partial_z+\frac{\mathrm d \psi}{\mathrm dz}(z)\cdot\partial_y\}
           \longrightarrow\operatorname{span}\{\partial_z\},
\]
as does the smallest Ricci eigenline $L$.
The graph lies strictly inside the spectral gap region after shrinking
fixed margins. This proves the local assertions.

Applied to the escaping sequences above, these estimates show that
there is a compact set outside which all zeros have rank-three transverse derivative
and satisfy $\mu_2-\mu_1>\mathrm{R}/4$. Otherwise a sequence of failures would
converge to the normalized tip, where both conclusions hold. Thus all
zeros outside this compact set are smooth curves and stay away from the boundary of
$\mathcal U$. Any limit point in $M\setminus K$ still has gap at least
$\mathrm{R}/4$ and still satisfies $P_E\nabla \mathrm{R}=0$. Hence the zero set is
relatively closed there. Local uniqueness of the graph gives an
embedded submanifold and agreement between overlapping charts.
This proves the stated conclusion.
\end{proof}

A point $p\in\Gamma\setminus K$, with $K$ as in
Proposition~\ref{prop:tipcurves}, is called an
\textbf{intrinsic tip center}.
\begin{cor}\label{cor:tipnaturality}
Let $\Gamma,L,E,K$ be as in Proposition~\ref{prop:tipcurves}.
The definition of $\Gamma$ is invariant under isometry and constant
metric rescaling. Put
\[
 K_t:=\Phi_t^{-1}(K),\qquad \Gamma_t:=\Phi_t^{-1}(\Gamma).
\]
The same construction applied to $g_t$ gives
\[
 \Gamma(g_t)=\Gamma_t,\qquad
 \Gamma_t\setminus K_t=\Phi_t^{-1}(\Gamma\setminus K).
\]
After enlarging $K$ to contain a critical point of $f$ if it exists,
each component of $\Gamma\setminus K$ has a unit tangent $V_0$ with
\[
 \mathrm df(V_0)>0,\qquad
 \mathrm d\mathrm R(V_0)=-2\mu_1\,\mathrm df(V_0)<0.
\]
Every $p\in\Gamma\setminus K$ is a nondegenerate local maximum
of scalar curvature on its regular potential level. No component
of $\Gamma\setminus K$ is a closed curve.
\end{cor}
\begin{proof}
Under $g\mapsto\lambda_{\rm sc}g$, $\lambda_{\rm sc}>0$, Ricci eigenvalues and scalar curvature both
scale by $\lambda_{\rm sc}^{-1}$, the eigenspaces and projection remain unchanged,
and $\nabla \mathrm{R}$ scales by $\lambda_{\rm sc}^{-2}$. The gap condition and the zero set
are therefore preserved. Diffeomorphism naturality gives the canonical
flow assertion.

On $\Gamma$, the soliton identity and commutation of $P_E$ with the
Ricci endomorphism give
\[
 0=P_E\nabla \mathrm{R}=-2\Ric^\sharp P_EX.
\]
Positive Ricci curvature implies $P_EX=0$, so $X\in L$ and
$\nabla \mathrm{R}=-2\mu_1 X$. Positive Ricci curvature also gives at most one
critical point of $f$, which can be included in $K$ if present.
The angle estimate makes the projection of $X$ onto $T\Gamma$
nonzero outside a further compact set. We normalize that projection to
orient the curve. Then $\mathrm{d}f(V_0)>0$ and the displayed scalar derivative
identity follows. A closed component would have strictly positive
integral of $\mathrm{d}f$ over one traversal, contradicting exactness.
On the regular level $\Sigma=\{f=f(p)\}$ through $p\in\Gamma$,
$T_p\Sigma=E_p$ and $\mathrm{d}\mathrm{R}|_{T_p\Sigma}=0$. The hypersurface Hessian
formula gives
\[
 \Hess_{\Sigma}\mathrm{R}|_p
   =\Hess_g\mathrm{R}|_{E_p}+2\mu_1\Ric|_{E_p}.
\]
In curvature units this tends to $-\frac29g_B$ along escaping tip
centers, so it is negative definite outside a further compact set.
The cancellation in this formula uses
$\mathrm{d}\mathrm{R}(X/|X|)=-2\mu_1|X|$ and the second fundamental form
$\Ric|_{T\Sigma}/|X|$; it does not require a uniform lower bound for
$|X|$. This proves the local maximum assertion.
The eigenline $L$ need only be close to $T\Gamma$; the proof does not
identify them.
This finishes the proof.
\end{proof}

We next construct a normal tube around these tip curves and its
nearest-point projection.

\begin{prop}\label{prop:tipnormal}
Let $\Gamma$ be as in Proposition~\ref{prop:tipcurves}.
For every fixed $R_{\rm tube}>0$ there is a compact set
$K_{R_{\rm tube}}$ such that, with
$\Gamma_{R_{\rm tube}}:=\Gamma\setminus K_{R_{\rm tube}}$, the map
\[
 \mathcal E(p,v):=\exp_p(v),\qquad \mathrm{where}\qquad
 p\in\Gamma_{R_{\rm tube}},\quad v\perp T_p\Gamma,\quad
 |v|<R_{\rm tube}\mathrm R(p)^{-1/2},
\]
is injective and nonsingular. Its image $\mathcal T_0(R_{\rm tube})$
is called an \textbf{intrinsic tube}. It has a smooth projection
$\pi_0(\mathcal E(p,v))=p$, and
\[
 d_g(\mathcal E(p,v),\Gamma)=|v|,\qquad
 \pi_0(x)=\operatorname*{argmin}_{q\in\Gamma}d_g(x,q).
\]
Here $\operatorname*{argmin}_{p\in\Gamma}d_g(x,p)$ denotes
the point $p$ where the distance attains its minimum, rather than
the minimum distance itself; Proposition~\ref{prop:tipnormal} asserts uniqueness on the stated tube.

The minimizer is unique. In curvature-normalized charts centered at
$p_i\in\Gamma_{R_{\rm tube}}$, $p_i\to\infty$, these maps converge
smoothly to the normal exponential map and projection of the
Bryant-product tip line. If $x_i=\mathcal E(p_i,v_i)$ with
$|v_i|<R_{\rm tube}\mathrm R(p_i)^{-1/2}$, then $x_i\to\infty$.
\end{prop}
\begin{proof}
At an escaping sequence of bases, the curvature-normalized metrics
converge to a tip-centered Bryant product, and the tip curves converge
smoothly to its tip line. The normal exponential map of that line is
\[
 (z,v)\longmapsto(\exp_o^{g_B}v,z).
\]
It is a diffeomorphism: polar coordinates of the rotational Bryant
metric are globally regular away from the tip, with $\omega(s)>0$ for
$s>0$, and the exponential map is smooth and nonsingular at the tip.
On any fixed compact subset its differential has a positive lower
singular-value bound. Smooth dependence of geodesics and Jacobi fields
therefore gives a compact exclusion, depending on $R_{\rm tube}$,
outside which $\mathcal E$ has no singular differential for
$|v|\le R_{\rm tube}\mathrm R(p)^{-1/2}$.

To prove injectivity, suppose otherwise that two different normal data
$(p_i,v_i),(q_i,w_i)$ with bases escaping satisfy
$\exp_{p_i}v_i=\exp_{q_i}w_i$. We write $r(p):=\mathrm{R}(p)^{-1/2}$ and order the
bases so that $r(q_i)\ge r(p_i)$. Then we obtain
\[
 d(p_i,q_i)\le R_{\rm tube}(r(p_i)+r(q_i))\le2R_{\rm tube} r(q_i).
\]
In the pointed limit normalized at $q_i$, scalar curvature is uniformly
bounded on that fixed ball. Hence $\mathrm{R}(p_i)/\mathrm{R}(q_i)$ is bounded, so the
scales are comparable. Any subsequential limit of $p_i$ is a zero of
the intrinsic projected-gradient section in the Bryant product, and
therefore lies on its tip line. In particular, $\mathrm{R}(p_i)/\mathrm{R}(q_i)\to1$. Both sets of normal data now lie in one fixed
converging normal bundle chart. Their limiting exponential images agree.
Model injectivity makes the two limiting data equal. The uniform local
inverse theorem near those data then contradicts their being distinct
for large $i$. This proves injectivity on the stated normal bundle. We work with slightly
larger fiber bounds throughout to retain interior margins.

If $x_i=\exp_{p_i}v_i$ with $p_i\to\infty$ and
$|v_i|\le R_{\rm tube} r(p_i)$, Proposition~\ref{prop:sublinearcurvscale} gives
\[
 d(o,x_i)\ge d(o,p_i)-R_{\rm tube} r(p_i)\longrightarrow\infty.
\]
It also shows that a fixed compact part of the tip set cannot supply
a closer point to $x_i$ for large $i$. Choose a fixed compact set $K_1$ with $K\subset\operatorname{int}K_1$.
Relative closedness from Proposition~\ref{prop:tipcurves} makes
$\Gamma\cap(M\setminus\operatorname{int}K_1)$ closed. A nearest
point $q_i$ in this set exists by completeness: $p_i$ is a candidate
for large $i$, and the competing points lie in a compact ball. Since
$d(q_i,x_i)\le|v_i|$, we have
$d(p_i,q_i)\le2R_{\rm tube}r(p_i)$.
Proposition~\ref{prop:sublinearcurvscale} now implies $q_i\to\infty$.
Thus $q_i$ is an interior point of the smooth curve, rather than a
point on the truncation boundary, and its minimizing geodesic to
$x_i$ is normal to $\Gamma$. Both bases escape; convergence on this fixed
ball and the intrinsic zero-set condition give
$\mathrm{R}(q_i)/\mathrm{R}(p_i)\to1$. The competing normal vector consequently has
length less than $2R_{\rm tube} r(q_i)$ for large $i$. We apply the injectivity result
with fiber bound $2R_{\rm tube}$ to conclude $q_i=p_i$ and that the two normal
vectors agree. A contradiction sequence proves the nearest-point assertion for all
$p\in\Gamma_{R_{\rm tube}}$. Smooth projection and convergence
follow from the local inverse theorem.
The assertion follows.
\end{proof}

Fix $R_{\rm tube}>0$ and $K_{R_{\rm tube}}$ from
Proposition~\ref{prop:tipnormal}. With
$\Gamma_{R_{\rm tube}}:=\Gamma\setminus K_{R_{\rm tube}}$, define
\begin{equation}\label{eq:intrinsic-tube-domain}
 \begin{aligned}
 \mathcal T_0(R_{\rm tube})
 &:=\{\exp_pv:p\in\Gamma_{R_{\rm tube}},\ v\perp T_p\Gamma,
                     \ |v|<R_{\rm tube}\mathrm R(p)^{-1/2}\},\\
 \mathcal T_t(R_{\rm tube})&:=\Phi_t^{-1}(\mathcal T_0(R_{\rm tube})).
 \end{aligned}
\end{equation}
The projections are
\[
 \pi_0:\mathcal T_0(R_{\rm tube})\to\Gamma_{R_{\rm tube}},\qquad
 \pi_t:=\Phi_t^{-1}\circ\pi_0\circ\Phi_t
 :\mathcal T_t(R_{\rm tube})\to\Phi_t^{-1}(\Gamma_{R_{\rm tube}}).
\]
Let $B(s)$ be the tip-zero potential of the tip-normalized Bryant
soliton. Define
\begin{equation}\label{eq:ambientbeta}
 \beta_0(x):=B\left(\sqrt{\mathrm R_g(\pi_0(x))}\,d_g(x,\Gamma)\right),
 \qquad\mathrm{where}\qquad x\in\mathcal T_0(R_{\rm tube}).
\end{equation}
When the chosen radius is fixed, we abbreviate
$\mathcal T_t:=\mathcal T_t(R_{\rm tube})$.
The notation $\beta,\pi$ means $\beta_0,\pi_0$ at time zero.
Since $\Phi_t:(M,g_t)\to(M,g)$ is an isometry, we can also define
\[
 \begin{aligned}
 \beta_t(x):=\beta_0(\Phi_t(x))=B\left(\sqrt{\mathrm R_{g_t}(\pi_t(x))}\,
                       d_{g_t}(x,\Gamma_t)\right),
 \qquad\mathrm{where}\qquad x\in\mathcal T_t(R_{\rm tube}).
 \end{aligned}
\]
This identity gives the intrinsic radial coordinate for $g_t$.
The functions are smooth at the tip centers because $B(s)$ is
smooth and even at $s=0$, and the squared normal vector norm is
smooth on the normal bundle. On fixed normalized product charts,
$\beta_t$ converges smoothly to the model Bryant potential.
A local line coordinate is signed arclength along $\Gamma$ from a
chosen center to $\pi_0(x)$, multiplied by the square root of the
center's scalar curvature.

For $b>0$, the notation $\{\beta_t<b\}$ means
$\{x\in\mathcal T_t(R_{\rm tube}):\beta_t(x)<b\}$; points outside
$\mathcal T_t(R_{\rm tube})$ do not belong to this set.
To include all radial ranges used in an argument, first fix their
largest bound $b_{\max}$ and a $\delta_b>0$, and choose
\[
 R_{\rm tube}>B^{-1}(b_{\max}+\delta_b).
\]
Then choose $K_{R_{\rm tube}}$ from
Proposition~\ref{prop:tipnormal}. In every normal fiber,
$\beta_t\le b_{\max}$ now implies
$\sqrt{\mathrm R_{g_t}(\pi_t(x))}\,d_{g_t}(x,\Gamma_t)
\le B^{-1}(b_{\max})<R_{\rm tube}$.
All later compact sets may be enlarged after this radius is fixed.
On the overlap of two such tubes the projections agree, since both
are the unique nearest-point projection onto $\Gamma$. Hence the
coordinates $\beta_t$ also agree there. Increasing the tube radius
therefore preserves these coordinates wherever both constructions
are defined, after taking a common compact exclusion.

\begin{prop}
\label{prop:corecover}
Fix
$0<\eta_{\rm cyl}<1/2$, with cylindrical neighborhoods understood
as in Definition~\ref{cylinder:def:neighborhood}.
There are $b_0<\infty$ and a compact set $K$ such that every
$x\notin K$ which is not $\eta_{\rm cyl}$-cylindrical lies in the
set $\mathcal T_0(R_{\rm tube})$ defined in \eqref{eq:intrinsic-tube-domain} and satisfies
\[
 \beta_0(x)<b_0.
\]
Choose the tube radius to contain this radial range with a fixed
margin. For every fixed $b_*>b_0$, all points outside $K$ and the
core $\{\beta_0<b_*\}$ are $\eta_{\rm cyl}$-cylindrical.
\end{prop}
\begin{proof}
On the exact Bryant product, points whose tip-normalized potential
tends to infinity converge at their own curvature scale to the round
shrinking $S^2\times\R^2$ cylinder on every fixed backward parabolic
domain. This is the usual cylindrical asymptotic of the Bryant model;
in the normalized radial coordinates its spherical coefficient tends
to $2(1-t)$ and its radial coefficient to one, for fixed $t\le0$,
using $\omega(s)^2\sim2s$, $\mathrm{R}(s)\sim s^{-1}$, and $\frac{\mathrm d B}{\mathrm ds}(s)\to1$.
Here the center time is zero; shifting normalized time by $-1$
gives the cylinder $(-2t)g_{\sphere}+\mathrm dz_1^2+\mathrm dz_2^2$
centered at $t=-1$ in Definition~\ref{cylinder:def:neighborhood}.
The differentiated Bryant asymptotics give smooth convergence on compact
backward intervals. We choose $b_0$ large enough, with strict comparison and radial margins,
that all such model points beyond a smaller threshold already satisfy
the prescribed cylindrical comparison.

If the assertion failed, take an escaping sequence of non-neck points.
The dimension reduction alternatives of \cite{CMZ} give a subsequential
cylindrical or Bryant-product limit at their curvature scale. A cylinder
limit contradicts the non-neck assumption, using a slightly better
comparison accuracy in the limit. In a Bryant-product limit the basepoint
must therefore have tip-normalized radial potential below $b_0$, by
the choice just made. The model tip is at bounded scaled distance,
and smooth convergence of the intrinsic zero set supplies nearby
actual tip centers. Normal-exponential convergence on the enlarged
fixed chart then puts the original basepoints inside the actual tube,
with $\beta_0<b_0$. This contradicts the asserted failure.
For the resulting tip centers $p_i$, fixed-chart convergence gives
\[
 C^{-1}\mathrm R(x_i)\le\mathrm R(p_i)\le C\mathrm R(x_i),\qquad
 d_g(p_i,x_i)\le C\mathrm R(x_i)^{-1/2}.
\]
If $p_i$ stayed in a compact set, these inequalities would bound
$\mathrm R(x_i)^{-1/2}$ and $d_g(p_i,x_i)$, contradicting
$x_i\to\infty$. Thus $p_i\to\infty$, and the exterior tube
results apply.
This proves the stated conclusion.
\end{proof}

\begin{cor}
\label{cor:outsidecorescalar}
Let $\Gamma$
and $\beta_0$ be the tip set and tube coordinate of
Proposition~\ref{prop:tipcurves} and \eqref{eq:ambientbeta}.
Fix a core threshold $b_*$ and write $\rho_* =\mathrm{R}_{B,\mathrm{pot}}(b_*)$ for the
scalar curvature at that potential value in the tip-normalized Bryant
model. Then
\[
 \limsup_{\substack{x\to\infty\\x\notin\{\beta_0<b_*\}}}\mathrm{R}(x)
 \le\rho_*.
\]
In particular, for a sufficiently large fixed threshold, outside the
core and a compact set one has $\mathrm{R}\le2\rho_*<1/2$ and
$|\nabla f|\ge1/\sqrt2$. 
\end{cor}
\begin{proof}
Suppose an escaping sequence outside the core has
$\mathrm{R}(x_i)\to q>\rho_*$. The normalized potentials
$F_i=\sqrt{\mathrm{R}(x_i)}(f-f(x_i))$ converge on the pointed limit, after a
subsequence, and satisfy
\[
 \Hess F_\infty=\sqrt q\,\Ric_\infty,\qquad
 |\nabla F_\infty|^2=1-q\mathrm{R}_\infty.
\]
A cylindrical limit is
impossible: integrate the tangential Hessian equation on its whole
round spherical factor. In the Bryant-product alternative, let $v:=B(s)$
be the basepoint's radial potential in tip units. The intrinsic tube
convergence shows $v\ge b_*$, since otherwise the original points
would lie inside the core. The tip scalar in basepoint-normalized
units is $1/\mathrm{R}_{B,\mathrm{pot}}(v)$. Evaluating the nonnegative squared gradient
identity there gives $q\le\mathrm{R}_{B,\mathrm{pot}}(v)\le\rho_*$, a contradiction.
This proves the limsup bound. We choose $b_*$ large enough that
$2\rho_*<1/2$ and enlarge the compact set; the normalized soliton
identity then gives the gradient bound.
This finishes the proof.
\end{proof}

This gives the exterior core/neck cover. We next prove that a neck comparison center entering the core lies in the restriction band, and verify its
assigned scale at the correct time.

By Corollary~\ref{cor:tipnaturality}, we have
$\Gamma_t:=\Gamma(g_t)=\Phi_t^{-1}(\Gamma(g))$. We use its intrinsic tube projection
$\pi_t$ and radial coordinate $\beta_t$. All exterior exclusions below
are imposed on original-soliton representatives under $\Phi_t$.
We fix $b_*:=2\Lambda$, put
$\rho_*:=\mathrm{R}_{B,\mathrm{pot}}(b_*)$, and define
\begin{equation}
 \mathcal C_t:=\{y\in\mathcal T_t:\beta_t(y)<2\Lambda\}
              =\Phi_t^{-1}(\mathcal C_0).
\end{equation}
An exterior point
outside $\{\beta_t<b_*\}$ is a neck center. A point inside that core
is replaced, when a neck comparison center is requested there, by the cap center
$\pi_t(y)$. Assign
\[
 Q_N(y,t)=\mathrm{R}_{g_t}(y),\qquad Q_C(p,t)=\rho_*\mathrm{R}_{g_t}(p)
 \quad(p\in\Gamma_t).
\]
These assignments are invariant under canonical pullback. The core
radius is chosen before the comparison compact set.

\begin{prop}
\label{prop:inputcomparison}
Let $\Gamma$
and $\beta_0$ be the tip set and tube coordinate of
Proposition~\ref{prop:tipcurves} and \eqref{eq:ambientbeta}.
For every fixed $L,T<\infty$ there are a compact set $K$ and
$C>0$ such that, for $x\notin K$ and $r_x:=\mathrm R(x)^{-1/2}$,
\[
 C^{-1}\mathrm{R}(x)\le \mathrm{R}_{g_t}(y)\le C\mathrm{R}(x),\qquad
 d_{g_t}(x,y)\le C r_x
\]
whenever $-Tr_x^2\le t\le0$ and $d_g(x,y)\le Lr_x$.
For every nonempty compact set $A\subset M$ there is a compact
set $K_A\supset K$ such that
\[
 x\notin K_A,\quad -Tr_x^2\le t\le0,\quad d_g(x,y)\le Lr_x
 \quad\Longrightarrow\quad\Phi_t(y)\notin A.
\]
\end{prop}
\begin{proof}
A sequence violating a bound has a subsequential smooth
curvature-normalized ancient limit by \cite[Theorem 1.2]{CMZ}.
On the fixed normalized time-zero ball and backward interval, the
limiting metric is smooth and comparable to its time-zero metric,
and scalar curvature is positive. In the non-Bryant case the
limits are the positive-scalar Bryant products and round cylinders.
Smooth convergence on a slightly enlarged ball therefore bounds the
scalar above and below and bounds $g_t$ relative to $g_0$. A minimizing
time-zero geodesic from $x$ to $y$ stays in that ball, so metric
comparison bounds its time-$t$ length as well. This contradiction
proves uniform constants. For a nonempty compact $A$, put
$m_A:=\min_A\mathrm R>0$ and
\[
 D_A:=C\sqrt{C/m_A}+TC/m_A,\qquad
 K_A:=K\cup\{x:d_g(x,A)\le D_A\}.
\]
This set is compact by completeness. Apply
Lemma~\ref{lem:compact-exclusion} with spatial constant $C$, time
constant $T$, and curvature constant $C$. If $\Phi_t(y)\in A$,
it gives $d_g(x,A)\le D_A$, proving the stated implication.
The assertion follows.
\end{proof}

We now transfer this model-scale comparison to the actual soliton.
This gives the controlled neck neighborhood needed at the
cap-to-neck transition.

\begin{prop}
\label{prop:actualtransition}
With the assigned scales just defined, fix a neck comparison radius $L\ge1$, $T>0$, and
$\sigma>0$. For sufficiently large $\Lambda$, and then outside a
sufficiently large compact set of base centers, the following holds.
Let $x$ be a neck center outside the intrinsic core at time zero,
$r_x:=\mathrm{R}(x)^{-1/2}$, and let
\[
 -Tr_x^2\le t\le0,\qquad d_g(x,y)\le Lr_x.
\]
If $y$ is outside the time-$t$ core, then
\[
 Q_N(y,t)\le(1+\sigma)Q_N(x,0).
\]
If $y$ is in the core, put $p:=\pi_t(y)$. Then, for a fixed model
constant $C_2$, we have
\[
 2\Lambda-C_2L\sqrt\Lambda\le\beta_t(y)<2\Lambda,
 \qquad
 2\Lambda\le\beta_0(x)\le2\Lambda+C_2L\sqrt\Lambda,
\]
and
\[
 Q_C(p,t)\le(1+\sigma)Q_N(x,0).
\]
The intrinsic coordinate of $x$ in the second case is defined in a
slightly larger fixed tube. In particular the entering point is in
the cap restriction band when $\Lambda$ is sufficiently large.
The constants may depend on $L,T,\Lambda,\sigma$.
\end{prop}
\begin{proof}
We prove the uniform exterior estimates by contradiction. If an
estimate failed outside every compact set, we could choose a sequence
of violating centers escaping to infinity. Proposition~\ref{prop:inputcomparison}
keeps the earlier-time centers in fixed normalized regions with positive
bounded scalar curvature and makes their original representatives escape.
If $y\in\mathcal C_t$ and $p:=\pi_t(y)$, then
\[
 \Phi_t(p)=\pi_0(\Phi_t(y)),\qquad
 d_g(\Phi_t(y),\Phi_t(p))
 \le B^{-1}(2\Lambda)\mathrm R(\Phi_t(p))^{-1/2}.
\]
The representatives $\Phi_t(y)$ escape by
Proposition~\ref{prop:inputcomparison}. If their projected
representatives $\Phi_t(p)$ stayed in a compact set, the displayed
bound would contradict that escape. Thus all compact exclusions
apply to these projected centers as well.
The tip chart comparison on $\beta_t\le2\Lambda$ bounds
$\mathrm{R}_{g_t}(y)/\mathrm{R}_{g_t}(p)$ above and below by positive constants depending
on $\Lambda$. Hence $p$ has scale comparable to $r_x$ and is at bounded
normalized time-$t$ distance from $x$.

We take the pointed ancient limit at the base point $(x,0)$, on a domain large
enough to contain these tip projections if present. An entering-core
sequence cannot have a cylinder limit: its projected tip would give a
positive-scalar point with Ricci gap bounded below by a positive
multiple of scalar curvature, whereas the cylinder has two zero Ricci
eigenvalues. Thus that case has a Bryant-product limit. If there is no
core entry, a cylinder limit is allowed, and its scalar at earlier
times is at most its time-zero scalar. It already gives the required
neck-to-neck scale inequality with a strict $\sigma$ margin.

It remains to consider a Bryant-product limit. In tip units let
$v:=B_0(x)$ and $u:=B_t(y)$. The base point is outside the core, so $v\ge2\Lambda$:
if the limiting radial value were smaller, intrinsic tube convergence
would put the actual base points in the core. The comparison ball is measured
at time zero. Thus $|\nabla B_0|\le1$ and $B_t\ge B_0$ give
\[
 u\ge B_0(y)\ge v-L/\sqrt{\mathrm{R}_{B,\mathrm{pot}}(v)}
                      \ge v-C L\sqrt v.
\]
For a neck comparison center, intrinsic tube convergence likewise gives $u\ge2\Lambda$
whenever this is relevant; in any case the last inequality and the
Bryant logarithmic derivative bound give
$\mathrm{R}_{B,\mathrm{pot}}(u)/\mathrm{R}_{B,\mathrm{pot}}(v)\le1+C L/\sqrt\Lambda$.
For a core comparison center, $u\le2\Lambda$, so the same quadratic inequality as
in Proposition~\ref{prop:matchedscale} gives
\[
 2\Lambda\le v\le2\Lambda+C L\sqrt\Lambda,
 \qquad 2\Lambda-C L\sqrt\Lambda\le u\le2\Lambda.
\]
It also gives $\rho_*/\mathrm{R}_{B,\mathrm{pot}}(v)\le1+C L/\sqrt\Lambda$.
The projected tip has model tip scalar one at every time, so this is
exactly the assigned cap-to-neck scale ratio.

We choose $\Lambda$ so that these model errors are smaller than
$\sigma/2$, and use a larger constant $C_2$ in the band bounds to keep
strict margins. Smooth convergence of the metrics, intrinsic tip
curves and tube coordinates then contradicts every proposed exterior
failure. This establishes the actual inequalities. Only backward-time
convergence centered at $(x,0)$ was used; no future-time convergence
from an earlier cap tip was assumed.
This proves the stated conclusion.
\end{proof}

The analogous control is needed for the full cap comparison domain.
The next result compares its actual scale with the model scale
and retains the prescribed margins.

\begin{prop}
\label{prop:actualcapscale}
Use the core and assigned scales of
Proposition~\ref{prop:actualtransition}. Fix $\Lambda,L_C,T_C$ and
$\sigma>0$. There is a compact set $K$ such that for every
$p\in\Gamma\setminus K$, with $r_p:=\mathrm R(p)^{-1/2}$, each point $y$ with
\[
 d_g(p,y)\le L_Cr_p,\qquad -T_Cr_p^2\le t\le0
\]
satisfies
\[
 \begin{cases}
 Q_N(y,t)\le(1+\sigma)Q_C(p,0),&y\notin\mathcal C_t,\\
 Q_C(\pi_t(y),t)\le(1+\sigma)Q_C(p,0),&y\in\mathcal C_t.
 \end{cases}
\]
The compact set depends only on the fixed domain sizes,
$\Lambda,\sigma$, and the soliton.
\end{prop}
\begin{proof}
Suppose the assertion fails along a sequence of tips escaping to
infinity. We normalize at each tip and pass to a pointed ancient
limit. This limit is the tip-normalized Bryant product.
Required points stay in a fixed normalized region by
Proposition~\ref{prop:inputcomparison}. Core projections also stay in
a fixed enlarged region by the tube bounds and scalar comparability,
as in the previous proof. A limiting cap center lies on the product
tip line and has scalar one, so its assigned scale is $\rho_*$.
A limiting neck comparison center has radial value at least $2\Lambda$, hence
scalar at most $\rho_*$. Both cases contradict a failure with the
fixed positive margin $\sigma$. The domain sizes were fixed before
selecting the compact exclusion; the comparison error may therefore
be made small relative to $\sigma\rho_*$.
This finishes the proof.
\end{proof}

We next place all centers needed by the finite cap construction in
one comparison domain. This includes intermediate centers, all choices
at each stage, and replacements of core centers by their tip
projections. We choose the domain before selecting its comparison
accuracy. Corollary~\ref{cor:compressedcapinputs} will then compare every
required scale directly with the final cap scale.

\begin{prop}\label{prop:finiteinputchart}
Fix $\Lambda>0$, an integer $J\ge0$, and a finite family
$\mathcal F=\{(L_\ell,T_\ell)\}$ of positive pairs.
Use $\mathscr R_{\mathcal F}^{\rm par}$ from
Definition~\ref{def:recurrence-reachability}. We also allow an edge
from $(x,s)$ to $(\pi_t(y),t)$ whenever
$(x,s)\,\mathscr R_{\mathcal F}^{\rm par}\,(y,t)$ and
$y\in\mathcal C_t$.
\begin{enumerate}[label=\textup{(\roman*)},ref=\textup{(\roman*)}]
\item There are $D_J,T_J<\infty$ and a compact set $K$ such that every
center reachable in at most $J$ steps from $(p,0)$,
$p\in\Gamma\setminus K$, lies in
\[
 B_{g(0)}(p,D_Jr_p)\times[-T_Jr_p^2,0],\qquad\mathrm{where}\qquad
 r_p:=\mathrm R(p)^{-1/2}.
\]

\item For any fixed $0<\eta_{\rm cyl}<1/2$, we can choose a single
Bryant-product comparison chart with spatial image $\mathcal O_p$
and time interval $[-S_-r_p^2,S_+r_p^2]$ such that
\[
 \overline B_{g(0)}(p,(D_J+1)r_p)\subset\mathcal O_p,
 \qquad S_->T_J+1,\qquad S_+>0.
\]
For each $\eta_{\rm cyl}$-cylindrical center $(x,s)$ in this family,
let $\Psi_{x,s}$ be its map from
Definition~\ref{cylinder:def:neighborhood}. The chart also contains
\[
 \begin{aligned}
 &\Psi_{x,s}\bigl(\overline B_{\bar g(-1)}
                  (\bar x,\eta_{\rm cyl}^{-1})\bigr)\subset\mathcal O_p,\\
 &[s+(1-\eta_{\rm cyl}^{-1})r_{x,s}^2,\ s+r_{x,s}^2]
                   \subset[-S_-r_p^2,S_+r_p^2],
 \end{aligned}
\]
with fixed spatial and time margins. In particular it contains the
cylindrical comparison interval ending at
$s+(1-\eta_{\rm cyl})r_{x,s}^2>s$.
The chart sizes and compact set may depend on
$J,\Lambda,\mathcal F,\eta_{\rm cyl}$ and the prescribed comparison
accuracy and derivative order, but not on $p$.
\end{enumerate}
\end{prop}
\begin{proof}
We first prove (i) in the tip-normalized Bryant product. We let $L$ and
$T$ bound the fixed spatial radii and backward durations in the local
rules, respectively.
We start with a compact region containing the base point and define
\[
 \mathcal K_k:=\overline{B_{g_0}(p,D_k)}\times[-T_k,0],
 \qquad m_k=\min_{\mathcal K_k}\mathrm{R}>0.
\]
The positivity follows from completeness and positivity of scalar
curvature on compact model regions. A center in $\mathcal K_k$ has curvature length at most
$m_k^{-1/2}$. Therefore every comparison point required at this center
satisfies
\[
 d_{g_0}(p,y)\le D_k+L m_k^{-1/2},\qquad
 t_y\ge-T_k-Tm_k^{-1}.
\]
Here $g_s\ge g_0$ for $s\le0$ converts a parent-time spatial ball
into a time-zero ball of no larger radius.

If such a required point lies in the core, replacing it by its tip projection
adds at most
$C_\Lambda$ to its distance in time-zero units. Indeed, the model
core has fixed normal radius in tip units, all tips have scalar one,
and $d_{g_0}\le d_{g_{t_y}}$. Consequently, the recursion
\begin{equation}
 D_{k+1}=D_k+L m_k^{-1/2}+C_\Lambda+1,\qquad
 T_{k+1}=T_k+Tm_k^{-1}+1
\end{equation}
contains every child, with strict spatial and temporal margins.
For each $k<J$, the new compact region has $m_{k+1}>0$, so the
next enlargement is finite. Induction gives finite $D_J,T_J$
simultaneously for all centers required at the next stage. No positive lower
bound for $m_k$ independent of $k$ is required or claimed.

For actual escaping tip centers, the pointed ancient flows converge
to the tip-normalized Bryant product on any such fixed enlarged domain.
Inductively the parent scalar lower bounds and the spatial/time margins
persist. A minimizing path to a point in an allowed spatial ball cannot
leave the enlarged comparison image: its first exit would exceed the
reserved distance margin. Core projections remain in a further bounded
region by the intrinsic tube bound and scalar comparability. Their
original-soliton centers escape, as in
Proposition~\ref{prop:actualtransition}, so the exterior tube estimates
apply. Smooth convergence of these projections preserves the model
containments. Failure of containment along a sequence at one of the first
$J$ steps would therefore contradict the corresponding strict model
inclusion. This proves (i). 

We now prove (ii). For the additional cylindrical charts, let $m>0$ be a lower bound
for the model scalar curvature on the already selected compact family
of required centers. A neck center there has curvature length at most
$m^{-1/2}$ in parent units. Its $\eta_{\rm cyl}$-chart therefore needs
spatial radius at most $C m^{-1/2}(\eta_{\rm cyl}^{-1}+1)$ and a time
extension of size at most $m^{-1}(1+\eta_{\rm cyl}^{-1})$.
These are finite constants. Enlarge the parent product chart by these
amounts and fixed margins. For its future-time comparison, the
scale-invariant derivative bounds of
\cite[Appendix~A, Corollary~A.9]{CMZ} imply
\begin{equation}\label{cylinder:eq:forward-continuation}
 |\partial_t\mathrm R|\le C\mathrm R^2,\qquad
 \mathrm R(y,t)\le
 \frac{\mathrm R(y,u)}{1-C\mathrm R(y,u)(t-u)}
 \quad(t\ge u),
\end{equation}
while the denominator is positive. Together with the higher curvature
derivative bounds and $\partial_tg=-2\Ric$, this gives smooth
compactness on a short forward interval in the same spatial charts.
Uniqueness identifies the continued limit with the model. Iterating
extends convergence across any fixed interval on which that model has
bounded curvature. The Bryant product is eternal with bounded
curvature; the cylinder has $\mathrm R=(-t)^{-1}\le\eta_{\rm cyl}^{-1}$
on $[-1,-\eta_{\rm cyl}]$. Thus the required future portions are
covered in finitely many steps.
This enlargement supplies geometry only; it adds no rotation triple
centers and does not increase the depth $J$ of the construction. This proves (ii) and thus we have completed the proof.
\end{proof}

Local rotation triples must be measured on whole leaves to preserve their
normalization. We now arrange the exterior foliation and its
leafwise margins for that purpose.

\begin{prop}
\label{prop:exteriorleaves}
Let $\Gamma$
and $\beta_0$ be the tip set and tube coordinate of
Proposition~\ref{prop:tipcurves} and \eqref{eq:ambientbeta}.
Fix the measurement domains $\mathcal P_\rho$ and enlargements
$\mathcal P_\rho^+$ of Corollary~\ref{cor:intrinsicneckcontraction},
with $\rho\ge100$ and $\ell_{\rm neck}\ge10\sqrt{\rho+1}$. These constants
are chosen before $\Lambda$ and the compact exclusion.
For sufficiently large $\Lambda$, and outside a compact set, the
neck neighborhoods outside $\{\beta_0<\Lambda/2\}$ carry a common
local foliation by embedded CMC two-spheres inside regular potential
levels. Here the construction is made on larger neighborhoods,
including $\{\beta_0>\Lambda/4\}$.
For each fixed cap chart size, there is a compact set $K$ such that on every such chart centered at $p\in\Gamma\setminus K$ this foliation agrees with
the collar construction of Proposition~\ref{prop:alltiltcollar}.
On overlaps the whole spheres coincide as subsets. The assertion holds smoothly at all times in
any prescribed bounded normalized backward interval and is invariant
under canonical pullback. 
\end{prop}
\begin{proof}
Fix all measurement sizes and a finite normalized backward duration.
We first fix a small cylindrical parameter $\eta_{\rm cyl}$.
Choose it to satisfy \eqref{cylinder:eq:fixed-domain-threshold} for
all the prescribed measurement domains, margins, and derivative orders.
The required geometric neighborhood is then the full pointed chart
of Definition~\ref{cylinder:def:neighborhood}, with radius
$\eta_{\rm cyl}^{-1}$ and time interval
$[-\eta_{\rm cyl}^{-1},-\eta_{\rm cyl}]$.
All these sizes are fixed before the exterior cutoff is selected.
Proposition~\ref{prop:inputcomparison} and
Lemma~\ref{lem:compact-exclusion} allow us to enlarge the compact
exclusion so that every enlarged neighborhood has original-soliton
representatives outside the compact sets used for the neck and
potential comparisons. By
Proposition~\ref{prop:corecover}, increasing $\Lambda$ makes every
exterior point outside $\{\beta_0<\Lambda/4\}$ a neck of this quality.
By Corollary~\ref{cor:outsidecorescalar}, increase $\Lambda$ again
so the potential gradient is uniformly bounded below there.

In a normalized neck chart the potential is close to a linear
coordinate in the flat factor. For completeness, this follows from
the soliton identities: integration of $\Hess F=\lambda\Ric$ over
a whole almost totally geodesic spherical factor gives
$\lambda\le C\eta$, since its tangential Ricci trace has a positive
lower bound and $|\nabla F|\le1$. The Hamilton identity then fixes
the length of the flat gradient to $1+O(\eta)$, and the Hessian
equation gives the remaining fixed derivative bounds.
Thus, the potential level is a controlled perturbation of
$S^2\times\mathbb R$. We prescribe the mean of a normal graph over a
coordinate sphere and solve the CMC equation with its constant mean
curvature as an additional unknown. On mean-zero functions the
linearization at the product is the spherical Laplacian, which is
invertible. The implicit function theorem, with the mean parameter
allowed to vary, gives the full local family and a positive normal
lapse. This is the neck version of the collar construction above.

We explain why charts give the same family. We reserve margins large
enough that a sphere through an overlap lies wholly in both
enlarged charts. This is possible with a fixed margin: near-round
spheres have uniformly bounded diameter in curvature units and their
curvature scales are comparable throughout the overlap. Coordinate
changes between two such enlarged neck charts are close, after
normalization, to product isometries. One way to verify this assertion
is compactness as the common neck threshold tends to zero. Metric
comparison bounds their first derivatives and the connection
transformation formula bounds higher derivatives. A limiting change
of coordinates is a local isometry of the round product. Its positive
Ricci eigenspace is the angular two-plane distribution, so it sends
spherical slices to spherical slices. A whole slice maps onto a whole
slice: the induced local diffeomorphism of $S^2$ is a covering and
$S^2$ is simply connected. Consequently, at a sufficiently small
fixed threshold, a leaf of either chart is a small normal graph in
the other chart inside the same actual potential level.

We vary the reference sphere parameter in the second chart until the
graph has mean zero. The derivative of its mean with respect to
that parameter is bounded away from zero, by the positive lapse.
Both the transferred leaf and the second chart's leaf then solve
the same mean-zero graph equation, with an unknown constant mean
curvature. Local uniqueness identifies them. Thus, the two constructions give the same leaves, although their
parameter labels may differ.
It also proves compatibility on triple overlaps, so the local
families define a foliation on their union.

For a cap collar with $\beta_0\ge\Lambda/2$, choose $\Lambda$ so its
Bryant model is within the same neck threshold in the required
normalized neighborhoods. Its potential level CMC spheres are the
angular spheres. The actual cap family is a small perturbation of
these, uniformly over the allowed potential tilts by
Proposition~\ref{prop:alltiltcollar}. It is therefore a small graph
in a neck chart on each common collar, and the identical uniqueness
argument applies. Enlarging the chart before choosing its comparison
error preserves all whole-leaf margins. Backward time moves model
collar points radially outward; on a fixed actual chart the reserved
gap between $\Lambda/2$ and $\Lambda/4$ persists.

The implicit equations depend smoothly on time and the transverse
parameters. The same comparison and uniqueness argument applies on
each fixed time interval. Isometries, constant metric scaling and
additive changes of the potential preserve both the CMC condition
and this uniquely determined local family. In particular the
canonical pullbacks identify the resulting spheres exactly.
The assertion follows.
\end{proof}

We can consequently define chart-independent cap domains.
For $p\in\Gamma$ outside a sufficiently large compact set, put
$r_p:=\mathrm{R}(p)^{-1/2}$. On a sufficiently long
local arc of $\Gamma$ about $p$, let $z_p(x)$ be the signed arclength
from $p$ to $\pi_0(x)$ divided by $r_p$, oriented by
$\mathrm df(V_0)>0$ as in Corollary~\ref{cor:tipnaturality}.
For the constants of Proposition~\ref{prop:productdomains}, set
\begin{equation}\label{eq:intrinsic-cap-measurement-domains}
 D_p:=\{\beta_0\le B,\ |z_p|\le H\},\qquad
 C_p:=\{\Lambda/2\le\beta_0\le2J\Lambda,\ |z_p|\le H\}.
\end{equation}
These are fixed spatial sets in $\mathcal T_0(R_{\rm tube})$, centered at $p$,
where $R_{\rm tube}$ is chosen large enough to contain their enlargements.
Their comparison enlargements and spacetime sets are
\begin{equation}\label{eq:intrinsic-cap-domain-margins}
 \begin{aligned}
 D_p^+&:=\{\beta_0<B+\Lambda/8,\ |z_p|<H+1\},\\
 C_p^+&:=\{3\Lambda/8<\beta_0<2J\Lambda+\Lambda/8,
                                    \ |z_p|<H+1\},\\
 \mathcal Q_p&:=D_p\times[-Hr_p^2,0],\qquad
 \mathcal Q_{p,\rm collar}:=C_p\times[-Hr_p^2,0],\\
 \mathcal Q_p^+&:=D_p^+\times[-(H+1)r_p^2,0].
 \end{aligned}
\end{equation}

In tip curvature normalization, write
$\widehat g_p(\tau):=r_p^{-2}g(r_p^2\tau)$ and
$F_p(x,\tau):=r_p^{-1}(f(\Phi_{r_p^2\tau}(x))+r_p^2\tau)$,
up to an additive constant. The corresponding normalized spacetime sets are
\[
 \widehat{\mathcal Q}_p:=D_p\times[-H,0],\qquad
 \widehat{\mathcal Q}_{p,\rm collar}:=C_p\times[-H,0].
\]
For a time-independent triple $U=(U_1, U_2, U_3)$ defined
near $D_p^+$, let
\begin{equation}\label{eq:intrinsic-cap-error}
 \begin{aligned}
 \mathcal E_C(U;p):={}&\sup_{\widehat{\mathcal Q}_p}
 \left[\sum_{i=1}^3\left(
 \sum_{j=0}^2|\widehat\nabla^j\Lie_{U_i}\widehat g_p|^2
                         +|U_iF_p|^2\right)\right]^{1/2}\\
 &+\sup_{\widehat{\mathcal Q}_{p,\rm collar}}
 \left[\sum_{i=1}^3\left(
 |\langle U_i,\widehat n_f\rangle|^2
 +|\langle U_i,\widehat\nu\rangle|^2\right)\right]^{1/2}
 +\sup_{\tau,S}|Q_S(U)-I|_{\mathrm F},\\
 e_C(p):={}&\inf_U\mathcal E_C(U;p).
 \end{aligned}
\end{equation}
The map $(x,\tau)\mapsto(x,r_p^2\tau)$ identifies these normalized
sets with $\mathcal Q_p,\mathcal Q_{p,\rm collar}$.
All norms, unit normals, and Gram matrices in
\eqref{eq:intrinsic-cap-error} use $\widehat g_p(\tau)$, and
$\widehat\nabla$ is its spatial connection. The last supremum is over whole collar CMC
spheres at every $\tau\in[-H,0]$. The neck infimum $e_N(x):=e_N(x,0)$
is the infimum of the three rotation-error bounds of
Definition~\ref{cylinder:def:symmetry} on $\mathcal P_\rho(x,0)$,
as specified in Corollary~\ref{cor:intrinsicneckcontraction}.
At final time $t_*$, canonical pullback gives
\[
 e_C(p,t_*):=e_C(\Phi_{t_*}(p)),\qquad
 e_N(x,t_*)=e_N(\Phi_{t_*}(x),0),\qquad \mathrm{where}\qquad p\in\Gamma_{t_*}.
\]
The spatial map $\Phi_{t_*}^{-1}$ and time shift
$\sigma\mapsto t_*+\sigma$ transport the measurement domains and
admissible time-independent triples; thus these identities preserve
the full error functionals.

\begin{prop}
\label{prop:intrinsicwitnesses}
Let $\Gamma$
and $\beta_0$ be the tip set and tube coordinate of
Proposition~\ref{prop:tipcurves} and \eqref{eq:ambientbeta}.
Fix the neck domains of Corollary~\ref{cor:intrinsicneckcontraction}
and choose $\Lambda,J,B,H$ by Proposition~\ref{prop:productdomains}.
The cap domains are $D_p,C_p$ displayed above, with time interval
$[-Hr_p^2,0]$. The infima $e_C(p)$ and $e_N(x)$ are defined in
\eqref{eq:intrinsic-cap-error} and Corollary~\ref{cor:intrinsicneckcontraction}, respectively.
Both use the same CMC spheres of Proposition~\ref{prop:exteriorleaves},
and times other than zero are defined by canonical pullback.
\begin{enumerate}[label=\textup{(\roman*)},ref=\textup{(\roman*)}]
\item There is a compact set $K\subset M$ outside which these infima are finite and satisfy
\[
 e_C(p)\longrightarrow0\quad(p\in\Gamma,\ p\to\infty),
 \qquad
 e_N(x)\longrightarrow0\quad(x\to\infty\text{ outside the core}).
\]
\item There is a fixed $C_0<\infty$ such that a cap rotation triple centered at
$p=\pi_{t_*}(x)$ restricts to a neck rotation triple at $(x,t_*)$ with error
at most $C_0$ times its cap error whenever
$3\Lambda/4\le\beta_{t_*}(x)\le2\Lambda$ and the whole required
domains satisfy
\[
 \mathcal P_\rho^+(x,t_*)\cup\mathcal Q_{p,t_*}^+
 \subset\{(y,t):\Phi_t(y)\notin K\}.
\]
Here $\mathcal Q_{p,t_*}^+$ is the image of
$\mathcal Q_{\Phi_{t_*}(p)}^+$ under
$(y,\sigma)\mapsto(\Phi_{t_*}^{-1}(y),t_*+\sigma)$.
In particular, $e_N(x,t_*)\le C_0e_C(p,t_*)$
in this band. The core-box containments and the restriction of the constructed triples in Proposition~\ref{prop:productdomains} hold
for this same intrinsic class.
\end{enumerate}
\end{prop}
\begin{proof}
We first prove (i) by transferring exact rotations from the limiting models to obtain
local rotation triples with errors tending to zero. We then compare the
intrinsic domains and common leaves to prove the restriction estimate.

For escaping cap centers, the tip-normalized pointed flows converge
on any fixed enlarged domain to the Bryant product. By the intrinsic
tube convergence, $\beta_0$ and $z_p$ converge smoothly to the model
radial potential and line coordinate. The normalized soliton
potentials converge, after a subsequence and additive normalization,
to $\lambda B_t+az+\lambda t$, with
$\lambda^2+a^2=1$. The standard rotational fields are
time-independent in these model coordinates, have zero metric Lie derivative and
potential pairing, and are tangent to the angular spheres. Their
area-squared Gram matrices are exactly $I$ after one fixed constant
normalization, at every radius and time. We transfer these fields by a
time-independent pointed comparison map. Smooth convergence and
Proposition~\ref{prop:exteriorleaves} make all their errors
tend to zero on the actual intrinsic domains. Fixed enlarged domains
and strict margins ensure that this assertion includes all times and whole leaves, not just a central sphere.

The same contradiction argument proves the neck assertion. For a
prescribed cylindrical quality $\eta_{\rm cyl}$, first fix its radius
$\eta_{\rm cyl}^{-1}$, time interval
$[-\eta_{\rm cyl}^{-1},-\eta_{\rm cyl}]$, and derivative order
$\lfloor\eta_{\rm cyl}^{-1}\rfloor$. If the pointed limit is cylindrical,
\eqref{cylinder:eq:forward-continuation} supplies convergence on this
entire interval, and excluding a larger compact set gives the required
$\eta_{\rm cyl}$-cylindrical chart. The restriction rule
\eqref{cylinder:eq:fixed-domain-threshold} then supplies every fixed
measurement domain used in Theorem~\ref{cylinder:improvement}.
An escaping neck sequence has a cylindrical or Bryant-product limit.
In the latter case its basepoint is outside the fixed core; choosing
$\Lambda$ large relative to the neck size keeps the frozen rotation triple
domain away from the tip, including its whole-leaf margins. On either
limit there is an exact rotational triple with the required Gram
normalization, tangent to the potential level CMC spheres. For the
cylinder, the normalized limiting potential is linear in the flat
factor. For the Bryant product it is rotationally invariant by the
limiting soliton identities. Thus, the transferred triples have all
errors tending to zero. Uniform compactness rules out a
sequence of infima bounded below by a positive number. This proves (i). 

To prove (ii), recenter the model at the projection of $x$ onto
the tip line at time $t_*$. This is precisely the limit of
$\pi_{t_*}(x)$, rather than a tip chosen on its potential level.
After this canonical recentering, $\beta_{t_*}$ and normalized
arclength converge to the frozen model coordinates at time $t_*$. Every model inclusion in
Proposition~\ref{prop:productdomains} has a strict margin, so it
continues to hold for these intrinsic domains on sufficiently close
charts. The common actual leaves agree as sets by
Proposition~\ref{prop:exteriorleaves}. Comparing the same fields on
these sets therefore introduces no additive geometric error.
The curvature length ratio is comparable to $\sqrt\Lambda$;
metric Lie derivative derivatives and transverse components rescale as in that
proposition. On the collar, conversion of $UF$ to the unit-potential
pairing has a uniformly bounded denominator. The area-squared Gram
matrix is unchanged by the homothety. These facts give a fixed
multiplicative constant $C_0$. Taking approximate minimizers and
letting their excess error tend to zero proves the infimum inequality.

The identical strict-margin transfer applies to the core boxes and to the fixed domain enlargement for the constructed triples. No new size is introduced, and no potential gradient denominator is used
in the core. This proves the remaining assertions.
This proves (ii) and completes the proof.
\end{proof}

The intrinsic domains are now available. We first obtain radially
weighted neck triples for a cap, then prove that one fixed class is preserved by the combined cap and neck improvements.
\begin{prop}
\label{prop:radialtree}
Let $\Gamma$ and $\beta_t$ be the tip set and radial coordinates of
Proposition~\ref{prop:tipcurves} and \eqref{eq:ambientbeta}.
Fix $\Lambda>0$ and comparison sizes $L,T>0$. A recurrence step
$(x,s)\leadsto(y,t)$ means
\[
 s-T r_{x,s}^2\le t\le s,\qquad
 d_{g(s)}(x,y)\le Lr_{x,s},\qquad\mathrm{where}\qquad
 r_{x,s}:=\mathrm R_{g(s)}(x)^{-1/2}.
\]
There is $D=D(L,T)<\infty$ with the following properties.
\begin{enumerate}[label=\textup{(\roman*)},ref=\textup{(\roman*)}]
\item Fix a bounded tip-normalized spacetime region about $(p,0)$,
with a finite time interval contained in $(-\infty,0]$,
and an integer $J\ge0$. There are an enlarged tube and a compact
set $K$, depending on this region and $J$, such that for every
$p\in\Gamma\setminus K$, all chains of at most $J$ recurrence
steps starting in this region lie in the tube. Every step satisfies
\begin{equation}\label{eq:radialstep}
 \sqrt{\beta_t(y)+1}\ge\sqrt{\beta_s(x)+1}-D.
\end{equation}
The sizes $L,T$ may be enlarged in advance to include the required
rotation triple domains and whole CMC leaves.
\item Fix $N\ge1$, $c_N\ge1$, and $0<\theta<1$. There is
$B_N>2\Lambda$, independent of the heat radius $E$, such that the
following holds. For each fixed $E\ge B_N/2$, choose the region in
(i) to contain the image of $Q_E^+$ from
\eqref{eq:cap-heat-enlargement} in a tip-normalized cap chart, and put
\[
 J_E:=\left\lceil\frac{N\log(2E+c_N)}{|\log\theta|}\right\rceil.
\]
After choosing the enlarged tube and compact set $K_E$ for this
region and depth $J_E$, let $p\in\Gamma\setminus K_E$ and let
$(x,s)$ be any initial spacetime center in this image of $Q_E^+$ with
$v:=\beta_s(x)\in[B_N,2E]$. For every chain
\[
 (x_0,t_0)=(x,s)\leadsto(x_1,t_1)\leadsto\cdots
                      \leadsto(x_j,t_j),\qquad 0\le j\le m(v),
\]
where
\begin{equation}
 m(v):=\left\lceil\frac{N\log(v+c_N)}{|\log\theta|}\right\rceil,
\end{equation}
all vertices satisfy
\[
 \beta_{t_j}(x_j)\ge(v-1)/2>2\Lambda.
\]
\end{enumerate}
\end{prop}
\begin{proof}
We first prove (i). Fix the initial spacetime region and depth.
Proposition~\ref{prop:finiteinputchart}, with an initial step covering
that region, encloses all allowed chains and their measurement
margins in one fixed enlarged tip chart. Choose its tube radius
before its compact exclusion. We now prove the estimate on every
step in this chart.

In the tip-normalized Bryant product, put $v:=B_s(x)$ and $u:=B_t(y)$.
The function $B_{s}$ is $1$-Lipschitz for $g_s$, and backward time
increases $B_t$ at a fixed spatial point. The global model bound
$\mathrm{R}_{B,\mathrm{pot}}(v)\ge c_1/(v+1)$ therefore gives
\[
 u\ge B_{s}(y)\ge v-C_1L\sqrt{v+1}.
\]
We choose $D>C_1L+2$. On each fixed enlarged comparison domain, the
metrics, positive scalar curvatures, intrinsic projections and radial
coordinates converge smoothly. Consequently, after enlarging a compact set $K$ depending on this
fixed family, every chart centered at $p\in\Gamma\setminus K$ satisfies
\[
 \beta_t(y)\ge\beta_s(x)-D\sqrt{\beta_s(x)+1}.
\]
For the rest of this calculation put $v:=\beta_s(x)$.
This transfer is uniform for the entire compact family of allowed
parent and child spacetime points. The strict increase of the
constant absorbs perturbations, including the small perturbation of
backward radial monotonicity. No uniform comparison threshold over
unbounded $E$ or depth is claimed.

If $w=\sqrt{v+1}\ge D$, then
$w^2-Dw\ge(w-D)^2$, so the last inequality implies
\eqref{eq:radialstep}. If $w<D$, that equation holds automatically.
This proves (i).

We next prove (ii). Iteration gives, for every $0\le j\le m(v)$,
\[
 \sqrt{\beta_{t_j}(x_j)+1}\ge\sqrt{v+1}-Dj.
\]
Since $\log(v+c_N)/\sqrt{v+1}\to0$, choose $B_N$ so large that
\[
 D\left(\frac{N\log(v+c_N)}{|\log\theta|}+1\right)
 \le(1-2^{-1/2})\sqrt{v+1},\qquad
 (v-1)/2>2\Lambda+1
\]
for every $v\ge B_N$. These choices are independent of $E$.
They imply the asserted lower bound at all vertices, so no core
replacement occurs in these trees.

For the uniform actual assertion, fix $E$ first. All initial heat grid
centers lie in a fixed spacetime region about the cap center, and
$m(v)$ is bounded by the finite integer
$\lceil N\log(2E+c_N)/|\log\theta|\rceil$. We apply
Proposition~\ref{prop:finiteinputchart}, including one initial
cap-to-grid step, to enclose all possible descendants with margins.
It supplies one enlarged domain before the comparison error is chosen.
The preceding uniform transfer then applies throughout this domain.
The same reasoning works with fixed enlargements of $2E$ and the
grid's spatial and time bounds. Thus, the tree is not required to stay
inside the original heat box.
This proves (ii) and completes the proof.
\end{proof}

\begin{cor}
\label{cor:weightedneckinputs}
Use the soliton, cap heat region, and recurrence steps of
Proposition~\ref{prop:radialtree}. Fix $N\ge1$, $c_N\ge1$,
$0<\theta<1$, and $E\ge B_N/2$. Choose $L,T$ so that each step
contains the comparison set $\mathcal I(x,s)$ of
\eqref{cylinder:eq:neck-input-set}. Let $\mathcal V_E$ consist of
all initial spacetime centers $(x,s)$ in the cap heat region with $B_N\le\beta_s(x)\le2E$
and every point reachable from each such initial center by at most
$m(\beta_s(x))$ steps. We also include the initial centers with
$2\Lambda\le\beta_s(x)<B_N$. Here $B_N$ and $m$ are as in
Proposition~\ref{prop:radialtree}.
Assume $e_N\le\epsilon$ on $\mathcal V_E$, where $\epsilon>0$,
and, at each nonterminal vertex $(x,s)$ of these chains, assume
\[
 e_N(x,s)\le\theta\sup_{(y,t)\in\mathcal I(x,s)}e_N(y,t).
\]
All these domains and their whole CMC leaves are those of
Corollary~\ref{cor:intrinsicneckcontraction}. Choose the compact set
$K$ large enough for Proposition~\ref{prop:radialtree} and the common
CMC foliation, and assume that every required spacetime domain
is contained in $\{(y,t):\Phi_t(y)\notin K\}$.
Then every initial center $(x,s)$ with $v:=\beta_s(x)\in[2\Lambda,2E]$
admits a time-independent rotation triple on
$\mathcal P_\rho(x,s)$ satisfying all three inequalities of
Definition~\ref{cylinder:def:symmetry} with error at most
\begin{equation}
 C_N\epsilon(v+c_N)^{-N}.
\end{equation}
The constant $C_N$ is independent of $E$ and the initial center; the enlarged
comparison chart and compact exclusion may depend on $E$.

For the core boxes $\mathcal B_{0\ell m}$ of
Proposition~\ref{prop:productgrid}, the centers
$(q_\alpha,t_\alpha)$ lie on $\Gamma_{t_\alpha}$, so their radial
values are $v_\alpha:=\beta_{t_\alpha}(q_\alpha)=0$.
Let $U_\alpha$ be restrictions of cap triples to the field domains
$\mathcal W_\alpha\times T_\alpha$ and whole-leaf collars
$\mathcal C_\alpha^+$ of that proposition. Put
$r_\alpha:=\mathrm R(q_\alpha,t_\alpha)^{-1/2}$, and let $e_\alpha$
be the least common bound in
\begin{align*}
 \sup_{\mathcal W_\alpha\times T_\alpha}
 \left(\sum_{i=1}^3\sum_{j=0}^2
 r_\alpha^{2j}|\nabla^j\Lie_{U_{\alpha,i}}g|^2\right)^{1/2}
 &\le e_\alpha,\\
 \sup_{\mathcal C_\alpha^+}r_\alpha^{-1}
 \left(\sum_{i=1}^3\bigl(
 |\langle U_{\alpha,i},n_f\rangle|^2
 +|\langle U_{\alpha,i},\nu\rangle|^2\bigr)\right)^{1/2}
 &\le e_\alpha,\\
 \sup_{t,S}|Q_S(U_\alpha)-I|_{\mathrm F}&\le e_\alpha.
 \end{align*}
The last supremum is over whole spheres in $\mathcal C_\alpha^+$.
If these cap restrictions satisfy $e_\alpha\le A_{\rm grid}\epsilon$,
with $A_{\rm grid}$ independent of $E$, then they also satisfy
\[
 e_\alpha\le C_N\epsilon(v_\alpha+c_N)^{-N}
                  =C_N\epsilon c_N^{-N},\qquad\mathrm{where}\qquad
 C_N:=2+(B_N+c_N)^N+A_{\rm grid}c_N^N.
\]
\end{cor}
Note that without the core-box assertion, one may take
$C_N=2+(B_N+c_N)^N$.
\begin{proof}
For $v\ge B_N$, apply the assumed inequality down the full recurrence tree
of depth $m(v)$. All its vertices are necks by
Proposition~\ref{prop:radialtree}. Repeated suprema give
\[
 e_N(x,s)\le\theta^{m(v)}\epsilon
              \le\epsilon(v+c_N)^{-N}.
\]
The estimate uses the bound on every terminal choice, not a selected
recurrence path. For $2\Lambda\le v<B_N$, the original bound gives
$e_N\le(B_N+c_N)^N\epsilon(v+c_N)^{-N}$ without an improvement step.
We choose actual rotation triples within the positive tolerance
$\epsilon(v+c_N)^{-N}$ of each infimum. This only enlarges $C_N$ by a
fixed amount and requires no attainment of the infimum.
For a core box, $v_\alpha=0$ and
$A_{\rm grid}\epsilon\le C_N\epsilon c_N^{-N}$ by the stated
choice of $C_N$. Thus the same weighted estimate holds for all
three displayed rotation errors on its field domain and collar. The number of grid boxes
and iteration depth may depend on $E$; this affects the comparison
domain and its accuracy, not the displayed weight constant.
This proves the stated conclusion.
\end{proof}

\begin{cor}
\label{cor:compressedcapinputs}
Use the assigned scales $Q_N,Q_C$ of
Proposition~\ref{prop:actualtransition} and the weighted construction
of Corollary~\ref{cor:weightedneckinputs}.
Fix $\Lambda,N,E$, the local comparison sizes, and $\sigma>0$.
For a cap center $p\in\Gamma$ at time zero, let
$\mathcal A_E(p)$ consist of all neck and cap centers used in this
construction: the intermediate and terminal recurrence vertices,
the direct grid comparison centers, and the cap centers obtained
by the prescribed core replacement.
There are $L_*,T_*<\infty$ and a compact set $K\subset M$ such that,
for every $p\in\Gamma\setminus K$ and $(z,t)\in\mathcal A_E(p)$,
\[
 \begin{gathered}
 Q_{\rm req}(z,t)\le(1+\sigma)Q_C(p,0),\\
 d_g(p,z)\le L_*Q_C(p,0)^{-1/2},\qquad
 -T_*Q_C(p,0)^{-1}\le t\le0.
 \end{gathered}
\]
Here $Q_{\rm req}=Q_N$ at neck centers and $Q_{\rm req}=Q_C$
at cap centers. The constants and $K$ may depend on the fixed
parameters, but not on $p$.
\end{cor}
\begin{proof}
The radial construction has fixed finite depth once $E,N$ are fixed.
Proposition~\ref{prop:finiteinputchart} contains every intermediate
and terminal center in one frozen spatial ball and backward interval
about $p$, with finite sizes. We enlarge this fixed domain to include
the direct grid comparison centers as well. We apply
Proposition~\ref{prop:actualcapscale} directly to that entire ball
and interval. Each required point, after the prescribed core replacement,
is either outside the core or is a cap center on $\Gamma_t$,
precisely the two alternatives
of that proposition. Converting to assigned curvature units at the base point changes
only fixed spatial and time constants. This proves the claim.
\end{proof}

We fix the requested decay exponent $p>0$, then choose $\sigma>0$
so that
\[
 \frac{\log2}{\log(1+2\sigma)}>p.
\]
Next fix the neck measurement sizes $\rho,\ell_{\rm neck}$ and the
resulting input size $L$ of Corollary~\ref{cor:intrinsicneckcontraction},
with gain $\theta=1/2$. Then choose the core threshold $\Lambda$ as in
Proposition~\ref{prop:actualtransition}, the cap measurement margins
\eqref{eq:intrinsic-cap-domain-margins}, and the restriction constant
$C_0$ of Proposition~\ref{prop:intrinsicwitnesses}.
The radial weight $N=16$ and gluing loss $q=4$ are fixed before
choosing $E$ large enough so that
\[
                       C_0CE^{-9}<\tfrac14.
\]
With $E$ fixed, place all recursive vertices, core replacements,
and direct grid comparison centers in one finite comparison chart.
We choose $\eta$ only on this entire chart, requiring
\[
                       C_0C_E\eta<\tfrac14,
\]
and then exclude a sufficiently large compact set so that every
comparison estimate holds there.
Here $C$ is independent of $E$, whereas $C_E$ and the excluded
compact set may depend on $E$.

\begin{prop}
\label{prop:typedclosure}
Fix $o\in M$ and a decay exponent
$p>0$. There are a core threshold $2\Lambda$, fixed neck and cap
measurement domains $\mathcal P_\rho(x,0)$ and
$\mathcal Q_q=D_q\times[-Hr_q^2,0]$, together with the collar
$\mathcal Q_{q,\rm collar}$ of
\eqref{eq:intrinsic-cap-measurement-domains}--\eqref{eq:intrinsic-cap-domain-margins},
a compact set $K_p\subset M$, and a constant
$C_p<\infty$ such that
\begin{enumerate}
\item For every $x\in M\setminus(K_p\cup\mathcal C_0)$,
\[
 e_N(x,0)\le C_p(1+d_g(o,x))^{-p}.
\]
\item For every cap center $q\in\Gamma\setminus K_p$,
\[
 e_C(q)\le C_p(1+d_g(o,q))^{-p}.
\]
\end{enumerate}
Here $\Gamma$ and $\beta_0$ are given by
Proposition~\ref{prop:tipcurves} and \eqref{eq:ambientbeta}, and
\[
 \mathcal C_0:=\{x\in\mathcal T_0:\beta_0(x)<2\Lambda\}.
\]
The neck error $e_N(x,0)$ is the infimum on
$\mathcal P_\rho(x,0)$ defined in
Corollary~\ref{cor:intrinsicneckcontraction}. The cap error $e_C(q)$
is the infimum \eqref{eq:intrinsic-cap-error} on
$\mathcal Q_q=D_q\times[-Hr_q^2,0]$, with normal and Gram errors
on $\mathcal Q_{q,\rm collar}$ and whole collar spheres.
Both use the common CMC family of
Proposition~\ref{prop:exteriorleaves}.
All measurement sizes and $\Lambda$ are chosen once for the given
$p$, independently of $x,q$; $C_p$ and $K_p$ may depend on these
choices and on $o$.
\end{prop}
\begin{proof}
We choose $\sigma>0$ so small that
$\log 2/\log(1+2\sigma)>p$. We fix Corollary~\ref{cor:intrinsicneckcontraction} with gain
$\theta=1/2$ and then choose $\Lambda$ so that Proposition~\ref{prop:actualtransition} holds with
scale factor $1+\sigma$. We increase $\Lambda$ for all the leaf,
restriction and derivative requirements already specified. Let
$C_0\ge1$ be the fixed restriction constant. We weight the neck
and cap errors by setting
\[
 \mathscr E_N:=e_N,\qquad \mathscr E_C:=C_0e_C.
\]
For a neck improvement, Proposition~\ref{prop:actualtransition}
places every comparison point entering $\mathcal C_t$ in
$3\Lambda/4\le\beta_t<2\Lambda$, after increasing $\Lambda$.
Proposition~\ref{prop:intrinsicwitnesses} supplies its neck triple
from the projected cap, with error at most $\mathscr E_C$.
Comparison points outside $\mathcal C_t$ use $\mathscr E_N$.
The cap construction also uses its direct core boxes; their cap
restrictions are covered by the core-box assertion of
Corollary~\ref{cor:weightedneckinputs}.

We first assume that every comparison point required by the full local construction has
weighted symmetry error at most $\epsilon>0$. The neck rule gives an improved symmetry error at most $\epsilon/2$. For a cap center, fix $N=16$ and
apply Corollary~\ref{cor:weightedneckinputs} with the fixed shift
$c_N=c_{16}$. To provide fields on $Q_E^+$, use the weighted
construction and Proposition~\ref{prop:productgrid} at size $2E$,
as explained at the end of that proposition. Its annular boxes use
the weighted neck triples, and its core boxes use cap restrictions
with the three error bounds displayed in
Corollary~\ref{cor:weightedneckinputs}. All constants in these weighted estimates are
fixed before $E$. Proposition~\ref{prop:loworderalignment} supplies
the $C^3$ overlap comparisons from the low-order errors.
The controlled core paths and common whole spheres lie in the
domains of Proposition~\ref{prop:intrinsicwitnesses}.

The remaining size and potential pairing hypotheses of the product
gluing estimate follow in the same units. At a neck center $(x,t)$, put $v:=\beta_t(x)$.
Its curvature length in tip units is comparable to $\sqrt{v+1}$. The bounded rotation triple estimate in the proof of
Proposition~\ref{prop:loworderalignment} therefore gives field size
at most $C\sqrt E$ and the stated bounds through three spatial
derivatives. The unit-potential-normal bound gives
$|UF|\le C\sqrt{v+1}$ times the local error. Core rotation triples have
the $UF$ bound directly. Equation~\eqref{eq:productgrid-model-closeness} supplies the
required $C^3$ comparison with one fixed area-normalized model
triple after a constant orthogonal change of basis.

Divide the area-normalized glued fields by $\krot=\sqrt{6\pi}$
to match the bracket-one rotations $J_i$ in
Corollary~\ref{cor:alltiltheat}. This changes only fixed constants.
That corollary restores Gram normalization by its final constant
matrix and gives a cap triple with error
\[
 e_C\le(CE^{-9}+C_E\eta)\epsilon.
\]
Here fixed $\Lambda,N$ and fixed-domain constants may affect $C$,
but it is independent of $E$. We choose $E$ so that
$C_0CE^{-9}<1/4$, and then choose the common comparison and error
threshold so that $C_0C_E\eta<1/4$. The constructed triple restricts
to exactly the intrinsic class by
Proposition~\ref{prop:intrinsicwitnesses}. Thus its weighted symmetry error is
at most $\epsilon/2$. Approximate minimizers can be chosen with a
common positive excess error and that excess sent to zero. The same
argument with arbitrarily small positive $\epsilon$ handles zero
supremum. There is no need to assume an infimum is attained.

All internal cap iterations have finite depth after $E$ is fixed.
Corollary~\ref{cor:compressedcapinputs} compares each terminal required point
directly with the cap center with factor $1+\sigma$. Ordinary neck centers have the same factor by
Proposition~\ref{prop:actualtransition}. All these edges have finite
spatial and backward time bounds in assigned units; enlarging those
bounds includes direct grid comparison centers and every intermediate vertex,
as well as terminal neck comparison centers. Thus, the baseline error bound used
in the weighted iteration applies wherever it was assumed.
The spatial bounds just obtained are in the frozen final-time
metric. Proposition~\ref{prop:inputcomparison}, using
$\rho_*\mathrm R\le Q\le\mathrm R$, converts them into
$d_{g_t}(x,y)\le L_*Q(x)^{-1/2}$ with fixed enlarged constants.
It leaves the assigned-scale ratio $1+\sigma$ unchanged.
Let $\mathcal Z(x)$ consist of the original-soliton representatives
of the comparison centers in one full local improvement at $x$,
retaining their neck or cap type. With the preceding error weights,
the local improvement gives
\[
 \mathscr E(x)\le\tfrac12\sup_{z\in\mathcal Z(x)}\mathscr E(z),\qquad
 Q(z)\le(1+\sigma)Q(x),
\]
with one edge for the entire cap construction. The state space here is
\[
 \mathscr X:=(M\setminus\mathcal C_0)\sqcup\Gamma,
\]
with the two copies retaining their neck or cap type. On its first
copy $\mathscr E=e_N$ and $Q=\mathrm R$; on its second copy
$\mathscr E=C_0e_C$ and $Q=\rho_*\mathrm R$.
Distances and compact exclusions use the spatial representative
in $M$. Every spacetime center is first mapped to its original-soliton
representative by $\Phi_t$, preserving its type and error by canonical
pullback. Thus no neck error is evaluated at an unprojected point
inside the core.

By Proposition~\ref{prop:intrinsicwitnesses}, both weighted symmetry
errors are below the required threshold outside a compact set.
Extend these errors by arbitrary finite values inside that compact
set; only exterior inequalities will be used. On compact sets the
assigned scales have positive lower bounds, since
$Q\ge\rho_*\mathrm{R}>0$. We apply the assigned-scale version of
Corollary~\ref{cor:finite-exclusion} and
Proposition~\ref{prop:polynomialrecurrence} with
$\bar C=1+2\sigma$. It gives exponent
$\log 2/\log(1+2\sigma)>p$. Its spatial and time constants affect
the prefactor and starting compact set, not this exponent.
For clarity, the order of the choices in this proof is
\begin{equation}
 \begin{gathered}
 p\ \longrightarrow\ \sigma\ \longrightarrow\
 \text{fixed neck rule}\ \longrightarrow\
 (\Lambda,J,B,H,C_0)\ \longrightarrow\ (N=16,C_N),\\
 (N,C_N)\ \longrightarrow\ E\ \longrightarrow\
 \text{enlarged finite family of required centers}\ \longrightarrow\ (\eta,K).
 \end{gathered}
\end{equation}
In particular, the constants in $CE^{-9}$ are fixed before $E$,
whereas $C_E$ and the required comparison compact set may depend on
$E$. The cylindrical quality $\eta_{\rm cyl}:=\eta_{\rm neck}$ is fixed
at the neck-rule step by \eqref{cylinder:eq:neck-threshold}, for the
chosen $\rho,\ell_{\rm neck}$ and $\theta=1/2$, before $\Lambda$ is
chosen.
It is distinct from the later comparison tolerance $\eta$ on the
finite cap construction. Proposition~\ref{prop:finiteinputchart}
incorporates the full $\eta_{\rm cyl}$-charts into the parent geometric
domain. Choosing the final tolerance and excluded compact set then
preserves these fixed inclusions without changing the error gain or
the polynomial exponent. The compact-set recurrence yields a threshold sequence with
\[
 D_j\le D_*\bar C^j,\qquad
 \sup_{\substack{x\in\mathscr X\\x\notin K_j}}\mathscr E(x)
       \le 2^{-j}\epsilon_*,\qquad
 \bar C:=1+2\sigma.
\]
For either type of center $x\in\mathscr X$, put $r:=d_g(o,x)$.
For large $r$, take an integer $j$ with
$D_*\bar C^j<r\le D_*\bar C^{j+1}$. Since $K_j\subset B(o,D_j)$,
the point is outside $K_j$, and
\[
 \mathscr E(x)\le C r^{-\log2/\log\bar C}\le C_p(1+r)^{-p}.
\]
This proves both stated inequalities and completes the proof.
\end{proof}

Proposition~\ref{prop:typedclosure} supplies the initial error bounds
for Section~\ref{sec:globalization}. Proposition~\ref{prop:localrotationatlas}
uses them with a larger decay exponent and one additional local
improvement to obtain derivative estimates and constant overlap
comparisons for the local rotation triples.
Proposition~\ref{prop:curvaturelength} transfers the decay bounds to
all centers used in this additional improvement. These estimates
verify the hypotheses of Proposition~\ref{prop:rotationconnection};
Corollary~\ref{cor:rotationglobalization} then produces the global fields
whose estimates are stated in Proposition~\ref{prop:globalrotations}.

\section{Global alignment of the local rotations}\label{sec:globalization}
In this section, we construct three global approximate rotations from the local
cap and neck triples (Proposition~\ref{prop:globalrotations}). The main issue is to make their choices of basis
compatible across all overlaps.

We first construct a convex auxiliary exhaustion for the geometric
transport. The soliton potential is not yet known to be proper, so
its levels cannot be used for this construction. By Greene--Wu
\cite[Theorem 1(a)]{GW}, positive sectional curvature supplies a smooth
proper function $\psi\ge0$ with positive-definite Hessian and bounded
gradient. We normalize its unique minimum to $\psi(o)=0$, and write
$Y:=\nabla\psi$, $|Y|\le L$. All exterior bundle and frame constructions in this section take
place on auxiliary superlevels and bands
\begin{equation}
 \mathcal E_a:=\{\psi>a\},\qquad
 \mathcal A_{a,b}:=\{a\le\psi\le b\},\qquad
 \Sigma_a:=\{\psi=a\},\qquad\mathrm{where}\qquad b>a>0.
\end{equation}
A smaller exterior means $\mathcal E_b\subset\mathcal E_a$ for a
larger fixed $b$. Its threshold may depend on the requested decay
and derivative orders.

In this section a triple is viewed as the column map
$U_i:\mathbb R^3\to TM$, $U_iv=\sum_\alpha v_\alpha U_{i,\alpha}$.
Changes of basis act on the right:
$(U_iB)_\alpha=\sum_\beta U_{i,\beta}B_{\beta\alpha}$.
If Sections~\ref{sec:cylinder}--\ref{sec:capimprovement} give
$U_{i,\alpha}=\sum_\beta P_{\alpha\beta}U_{j,\beta}$, this convention
reads $U_i=U_jP^{\mathsf T}$; hence $O^{ij}:=P^{\mathsf T}$ below.

\begin{prop}
\label{prop:auxiliaryflow}
Let $(M^n,g)$ be a complete connected noncompact Riemannian manifold with positive sectional curvature. Let $\psi\ge0$ be a smooth proper function
with $\Hess\psi>0$, unique minimum $\psi(o)=0$, and
$Y=\nabla\psi$ satisfying $|Y|\le L$.
\begin{enumerate}[label=\textup{(\roman*)},ref=\textup{(\roman*)}]
\item There are $a_0>0$ and $C>1$ such that
\[
 C^{-1}d(o,x)-C\le\psi(x)\le Ld(o,x),\qquad
 |Y(x)|\ge C^{-1}\quad\text{on }\{\psi\ge a_0\}.
\]
\item Let $\Psi_t$ be the flow of $-Y$, and put
$C_a:=\{\psi\le a\}$, $\Sigma_a:=\{\psi=a\}$. For each $a\ge a_0$
there is a retraction $P_a:M\to C_a$, equal to the identity on
$C_a$, that is nonexpanding for the ambient distance. On $\psi>a$
define
\[
 \tau_a(x):=\inf\{t\ge0:\psi(\Psi_t(x))=a\},\qquad
 P_a(x):=\Psi_{\tau_a(x)}(x).
\]
Thus $\psi(\Psi_t(x))>a$ for $0\le t<\tau_a(x)$ and
$\psi(P_a(x))=a$. Strict decrease of $\psi$ along this trajectory
makes the intersection unique.
Its hitting time $\tau_a$ is smooth on $\{\psi>a\}$ and satisfies
\[
 0\le\tau_a(x)\le C^2(\psi(x)-a),\qquad
 |\nabla \tau_a|\le C,\qquad |DP_a|\le1.
\]
The extension of $P_a$ across $\Sigma_a$ is locally Lipschitz;
$|DP_a|\le1$ there is understood almost everywhere. Set
$\tau_a=0$ on $\Sigma_a$.
\item On $a\le\psi\le b$, the homotopy
$H(s,x)=\Psi_{s\tau_a(x)}x$, $0\le s\le1$, stays in that annulus
and obeys
\[
 |D_xH|\le1+LC,\qquad
 |\partial_sH|\le LC^2(b-a).
\]
Every loop $\gamma$ in $\{a\le\psi\le b\}$ is homotopic there to
$P_a\gamma\subset\Sigma_a$ through a cylinder of parametrized area
at most $C(b-a)\operatorname{Length}(\gamma)$. The constants are
independent of $a,b$ once $a\ge a_0$.
\end{enumerate}
\end{prop}

\begin{proof}
We first prove (i). On a fixed small ball about $o$, the Hessian of $\psi$ has a
positive lower bound. Integrating along minimizing radial geodesics
gives a uniform positive lower bound for the outward derivative at
that ball's boundary. Convexity preserves it farther out. This proves
both the linear lower bound and the gradient lower bound outside
the ball. We choose $a_0$ larger than the maximum of $\psi$ on that
ball. The upper bound follows from $|Y|\le L$. This proves (i). 

We now prove (ii). The bounded field $Y$ is complete. If $v$ is a variational vector
along its inward flow, then we obtain
\[
 \frac{\mathrm{d}}{\mathrm{d}t}|D\Psi_t(v)|^2
 =-2\Hess\psi(D\Psi_t(v),D\Psi_t(v))\le0.
\]
Consequently $|D\Psi_t|\le1$. While $\psi\ge a$, its inward
derivative is $-|Y|^2\le-C^{-2}$, proving finite hitting and the time
bound. The implicit function theorem at the regular target level gives
\[
 \mathrm{d}\tau_a(v)=\frac{\mathrm{d}\psi_{P_a(x)}(D\Psi_{\tau_a}v)}{|Y(P_a(x))|^2},
 \qquad
 DP_a(v)=\left(I-\frac{Y\otimes Y^\flat}{|Y|^2}\right)_{P_a(x)}
                                  D\Psi_{\tau_a}v.
\]
The last factor is contracting and the first is orthogonal projection.
This proves the differential estimates. Extending $P_a$ by the
identity inside $C_a$ gives a locally Lipschitz map across its smooth
boundary, with derivative norm at most one on both sides. Applying
this bound along minimizing geodesics proves the global distance
inequality; no convexity of the exterior domain is needed. This proves (ii).

For (iii), differentiation at a fixed $s$ gives
\[
 D_xH=D\Psi_{s\tau_a}-sY(H)\otimes \mathrm{d}\tau_a,
 \qquad \partial_sH=-\tau_a Y(H).
\]
These give the stated bounds. Monotonicity of $\psi$ along each
trajectory keeps its values between $a$ and $b$. Integrating
$|\partial_sH|\,|D_xH(\frac{\mathrm d\gamma}{\mathrm dt})|$
over a parameter $t$ of the loop proves the area estimate.
Piecewise smooth or Lipschitz loops follow by approximation.
This proves (iii) and completes the proof.
\end{proof}

The exterior gauge construction uses polynomial area bounds for the
auxiliary levels and their outward transport. We establish these bounds
before constructing the limiting connection.

\begin{prop}
\label{prop:levelgeometry}
Let $(M^n,g,\psi,o)$ satisfy Proposition~\ref{prop:auxiliaryflow}.
Assume fixed-radius ambient normal charts with uniformly comparable
metrics, uniformly bounded curvature, and a positive lower bound
for the injectivity radius, as in
Lemma~\ref{lem:boundedlocalgeometry}. After increasing $a_0$, there are constants $r_0,v_0,C>0$
independent of $a\ge a_0$ such that
\[
 \Area(\Sigma_a)\le C a^{n-1},\qquad
 \operatorname{Vol}_{\Sigma_a}
       (B_{\Sigma_a}(x,r_0))\ge v_0\quad(x\in\Sigma_a),
 \qquad \operatorname{diam}_{\Sigma_a}(\Sigma_a)\le C a^{n-1}.
\]
Each level has a covering by at most $Ca^{n-1}$ intrinsic balls of
radius $2r_0$. Moreover every point of every such level has a graph
chart of fixed size with uniformly bounded slope in ambient normal
coordinates.  The levels are diffeomorphic to $S^{n-1}$ and have positive
intrinsic sectional curvature.
\end{prop}
\begin{proof}
We write ${\mathscr S}(a):=\Area(\Sigma_a)$. The outward flow of
$Y/|Y|^2$ identifies any two positive levels through their compact
intervening band. Its tangential area variation is
\[
 \operatorname{div}_{\Sigma_a}\frac{Y}{|Y|^2}
       =\frac{\operatorname{tr}_{T\Sigma_a}\Hess\psi}{|Y|^2}>0.
\]
Consequently ${\mathscr S}$ is increasing. Coarea, $|Y|\le L$, and the linear
lower bound for $\psi$ give
\[
 \frac{a}{L}{\mathscr S}(a)
 \le\int_a^{2a}\frac{{\mathscr S}(t)}L\,\mathrm{d}t
 \le\operatorname{Vol}(C_{2a})
 \le\operatorname{Vol}(B(o,Ca))\le C'a^n.
\]
The last inequality is Bishop--Gromov volume comparison, since
$\Ric\ge0$. This proves the area bound.

We prove this by a local chart argument that requires no upper bound
for $\Hess\psi$. The proof of Proposition~\ref{prop:auxiliaryflow}
gives a constant $c_1>0$ such that the derivative of $\psi$ along
any minimizing geodesic from $o$ is at least $c_1$ after it leaves
a fixed ball. We fix $\delta>0$ smaller than a fixed fraction of the
uniform ambient injectivity radius. For $x\in\Sigma_a$, choose $p$
a distance $\delta$ before $x$ on such a radial minimizing geodesic.
Increasing $a_0$ makes this entire last segment lie outside that ball.
Thus $\psi(x)-\psi(p)\ge c_1\delta$. For
$q\in B(x,\epsilon)$, with $\epsilon\ll\delta$ fixed, convexity on
the unique geodesic from $p$ to $q$ yields
\[
 \mathrm{d}\psi_q(\nabla d_p)
 \ge\frac{\psi(q)-\psi(p)}{d(p,q)}
 \ge\frac{c_1\delta-L\epsilon}{\delta+\epsilon}.
\]
We choose $\epsilon$ small enough that this is at least $c_1/2$.
The distance function $d_p$ is smooth in this neighborhood, and its
Hessian is uniformly bounded there: its distance from $p$ stays
between $\delta/2$ and $2\delta$, within the uniform normal radius.
The ambient bounded geometry therefore allows normal coordinates
$(u,t)\in\mathbb R^{n-1}\times\mathbb R$ centered at $x$ with the
$t$ direction equal to $\nabla d_p(x)$, and, after another uniform
shrinking,
\[
 \partial_t\psi\ge c_1/4=:k>0,\qquad
 |\partial_u\psi|\le C L.
\]
Indeed, the coordinate $t$ direction and $\nabla d_p$ differ by
$O(\epsilon)$ with uniform constant, and $|\nabla \psi|\le L$.
We take a small fixed coordinate height $h$ and then a horizontal radius
$\rho<kh/(2CL)$, with the entire cylinder in this chart. At its
upper and lower faces, respectively,
\[
 \psi(u,h)>a,\qquad \psi(u,-h)<a\qquad (|u|<\rho).
\]
Strict monotonicity in $t$ and the implicit function theorem give
exactly one graph $t=t(u)$ of $\Sigma_a$ in this cylinder, with
$t(0)=0$ and $|Dt|\le CL/k$. All radii and constants are independent
of $x$ and $a$.

The ambient chart metric and this slope estimate bound the induced
graph metric above and below by fixed multiples of the Euclidean
metric. In particular, a horizontal disk of a sufficiently small
fixed radius maps into an intrinsic ball of radius $r_0$, and its
induced volume is bounded below by a fixed $v_0>0$. This proves the
claimed intrinsic ball estimate, reducing $r_0$ once if needed.
The same argument gives such a lower bound for any smaller fixed
radius used below.

The unique minimum of $\psi$ is nondegenerate. Its sufficiently small
positive levels are spheres by Morse coordinates. The normalized
outward gradient flow across compact regular bands proves that every
positive level is a sphere, in particular connected. On a minimizing
intrinsic geodesic, centers spaced $3r_0$ apart have disjoint intrinsic
$r_0$ balls. Their volume lower bounds and the area estimate bound
the length by $Ca^{n-1}$. This proves the diameter estimate. A maximal
$2r_0$-separated set likewise has disjoint $r_0$ balls and cardinality
at most $Ca^{n-1}$; its $2r_0$ balls cover the level. Finally its
second fundamental form with outward normal is
$\Hess\psi|_{T\Sigma_a}/|Y|$, which is positive definite. The Gauss
equation and ambient positive sectional curvature give positive
intrinsic sectional curvature.
This finishes the proof.
\end{proof}

\begin{prop}
\label{prop:transportedfillings}
Assume the hypotheses of Propositions~\ref{prop:auxiliaryflow}
and~\ref{prop:levelgeometry}, with $n\ge3$. Fix $a_*\ge a_0$.
For $a\ge a_*$, define
\[
 G_a:=\bigl(P_{a_*}|_{\Sigma_a}\bigr)^{-1}:
                  \Sigma_{a_*}\longrightarrow\Sigma_a.
\]
Equivalently, $G_a(x)$ is the value at level $a$ of the solution of
$\partial_aG_a=Y(G_a)/|Y(G_a)|^2$, $G_{a_*}(x)=x$.
Set
\[
 g_{\Sigma_a}:=g|_{T\Sigma_a},\qquad h_a:=G_a^*g_{\Sigma_a}.
\]
Then $h_b\ge h_a$ for $b\ge a\ge a_*$, and
\[
 \int_{\Sigma_{a_*}}J_k(\mathrm dG_a)\,\mathrm d\mu_{h_{a_*}}
 \le C(1+a)^{n-1},\qquad 1\le k\le n-1.
\]
Here $J_k(\mathrm dG_a)(x)$ is the largest factor by which
$\mathrm dG_a$ expands a $k$-dimensional volume, i.e., 
\[
 J_k(\mathrm dG_a)(x):=
 \sup_{\substack{v_1,\ldots,v_k\in T_x\Sigma_{a_*}\\
                         h_{a_*}(v_i,v_j)=\delta_{ij}}}
 \left(\det\!\left[g_{\Sigma_a}(\mathrm dG_a v_i,
                         \mathrm dG_a v_j)\right]_{i,j=1}^k\right)^{1/2}.
\]
It is the operator norm of $\bigwedge^k\mathrm dG_a$, using the
metrics induced by $h_{a_*}$ and $g_{\Sigma_a}$ on the two exterior
powers. If $\sigma_1\ge\cdots\ge\sigma_{n-1}>0$ are the singular
values in these metrics, then $J_k=\sigma_1\cdots\sigma_k$.
\end{prop}
\begin{proof}
The projection $P_a|_{\Sigma_b}$ is inverse to the outward level
flow from $\Sigma_a$ to $\Sigma_b$. Its differential has norm at
most one by Proposition~\ref{prop:auxiliaryflow}. Thus, the outward
map expands every tangent vector and $h_b\ge h_a$. In particular,
all singular values of $\mathrm{d}G_a$, measured from $h_{a_*}$ to the
induced metric on $\Sigma_a$, are at least one. The product of the
largest $k$ singular values is therefore at most their full product:
\[
 J_k(\mathrm{d}G_a)\le J_{n-1}(\mathrm{d}G_a),\qquad
 \int_{\Sigma_{a_*}}J_{n-1}(\mathrm{d}G_a)\,\mathrm{d}\mu_{h_{a_*}}
       =\Area(\Sigma_a)\le C(1+a)^{n-1}.
\]
Here the area bound is Proposition~\ref{prop:levelgeometry}.
This proves the integral estimate; it does not assert a uniform
pointwise bound for $\mathrm{d}G_a$.

\end{proof}

We next use the transported area bounds to construct a limiting
connection in a radial gauge.

\begin{prop}
\label{prop:radialconnectionlimit}
Assume the geometric hypotheses and use $G_a,h_a,\Sigma_{a_*}$ of
Proposition~\ref{prop:transportedfillings}. Let $D$ be a smooth
metric connection on a bundle over $\{\psi\ge a_*\}$ with
\[
 |\mathcal F_D|\le C_p(1+\psi)^{-p},\qquad p>n.
\]
Here $\mathcal F_D$ is the curvature of $D$, the endomorphism-valued
two-form defined by
\[
 \mathcal F_D(V,W):=D_VD_W-D_WD_V-D_{[V,W]}.
\]
Identify the bundles over $\Sigma_a$ by parallel transport along
outward level trajectories, and let $D_a$ be the pulled-back
tangential connection on $\Sigma_{a_*}$.
There is a connection $D_\infty$, whose difference from a fixed
smooth reference connection is in $L^1$, satisfying, for every $a\ge a_*$,
\[
 \|D_\infty-D_a\|_{L^1(h_{a_*})}
 +\|D_\infty-D_a\|_{L^1(h_a)}\le C(1+a)^{n-p}.
\]
Each norm uses the indicated metric for both covectors and volume,
and the fixed fiber metric.
\end{prop}
\begin{proof}
We pull the connections back to a fixed level by outward parallel
transport. Integrating the curvature estimate along the flow gives
the limiting connection and its quantitative tail bound.

In this parallel gauge the connection has no component in the
level direction. The curvature identity consequently gives, for a
fixed reference tangent vector $v$,
\[
 (\partial_sD_s)(v)
   =\mathcal F(\partial_sG_s,\mathrm{d}G_s v),
\]
with the fiber endomorphism pulled back by the same isometric
transport. For a fixed $a$, set
$G_{a\to s}:=G_s\circ G_a^{-1}:\Sigma_a\to\Sigma_s$ and denote
by $D_s^{(a)}$ the connection pulled back to the moving reference
level $\Sigma_a$. The speed of this map is bounded by the constant $C$ in Proposition~\ref{prop:auxiliaryflow}.
With the induced metric $g_{\Sigma_a}$ we obtain
\[
 \int_{\Sigma_a}|\partial_sD_s^{(a)}|_{g_{\Sigma_a}}
                        \,\mathrm d\mu_{g_{\Sigma_a}}
 \le C(1+s)^{-p}\int_{\Sigma_a}J_1(\mathrm dG_{a\to s})
                        \,\mathrm d\mu_{g_{\Sigma_a}}
 \le C(1+s)^{n-1-p}.
\]
Pulling this estimate back by $G_a$ is exactly the corresponding
$L^1(h_a)$ estimate on $\Sigma_{a_*}$.
For the last inequality, every singular value of
$\mathrm dG_{a\to s}$ is at least one, so
\[
 J_1(\mathrm dG_{a\to s})\le J_{n-1}(\mathrm dG_{a\to s}),\qquad
 \int_{\Sigma_a}J_{n-1}(\mathrm dG_{a\to s})\,\mathrm d\mu_{g_{\Sigma_a}}
 =\Area(\Sigma_s)\le C(1+s)^{n-1}.
\]
This is the change of variables estimate of
Proposition~\ref{prop:transportedfillings}. All constants are uniform
in $s\ge a\ge a_*$. Integration from $a$ to infinity proves that
$D_s$ is Cauchy in $L^1(h_a)$ with the claimed tail bound. Starting
with $a=a_*$ constructs a single limit. For every other fixed $a$,
the two smooth metrics on the compact reference sphere are equivalent,
so this is the same limit, proving the second assertion as well.
Here an $L^1$ connection means a fixed smooth connection plus an
$L^1$ endomorphism-valued one-form; this is invariant under smooth
changes of reference trivialization.
The assertion follows.
\end{proof}

The pulled-back curvature tends to zero in $L^1$, since its norm is
at most $C(1+a)^{-p}J_2(\mathrm{d}G_a)$ and its integral is
$O((1+a)^{n-1-p})$. This does not allow passage to a limit in
${\mathsf a}_a\wedge {\mathsf a}_a$ from $L^1$ convergence alone.
Here $\mathsf a_a$ is the matrix-valued connection one-form, and
\[
 (\mathsf a_a\wedge\mathsf a_a)(V,W)
 :=\mathsf a_a(V)\mathsf a_a(W)-\mathsf a_a(W)\mathsf a_a(V),
 \qquad \mathcal F_{D_a}=\mathrm d\mathsf a_a+\mathsf a_a\wedge\mathsf a_a.
\]
 The next lemma obtains
a weak parallel frame directly, avoiding that quadratic limit.

\begin{lem}
\label{lem:Loneparallel}
Let $E$ be a metric vector bundle over a compact simply connected
smooth manifold $N$. Suppose smooth metric connections $D_j$ converge
to a connection $D_\infty$ in $L^1$, measured relative to a fixed
smooth reference connection and metric, and
\[
 \|\mathcal F_{D_j}\|_{L^1(N)}\longrightarrow0.
\]
Then $E$ has a global orthonormal frame in $W^{1,1}$ whose columns
$e_\alpha$ satisfy $D_\infty e_\alpha=0$ almost everywhere.
\end{lem}
\begin{proof}
The main point is to construct a parallel frame without taking a
limit of the quadratic connection term. We use radial gauges, average
the curvature estimate over their centers, and pass to the limit in
the first-order frame equation.

We work first in a coordinate ball represented by a Euclidean convex
ball $B$, with a fixed smooth orthonormal trivialization. All
coordinate and geometric norms are uniformly equivalent on this
fixed chart. We write ${\mathsf a}_j$ for the connection matrix. For any
$y\in B$, parallel transport a fixed orthonormal basis at $y$ along
straight segments to obtain a smooth orthogonal frame matrix
$g_j^y(x)$. Its connection matrix $B_j^y$ is in radial gauge and
satisfies
\[
 |B_j^y(x)|\le C|x-y|\int_0^1
       t|\mathcal F_{D_j}(y+t(x-y))|\,\mathrm{d}t.
\]
Indeed, the radial component vanishes, and integrating the curvature
identity along the segment gives this formula with curvature conjugated
by orthogonal transports, which preserve its norm.

Figure~\ref{fig:radialgauge} illustrates the radial construction and
the averaging used to choose the centers.
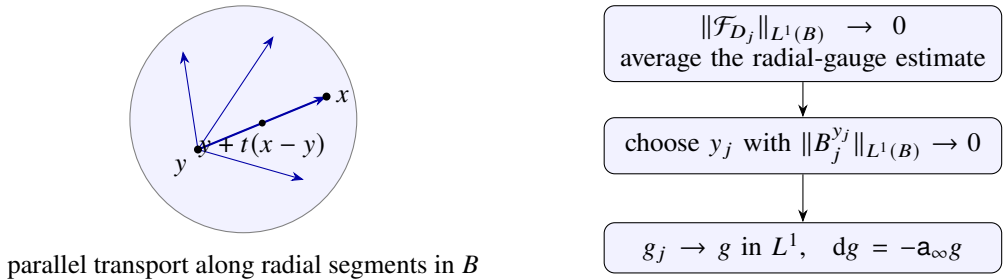
\begin{figure}[!htbp]
\centering
\begin{tikzpicture}[>=Stealth,font=\small]
\begin{scope}
 \fill[blue!5] (0,0) circle (1.5);
 \draw[gray] (0,0) circle (1.5);
 \coordinate (y) at (-.6,-.4);
 \foreach \x/\v in {-.8/.9,.4/1.1,1.1/.3,.8/-.8} {
  \draw[->,blue!65!black] (y)--(\x,\v);
 }
 \draw[thick,blue!65!black] (y)--(1.1,.3);
 \fill (y) circle (1.5pt) node[below left] {$y$};
 \fill (1.1,.3) circle (1.5pt) node[right] {$x$};
 \fill (.25,-.05) circle (1.3pt) node[below] {$y+t(x-y)$};
 \node at (0,-1.95) {parallel transport along radial segments in $B$};
\end{scope}
\begin{scope}[xshift=7.4cm]
 \node[draw,rounded corners,align=center,text width=5cm,fill=blue!5] (a) at (0,1)
 {$\|\mathcal F_{D_j}\|_{L^1(B)}\to0$\\average the radial-gauge estimate};
 \node[draw,rounded corners,align=center,text width=5cm,fill=blue!5] (b) at (0,-.35)
 {choose $y_j$ with $\|B_j^{y_j}\|_{L^1(B)}\to0$};
 \node[draw,rounded corners,align=center,text width=5cm,fill=blue!5] (c) at (0,-1.7)
 {$g_j\to g$ in $L^1$,\quad $\mathrm{d}g=-{\mathsf a}_\infty g$};
 \draw[->] (a)--(b);\draw[->] (b)--(c);
\end{scope}
\end{tikzpicture}
\caption{The local construction in Lemma~\ref{lem:Loneparallel}.
The left panel is drawn in coordinate space. Curvature is integrated
along each segment and then averaged over both endpoints.
The chosen gauges converge strongly in $L^1$; their derivative
identity supplies the weak parallel frame.}
\label{fig:radialgauge}
\end{figure}

Average over both $x,y\in B$. For $0\le t\le1/2$ change variables
from $y$ to $z=(1-t)y+tx$, with $x$ fixed; the Jacobian factor is
at most $2^{\dim N}$. For $1/2\le t\le1$ instead change variables
from $x$ to $z$, with $y$ fixed, with the same bound. Convexity keeps
$z$ in $B$, and $|x-y|\le\operatorname{diam}B$. Consequently, we have
\[
 \int_B\int_B|B_j^y(x)|\,\mathrm{d}x\,\mathrm{d}y
       \le C\operatorname{Vol}(B)\,
                     \|\mathcal F_{D_j}\|_{L^1(B)}.
\]
We choose $y_j$ for which the inner integral is no larger than this
average, and abbreviate $g_j=g_j^{y_j}$, $B_j=B_j^{y_j}$. Thus
$\|B_j\|_{L^1(B)}\to0$. With the frame convention just chosen,
\[
 \mathrm{d}g_j=g_jB_j-{\mathsf a}_jg_j.
\]
In the fixed local trivialization this gives the quantitative bound
\[
 \|\mathrm{d}g_j\|_{L^1(B)}
       \le C\bigl(\|B_j\|_{L^1(B)}+\|{\mathsf a}_j\|_{L^1(B)}\bigr)
       \le C.
\]
The compactness statement used next is strong $L^1$ compactness of
bounded $W^{1,1}$ families on a smaller ball, or on a ball with its
standard extension operator. It does not assert weak compactness of
the derivatives in $L^1$.

The matrices $g_j$ are uniformly bounded, and this identity bounds
their $W^{1,1}$ norms. Compactness of the embedding
$W^{1,1}(B)\hookrightarrow L^1(B)$ gives a subsequence converging
in $L^1$ and almost everywhere to an orthogonal matrix $g$.
Since ${\mathsf a}_j\to {\mathsf a}_\infty$ in $L^1$ and $g_j$ are uniformly bounded, we have
\[
 {\mathsf a}_jg_j\longrightarrow {\mathsf a}_\infty g\quad\hbox{in }L^1(B).
\]
Indeed, pointwise matrix multiplication gives
\[
 \|{\mathsf a}_jg_j-{\mathsf a}_\infty g\|_{L^1}
 \le C\|{\mathsf a}_j-{\mathsf a}_\infty\|_{L^1}
                 +\int_B |{\mathsf a}_\infty|\,|g_j-g|\,\mathrm d x\longrightarrow0.
\]

For the term with ${\mathsf a}_\infty$ this is dominated convergence; the
remaining term is bounded by $C\|{\mathsf a}_j-{\mathsf a}_\infty\|_{L^1}$.
Also, $g_jB_j\to0$ in $L^1$. Passing to distributional derivatives
therefore gives $\mathrm{d}g=-{\mathsf a}_\infty g\in L^1$. In particular, $g$ belongs
to $W^{1,1}$ and its columns are weakly parallel. This argument
avoids any passage to a limit in ${\mathsf a}_j\wedge {\mathsf a}_j$.

We take a finite good cover by such small coordinate balls and perform
this construction on every ball. On a connected overlap, two weak
parallel orthonormal frames differ by an orthogonal matrix $C(x)$.
The product rule for bounded $W^{1,1}$ functions and metric
compatibility give $\mathrm{d}C=0$ almost everywhere. Hence $C$ is constant
on that overlap. These constant transition matrices obey the cocycle
identities. They define a local system; its transport is invariant
under homotopies because the constants agree on triple overlaps.
Since $N$ is simply connected, transport from one fixed chart gives
constant changes of basis making all the local frames agree.
The resulting global frame is orthonormal almost everywhere,
belongs to $W^{1,1}$, and is parallel for $D_\infty$.
This proves the stated conclusion.
\end{proof}

\begin{cor}
\label{cor:weakexteriorframe}
Use the complete positively curved manifold, convex exhaustion
$\psi$, base point $o$, and regular levels $\Sigma_a$ of
Proposition~\ref{prop:auxiliaryflow}.
Assume $n\ge3$ and the ambient local geometry of
Proposition~\ref{prop:levelgeometry}. Let $E$ be a metric vector bundle over
$\{\psi\ge a_*\}$ with a smooth metric connection $D$ and
\[
 |\mathcal F_D|\le C_p(1+\psi)^{-p},\qquad p>n+1.
\]
Then $E$ has an orthonormal frame $e$ in $W^{1,1}_{\mathrm{loc}}$
on this exterior, parallel along outward level trajectories, with
\[
 \int_{\Sigma_a}|De|\,\mathrm{d}\mu_{\Sigma_a}
       \le C(1+a)^{n-p},\qquad
 \int_{\{\psi\ge a\}}|De|\,\mathrm{d}\mu_g
       \le C(1+a)^{n+1-p},
\]
for any $a\ge a_*$; $De$ is the weak covariant derivative and the
level estimate holds for almost every $a$ for its Sobolev traces.
Constants are independent of $a$.
\end{cor}
\begin{proof}
We apply Proposition~\ref{prop:radialconnectionlimit} on the reference
sphere $\Sigma_{a_*}$. The pulled-back curvature tends to zero in
$L^1$ by the $J_2$ estimate there. This sphere is simply connected
for $n\ge3$, so Lemma~\ref{lem:Loneparallel} supplies a global
weak parallel frame for $D_\infty$. We extend it by parallel transport
along the level trajectories, making it independent of the level
parameter in the outward parallel gauge. It is locally $W^{1,1}$:
on each compact band the flow and the connection are smooth, while
the reference frame is $W^{1,1}$ and bounded. Its radial covariant
derivative vanishes. Its tangential derivative at level $a$ is
\[
 D_a e=(D_a-D_\infty)e.
\]
Since the frame is orthonormal, the $L^1(h_a)$ estimate of
Proposition~\ref{prop:radialconnectionlimit} proves the first bound.
The gradient of $\psi$ is normal to the levels, so this also controls
the full covariant derivative. Coarea and $|\nabla\psi|\ge C^{-1}$ give
\[
 \int_{\{\psi\ge a\}}|De|\,\mathrm{d}\mu_g
 \le C\int_a^\infty
                 C(1+s)^{n-p}\,\mathrm{d}s
 \le C'(1+a)^{n+1-p}.
\]
This proves the tail bound. No pointwise bound for the inverse
projection's differential has been used.
This finishes the proof.
\end{proof}

\begin{prop}
\label{prop:smoothexteriorgauge}
Use the complete positively curved manifold, convex exhaustion
$\psi$, base point $o$, and regular levels $\Sigma_a$ of
Proposition~\ref{prop:auxiliaryflow}.
Assume $n\ge3$ and fixed-radius ambient charts with uniform metric,
inverse metric, and connection coefficient bounds through order $m+3$,
and positive injectivity radius, as supplied by
Lemma~\ref{lem:boundedlocalgeometry} in the soliton application. Fix $m\ge0$ and $p>n+1$. Suppose a smooth metric connection on an exterior
metric vector bundle satisfies
\[
 |D^j\mathcal F_D|(x)\le C_j(1+\psi(x))^{-p},
                  \qquad 0\le j\le m+1,
\]
where derivatives include the ambient Levi--Civita connection on
covariant indices. There are $a_1\ge a_*$ and constants $C'_0,\ldots,C'_m$, depending
on these geometric and curvature bounds, $m$, and $p$, such that on
$\{\psi>a_1\}$ there is a smooth
orthonormal frame $\widetilde e$ such that
\[
 |D^j(D\widetilde e)|(x)
       \le C'_j(1+\psi(x))^{-(p-n-1)},\qquad 0\le j\le m.
\]
\end{prop}
\begin{proof}
Smooth the weak parallel frame in fixed local charts. We then
project the smoothed frame back to the orthonormal frame bundle and
control its covariant derivatives.

We set $q:=p-n-1>0$ and take the weak orthonormal frame of
Corollary~\ref{cor:weakexteriorframe}. We choose a fixed small radius
$\rho$ below the uniform normal radius. Since $\psi$ is Lipschitz,
its weight is comparable throughout every $8\rho$ ball sufficiently
far out. The tail estimate in that corollary implies, on any such
ball centered at $x$,
\[
 \int_{B(x,8\rho)}|De|\,\mathrm{d}\mu_g
                    \le C(1+\psi(x))^{-q}.
\]
We use normal coordinates and a radial parallel orthonormal bundle
frame on this ball. The radial-gauge curvature formula gives its
connection matrix ${\mathsf a}$ the bounds
\[
 \|{\mathsf a}\|_{C^{m+1}}
                 \le C(1+\psi(x))^{-p}.
\]
These finite-order bounds follow by differentiating the radial
transport equation and curvature formula inductively. The first
bound is the integral of curvature along radial segments. At higher
orders, coordinate derivatives of curvature differ from covariant
ones by the already controlled connection coefficients, while the
ambient coordinate coefficients have uniform bounds. Thus, the
induction has uniform constants on this fixed-radius ball.

We represent the weak frame by an orthogonal matrix-valued function
$u$. The identity $\mathrm{d}u=De-\mathsf a\,u$ gives
\[
 \|\mathrm{d}u\|_{L^1(B(x,8\rho))}
                    \le C(1+\psi(x))^{-q}.
\]
The Euclidean $L^1$ Poincar\'e inequality in the normal coordinate
ball, with uniform metric comparisons, bounds $u$ in $L^1$ from
its mean by the same quantity. Some almost-everywhere value of $u$
is within a fixed multiple of this bound from the mean. Denoting
that orthogonal matrix by $Q_x$, we obtain
\[
 \|u-Q_x\|_{L^1(B(x,8\rho))}
                    \le C(1+\psi(x))^{-q}.
\]
Mollify $u$ at a fixed coordinate radius smaller than $\rho$ and
restrict to $B(x,4\rho)$. The resulting smooth local frame $s_x$
need not be orthonormal, but its coordinate matrix is within
$C(1+\psi(x))^{-q}$ of $Q_x$ in $C^{m+1}$. Hence it is uniformly
bounded, its Gram matrix differs from $I$ by this amount, and
\[
 |D^j(Ds_x)|\le C(1+\psi(x))^{-q},\qquad 0\le j\le m.
\]
It also differs from $e$ in $L^1$ on this smaller ball by the same
bound. Constants incorporate the fixed mollifier radius.

We choose a uniformly locally finite normal ball cover with a smooth
partition $\{\chi_i\}$ supported in $B(x_i,2\rho)$ and with uniform
bounds on its derivatives through order $m+1$. Such a cover and
partition follow from the ambient injectivity and curvature bounds.
We use the local frames $s_i$ on $B(x_i,4\rho)$. If two partition
supports meet at $z$, both frames are defined on $B(z,\rho)$.
There their difference has $L^1$ norm $O((1+\psi(z))^{-q})$ and
covariant derivative of the same order in the supremum norm.
Integrating along short geodesics in this ball, and using its
uniform volume lower bound, bounds their difference at $z$ by
$C(1+\psi(z))^{-q}$. Their higher covariant derivatives have the
same bound by the estimates already proved for each frame.

We set $s:=\sum_i\chi_i s_i$. In differentiated partition terms subtract
one active $s_i$, using $\sum_i \mathrm{d}\chi_i=0$ and the analogous
identities for higher derivatives. The overlap estimates give
\[
 |D^j(Ds)|\le C(1+\psi)^{-q},\qquad
 |s^*s-I|\le C(1+\psi)^{-q}.
\]
Choose $a_1$ so large that the preceding estimate gives
$|s^*s-I|<1/2$ on $\{\psi>a_1\}$. Then $s^*s\ge\tfrac12I$ there. The single pointwise polar normalization
$\widetilde e=s(s^*s)^{-1/2}$ is smooth and orthonormal.
Differentiating this formula, using metric compatibility and the
previous bounds, proves all the claimed covariant derivative
estimates. Here a varying normalization is allowed: its derivatives
are bounded by the same small error, rather than being discarded.
The assertion follows.
\end{proof}

We now apply Proposition~\ref{prop:smoothexteriorgauge} to the
local approximate rotations. Assume from now on that $(M^4,g,f)$
satisfies the hypotheses of Theorem~\ref{thm:main} and
\eqref{eq:normalization} and is not the four-dimensional Bryant soliton. Proposition~\ref{prop:curvaturelength}
first transfers their polynomial error bounds to all comparison
centers. Proposition~\ref{prop:rotationconnection} then constructs a
rank-three bundle determined by these pairs, with a metric connection
whose curvature has polynomial decay. The smooth orthonormal frame
provided above yields global vector fields by taking its vector field
components; Proposition~\ref{prop:globalrotations} records their
estimates.

\begin{prop}
\label{prop:curvaturelength}
Fix $o\in M$. There is $C<\infty$ such that
\[
 |\nabla \mathrm{R}|\le C \mathrm{R}^{3/2},\qquad
 \mathrm{R}(x)^{-1/2}\le C(1+d(o,x)),\qquad
 \frac{\mathrm{R}(x)^{-1/2}}{d(o,x)}\longrightarrow0.
\]
Fix $L,T,C_0>0$. Put $r_x:=\mathrm R(x)^{-1/2}$. Suppose that $(y,t)$ satisfies
\[
 -Tr_x^2\le t\le0,\qquad d_g(x,y)\le Lr_x,\qquad
 C_0^{-1}\le r_x^2\mathrm R_{g_t}(y)\le C_0.
\]
Then for $z:=\Phi_t(y)$, we have
\[
 1+d(o,x)\le C(1+d(o,z))^2.
\]
The constant in this distance estimate depends only on
$(M,g,f)$, $o$, and $L,T,C_0$, and is independent of $x,y,t$.
For an intrinsic cap center $p\in\Gamma$, take $x=p$ and
$r_p:=\mathrm R(p)^{-1/2}$. Thus every
\[
 (y,t)\in B_g(p,Lr_p)\times[-Tr_p^2,0],\qquad
 C_0^{-1}\le r_p^2\mathrm R_{g_t}(y)\le C_0,
\]
satisfies $1+d(o,p)\le C(1+d(o,\Phi_t(y)))^2$ with the same
uniform constant. On each fixed enlarged Bryant-product chart,
Proposition~\ref{prop:inputcomparison} supplies such $C_0$
uniformly for $p$ outside a compact set.
\end{prop}
\begin{proof}
If the first inequality failed outside compact sets, an escaping
sequence with unbounded $|\nabla \mathrm{R}|/\mathrm{R}^{3/2}$ would contradict smooth
curvature scale dimension reduction: every subsequential normalized
limit is a smooth $B^3\times\mathbb R$ product or shrinking cylinder with
base scalar one. On a compact set the ratio is bounded because
$\mathrm{R}>0$. Hence $|\nabla \mathrm{R}^{-1/2}|\le C/2$ globally. Integrating along
a minimizing geodesic from $o$ proves the second inequality.

The strict asymptotic statement is Proposition~\ref{prop:sublinearcurvscale}. For the distance estimate put $r:=\mathrm R(x)^{-1/2}$.
The upper scalar bound at the given point yields
\[
 \mathrm R(z)=\mathrm R_{g_t}(y)\le C_0r^{-2},\qquad
 r\le\sqrt{C_0}\,\mathrm R(z)^{-1/2}\le C(1+d(o,z)).
\]
Since $d_g(x,y)\le Lr$, $|t|\le Tr^2$, and $|X|\le1$,
we obtain
\[
 d(o,x)\le d(o,z)+Lr+Tr^2
             \le C'(1+d(o,z))^2.
\]
Taking $x=p\in\Gamma$ and $r=r_p$ proves the stated cap center case.
This proves the stated conclusion.
\end{proof}

We apply the connection construction to the pairs
$(U,\operatorname{skew}\nabla U)$ formed from the local approximate
rotations. The resulting rank-three bundle records their values and
skew-adjoint covariant derivatives.

\begin{prop}
\label{prop:rotationconnection}
Set $w(x):=1+d_g(o,x)$ and
$\mathcal J:=TM\oplus\mathfrak{so}(TM)$, with its natural metric
and connection. Fix $r_0>0$ below the uniform normal radius and a
cover
\[
 \mathcal B_i:=B(x_i,r_0),\qquad
 \mathcal B_i^+:=B(x_i,4r_0)
\]
with bounded overlap multiplicity and a subordinate partition of
unity with uniformly bounded derivatives through order eight.
Let $U_i=(U_{i,1},U_{i,2},U_{i,3})$ be smooth on
$\mathcal B_i^+$, and define
\[
 T_i v:=\left(\sum_{\alpha=1}^3v_\alpha U_{i,\alpha},\
       \operatorname{skew}\nabla\sum_{\alpha=1}^3v_\alpha U_{i,\alpha}\right),
 \qquad v\in\mathbb R^3.
\]
Here $\nabla U$ is the endomorphism $v\mapsto\nabla_vU$, and
\[
 \operatorname{skew}\nabla U
 :=\tfrac12\bigl(\nabla U-(\nabla U)^*\bigr)
                       \in\mathfrak{so}(TM).
\]
Thus the second component of $T_i v$ is the skew-adjoint
part of the covariant derivative of $\sum_\alpha v_\alpha U_{i,\alpha}$.
Suppose, for fixed $b\ge0$, we have
\[
 \|T_i\|_{C^8(\mathcal B_i^+)}\le Cw(x_i)^b,\qquad
 |T_i(x)v|\ge C^{-1}w(x_i)^{-b}|v|
 \quad(x\in\mathcal B_i^+),
\]
and on every nonempty enlarged overlap there is a constant
$O^{ij}\in O(3)$ with
\[
 \|T_i-T_jO^{ij}\|_{C^8(\mathcal B_i^+\cap\mathcal B_j^+)}
 \le Cw(x_i)^{-N}.
\]
Then there is a finite $L=L(b)$, independent of $N$, such that the following hold for sufficiently large $N$ outside a compact set.
\begin{enumerate}[label=\textup{(\roman*)},ref=\textup{(\roman*)}]
\item There are a smooth
rank-three subbundle $E\subset\mathcal J$, a metric $h$ on $E$, and
local $h$-orthonormal frames $V_i$ with
\[
 \|V_i-T_i\|_{C^5(\mathcal B_i)}\le Cw(x_i)^{-N+L}.
\]
\item There is a metric connection $D$ on $(E,h)$ such that
\[
 |D^j\mathcal F_D|(x)\le C_jw(x)^{-N+L},\quad 0\le j\le3.
\]
In frame $V_i$, the connection matrix has $C^4(\mathcal B_i)$ norm
at most $Cw(x_i)^{-N+L}$. Frame norms use the ambient bundle metric;
curvature norms use $h$ and $g$. Constants depend on the fixed atlas
bounds and $b$.
\end{enumerate}
\end{prop}
\begin{proof}
We first prove (i). We organize the local first-order data into a rank-three bundle and
compare the resulting local connections on overlaps. The small errors
of the approximate rotations will control the curvature of the
induced metric connection.

All estimates below are local on the fixed enlarged balls. Their
weights are comparable, the number of neighbors is uniformly
bounded, and each operation uses only a fixed number of derivatives.
Here $T_i^*$ denotes the adjoint for the Euclidean metric on
$\mathbb R^3$ and the natural direct-sum metric on $\mathcal J$.
We set
\[
 P_i:=T_i(T_i^*T_i)^{-1}T_i^*,
\]
the orthogonal projection onto $\operatorname{im}T_i$ in this
ambient metric. The metric $h$ constructed below is a separate
metric on $E$; its orthonormal frames are not obtained by
normalizing $T_i$ in the ambient metric.
The lower singular-value bound gives
$\|(T_i^*T_i)^{-1}\|\le Cw^{2b}$. For an invertible matrix field $G=T_i^*T_i$, the first two
derivative identities are
\[
 \nabla G^{-1}=-G^{-1}(\nabla G)G^{-1},
\]
\begin{align*}
 \nabla_i\nabla_jG^{-1}={}&
 G^{-1}(\nabla_iG)G^{-1}(\nabla_jG)G^{-1}\\
 &+G^{-1}(\nabla_jG)G^{-1}(\nabla_iG)G^{-1}
       -G^{-1}(\nabla_i\nabla_jG)G^{-1}.
\end{align*}
In general $\nabla^jG^{-1}$ is a finite sum of products
\[
 G^{-1}(\nabla^{i_1}G)G^{-1}\cdots
               (\nabla^{i_\ell}G)G^{-1},
 \qquad i_\nu\ge1,\quad \sum i_\nu=j.
\]
Consequently, the power of $w$ lost at any fixed derivative order is
independent of the requested decay exponent $N$.

Differentiating the inverse formula
repeatedly proves polynomial bounds through order eight, with
exponents depending only on $b$ and the order. The same formula,
applied to $T_i$ and $T_jO^{ij}$ and using invariance under constant
orthogonal changes, gives
\[
 \|P_i-P_j\|_{C^8}\le Cw^{-N+L_1(b)}.
\]
No small cocycle identity is needed for this estimate.

Let $P^0:=\sum_i\chi_iP_i$. It is self-adjoint and, on each active
chart, differs from $P_i$ by $O_{C^8}(w^{-N+L_1})$. For $N>L_1$
and sufficiently far out, its spectrum is in small disjoint
neighborhoods of zero and one. Let $P$ be the spectral projection
onto the latter cluster and put $E:=\operatorname{im}P$. This has
rank three. A contour formula for $P$, with a fixed contour around
one, and differentiation of the resolvent show that $P$ has polynomial
$C^8$ bounds and
\[
 \|P-P_i\|_{C^8}\le Cw^{-N+L_2(b)}.
\]
To retain the small factor, we subtract the two resolvent formulas.
Every resulting term contains $P^0-P_i$ or one of
its derivatives, and all remaining factors have polynomial bounds.
We set $S_i:=PT_i$. These are frames of $E$, close to $T_i$ with a fixed
polynomial loss, and
$S_i-S_jO^{ij}=P(T_i-T_jO^{ij})$ has the same form of estimate.

We define $h_i$ on $E$ by declaring $S_i$ orthonormal, and set
$h:=\sum_i\chi_i h_i$. The matrices
$C^{ji}=S_j^{-1}S_i$ satisfy
\[
 C^{ji}=O^{ij}+O_{C^8}(w^{-N+L_3(b)}).
\]
Here the left inverse of $S_j$ and its derivatives have polynomial
bounds by the singular-value estimate, which persists after the
small projection error. Consequently, the matrix of $h$ in the
frame $S_i$ is
\[
 H_i=h(S_i,S_i)=I+O_{C^8}(w^{-N+L_4(b)}).
\]
Derivatives of the partition introduce no nondecaying term, because
the constant contribution is exactly $I\sum_j\chi_j=I$.
We define $V_i:=S_iH_i^{-1/2}$. These frames are $h$-orthonormal and
$V_i-T_i=O_{C^8}(w^{-N+L_5(b)})$ after another fixed loss. Their
exact transition matrices $\tau^{ji}$ are orthogonal and obey
\[
 V_i=V_j\tau^{ji},\qquad
 \tau^{ji}=O^{ij}+O_{C^8}(w^{-N+L_6(b)}).
\]
This choice of metric is essential: directly orthonormalizing $T_i$
with the natural metric on $TM\oplus\mathfrak{so}(TM)$ would not give a change close to the
identity and could alter the vector field metric Lie derivative at leading order. This proves (i). 

We now prove (ii). Let $D_i$ be the flat metric connection declaring $V_i$ parallel
on its chart, and form the affine combination $D:=\sum_i\chi_iD_i$.
It is a globally defined metric connection because the coefficients
sum to one. In the frame $V_i$, each $D_j-D_i$ is expressed by a
matrix $(\tau^{ji})^{-1}\mathrm{d}\tau^{ji}$, up to the frame convention.
It and its derivatives are small, since $O^{ij}$ is constant.
For example, with $V_i=V_j\tau^{ji}$ and column-frame
conventions, the connection matrix of $D_j$ in the frame $V_i$ is
$(\tau^{ji})^{-1}\mathrm{d}\tau^{ji}$. Thus, we have
\[
 {\mathsf a}_i=\sum_j\chi_j(\tau^{ji})^{-1}\mathrm{d}\tau^{ji},\qquad
 \mathrm{d}O^{ij}=0.
\]
Every term of $\nabla^k\mathsf a_i$ contains at least one differentiated
factor $\tau^{ji}-O^{ij}$. In particular there is no term of the
size of $\mathrm{d}\chi_j$ alone. This cancellation gives the required smallness of the connection
curvature.

The partition bounds therefore give a $C^4$ connection matrix
estimate of the required form. In this orthonormal frame,
$\mathcal F_D=\mathrm{d}{\mathsf a}_i+{\mathsf a}_i\wedge {\mathsf a}_i$, and its first three covariant
derivatives are controlled by that estimate and the ambient local
geometry. Taking the maximum of the finitely many power losses
above gives $L(b)$. The use of eight derivatives leaves margins
beyond the stated five and four derivatives. All smallness conditions
are met by choosing $N$ larger than these fixed losses and then
increasing the compact exclusion.
This proves (ii) and completes the proof.
\end{proof}

\begin{cor}
\label{cor:rotationglobalization}
Assume the hypotheses of Proposition~\ref{prop:rotationconnection} and
the exhaustion hypotheses of Proposition~\ref{prop:auxiliaryflow},
with $n\ge3$. Fix $q>0$. For sufficiently large $N$ there is
$a_1$ such that on $\{\psi>a_1\}$ the bundle $(E,h)$ has a smooth
global orthonormal frame, and there are smooth fields
$U=(U_1,U_2,U_3)$ satisfying, on each
$\mathcal B_i\cap\{\psi>a_1\}$, for some $B_i:\mathcal B_i\cap\{\psi>a_1\}\to O(3)$, we have
\[
 \|U-U_iB_i\|_{C^3}\le Cw(x_i)^{-q},\qquad
 \sum_{j=1}^3|\nabla^jB_i|\le Cw(x_i)^{-q}.
\]
Here $(U_iB_i)_\alpha:=\sum_\beta U_{i,\beta}(B_i)_{\beta\alpha}$.
There is a fixed exponent $b_1<\infty$, depending only on the
local polynomial growth exponent $b$ and the fixed geometry, such that
\[
 \sum_{\alpha=1}^3\bigl(|U_\alpha(x)|+|\nabla U_\alpha(x)|\bigr)
 \le C(1+d_g(o,x))^{b_1}\qquad(\psi(x)>a_1).
\]
The choice of $N,a_1$ may depend on $q$ and the fixed exponent $b$.
\end{cor}
\begin{proof}
Since $\psi$ and $w$ are comparable outside a compact set, apply
Proposition~\ref{prop:smoothexteriorgauge} with $m:=2$ and
$p=N-L(b)$. It gives a global $h$-orthonormal frame $e$ with
covariant derivatives through order three bounded, after the first,
by $Cw^{-(N-L(b)-n-1)}$. We write $e:=V_iB_i$ on a chart. Both frames
are orthonormal, so $B_i\in O(3)$. The connection matrix in $V_i$
is $O_{C^4}(w^{-N+L(b)})$. The equation
$De=V_i(\mathrm{d}B_i+{\mathsf a}_iB_i)$ and its differentiated versions bound the
first three ordinary covariant derivatives of $B_i$ by the sum of
these two decay rates.

Let $U$ be the first component of $e\in E\subset\mathcal J$.
Projection $\mathcal J\to TM$ is parallel for the ambient direct-sum
connection. Since $V_i-T_i$ is small with its first five derivatives, we have
\[
 U=U_i B_i+O_{C^3}(w^{-N+L'(b)}).
\]
Products with the bounded polynomial derivatives of $T_i$ lose only
another fixed power. We choose $N$ larger than the resulting losses,
$n+1$, and the requested $q$. This proves all assertions. The
orthogonality of $B_i$ and the local polynomial bounds also prove
the size and first derivative bounds for $U$.
This finishes the proof.
\end{proof}

The preceding construction starts from approximate overlap matrices
and produces a bundle with exact transitions. We next verify its
hypotheses for the local rotation atlas and establish the global
Gram, bracket, and transverse estimates.

\begin{prop}
\label{prop:localrotationatlas}
Fix $o\in M$
and $w(x):=1+d_g(o,x)$. For every sufficiently large $N$ and every
$q>0$, outside a compact set, there is an atlas as in
Proposition~\ref{prop:rotationconnection}, with $b=1$ and
$O^{ij}\in SO(3)$, whose local triples satisfy
\[
 \sum_{j=0}^8\bigl(|\nabla^j\Lie_{U_{i,\alpha}}g|
                  +|\nabla^j(U_{i,\alpha}f)|\bigr)\le Cw^{-N}
 \quad(\alpha=1,2,3)
\]
on the enlarged fixed domains. On the balls, for cyclic
$(\alpha,\beta,\gamma)$, we have
\[
 \sum_{j=0}^3\left|\nabla^j
 ([U_{i,\alpha},U_{i,\beta}]-\krot U_{i,\gamma})\right|
 \le C_qw^{-q}.
\]
On every whole collar sphere $S$, we have
\[
 \sum_{\alpha=1}^3\bigl(|\langle U_{i,\alpha},n_f\rangle|
                    +|\langle U_{i,\alpha},\nu\rangle|\bigr)
 +|Q_S(U_i)-I|_{\mathrm F}\le C_qw^{-q}.
\]
The weight $w$ can be evaluated at any point of the sphere.
Constants, the core width, and the compact exclusion may depend on
$N,q$.
\end{prop}
\begin{proof}
We fix the enlarged domain sizes and collar margins first. In particular,
choose cap domains with a strict margin around the core
$\{\beta\le2\Lambda\}$ and around the collar needed for comparisons.
We choose neck fixed domains with margins containing every whole leaf
through a smaller neighborhood. The domain constructions
and local improvement proofs allow these fixed enlargements before
the heat parameter and comparison accuracy are selected.

Choose a starting exponent $P$ larger than twice the desired
exponent plus all fixed losses from rescaling. By
Proposition~\ref{prop:curvaturelength}, every required center $z$
in the additional local construction at $x$ satisfies
\[
 1+d(o,z)\ge C^{-1}(1+d(o,x))^{1/2},\qquad
 (1+d(o,z))^{-P}\le C(1+d(o,x))^{-P/2}.
\]
Apply Proposition~\ref{prop:typedclosure} with this $P$, then
perform one additional local improvement at each desired center,
retaining the higher-order conclusion estimates. For a neck these are
provided by the derivative estimates in the cylindrical proof;
for a cap use Corollary~\ref{cor:alltiltheat}, which retains eight
metric Lie derivative and potential pairing derivatives. The displayed distance comparison applies to every intermediate
center and cap replacement in this fixed finite construction.

The heat-corrected fields are close, through the required derivatives, to
one normalized model rotation triple on their fixed domains. They
are time-independent. Their constant Gram normalization and the
constant structure-coefficient argument give
\[
 [U_{i,\alpha},U_{i,\beta}]
       =\krot\epsilon_{\alpha\beta\gamma}U_{i,\gamma}
             +\text{small error}.
\]
For caps, the last assertion follows from the cap version of
Claim~\ref{cylinder:lem:normalized-brackets}: its structure coefficients
are constant throughout the connected cap and its Gram normalization
is measured on the collar. Retaining three derivatives in this
argument uses fewer derivatives than the available derivative reserve.
All normal estimates are imposed only on the regular collars.

We explain the scale conversion. If $r_i$ is the neck curvature
length or the cap tip length, model closeness gives
$|\nabla^jU_i|\le C_j r_i^{1-j}$ for $0\le j\le9$ on the smaller domains. The local potential has the form $F=r_i^{-1}f$ plus a
spatial constant at the time. Hence converting $U_iF$ to
$U_if$ costs at most one factor $r_i$. Vector errors also cost at
most this factor, whereas metric Lie derivative derivatives have factors $r_i^{-j}$.
Since $r_i\ge1$ and $r_i\le Cw$, these are fixed polynomial losses.
We increase the initial exponent to absorb them. Every fixed spatial domain has diameter $O(r_i)=o(d(o,x_i))$ by
Proposition~\ref{prop:curvaturelength}; consequently its distances
from $o$ are comparable to the center's distance. For core centers,
use the tip center in this observation. This proves the asserted
absolute polynomial estimates throughout the fixed domains.

We take a fixed small-radius normal ball cover of the remaining exterior,
with fixed further enlargements and bounded overlap. At a core point
use a cap triple based at its intrinsic projection on $\Gamma$;
elsewhere use a neck triple. All enlarged small balls fit in their
chosen domains. The radius can be fixed because $r_i\ge1$ and all
normalized spatial margins are fixed and positive.

There are three overlap cases. For two core charts the intrinsic
projection is uniformly controlled in tip units, and nearby centers
have comparable scales. Their cap domains therefore contain a common
regular collar, a whole leaf, and paths of bounded normalized length
from that collar to the small-ball overlap. We choose the small cover
radius within the reserved margins. For two neck charts the whole-leaf
margins give a common whole sphere and its open collar. For a core/neck
pair the centers lie within the fixed small cover radius of
$\beta=2\Lambda$. The function $\beta$ has bounded derivative in
comparison units there, so both neighborhoods lie in the interior
of the cap collar $\{\Lambda/2<\beta<2J\Lambda\}$ after reserving
margins. They have the same whole leaves by
Proposition~\ref{prop:exteriorleaves}. Thus in every case
Proposition~\ref{prop:loworderalignment} supplies one constant
$O^{ij}\in O(3)$ with small $C^3$ difference, including at core
points by transport from the common collar.

For this same constant matrix, set $E:=U_i-U_jO^{ij}$ and
$h_\alpha:=\Lie_{E_\alpha}g$ for each component. Each $h_\alpha$
is small through order eight. Starting from the already small
$C^3$ norm of $E_\alpha$, differentiate
\[
 \nabla^2E_\alpha=\Rm*E_\alpha+\mathcal A(\nabla h_\alpha),
\]
where $\mathcal A$ is the fixed index permutation in the Killing
transport identity. Induction through seven further derivatives of
this identity uses at most eight derivatives of $h_\alpha$ and gives the
field difference through order nine. The corresponding
pairs $(U,\operatorname{skew}\nabla U)$ differ through order eight. Converting scales
and using $r_i\le Cw$ costs only fixed powers already allowed by the
choice of exponent. This proves the overlap hypothesis for $T_i$.

The $C^8$ size bound for $T_i$ holds with $b=1$. For the lower bound, in each
normalized model the map from the three rotational coefficients to
value and skew derivative is injective at every point: a Killing
field with vanishing value and covariant derivative at one point is
zero on the connected model. Its least
singular value is uniformly positive on the compact normalized
measurement sets, including the cap axis. Model closeness preserves
this bound. In the original metric it reads
\[
 r_i^{-2}|U_i v|^2+
          |\operatorname{skew}\nabla(U_i v)|^2\ge C^{-1}|v|^2.
\]
Since $r_i\ge1$, it implies the required unweighted lower bound for $T_i$,
in fact with no loss in $w$.

Finally the overlap matrices preserve orientation. The local triples
have the standard bracket convention with coefficient $\krot$.
An orthogonal change of basis transforms this coefficient by
$\det O^{ij}$. We compare the bracket identities using the $C^1$
overlap estimate, then use the nondegenerate whole-leaf Gram matrix
to bound their constant coefficient discrepancy. This gives
$|1-\det O^{ij}|\le C\epsilon$, in comparison units. For the
chosen small threshold, the determinant cannot be $-1$. Thus all
$O^{ij}$ belong to $SO(3)$, as asserted.
The assertion follows.
\end{proof}

The bases of the local triples can now be aligned and the fields interpolated.
The following proposition records the resulting global approximate
rotations and the decay of their errors.

Every fixed polynomial loss can be absorbed by increasing the starting
decay exponent: if $e\le C(1+r)^{-P}$, then
$(1+r)^me\le C(1+r)^{-p}$ for $P\ge p+m$.
Comparable distance weights have the same decay exponent, with a
changed constant. We use this convention below.

\begin{prop}
\label{prop:globalrotations}
Fix $o\in M$
and $r(x):=d_g(o,x)$. For every $p>0$ there are a compact set
$K_p$, a constant $C_p$, and smooth fields
$\mathcal U=(\mathcal U_1,\mathcal U_2,\mathcal U_3)$ on
$M\setminus K_p$ such that, for $i=1,2,3$, we have
\[
 |\mathcal U_i|+|\nabla\mathcal U_i|\le C_p(1+r),\qquad
 \sum_{j=0}^2\bigl(|\nabla^j\Lie_{\mathcal U_i}g|
                  +|\nabla^j(\mathcal U_if)|\bigr)
 \le C_p(1+r)^{-p}.
\]
For cyclic $(i,j,k)$,
$|[\mathcal U_i,\mathcal U_j]-\mathcal U_k|\to0$ uniformly as
$r\to\infty$. On the common CMC spheres of
Proposition~\ref{prop:exteriorleaves}, we have
\[
 |Q_S(\mathcal U)-\krot^{-2}I|_{\mathrm F}
 +\sup_S\sum_{i=1}^3
   (|\langle\mathcal U_i,n_f\rangle|
    +|\langle\mathcal U_i,\nu\rangle|)\longrightarrow0
 \quad\text{as }\inf_S r\to\infty.
\]
And there is a sequence of such whole spheres with $\inf_Sr\to\infty$.
\end{prop}
\begin{proof}
We apply Proposition~\ref{prop:localrotationatlas}, then
Proposition~\ref{prop:rotationconnection} and
Corollary~\ref{cor:rotationglobalization}. The preceding polynomial-loss convention
permits the initial exponent to absorb every fixed polynomial loss
while retaining the prescribed final decay exponent. The oriented
local frames define an orientation of $E$: their exact transitions
are close to $SO(3)$ and therefore have positive determinant.
We choose the global smooth orthonormal frame with this orientation;
if necessary change one column by a constant sign. In its local
representations the matrices $B_i$ then belong to $SO(3)$.
Write $w:=1+r$. Choose the local error exponent $N>p$ and the frame comparison
exponent $q>p+4$, increasing the initial exponents to absorb the
fixed polynomial losses above. Write the resulting evaluated fields
as $U=U_iB_i+E_i$, with $\|E_i\|_{C^3}+\|\nabla B_i\|_{C^2}
\le Cw^{-q}$.

The errors $E_i$ are small through three derivatives, and the
first three derivatives of $B_i$ have any prescribed polynomial
decay. The local fields have at most linear size, while their
positive-order derivatives are uniformly bounded. For each component the product identities are
\[
 \begin{aligned}
 \Lie_{U_a}g
 &=\sum_\alpha B_{i,\alpha a}\Lie_{U_{i,\alpha}}g
   +2\sum_\alpha\mathrm dB_{i,\alpha a}\odot U_{i,\alpha}^{\flat}
   +\Lie_{E_{i,a}}g,\\
 U_af&=\sum_\alpha B_{i,\alpha a}(U_{i,\alpha}f)+E_{i,a}f.
 \end{aligned}
\]
Thus their $C^2$ norms on each fixed ball are at most
$C(w^{-N}+w^{1-q})\le C_pw^{-p}$. For the potential term use the local bounds on $U_if$
and the uniform fixed derivatives of $f$, which follow from the
soliton equation and ambient curvature derivative bounds. The same
product rule gives the stated linear size bound. Since a constant
$SO(3)$ matrix preserves the cyclic bracket convention, all bracket
errors contain a local bracket defect, a derivative of $B_i$, or an
error $E_i$ and its first derivative. Their sum is bounded by $C(w^{-q}+w^{2-q})$, which tends to zero. The local normal estimates
likewise imply absolute decay of the two normal components of $U$.

We verify the whole-leaf normalization; pointwise rank two of the
vector fields alone is not sufficient for this step. We fix a sphere $S$ and one local constructed triple $U_0=(U_{0,1},U_{0,2},U_{0,3})$ defined on its whole
leaf neighborhood. Every atlas chart meeting $S$ can be compared
directly to this triple on an enlarged common whole-leaf collar,
with controlled normalized paths if it is a cap chart. The same
alignment argument gives a constant $O_i\in SO(3)$, with all needed
small derivatives of $U_i-U_0O_i$ near its part of $S$.
Define $T_0$ from $U_0$ and $T$ from the global fields $U$ by
the same formula $v\mapsto(Uv,\operatorname{skew}\nabla(Uv))$ as for $T_i$. On the leaf
neighborhood set
\[
 C:=(T_0^*T_0)^{-1}T_0^*T.
\]
The lower singular-value bound makes this well-defined. The local
representations and the small derivatives of $B_i$ imply
$C=O_iB_i+\text{small error}$, with small first derivative. Indeed
forming $(U_iB_i,\operatorname{skew}\nabla(U_iB_i))$ adds terms containing $\mathrm{d}B_i$;
these have already been made arbitrarily small. Differentiating
the displayed left inverse costs only fixed polynomial factors.
After increasing the starting exponents for the fixed polynomial
factors, these estimates give, with $w_S:=1+\inf_S r$,
\[
 \operatorname{dist}(C,SO(3))+|\nabla C|\le Cw_S^{-q},\qquad
 |U-U_0C|\le Cw_S^{-q}\quad\text{on }S.
\]

The sphere has diameter $O(r_S)$ and scalar comparable to
$r_S^{-2}$. Its minimum and maximum distances to $o$ are comparable
outside a compact set. For example the curvature length estimate
at a minimum-distance point gives
$r_S\le C(1+\min_S r)$, and adding the diameter proves this comparison.
Choose $B_S\in SO(3)$ closest to $C$ at one point of $S$.
Integrating $\nabla C$ along its intrinsic geodesics gives
$\sup_S|C-B_S|\le Cw_S^{1-q}$, since
$\operatorname{diam}(S)\le Cr_S\le Cw_S$.
Consequently, with $\delta_S:=\sup_S|U-U_0B_S|$,
\[
 \delta_S\le Cw_S^{2-q}\longrightarrow0,\qquad
 |Q_S(U)-I|_{\mathrm F}
 \le |Q_S(U_0)-I|_{\mathrm F}
       +C\left(\frac{\delta_S}{r_S}
                     +\frac{\delta_S^2}{r_S^2}\right)\longrightarrow0.
\]
The last estimate uses $|U_0|\le Cr_S$,
$\Area(S)\asymp r_S^2$, and $r_S\ge C^{-1}>0$.

Escaping spheres exist as follows. If the tip set has an
unbounded sequence, take a fixed collar sphere in the cap neighborhood at
each such tip. Its distance from the tip is $O(r_{\rm tip})$, which
is little-oh of the tip distance. These spheres escape. If the tip
set is bounded, Proposition~\ref{prop:corecover} gives a compact
set outside which every point is a neck center, and their whole neck spheres escape by
the same curvature length estimate. The empty-tip case is included.
Finally set $\mathcal U_\alpha:=U_\alpha/\krot$. This gives bracket
coefficient one and Gram limit $\krot^{-2}I$, while preserving all
other conclusions.
This proves the stated conclusion.
\end{proof}

\section{Global correction and exact symmetry}\label{sec:exact}
In this section, we correct the global approximate rotations to
Killing fields and complete the symmetry proof.
Proposition~\ref{prop:globalcorrection} solves the drift equation,
and Corollary~\ref{cor:globalKilling} uses integrable tensor rigidity
to obtain exact rotations. We then construct the angular action in
Proposition~\ref{prop:orbitcompletion}, establish the critical point
and properness of the potential in
Proposition~\ref{reflcrit:critical}, and determine the metric blow-down
in Theorem~\ref{lem:metric-blowdown}.

\subsection{Global drift correction}
We write $r(x):=d(o,x)$ and $\mathcal D:=\Delta+\nabla_X$, with the rough
Laplacian on vector fields and the scalar Laplacian on functions.
The analytic equations in this subsection are on all of $M$.
In the application, an approximate field defined on $\mathcal E_a$
is extended by $\widetilde U=\chi(\psi)U$, where
$\chi=0$ on $(-\infty,a]$ and $\chi=1$ on $[b,\infty)$ for a
fixed $b>a$. It is set to zero inside $\{\psi\le a\}$; derivatives
of the cutoff are supported in the compact band
$\mathcal A_{a,b}$, where
\[
 \mathcal E_a:=\{\psi>a\},\qquad
 \mathcal A_{a,b}:=\{a\le\psi\le b\}.
\]
The comparison with the original field is made
on $\mathcal E_b$. If a Dirichlet exhaustion is used somewhere,
we always mean to choose smooth relatively compact domains $\Omega_j$ with
\[
 \overline B_g(o,2^j)\subset\Omega_j
       \Subset B_g(o,2^{j+1}),\qquad j\to\infty;
\]
interior limits solve the equation on $M$.

\begin{prop}
\label{prop:positive-resolvent-kernel}
Let $(M^n,g)$ be complete and connected, with $n\ge3$,
$\Ric\ge0$, and $v_0:=\inf_x\operatorname{Vol}B(x,1)>0$.
The positive Green kernel $G$ of $-\Delta+1/4$ satisfies
\[
 G(x,y)\le C(n,v_0)
 \begin{cases}
 d(x,y)^{2-n},&0<d(x,y)\le1,\\
 (1+d(x,y))e^{-d(x,y)/2},&d(x,y)\ge1.
 \end{cases}
\]
For a metric vector bundle with a smooth metric connection and a
smooth scalar potential $V\ge1/4$, the operator norm of the Green
kernel of $-\Delta+V$ satisfies the same estimate. The constant is independent of
$V$ and the bundle rank.
\end{prop}

\begin{proof}
Since the result is quite standard, we only sketch the heat kernel argument, using the Li--Yau Harnack inequality
\cite{LiYau} and the Davies--Gaffney estimate
\cite[Section 3]{CS}. Write $d=d(x,y)>0$ and
$\ell:=\min\{1,\sqrt t,t/d\}$. Harnack in each variable, followed by
averaging on $B(x,\ell)\times B(y,\ell)$ and Davies--Gaffney, gives
\[
 p_t(x,y)\le C\ell^{-n}
 \exp\!\left(-\frac{(d-2\ell)_+^2}{4(t+2\ell^2)}\right).
\]
Here Bishop--Gromov and the unit-ball lower bound give
$\operatorname{Vol}B(x,\ell)\ge v_0\ell^n$.
Since $\ell^2\le t$ and $d\ell\le t$, the exponent differs from
$-d^2/(4t)$ by at most $3/2$. Thus
\begin{equation}\label{eq:exactgaussian}
 p_t(x,y)\le C\max\{1,t^{-n/2},(d/t)^n\}e^{-d^2/(4t)}.
\end{equation}
Integrate $e^{-t/4}p_t$ over $t>0$. For $d\le1$, splitting at
$d^2$ and $1$ gives $G(x,y)\le Cd^{2-n}$. For $d\ge1$, use
\[
 \frac t4+\frac{d^2}{4t}=\frac d2+\frac{(t-d)^2}{4t}.
\]
On $[d/2,2d]$ the prefactor in \eqref{eq:exactgaussian} is bounded,
so its contribution is at most $Cd e^{-d/2}$. On the two remaining
intervals, $(t-d)^2/(4t)\ge t/16$ or $d^2/(16t)$, respectively;
integration gives at most $C(1+d)e^{-d/2}$.
For a metric connection, the norm Kato inequality and $V\ge1/4$
dominate the bundle heat semigroup by $e^{-t/4}e^{t\Delta}$ on
scalar norms. Integration gives the same Green-kernel bound,
independently of $V$ and the bundle rank.
\end{proof}

We now apply this Green-kernel estimate to the drift operator. The
soliton identities allow us to conjugate it to a connection operator
with a positive scalar potential. Throughout this subsection, let
$n\ge3$ and let $(M^n,g,f)$ be a complete gradient steady soliton
that is $\kappa$-noncollapsed at all scales, with positive curvature
operator and
\[
 \Ric=\Hess f,\qquad \mathrm R+|\nabla f|^2=1,\qquad
 \inf_x\operatorname{Vol}B(x,1)>0.
\]
Lemma~\ref{lem:boundedlocalgeometry} verifies all the uniform local
coefficient, injectivity radius, and volume hypotheses we will use.

\begin{prop}
\label{prop:globalcorrection}
Fix $o\in M$ with $\mathrm R(o)>0$, and put
$r(x):=d_g(o,x)$, $X:=\nabla f$, and
$\mathcal D:=\Delta+\nabla_X$. If a smooth vector field $Q$ satisfies
\[
 |Q(x)|\le C_Q(1+r(x))^{-p},\qquad p>n+1,
\]
then $\mathcal DV=Q$ has a unique smooth solution tending uniformly
to zero at infinity, and
\begin{equation}\label{eq:correctiondecay}
 |V(x)|+|\nabla V(x)|\le C(1+r(x))^{-(p-n-1)}.
\end{equation}
If we also assume
$|\nabla Q|\le C(1+r)^{-p}$, the same holds for $|\nabla^2V|$.
Constants depend only on the stated geometric and source bounds.
\end{prop}
\begin{proof}
We conjugate the drift operator to a positive connection operator.
We then estimate its inverse separately for sources near the base point
and sources in the exterior. This yields the claimed polynomial decay.

For any smooth scalar $\varphi$ and vector field $W$, the product
rule for the metric connection gives
\[
 \mathcal D(\varphi W)
 =\varphi\Delta W+2\nabla_{\nabla\varphi}W
     +\varphi\nabla_XW+(\Delta\varphi+X\varphi)W.
\]
With $\varphi=e^{-f/2}$, the first-order terms in $W$ cancel because
$2\nabla\varphi=-\varphi X$. Moreover, we have
\[
 \frac{\Delta\varphi+X\varphi}{\varphi}
   =-\frac12\Delta f-\frac14|X|^2
   =-\frac14(1+\mathrm{R}).
\]
Thus, the exact conjugation is
\begin{equation}\label{eq:driftconjugation}
 \mathcal D(e^{-f/2}W)
 =-e^{-f/2}\mathcal H W,
 \qquad \mathcal H:=-\Delta+\frac{1+\mathrm{R}}{4}.
\end{equation}
These equalities use $\Delta f=\mathrm{R}$ and $|\nabla f|^2=1-\mathrm{R}$. The connection
operator $\mathcal H$ is self-adjoint on unweighted $L^2$ and is
bounded below by $1/4$. We denote its inverse kernel by $G_{\mathcal H}$.
The norm kernel for the conjugated inverse is at most
\[
 e^{(f(y)-f(x))/2}|G_{\mathcal H}(x,y)|.
\]
Since $|\nabla f|\le1$, Proposition~\ref{prop:positive-resolvent-kernel} bounds this
by $C d(x,y)^{2-n}$ for $d(x,y)\le1$ and by $C(1+d(x,y))$ otherwise.

An additional estimate is essential for sources near the fixed base.
We set $a_o:=|X(o)|=\sqrt{1-\mathrm{R}(o)}<1$. Convexity of $f$ along a minimizing
geodesic from $o$, and its Lipschitz bound, give
\[
 f(x)\ge f(o)-a_o r(x),\qquad f(y)\le f(o)+r(y).
\]
We write $r:=r(x)$ and $s:=r(y)$ and choose
$\kappa_o=(1-a_o)/4>0$. If $s\le\kappa_o r$, then we obtain
\begin{equation}\label{eq:innergreensource}
 \frac{f(y)-f(x)-d(x,y)}2
 \le-\frac{1-a_o}{2}r+s
 \le-\frac{1-a_o}{4}r.
\end{equation}
This is where positive scalar curvature at a single point is used.

For a compactly supported smooth cutoff $Q_j$ of $Q$ with
$|Q_j|\le|Q|$, define
\[
 V_j(x):=-e^{-f(x)/2}\int_M G_{\mathcal H}(x,y)
                              e^{f(y)/2}Q_j(y)\,\mathrm{d}\mu(y).
\]
This formula is legitimate for compact sources by the $L^2$
resolvent, and \eqref{eq:driftconjugation} gives
$\mathcal D V_j=Q_j$. We next bound it independently of $j$.
For large $r$, sources with $s\le\kappa_o r$ have $d(x,y)\ge1$.
Their contribution is at most
$C(1+r)e^{-(1-a_o)r/4}\int_M|Q|\,\mathrm d\mu_g$, by
\eqref{eq:innergreensource}; the integral is finite since $p>n$.
For $s>\kappa_o r$ and $d(x,y)\ge1$, the coarse kernel bound gives
\[
 \int_{s>\kappa_o r}C(1+r+s)(1+s)^{-p}\,\mathrm{d}\mu(y)
 \le C(1+r)^{n+1-p}.
\]
To see the power explicitly, decompose the region into annuli
$s\asymp 2^j\kappa_o r$, $j\ge0$. Bishop--Gromov gives volume
at most $C(2^jr)^n$ on each annulus, while
$1+r+s\le C(2^jr)$ there, with constants depending on $\kappa_o$.
Hence the integral is bounded by
\[
 C r^{n+1-p}\sum_{j=0}^{\infty}2^{j(n+1-p)}
      \le C' r^{n+1-p}.
\]
This is the loss of $n+1$ powers in the inverse estimate; the series
converges precisely for the indicated sufficient condition $p>n+1$.
The remaining unit ball about $x$ contributes at most $C(1+r)^{-p}$,
because $d^{2-n}$ is uniformly locally integrable under the same
volume upper bound. These estimates prove \eqref{eq:correctiondecay}
uniformly for $V_j$.

We choose cutoffs equal to one on an exhausting sequence of compact
sets. Local elliptic estimates and the preceding uniform bound give
a subsequence converging locally to a smooth solution $V$ with the
same bound. This avoids applying the unweighted $L^2$ resolvent
directly to the possibly non-$L^2$ source $e^{f/2}Q$.
Uniform unit-scale $W^{2,q}$ estimates with $q>n$ give the first
derivative bound; local Schauder estimates using the first derivative
of $Q$ give the second. The weight $1+r$ is comparable on unit balls.
For uniqueness, the difference $Z$ of two such solutions satisfies
$\mathcal D|Z|^2=2|\nabla Z|^2\ge0$ and tends to zero. If it were
nonzero, its positive maximum would be attained in a compact set,
contradicting the strong maximum principle.
This finishes the proof.
\end{proof}

We next prove a vanishing result for integrable solutions of the
stationary drift--Lichnerowicz equation. This will apply to the metric Lie derivative
of each corrected field. Let $o$ be a fixed base point and set $r=d_g(o,\cdot)$.
\begin{prop}
\label{prop:distanceLone}
Let $h$ be a smooth symmetric
covariant two-tensor satisfying
\[
 (\Delta+\nabla_X)h+2\Rm(h)=0,
 \qquad h\in L^1(\mathrm{d}\mu).
\]
Then $h\equiv0$. In particular, this applies when
$|h|\le C(1+r)^{-s}$ for some $s>n$.
\end{prop}
\begin{proof}
We apply the norm inequality for the stationary tensor equation and
then integrate it against distance cutoffs. Integrability removes
the boundary terms, while positive sectional curvature forces the
tensor to vanish.

We set $q:=\lambda_{\max}(\Ric)$. In an orthonormal frame
diagonalizing $h$, nonnegative sectional curvature gives
\[
 \langle\Rm(h),h\rangle
 =\sum_{i\ne j}K_{ij}h_i h_j
 \le\sum_i\Ric_{ii}h_i^2\le q|h|^2.
\]
To justify the inequality at the zeros of $h$, put
$v_\epsilon:=(|h|^2+\epsilon^2)^{1/2}$. The tensor equation and the scalar
chain rule give
\begin{align*}
 \mathcal D v_\epsilon
 &=\frac{|\nabla h|^2-2\langle\Rm(h),h\rangle}{v_\epsilon}
                         -\frac{|\nabla|h|^2|^2}{4v_\epsilon^3}\\
 &\ge-\frac{2q|h|^2}{v_\epsilon}\ge-2qv_\epsilon.
\end{align*}
The first inequality uses
$|\nabla|h|^2|^2\le4|h|^2|\nabla h|^2\le4v_\epsilon^2|\nabla h|^2$.
Testing against nonnegative compactly supported smooth functions
and taking $\epsilon\downarrow0$ yields
$\mathcal D|h|\ge-2q|h|$ distributionally.
Moreover, $b=\mathrm{R}-2q>0$: in a Ricci eigenframe whose first vector
corresponds to $q$,
\[
 \mathrm{R}-2q=2\sum_{2\le i<j\le n}K_{ij}>0.
\]

Positive sectional curvature and the normalization give
$|\Rm|\le C_n\mathrm R\le C_n$. Shi's estimates for the complete
canonical flow and $\kappa$-noncollapsing give uniformly bounded
curvature derivatives and a positive uniform injectivity radius.
We construct distance cutoffs instead of potential cutoffs. A
uniformly locally finite covering by small normal balls admits a
partition of unity $\phi_i$ with uniformly bounded first two
derivatives. If $x_i$ are the centers, set
$\rho:=1+\sum_i\phi_i r(x_i)$. Then $|\rho-(1+r)|\le C$ and
$|\nabla\rho|+|\nabla^2\rho|\le C$: in differentiated sums subtract
any one nearby center value and use the uniformly bounded differences
$|r(x_i)-r(x_j)|$. We add a constant if necessary to make $\rho\ge1$.
It is a smooth proper function, constructed independently of $f$.
For a fixed smooth cutoff $\eta$, equal to one on $[0,1]$ and zero
on $[2,\infty)$, put $\chi_s:=\eta(\rho/s)$. Then we obtain
\[
 |\Delta\chi_s-X\chi_s|\le C/s\quad(s\ge1).
\]
Since $\operatorname{div}X=\mathrm{R}$, the unweighted scalar adjoint is
$\mathcal D^*=\Delta-X-\mathrm{R}$. Testing the Kato inequality with
$\chi_s$ therefore gives
\[
 \int_M\chi_s b|h|\,\mathrm{d}\mu
 \le\int_M|h|(\Delta\chi_s-X\chi_s)\,\mathrm{d}\mu\longrightarrow0.
\]
Dominated convergence applies on the left, since $0<b\le1$ and
$h\in L^1$. Thus $b|h|=0$, which proves the result. The polynomial
criterion follows from the volume upper bound by dyadic annuli.
The assertion follows.
\end{proof}

\begin{cor}
\label{cor:globalKilling}
Let $(M^n,g,f)$ satisfy the geometric assumptions of
Proposition~\ref{prop:distanceLone}. Let $U=(U_1,U_2,U_3)$ be smooth outside a compact set, with
\[
 |U_i|+|\nabla U_i|\le C(1+r)^m,\qquad
 \sum_{j=0}^2\bigl(|\nabla^j\Lie_{U_i}g|
                  +|\nabla^j(U_if)|\bigr)\le C(1+r)^{-p},
\]
where $m\ge0$, $p>\max\{2n+1,m+n+1\}$, and $i=1,2,3$.
Then there are global Killing fields $W_1,W_2,W_3$ with $W_if=0$ and
\[
 |W_i-U_i|+|\nabla(W_i-U_i)|\le C(1+r)^{-(p-n-1)}
\]
on that exterior. If $[U_i,U_j]-U_k\to0$ uniformly at infinity
for cyclic $(i,j,k)$, then $[W_i,W_j]=W_k$.
If this bracket limit holds and whole embedded spheres $S_\ell$ satisfy
\[
\lim_{\ell\to\infty} \inf_{S_\ell}r= \infty,\qquad \Area(S_\ell)\ge a_*>0,\qquad
 \Area(S_\ell)^{-2}\int_{S_\ell}
       \left|\sum_{i=1}^3v_iU_i\right|^2\,\mathrm d\mu
       \ge q_*|v|^2\quad(v\in\mathbb R^3)
\]
for some $q_*>0$, the fields $W_i$ are linearly independent over
constants and generate a nontrivial isometric $SU(2)$ action.
\end{cor}
\begin{proof}
We extend the $U_a$ over the compact core by a smooth cutoff. The identity
\[
 (\mathcal D U_a)^\flat
 =\operatorname{div}(\Lie_{U_a}g)
  -\tfrac12\mathrm{d}\operatorname{tr}(\Lie_{U_a}g)
  +(\Lie_{U_a}g)(X,\cdot)-d(U_af)
\]
shows that $Q_a=\mathcal D U_a$ and its first derivative have the
required decay. Proposition~\ref{prop:globalcorrection} solves
$\mathcal D V_a=Q_a$ with the stated decay through two derivatives.
We set $W_a:=U_a-V_a$. We give the details of the tensor identity used for
these corrected fields. For $h:=\Lie_Wg$, define the Bianchi vector
\[
 \mathcal B(h)=\bigl(\operatorname{div}h-\tfrac12\mathrm{d}\tr h\bigr)^\sharp
                     =\Delta W+\Ric^\sharp W.
\]
Naturality and the Ricci variation formula give
\[
 2\Lie_W\Ric=2(D\Ric)_g[h]
                  =-\Delta_Lh+\Lie_{\mathcal B(h)}g.
\]
On the other hand $2\Ric=\Lie_Xg$ and the Lie derivative commutator
imply
\[
 2\Lie_W\Ric=\Lie_W\Lie_Xg
                   =\Lie_Xh-\Lie_{[X,W]}g.
\]
Equating these expressions yields
\[
 \Delta_Lh+\Lie_Xh=\Lie_{\mathcal B(h)+[X,W]}g
                             =\Lie_{\mathcal D W}g,
\]
since $\nabla_WX=\Ric^\sharp W$. The two algebraic Ricci terms
cancel as follows:
\begin{align*}
 (\Delta_Lh)_{ij}
 &=\Delta h_{ij}+2\Rm(h)_{ij}
              -\Ric_i{}^k h_{kj}-\Ric_j{}^k h_{ik},\\
 (\Lie_Xh)_{ij}
 &=\nabla_Xh_{ij}+\Ric_i{}^k h_{kj}+\Ric_j{}^k h_{ik}.
\end{align*}
Consequently, the exact identity is
\begin{equation}
 (\Delta+\nabla_X)\Lie_Wg+2\Rm(\Lie_Wg)
                         =\Lie_{\mathcal D W}g.
\end{equation}
Thus $h_a=\Lie_{W_a}g$ solves the homogeneous tensor equation
and decays with exponent $p-n-1>n$. We apply
Proposition~\ref{prop:distanceLone} to conclude $h_a=0$.

For any Killing field $K$, $\Hess(Kf)=\Lie_K\Ric=0$.
Its gradient is parallel, hence zero because $\Ric>0$.
Therefore $Kf$ is constant and the displayed vector identity gives
$\mathcal D K=0$. For our $W_a$ the constant is zero, since
$U_af\to0$, $V_a\to0$, and $|X|\le1$.
The bracket defect $Z=[W_a,W_b]-W_c$ is Killing. If
$s=p-n-1>m$, its exact expansion is
\[
 Z=[U_a,U_b]-U_c-[U_a,V_b]-[V_a,U_b]+[V_a,V_b]+V_c.
\]
The estimates for values and first derivatives imply
\[
 |Z|\le |[U_a,U_b]-U_c|
             +C(1+r)^{m-s}+C(1+r)^{-2s}+C(1+r)^{-s}
             \longrightarrow0.
\]
Since $\mathcal D|Z|^2=2|\nabla Z|^2\ge0$, the maximum principle
and decay at infinity imply $Z=0$. This uses ordinary compact
exhaustions, not compact potential sublevels.

If $\sum_a v_aW_a=0$ for a constant unit vector $v$, then
$\sum_a v_aU_a=\sum_a v_aV_a$, and hence
\[
 q_*\le\Area(S_\ell)^{-2}\int_{S_\ell}
       \left|\sum_a v_aU_a\right|^2\,\mathrm d\mu
 \le a_*^{-1}\sup_{S_\ell}\sum_a|V_a|^2\longrightarrow0,
\]
a contradiction. Thus the $W_a$ are linearly independent over constants.
Killing fields on a complete manifold are complete;
their standard Lie brackets therefore integrate to the simply
connected compact group $SU(2)$. No fixed point or critical point
of $f$ was used in this integration step.
This proves the stated conclusion.
\end{proof}

For the particular fields obtained from the actual atlas, the exact
fields in Corollary~\ref{cor:globalKilling} satisfy
\begin{equation}\label{eq:exact-generator-size}
 \left(\sum_{a=1}^3|W_a(x)|^2\right)^{1/2}
                       \le C \mathrm{R}(x)^{-1/2}.
\end{equation}
\begin{proof}[Proof of \eqref{eq:exact-generator-size}]
We write the evaluated approximate triple
on an exterior chart as $U:=U_iB_i+E_i$, as in
Corollary~\ref{cor:rotationglobalization}. Orthogonality gives the exact
pointwise identity
\[
 \sum_a|(U_iB_i)_a|^2=\sum_a|U_{ia}|^2.
\]
The model estimates in Proposition~\ref{prop:localrotationatlas} give
$|U_i|\le Cr_i$. On neck neighborhoods, the curvature length $r_i$ is
comparable to $\mathrm{R}(x)^{-1/2}$; on the fixed-width core neighborhoods,
scalar comparability gives the same comparison for the tip scale.
The core width and the fixed margins are fixed, so these constants
are uniform over the exterior cover. The error $E_i$ decays, the
constant bracket normalization changes only $C$, and the drift
correction satisfies $W-U=o(1)$. Since $\mathrm{R}\le1$, these errors are
absorbed in the right-hand side of \eqref{eq:exact-generator-size}.
Enlarging $C$ handles the remaining compact set. This proves \eqref{eq:exact-generator-size}
\end{proof}

\subsection{Exact Killing fields, the global group action, and the critical geometry}\label{sec:exact-action}

In this subsection, we construct the angular action and establish
the metric blow-down and critical point geometry. We first
obtain a fixed point from an averaged auxiliary convex exhaustion.
Its differential representation determines the action in
Proposition~\ref{prop:orbitcompletion}. Theorem~\ref{lem:metric-blowdown} then identifies the
metric limit.

\begin{prop}
\label{prop:equivariantsoul}
If a compact Lie group $G$ acts isometrically on a complete connected noncompact
manifold of strictly positive sectional curvature, then there is a fixed
point $s$ and a $G$-equivariant diffeomorphism $T_sM\to M$ for the
differential action on $T_sM$.
\end{prop}

\begin{proof}
We choose a smooth strictly convex Lipschitz exhaustion $\psi$ by
\cite[Theorem 1(a)]{GW}, and average it over $G$ with Haar probability
measure. The function
\[
 \bar\psi(x)=\int_G\psi(gx)\,\mathrm{d}g
\]
is smooth, invariant, Lipschitz, and has positive-definite Hessian.
It is proper: the linear lower bound proved in
Proposition~\ref{prop:auxiliaryflow} is uniform for $\psi(gx)$,
since the compact orbit $G\cdot o$ has bounded diameter. Thus
$\bar\psi$ has a unique nondegenerate minimum $s$, fixed by $G$.

Equivariant Morse coordinates identify a neighborhood of $s$ with
its tangent representation, with $\bar\psi-\bar\psi(s)=|v|^2$.
These coordinates can be obtained directly in equivariant normal
coordinates: write the Taylor expansion as $\langle {\mathsf H}(v)v,v\rangle$
with ${\mathsf H}(v)$ positive definite and equivariant, and use the local
diffeomorphism $v\mapsto {\mathsf H}(v)^{1/2}v$.
On a small punctured annulus choose the corresponding radial field
$Z$ with $Z\bar\psi=1$. Outside a slightly larger annulus use
$\nabla\bar\psi/|\nabla\bar\psi|^2$. An invariant cutoff interpolates
between them while preserving $Z\bar\psi=1$. Its flow between any
two positive levels exists because the intervening band is compact
and contains no critical point. Extending the local identification
by this flow gives the global equivariant diffeomorphism. This uses
an auxiliary convex exhaustion, not the soliton potential.
This finishes the proof.
\end{proof}

We next pass from exact Killing fields to a global group action.
The bracket relations and the isotropy representation determine
the orbit structure of that action.

\begin{prop}
\label{prop:orbitcompletion}
Let $(M^4,g,f)$ be a complete connected noncompact steady soliton with positive
sectional curvature and $\Ric=\Hess f$. Suppose linearly independent
complete Killing fields $W_1,W_2,W_3$ satisfy
\[
 [W_1,W_2]=W_3,\qquad [W_2,W_3]=W_1,\qquad [W_3,W_1]=W_2,
 \qquad W_af=0\quad(a=1,2,3).
\] Assume there are escaping points $x_j$ and
unit coefficient vectors $z_j\in\mathbb R^3$ such that
\begin{equation}\label{eq:almostisotropy}
 \left|\sum_a(z_j)_aW_a(x_j)\right|\longrightarrow0.
\end{equation}
Then their action factors through an effective isometric $SO(3)$
action with principal orbits $S^2$. It extends to an effective
isometric $O(3)$ action with the same principal orbits. 
\end{prop}
\begin{proof}
We first integrate the Killing fields to a compact group action.
A fixed point allows us to identify its linear representation, and
the orbit-rank information determines the angular action.

We integrate the Killing algebra to $SU(2)$ and apply
Proposition~\ref{prop:equivariantsoul}. Its differential at the fixed
point $s$ is nontrivial: an isometry fixing a point and acting
identically on its tangent space is the identity. The elementary
classification of real orthogonal $SU(2)$ representations of
dimension four leaves two possibilities: the irreducible quaternionic
representation, or the three-dimensional adjoint representation plus
a trivial line. The entirely trivial representation is excluded by
linear independence of the Killing fields.

Suppose the representation is quaternionic. Since $f$ is invariant,
$\nabla f(s)$ is a fixed vector in this representation, hence zero.
Thus this branch supplies, rather than assumes, a critical point.
Positive Ricci curvature makes it unique. It also makes $f$ proper:
choose a small ball about $s$ on which $\Ric\ge C^{-1}g$; integration
along minimizing geodesics gives a uniform positive radial derivative
at that ball's boundary, and convexity preserves this lower bound
farther out. Consequently $f(x)\ge c_1d(s,x)-C_1$.

We fix a regular level $\Sigma=\{f=f(s)+u\}$ with $u>0$. It is compact.
For large $j$, inward gradient flow takes $x_j$ to a point
$y_j\in\Sigma$ in finite time: its intervening compact potential
band has a positive minimum of $|\nabla f|$. Every Killing field $K$ in this algebra satisfies
\[
 [K,X]=\Lie_K(\nabla f)=\nabla(Kf)=0.
\]
Metric compatibility and $\nabla X=\Ric^\sharp$ consequently give
\[
 X(|K|^2)=2\langle\nabla_XK,K\rangle
        =2\langle\nabla_KX,K\rangle=2\Ric(K,K)\ge0.
\]
For $K_j=\sum_a(z_j)_aW_a$ this gives
$|K_j(y_j)|\le|K_j(x_j)|\to0$.
But global equivariant linearization identifies the quaternionic
action away from $s$ with the free left $SU(2)$ action on
$\mathbb H\setminus\{0\}$. Its infinitesimal action is injective
at every point of $\Sigma$. Compactness of $\Sigma\times S^2$
then gives a positive lower bound for $|\sum_a z_aW_a(y)|$,
contradicting the preceding limit.

The representation is consequently the adjoint representation plus
a trivial line. The central element $-1\in SU(2)$ has identity
differential at $s$, hence acts identically on $M$. The induced
$SO(3)$ action is effective and, by the global equivariant
diffeomorphism, has the standard model on
$\mathbb R^3\times\mathbb R$. We write its coordinates as
$(v,z)$, with $\rho:=|v|$ and $\vartheta:=v/|v|$ on $\rho>0$.
The isotropy $SO(2)$ fixes the two transverse directions and has
no fixed covector on the angular tangent plane. Invariance of the
metric therefore forces
\[
 g=\mathfrak a(\rho,z)\mathrm{d}\rho^2+2B(\rho,z)\mathrm{d}\rho\,\mathrm{d}z
       +C(\rho,z)\mathrm{d}z^2+D(\rho,z)g_{S^2}.
\]
It is invariant under $(v,z)\mapsto(-v,z)$ on this principal set.
This map is smooth on the full equivariant model, including its
axis, so continuity of the metric identity makes it a global
isometry. Together with $SO(3)$ it yields the standard effective
$O(3)$ action. Off the axis the orbits are precisely the spheres
$\{|v|=\rho,\ z=\mathrm{constant}\}$. In this representation branch
no critical point of $f$ has been asserted or used.
The assertion follows.
\end{proof}

\begin{cor}
\label{cor:globalfieldcriterion}
Let $(M^4,g,f)$ be a complete connected noncompact steady soliton
with positive sectional curvature and $\Ric=\Hess f$.
Let $W_1,W_2,W_3$ be linearly independent complete Killing fields with
$W_af=0$ and
\[
 [W_1,W_2]=W_3,\qquad [W_2,W_3]=W_1,\qquad [W_3,W_1]=W_2.
\]
Suppose smooth vector fields $U_a$ satisfy $|U_a-W_a|\to0$ at infinity.
Suppose
there are points $x_j\to\infty$ and orthonormal pairs
$n_{1,j},n_{2,j}\in T_{x_j}M$ such that
\[
 \sum_{a=1}^3\sum_{\ell=1}^2
       |\langle U_a(x_j),n_{\ell,j}\rangle|^2\longrightarrow0.
\]
 Then the $W_a$ generate the effective
isometric $SO(3)$ action with principal orbits $S^2$ of
Proposition~\ref{prop:orbitcompletion}, and this action extends to
an effective isometric $O(3)$ action preserving $f$.
\end{cor}
\begin{proof}
At each $x_j$ the tangential projection of the triple is a linear
map from $\mathbb R^3$ to a two-plane, so choose a unit vector $z_j$
in its kernel. The assumed transverse bounds give
$|\sum_a(z_j)_aU_a(x_j)|\to0$. Since $|U_a-W_a|\to0$ uniformly at infinity,
\eqref{eq:almostisotropy} holds for the corrected Killing triple.
We apply Proposition~\ref{prop:orbitcompletion}.
This proves the stated conclusion.
\end{proof}

We now apply these results to the approximate rotations constructed
in Section~\ref{sec:globalization}. This completes the angular part
of the main theorem.
\begin{proof}[Construction of the angular action in Theorem~\ref{thm:main}]
The four-dimensional Bryant soliton already has the required angular
action. In the other branch, we apply
Proposition~\ref{prop:globalrotations} with $p>9$, $n=4$, and
linear growth exponent $m=1$. Its fields satisfy all the hypotheses
of Corollary~\ref{cor:globalKilling}, including the escaping Gram
condition and cyclic brackets. Their two normal components satisfy
Corollary~\ref{cor:globalfieldcriterion}. The latter gives the effective
$O(3)$ action and the specified principal orbits, without assuming
scalar decay or a critical point.
This finishes the proof.
\end{proof}

Next, we study metric blow-downs with a fixed base
point. The sector classification uses the angular action constructed
above. We first establish the ray
criterion directly from scalar decay and the geometry of distance
spheres; these two implications are independent of the action.

\begin{lem}
\label{lem:decay-ray-blowdown}
Let $(M^n,g,f)$ be a complete connected noncompact steady gradient
soliton with $\Ric=\Hess f>0$ and $\mathrm{R}+|\nabla f|^2=1$.
If $\mathrm{R}\to0$ uniformly at infinity, then for every fixed $q\in M$,
\[
 (M,r^{-2}g,q)\longrightarrow([0,\infty),0)
                  \qquad(r\to\infty)
\]
in the pointed Gromov--Hausdorff sense as $r\to\infty$ without passing to a subsequence.
\end{lem}
\begin{proof}
We first choose a compact starting level for the transport argument.
Scalar decay makes $\mathrm{R}$ attain a positive maximum at a point $p$.
The identity $\mathrm{d}\mathrm{R}=-2\Ric(\nabla f,\cdot)$ gives $\nabla f(p)=0$.
We normalize $f(p)=0$. Strict convexity makes $p$ the unique minimum.
We choose $r_0,\mu>0$ with $\Ric\ge\mu g$ on $B(p,r_0)$.
Along every unit-speed minimizing geodesic from $p$, we have
\[
 \frac{\mathrm d}{\mathrm dr}(f\circ\gamma)(r_0)\ge\mu r_0,
 \qquad \mu r_0\,d(p,x)-C\le f(x)\le d(p,x).
\]
Thus $f$ is proper. A positive level $\Sigma=\{f=s_0\}$ is compact.
Let $\Psi_t$ be the complete forward flow of $X:=\nabla f$ and put
$\Sigma_t:=\Psi_t\Sigma$. These images need not be potential levels.
On these trajectories $|X|\ge C^{-1}>0$, since
$X(|X|^2)=2\Ric(X,X)>0$. They escape uniformly, and hence, we have
\[
 f(\Psi_t y)=s_0+\int_0^t(1-\mathrm{R}(\Psi_s y))\,\mathrm{d}s=t+o(t),
 \qquad d(p,\Psi_t y)=t+o(t)
\]
uniformly for $y\in\Sigma$. The distance estimate follows from
$f\le d(p,\cdot)$ and the speed bound $|X|\le1$.

We next show that the transported levels have sublinear diameter.
For fixed $L>0$ set
\[
 \delta_L(t):=\sup_{d(z,\Sigma_t)\le L}\mathrm{R}(z)\longrightarrow0,
 \qquad D(t)=\operatorname{diam}_g\Sigma_t.
\]
For $y,z\in\Sigma$ write $\ell(t)=d(\Psi_t y,\Psi_t z)$.
If $\ell(t)\ge2L$, choose a minimizing geodesic between its endpoints
and use a cutoff equal to one except on the two endpoint segments
of length $L$, where it rises or falls linearly. The summed index
forms of $n-1$ parallel normal fields give
\[
 \int_0^{\ell(t)}\phi^2\Ric(\dot\gamma,\dot\gamma)\,\mathrm{d}u
       \le\frac{2(n-1)}L,
 \qquad
 \int_0^{\ell(t)}(1-\phi^2)\Ric(\dot\gamma,\dot\gamma)\,\mathrm{d}u
       \le2L\delta_L(t).
\]
Since $\partial_t\Psi_t^*g=2\Ric_{\Psi_t^*g}$, upper first variation
and the same geodesic as a support at cut points imply
\[
 D^+\ell(t)\le\frac{2(n-1)}L+2L\delta_L(t)
                         \qquad(\ell(t)\ge2L).
\]
We integrate this bound for $\max\{\ell(t),2L\}$ before taking the
supremum over pairs. The result is
\[
 D(t)\le\max\{D(0),2L\}+\frac{2(n-1)t}L
                              +2L\int_0^t\delta_L(s)\,\mathrm{d}s.
\]
First let $t\to\infty$ with $L$ fixed, then let $L\to\infty$.
This proves $D(t)=o(t)$ without differentiating a supremum.

These estimates allow us to construct the approximations to the limiting ray.
Every point outside the compact sublevel $\{f\le s_0\}$ is uniquely
$\Psi_t y$ with $t\ge0$ and $y\in\Sigma$. The backward trajectory
reaches $\Sigma$ in a compact band without critical points.
For $t\ge s$, we have
\[
 |f(\Psi_t y)-f(\Psi_s z)|
 \le d(\Psi_t y,\Psi_s z)\le t-s+D(s).
\]
We set $e(t):=\sup_\Sigma|f(\Psi_t y)-t|=o(t)$.
A point in $B(p,R_0r)$ has $t\le C^2R_0r$, since
$f(\Psi_t y)\ge s_0+C^{-2}t$. Therefore for exterior points in that ball,
\[
 0\le d(x,x')-|f(x)-f(x')|
 \le2\sup_{0\le t\le C^2R_0r}e(t)
                    +\sup_{0\le t\le C^2R_0r}D(t)=o(r).
\]
Points in the compact sublevel add only a bounded error. The map
$x\mapsto f(x)/r$ on $\overline B(p,R_0r)$ thus has distortion $o(1)$
in the metric $r^{-1}d_g$. Its image lies in $[0,R_0]$ and is
$o(1)$-dense there: on a minimizing ray, the point at distance $R_0r$
has $f=R_0r-o(r)$, and continuity supplies the intervening values.
This gives the asserted pointed approximations.
A fixed change of base point has rescaled distance $r^{-1}d(p,q)\to0$.
The assertion follows.
\end{proof}

For the blow-down argument, we must relate a ray at large scale to
curvature decay on the original soliton. The following lemma
supplies the required implication.

\begin{lem}
\label{lem:ray-direct-converse}
Let $(M^n,g,f)$ be a complete connected noncompact gradient steady
soliton with
\[
 \Ric=\Hess f\ge0,\qquad \mathrm{R}+|\nabla f|^2=1.
\]
Suppose, for a fixed $p\in M$, that for large $r$, we have
\begin{equation}\label{eq:ray-sphere-collapse}
 D_p(r):=\operatorname{diam}_g\partial B(p,r)=o(r).
\end{equation}
Then, after setting $f(p)=0$, uniformly as $x\to\infty$, we have
\[
 \frac{f(x)}{d(p,x)}\longrightarrow1,\qquad
 |\nabla f(x)|\longrightarrow1,\qquad \mathrm{R}(x)\longrightarrow0.
\]
In particular this holds if $(M,r^{-2}g,p)$ converges to a ray as
$r\to\infty$. 
\end{lem}
\begin{proof}
Ray convergence implies \eqref{eq:ray-sphere-collapse}: on rescaled
balls of radius two, pointed Gromov--Hausdorff approximations send
every point at distance one from the base point to within a
vanishing error of the coordinate one. Their distortion then
bounds the diameter of that sphere by a vanishing error.

We choose a unit-speed minimizing ray $\sigma$ from $p$. Convexity
and the gradient bound give
\[
 \frac{\mathrm d^2}{\mathrm dt^2}(f\circ\sigma)\ge0,\qquad |\frac{\mathrm d}{\mathrm dt}(f\circ\sigma)|\le1.
\]
Thus $\frac{\mathrm d}{\mathrm dt}(f\circ\sigma)$ has a limit $a\in[-1,1]$ and
$f(\sigma(r))/r\to a$. For $d(p,x)=r$, the Lipschitz bound yields
\[
 \left|\frac{f(x)}r-\frac{f(\sigma(r))}r\right|
 \le\frac{D_p(r)}r\longrightarrow0.
\]
Consequently $f(x)=a\,d(p,x)+o(d(p,x))$ uniformly at infinity.

Next, put $m:=\inf_M \mathrm{R}\ge0$ and $V(r):=\Vol_gB(p,r)$.
The divergence theorem and coarea, interpreted by radial Lipschitz
cutoffs at nonsmooth distance spheres, give for almost every $r$
\[
 mV(r)\le\int_{B(p,r)}\mathrm{R}\,\mathrm d\mu_g
 =\int_{\partial B(p,r)}\langle\nabla f,\nu\rangle\,\mathrm d\mu_{\partial B(p,r)}
 \le \frac{\mathrm d V}{\mathrm dr}(r)\le\frac nr V(r).
\]
The last inequality is Bishop--Gromov comparison using $\Ric\ge0$.
Taking arbitrarily large such radii proves $m=0$, and hence, we have
\[
                         \sup_M|\nabla f|=1.
\]
We fix $q$ with $b=|\nabla f(q)|>0$, and let $\eta$ be the complete
unit-speed geodesic starting at $q$ in the direction $\nabla f(q)/b$.
It need not remain minimizing. Convexity along this geodesic gives
\[
 f(\eta(t))\ge f(q)+bt,\qquad
 d(p,\eta(t))\le d(p,q)+t\quad(t\ge0).
\]
The first inequality and $|\nabla f|\le1$ imply
$d(p,\eta(t))\to\infty$. For all sufficiently large $t$ the
numerator below is positive, so the uniform slope limit implies
\[
 a=\lim_{t\to\infty}\frac{f(\eta(t))}{d(p,\eta(t))}
 \ge\lim_{t\to\infty}\frac{f(q)+bt}{d(p,q)+t}=b.
\]
Taking the supremum over $q$ gives $a\ge1$, and thus, we conclude that $a=1$.

Finally, if $\gamma:[0,r]\to M$ is a unit-speed minimizing geodesic
from $p$ to $x$, convexity gives
\[
 f(x)=\int_0^r\frac{\mathrm d}{\mathrm ds}(f\circ\gamma)(s)\,\mathrm{d}s
 \le r\frac{\mathrm d}{\mathrm dr}(f\circ\gamma)(r)\le r|\nabla f(x)|.
\]
It follows that $f(x)/d(p,x)\le|\nabla f(x)|\le1$. The uniform
slope limit proves $|\nabla f|\to1$, and the Hamilton identity
then proves uniform scalar decay.
This proves the stated conclusion.
\end{proof}

The orbit estimates now allow the large-scale geometry to be read
from the meridian. We identify the metric blow-down and distinguish
a ray from a sector of positive angle. Recall that for $\alpha\in (0,\pi)$, $\mathcal C_\alpha$ is the closed Euclidean
sector of opening angle $\alpha$, with its intrinsic metric and
vertex $o$; $\mathcal C_0=[0,\infty)$ denotes the ray.

\begin{thm}
\label{lem:metric-blowdown}
Let $(M^4,g,f)$ be a complete connected gradient steady Ricci soliton
with positive curvature operator, $\Ric=\Hess f$, and
$\kappa$-noncollapsing at all scales. There is a number
$\alpha\in[0,\pi)$ such that, for every fixed $p\in M$,
\[
 (M,\lambda g,p)\longrightarrow(\mathcal C_\alpha,o)
 \qquad\text{as }\lambda\downarrow0
\]
in the pointed Gromov--Hausdorff sense, along the full family of
scales. Moreover, we have
\[
 \alpha=0\quad\Longleftrightarrow\quad \mathrm{R}(x)\longrightarrow0
                         \text{ as }x\to\infty.
\]
 In the non-Bryant branch, let
$S:=M^{O(2)}$ be a meridian plane for the standard subgroup fixing
one coordinate of $\mathbb R^3$ in the equivariant model. Then we have
\begin{equation}
 \alpha=\pi-\frac12\int_S K_S\,\mathrm{d}\mu_S.
\end{equation}
\end{thm}

\begin{proof}
We first reduce the blow-down problem to the meridian surface.
Its ideal boundary determines the limiting cone, and the angular
involution identifies the quotient with the stated sector.

If the soliton is the four-dimensional Bryant soliton, its radial form
$\mathrm{d}s^2+a(s)^2g_{S^3}$ has $a(s)=O(\sqrt{s})$ at infinity. Projection
to $s$ therefore gives the ray blow-down, and $\mathrm{R}\to0$. We henceforth
work in the non-Bryant branch.

We first show that the rotational orbits collapse at the blow-down scale, leaving the meridian as the large-scale model.
We choose $o$ on the fixed axis of the action and put $r(x):=d(o,x)$.
The exact Killing generators $W_a$ constructed above obey the
stronger size bound
\begin{equation}\label{eq:orbit-sublinear}
 \sum_a|W_a(x)|\le C \mathrm{R}(x)^{-1/2}=o(r(x)).
\end{equation}
The first inequality is \eqref{eq:exact-generator-size}, with an
equivalent norm on triples; Proposition~\ref{prop:curvaturelength}
gives the final little-oh estimate.

Let $Q:=M/O(3)$ with its quotient length metric and let $\pi:M\to Q$
be the quotient map. In the standard action the $SO(3)$ and $O(3)$
orbits coincide, so the disconnected component creates no additional
orbit points. Integrating the generators along bounded-length paths
in $SO(3)$ bounds each orbit's diameter by
$C\sup_{O(3)x}\sum_a|W_a|$. Since $o$ is fixed, both $r$ and $\mathrm{R}$
are constant on an orbit. Consequently, for every fixed $R_0$, we have
\[
 \sup_{x\in B(o,R_0r)}\operatorname{diam}_g(O(3)x)=o(r).
\]
Indeed, outside a fixed compact set the orbit diameter is at most
$\epsilon d(o,x)$ by \eqref{eq:orbit-sublinear}; inside it the
orbit diameters are uniformly bounded. Divide by $r$, let
$r\to\infty$, and then let $\epsilon\downarrow0$.
Also, we have
\[
 d_Q(\pi x,\pi y)\le d_g(x,y)
 \le d_Q(\pi x,\pi y)+\operatorname{diam}_g(O(3)y).
\]
The second inequality follows by minimizing $d_g(x,ky)$ over the
compact group. Quotient balls centered at $\pi o$ are images of
balls centered at $o$. Thus $M$ and $Q$ have the same pointed
metric blow-downs.

We next identify the cone of the meridian surface and relate it to the ambient blow-down.
In the equivariant model $\mathbb R^3\oplus\mathbb R$, let
$H\cong O(2)$ fix the first coordinate in $\mathbb R^3$. Its fixed
set $S:=M^H$ is a smooth totally geodesic surface diffeomorphic to
$\mathbb R^2$. We give $S$ its intrinsic induced Riemannian distance
$d_S$, not the restriction of the ambient distance. It is complete:
a geodesic initially tangent to $S$ remains fixed by $H$ by
uniqueness of geodesics, and extends for all time in the complete
ambient manifold. The Gauss equation gives $K_S=\sec_M(TS)>0$.
The angular inversion restricts to an
isometric involution $\tau$ of $S$ which reverses its orientation,
fixes the axis, and exchanges the two meridian half-planes. The
invariant metric has no mixed angular-horizontal terms, as shown
in Proposition~\ref{prop:orbitcompletion}; projecting lengths to
the meridian coordinates therefore identifies
\[
                         Q=S/\langle\tau\rangle
\]
as length spaces. More explicitly, on the principal set write
$g:=g_Q+w^2g_{S^2}$. Dropping the angular velocity decreases the
length of a path. Conversely, a path in the meridian half-plane
lifts with constant angular coordinate and the same length;
paths meeting the axis follow by piecewise lifting or approximation.
The metric on $S$ consists of two copies of $g_Q$ joined smoothly
along the axis. Folding between these copies preserves lengths.
Thus for $x,y\in S$ the quotient distance is exactly
\begin{equation}\label{eq:meridian-quotient-distance}
 d_Q(\pi x,\pi y)=\min\{d_S(x,y),d_S(x,\tau y)\}.
\end{equation}
No convexity of $S$ in the ambient distance is being assumed.

For completeness, the surface facts used here are the
Cohn--Vossen inequality and the description of the ideal boundary
of a complete plane of finite total curvature
\cite[\S1, Theorem~1.5]{ShioyaBoundary}. We set
\[
 c(S)=\int_S K_S\,\mathrm{d}\mu_S,\qquad \ell:=2\pi-c(S).
\]
They give $0<c(S)\le2\pi$: the ideal boundary is a point when
$\ell=0$, and otherwise is a circle of circumference
$\ell\in(0,2\pi)$. The unique metric blow-down of a complete
nonnegatively curved manifold is the Euclidean cone over its
ideal boundary; see \cite[Theorem~0.1]{MNO}. Thus, the blow-down of
$S$ is respectively a ray or the flat cone of total angle $\ell$.
The two ideal-boundary descriptions agree here: \cite[Theorem A1]{ShioyaBoundary} identifies the limiting distance between
unit-speed rays with $2\sin(d_\infty/2)$, and the circle above has
diameter $\ell/2<\pi$. Thus its metric is precisely the angular
metric determining the cone, without any ambiguity from truncation
at $\pi$.

This description is compatible with $\tau$. Indeed, the circular
order of the surface ideal boundary is defined by the order of
rays around its single end. The oriented arc lengths in that
construction are additive \cite[\S1]{ShioyaBoundary}; an
orientation-reversing isometry reverses this order and preserves
these lengths. For $\ell>0$, its induced map is therefore a
reflection of that circle, not the identity: three distinct ideal
points retain a circular order, which $\tau$ reverses.
To justify convergence of the quotients, use the ray-net
approximations to the cone in \cite[\S2.2]{MNO}, based at the fixed
point $o\in S$. We adjoin the $\tau$-image of each ray in every finite
net. The resulting nets are invariant, still approximate each fixed
rescaled ball, and their pairwise distances converge together with
the distances to the images under $\tau$. Formula
\eqref{eq:meridian-quotient-distance}, which takes the minimum of
these two distances, therefore gives convergent quotient nets.
This proves convergence to the quotient by the induced reflection,
not merely ordinary convergence of $S$. The quotient of a circle
of circumference $\ell$ by reflection is an interval of length
$\ell/2$, whose Euclidean cone is the required sector. When
$\ell=0$, the cone and its quotient are both rays. This proves
existence and uniqueness of the asserted blow-down, with
$\alpha=\ell/2\in[0,\pi)$. A different fixed base point has
rescaled distance $\sqrt\lambda\,d(o,p)\to0$ from $o$ and gives
the same pointed limit. The angular reflection $\tau$ fixes the
axis pointwise.

Finally, we distinguish the ray case from a sector of positive angle using the curvature decay criterion.
If $\mathrm{R}\to0$ uniformly, Lemma~\ref{lem:decay-ray-blowdown}
gives a ray, so $\alpha=0$.
Conversely, if $\alpha=0$, Lemma~\ref{lem:ray-direct-converse}
applies and gives uniform scalar decay directly. This proves the
equivalence without assuming or first constructing a critical point.
The angle formula follows from $\alpha=\ell/2$ and the definition
$\ell=2\pi-c(S)$ above.
This finishes the proof.
\end{proof}

\begin{rem}
\label{rem:ray-criterion-generality}
Hebbar--\v Se\v sum \cite[Theorems~1.6--1.7]{HS} establish related
implications under their stated hypotheses. Their Theorem~1.6 assumes
positive Ricci curvature, a critical point of the potential, a
noncollapsed canonical flow with nonnegative sectional curvature,
and the non-Bryant alternative. Under uniform scalar decay, they proved
\[
 \mathrm R(x)d_g(x,\Gamma)\longrightarrow0
 \quad\text{as }d_g(o,x)\to\infty,
\]
and that the asymptotic cone is a ray, where $\Gamma$ is the union
of their two edge curves. Their Theorem~1.7 assumes positive Ricci
curvature, a critical point, and nonnegative sectional curvature,
and proves that a ray asymptotic cone implies uniform scalar decay.
Their converse can also
be used here once a critical point is supplied. One alternative is
to argue that, if $\nabla f$ has no zero, both ends of any complete
gradient trajectory escape compact sets: a compact accumulation
point with nonzero gradient would force a fixed change of the
monotone function $f$ on each visit to a small neighborhood.
On the forward end $|\nabla f|$ is bounded below by its positive
initial value $b$. Points on the two ends at the same large distance
$r$ from the starting point then have potential difference at least
$br$, contradicting the Lipschitz bound for $f$ and
\eqref{eq:ray-sphere-collapse}. Thus a ray limit supplies the critical
point required in \cite[Theorem~1.7]{HS}. The direct proof in
Lemma~\ref{lem:ray-direct-converse} avoids this intermediate step
and already works with $\Ric\ge0$.
\end{rem}

Next, we prove that the potential has a unique critical point and is
proper, and then establish the quantitative geometry of its large
levels. Throughout the critical point argument we use the normalization
$\Ric=\Hess f$ and $\mathrm{R}+|\nabla f|^2=1$. The four-dimensional
Bryant soliton already has a critical point, so we always assume that $(M,g)$ is not the Bryant soliton.
We write $X:=\nabla f$. The action is the standard action on
$\mathbb R^3\oplus\mathbb R$, by Proposition~\ref{prop:orbitcompletion}.
Its fixed axis is denoted by $\mathcal A$; on its complement write
\[
                   g=h+u^2g_{S^2},\qquad W=u^2.
\]
Here $h$ is the two-dimensional quotient metric, $K$ is its
Gaussian curvature, and $u$ is the orbit radius, determined by
$\Area(O(3)\cdot x)=4\pi u(x)^2$. The field $X$ descends to
the quotient; $Xu$ denotes the derivative of $u$ in that direction. Invariance of
$f$ follows either from the construction or from averaging: the difference
of two soliton potentials has parallel gradient, which positive Ricci
curvature excludes. The constants introduced below may depend on the soliton.
We use $c$ for the scalar value of the potential and write $c:=f$
when using it as a level coordinate. Thus
\[
 \Sigma_c:=\{x\in M:f(x)=c\},\qquad c_i:=f(x_i).
\]
The potential itself remains the function $f$ on $M$; expressions
such as $F(z,c)$ and $\partial_cF$ refer to profiles indexed by the
level value. In the finite-infimum case, we instead
use $c=f-m$, where $m=\inf_M f$, and state this normalization explicitly.

The warped-product equations give
\begin{align}
 \Hess_h f=K h-2u^{-1}\Hess_hu,\qquad
 u\Delta_hu+|\nabla_hu|_h^2+uXu&=1,\label{reflcrit:quotient}
 \end{align}
 and
 \begin{align}
 (\Delta_h+X)W=2,\qquad
 \Hess_hu&<0,\qquad |\nabla_hu|_h\le1.\nonumber
\end{align}
The last two inequalities follow from positive radial--orbital and orbital
sectional curvature. The exact generators satisfy
\eqref{eq:exact-generator-size}; their sum of squared lengths is a fixed
positive multiple of $u^2$. Consequently, we have
\begin{equation}\label{reflcrit:Ru}
                           \mathrm{R}u^2\le C.
\end{equation}
We also use the pointwise estimates from \cite[Corollary A.9]{CMZ}
\begin{equation}\label{reflcrit:derivatives}
                   |\nabla^j\Rm|\le C_j\mathrm{R}^{1+j/2}.
\end{equation}

\begin{lem}
\label{reflcrit:foliation}
Suppose $X:=\nabla f$ has no zero. Use
$\mathcal A$ for the fixed axis and $m=\inf_M f\in[-\infty,\infty)$.
Every potential level is a complete connected positively curved rotational
$\mathbb R^3$ with one tip on $\mathcal A$. Every $X$-trajectory crosses
every value in $(m,\infty)$ once. For each $c_0\in(m,\infty)$,
the flow of $Z:=X/|X|^2$ defines a diffeomorphism
\[
 \Theta:\{f=c_0\}\times(m,\infty)\longrightarrow M,\qquad
 \partial_c\Theta=Z\circ\Theta,\quad \Theta(x,c_0)=x,\quad
 f(\Theta(x,c))=c.
\] If $\gamma:\mathbb R\to\mathcal A$ is parametrized by unit
arclength and oriented so that $F=f\circ\gamma$ increases, then
\[
 F(\mathbb R)=(m,\infty),\quad
 0\le a_-:=\lim_{s\to-\infty}\frac{\mathrm d F}{\mathrm ds}(s)
 <a_+:=\lim_{s\to\infty}\frac{\mathrm d F}{\mathrm ds}(s)\le1.
\]
The induced metric on each level is
$\mathrm dr^2+u(r)^2g_{S^2}$, where $r\ge0$ is intrinsic distance
from its unique tip. In these coordinates, $u(0)=0$, ${\partial_{r} u}(0)=1$,
${\partial_{r} u}>0$, and ${\partial_{r}^2 u}<0$ for $r>0$.
\end{lem}
\begin{proof}
The fixed axis is a complete properly embedded geodesic parametrized by
$\mathbb R$, by the standard equivariant model and completeness. Its
intrinsic unit-speed parametrization is not asserted to minimize globally.
Its invariant gradient is tangent to
it. Since $\frac{\mathrm d^2 F}{\mathrm ds^2}=\Ric(\frac{\mathrm d\gamma}{\mathrm ds},\frac{\mathrm d\gamma}{\mathrm ds})>0$ and $\frac{\mathrm d F}{\mathrm ds}$ never vanishes, orient it
so that $\frac{\mathrm d F}{\mathrm ds}>0$. The stated derivative limits follow from $\frac{\mathrm d F}{\mathrm ds}\le1$.

Both $X$ and $-X$ are complete since $|X|\le1$. Along the forward flow,
$Xf=|X|^2$ and $X|X|^2=2\Ric(X,X)>0$, so $f\to\infty$. Along a backward
trajectory $x(t)$ suppose the limiting potential $\ell$ were greater than
$m$. We choose $y$ with $f(y)<\ell-\delta$. Convexity along a minimizing
segment from $x(t)$ to $y$ gives
\[
 \langle X(x(t)),\exp_{x(t)}^{-1}y\rangle\le f(y)-f(x(t))\le-\delta.
\]
The upper first variation inequality, also at cut points, then gives
$D^+[d(x(t),y)^2/2]\le-\delta$, an impossibility for all $t\ge0$.
The same argument excludes a finite lower limit when $m=-\infty$.
Thus every trajectory has potential range $(m,\infty)$. Reparametrization
by $f$ gives the asserted smooth product, generated by $X/|X|^2$.

A regular level is a closed embedded hypersurface, hence complete in its
induced metric. Its second fundamental form is $\Ric/|X|$, so it has
positive sectional curvature. An axis-free connected level component would
have a complete circle or line as its one-dimensional quotient space. Its
warping radius would be positive and strictly concave on that circle or
line, which is impossible. Thus every component meets the axis. Strict
increase of $F$ allows only one such point, so the level is connected,
has one tip, and is a rotational $\mathbb R^3$. This also proves
$F(\mathbb R)=f(M)$. Finally, ${\partial_{r}^2 u}<0$ and smoothness at the tip give the
last assertions: if ${\partial_{r} u}$ became nonpositive, strict decrease would force
$u$ to vanish again at finite distance.
The assertion follows.
\end{proof}

The following lemma supplies the equivariant convergence and the
persistence of axis points used in the critical point argument.

\begin{lem}\label{reflcrit:equivariant-models}
Let $W_1,W_2,W_3$ be the exact Killing fields constructed in
Corollary~\ref{cor:globalKilling}, with the standard effective
$SO(3)$ action of Proposition~\ref{prop:orbitcompletion} and the
size bound \eqref{eq:exact-generator-size}. Let $x_i\to\infty$,
$r_i:=\mathrm R(x_i)^{-1/2}$, and suppose the curvature-normalized
canonical flows
\[
 (M,r_i^{-2}g(r_i^2t),x_i)\longrightarrow
 (M_\infty,g_\infty(t),x_\infty)
\]
converge smoothly to a Bryant product or a shrinking
$S^2\times\mathbb R^2$ cylinder. After passing to a subsequence,
the generators and their $SO(3)$ actions converge smoothly on
compact sets. Up to one constant $SO(3)$ change of basis of the
limiting generators, the limiting action is the standard rotation
action on the Bryant or spherical factor and is trivial on the
Euclidean factor.
Every point fixed by the limiting action is the smooth-chart limit
of actual axis points. In particular, a Bryant tip-line point has
nearby actual axis points, whereas a sequence of actual axis points
cannot have a cylindrical pointed limit.
\end{lem}
\begin{proof}
Put $g_i(t):=r_i^{-2}g(r_i^2t)$. We keep the same vector fields
$W_a$, without multiplying them by $r_i$: constant metric scaling
does not change their brackets or their group periods. Since
$W_af=0$, they commute with $X$ and generate isometries of every
$g_i(t)$. The bound \eqref{eq:exact-generator-size} becomes
\[
 \left(\sum_a|W_a|_{g_i(t)}^2\right)^{1/2}
       \le C\mathrm R_{g_i(t)}^{-1/2}.
\]
Scalar curvature is positive on each compact subset of either
limiting model, so this gives uniform bounds there. The Killing
equation implies $\Delta W_a+\Ric^\sharp W_a=0$ and the Killing
transport identity expresses $\nabla^2W_a$ as a curvature
contraction with $W_a$. Interior estimates and smooth metric
convergence therefore bound every spatial derivative on smaller
compact sets. A diagonal subsequence gives smooth limiting Killing
fields with the same brackets.

We also obtain convergence of the entire compact group action.
Scalar curvature is constant on each orbit. For a point in a fixed
compact comparison region, integrate the displayed generator bound
along group paths of uniformly bounded length in $SO(3)$. Its whole
orbit then stays within a uniformly bounded $g_i(0)$-distance of
that region. A larger pointed comparison domain contains all these
orbits. Smooth dependence for the generator flows, using finitely
many exponential coordinate neighborhoods of the compact group,
gives smooth convergence of the action maps, uniformly over the
group, on compact spatial sets. The group identities pass to the
limit, so the limiting fields generate an $SO(3)$ action on the
whole model. The same maps preserve the limiting canonical flow.

For later use, we prove persistence of fixed points, including when
the limiting action might be trivial. If $q_\infty$ is fixed and
$\chi_i$ are smooth convergence charts, the entire orbit of
$\chi_i(q_\infty)$ has diameter tending to zero. In a uniformly
small convex normal neighborhood its Haar-averaged squared-distance
function has a unique center $q_i$, tending to $q_\infty$ in these
charts. Invariance of the orbit and uniqueness of this center imply
that $q_i$ is fixed. More explicitly, both the center and each of
its group translates lie in the same convex neighborhood for large
$i$, and are critical points of the same strictly convex function.
Thus $q_i$ lies on the actual fixed axis $\mathcal A$.

The limiting action cannot be trivial. Otherwise its generators and
their first derivatives would vanish in the limit, and the preceding
center construction at $x_\infty$ would give axis points $q_i$ at
normalized distance tending to zero from $x_i$. At every actual
axis point, Proposition~\ref{prop:orbitcompletion} identifies the
differential representation with the adjoint representation on the
normal three-plane and the trivial representation on the axis.
For a bracket-one basis this gives, in the Hilbert--Schmidt norm,
\[
                  \sum_{a=1}^3|\nabla W_a(q_i)|_{g_i(0)}^2=6.
\]
The identity is unchanged by constant metric scaling or orthogonal
conjugation of the representation. It contradicts convergence of
these derivatives to zero.

Finally, an isometry of the Bryant product preserves its parallel
line distribution and its orthogonal Bryant factor. A connected
compact group acts trivially on the line. On the round cylinder,
the corresponding splitting is $S^2\times\mathbb R^2$; a compact
Euclidean action has a fixed point, and every homomorphism from
$SO(3)$ to $O(2)$ is trivial. Thus the limiting action is trivial
on the Euclidean factor in both cases. The connected isometry
group of the remaining Bryant or round spherical factor is
$SO(3)$. Its nontrivial homomorphism from $SO(3)$ is an automorphism,
and every such automorphism is inner. This proves the assertion
about the standard generators and the constant change of basis.
The Bryant fixed set is its tip line; the cylindrical action has no
fixed points. Persistence proves the first assertion about axis
points. Conversely, if $x_i$ are fixed, convergence of the actions
makes $x_\infty$ fixed, excluding the cylinder.
\end{proof}

\begin{lem}
\label{reflcrit:decay}
In the setting of Lemma~\ref{reflcrit:foliation}, set $L_+:=1-a_+^2$.
Assume $L_+=0$. For $c$ in the range of $f$, let
$p_c:=\gamma(F^{-1}(c))$ be the unique axis point on $\{f=c\}$.
Then, as $c\to\infty$, we have
\[
 \sup_{\{f\ge c\}}\mathrm{R}\longrightarrow0,\qquad
                 c\mathrm{R}(p_c)\longrightarrow\infty.
\]
\end{lem}
\begin{proof}
If the first assertion failed, take $x_i$ with $f(x_i)\to\infty$ and
$\mathrm{R}(x_i)\ge\epsilon>0$. The scalar-scale limit inherits a steady potential:
the centered potentials have uniformly bounded gradients and Hessian equal
to Ricci. A round $S^2\times\mathbb R^2$ cannot inherit this equation,
since restriction and tracing on its compact sphere would give a function
with strictly positive Laplacian everywhere. Thus, the limit is a Bryant
product. Lemma~\ref{reflcrit:equivariant-models} identifies the
limiting rotations with the standard Bryant rotations and supplies
actual axis points $y_i$ near a point on its tip line, with
\[
 d(x_i,y_i)\le C \mathrm{R}(x_i)^{-1/2},\qquad \mathrm{R}(y_i)\ge C^{-1}\mathrm{R}(x_i).
\]
Since $|\nabla f|\le1$, $f(y_i)\to\infty$. These are outgoing axis points,
contradicting $L_+=0$.

At an outgoing axis point every scalar-scale limit is based on the Bryant
tip line, by Lemma~\ref{reflcrit:equivariant-models}; its standard
cylindrical action has no fixed point. The time derivative of scalar curvature at the
Bryant tip is zero. Smooth convergence of the canonical flows gives
\[
 \frac{2\frac{\mathrm d F}{\mathrm ds}(s)^2\frac{\mathrm d^2 F}{\mathrm ds^2}(s)}{\mathrm{R}(\gamma(s))^2}\longrightarrow0.
\]
As $\bigl(\frac{\mathrm d}{\mathrm d s}F\bigr)^2\to1$, it follows that
\[
 \frac{\mathrm{d}}{\mathrm{d}c}\frac1{\mathrm{R}(\gamma(F^{-1}(c)))}
                   =\frac{2\frac{\mathrm d^2 F}{\mathrm ds^2}}{\mathrm{R}^2}\longrightarrow0.
\]
Integration proves $1/\mathrm{R}=o(c)$, as claimed.
This proves the stated conclusion.
\end{proof}

We begin the critical point argument by controlling the maximum
orbit radius on each potential level. This gives an envelope
for the subsequent barrier comparisons.

\begin{lem}
\label{reflcrit:envelope}
Assume that the potential function has no critical point.
If $m> -\infty$, put $c:=f-m$. Then every level has bounded orbit radius, and
\[
 U(c):=\sup_{\{f=m+c\}}u=\omega(B^{-1}(c)),
 \qquad |U(c)^2-2c|\le C\log(2+c),
\]
where $\mathrm{d}r^2+\omega(r)^2g_{S^2}$ is the Bryant metric and
$B(r)$ is its tip-zero potential. Fix $c>0$, let $p(c)$ be the
unique tip of $\Sigma_{m+c}:=\{f=m+c\}$, and let
\[
 x_i\in\Sigma_{m+c},\qquad
 r_i:=d_{\Sigma_{m+c}}(p(c),x_i)\longrightarrow\infty.
\]
Every such sequence has a subsequence for which
\begin{equation}
 (M,g,f-f(x_i),x_i)\longrightarrow
 (B^3\times\mathbb R,g_B+\mathrm{d}t^2,B(s)-c,(q_c,0)),
 \qquad d_{g_B}(o_B,q_c)=B^{-1}(c).
\end{equation}
Here $s=d_{g_B}(o_B,\cdot)$, $o_B$ is the Bryant tip. In particular, with $r$ the intrinsic level distance from its tip,
\begin{equation}\label{eq:finite-infimum-end}
 \lim_{r\to\infty}u(r,c)=U(c),\qquad
 \lim_{r\to\infty}|\nabla_hu|_h^2(r,c)=\Phi_B(U(c)).
\end{equation}
We call the pointed product convergence above and
\eqref{eq:finite-infimum-end} the \textbf{finite-infimum end description}.
These limits are sequentially uniform for $c$ in compact
subintervals of $(0,\infty)$; here the regular Bryant squared slope profile is defined by
$\Phi_B(\omega(r)):=\bigl(\frac{\mathrm d\omega}{\mathrm dr}(r)\bigr)^2$.
Conversely, the existence of a bounded-radius level forces $m> -\infty$.
\end{lem}

\begin{proof}
Suppose $m$ is finite and a level has points $x_i$ with $u(x_i)\to\infty$.
Let $y_i:=\gamma(s_i)$ be nearest axis points and $d_i:=d(x_i,y_i)$.
Since $u\le d_{\mathcal A}$, $d_i\to\infty$. The normal minimizing chords
have zero initial potential derivative and convex potential. Their endpoint
potential is fixed, so $F(s_i)$ is bounded above. If $s_i$ stayed in a
compact interval, a convergent subsequence of the chords would give a
constant-potential minimizing ray, contradicting $\Ric>0$. Thus
$s_i\to-\infty$. Finiteness of $m$ gives $\frac{\mathrm d F}{\mathrm ds}\to0$ there. The function
$\frac{\mathrm d^2 F}{\mathrm ds^2}>0$ is integrable and uniformly continuous by
Lemma~\ref{lem:boundedlocalgeometry}, hence tends to zero. Axis isotropy
therefore makes each transverse Ricci eigenvalue tend to $1/3$. A fixed
short initial piece of every late normal chord has Ricci curvature at least
$b>0$. Integrating $\Hess f=\Ric$ on that piece and using convexity on the
rest gives $f(x_i)-f(y_i)\ge bD(d_i-D)$ for fixed $b,D>0$. This contradicts
$f(x_i)-f(y_i)\le f(x_i)-m$. Thus every radius is bounded on each level.

On any bounded-radius level, $u\uparrow U>0$, ${\partial_{r} u}\downarrow0$, and the
intrinsic orbital sectional curvature tends to $1/U^2$. If ambient $\mathrm{R}$
tended to zero along a subsequence, $|X|\to1$ and $\mathrm{II}=\Ric/|X|$
would tend to zero. Gauss would contradict the positive limiting orbital
curvature. Thus, the scalar curvature has a positive lower bound at the end.
The dimension reduction in \cite[Theorems 1.2--1.3]{CMZ} gives a Bryant product at a nonvanishing scalar scale, and
therefore also an unscaled smooth limit with potential $B+\alpha t+b$.
The outward unit level tangents have $\mathrm{d}u(e_r)\to0$ and $\mathrm{d}f(e_r)=0$.
In the limit the radius has strictly positive Bryant radial derivative;
thus these tangents become $\pm\partial_t$, which implies $\alpha=0$.
The Hamilton normalization then fixes the Bryant tip scalar to one.

By Lemma~\ref{reflcrit:equivariant-models}, persist a point on the
limiting tip line to actual axis points $y_i$.
Their distance to the original level points is bounded, and their potentials
tend to $f(x_i)-B(r_*)$. They escape along the incoming axis, since the
outgoing axial potential tends to infinity. Hence $m$ is finite and
$m=f(x_i)-B(r_*)$. This proves both the converse and the precise formula
for $U(c)$, as well as the centered potential statement.

For the remainder estimate, differentiate in the Bryant
distance $r$; put $p:=\frac{\mathrm d B}{\mathrm dr}$, $z:=\omega\frac{\mathrm d \omega}{\mathrm dr}$, $e:=1-p$, $d:=z-1$.
The Bryant equations \eqref{eq:BryantODE} give
\[
 \frac{\mathrm d z}{\mathrm dr}=1-pz,\qquad \frac{\mathrm d p}{\mathrm dr}=1-p^2-2pz/\omega^2,
 \qquad \frac{\mathrm d e}{\mathrm dr}=-(1+p)e+2pz/\omega^2,\qquad \frac{\mathrm d d}{\mathrm dr}+pd=e.
\]
Since $p,z\to1$ and $\omega^2\asymp r$, variation of constants gives
$e=O(r^{-1})$ and $d=O(r^{-1})$. For example split each integral of the inhomogeneous term
at $r/2$; its earlier portion is exponentially small and its later portion
is $O(r^{-1})$. Thus $\frac{\mathrm d}{\mathrm dr}(\omega^2-2B)=2(z-p)=O(r^{-1})$, proving the logarithmic
remainder after integration.

On any fixed positive level both $\mathrm{R}$ and $|X|$ have positive infima, by
the identified Bryant end and compactness of the remainder. Along normal
trajectories $\mathrm{R}$ decreases and $|X|$ increases. Hence a positive level gives
a positive lower bound for the lapse on all higher levels and for
scalar curvature on all lower levels.
If $c_i\to c_\infty>0$, transport points from level $c_i$ to $c_\infty$; the transport
length is $o(1)$ by this lower bound for the lapse. Applying the fixed-level conclusion
proves all the sequentially uniform assertions.
This finishes the proof.
\end{proof}

The next estimate identifies a region where the orbit radius
remains large enough to apply the local geometric models.

\begin{lem}
\label{reflcrit:thick}
Assume that the potential function has no critical point and we use notation of Lemma~\ref{reflcrit:foliation}.
Let $L_+:=1-a_+^2$, with $a_+$ defined in
Lemma~\ref{reflcrit:foliation}, and let $\Ric_{\rm orb}$ denote the
ambient Ricci eigenvalue in a unit orbit direction. Suppose either
$m$ is finite, or $L_+=0$ and $m=-\infty$. Use $c:=f$ as the level coordinate; only sufficiently large
positive values of $c$ are used. For every fixed $d>0$ and every sequence $x_i$ with
$c_i:=f(x_i)\to\infty$ and $u(x_i)\ge d\sqrt{c_i}$, every
subsequential pointed limit of the canonical flow rescaled by
$\mathrm R(x_i)$ is the shrinking $S^2\times\mathbb R^2$
normalized to scalar curvature one at its base point. Uniformly at points $x$ with $c=f(x)\ge c_0$ and
$u(x)\ge d\sqrt c$, as $c_0\to\infty$, we have
\begin{align}
\mathrm{R}u^2\to2,\quad \Ric_{\rm orb}/\mathrm{R}\to\tfrac12,
 \quad XW\to2,\quad |\nabla_hu|_h^2\to0.	
\end{align}

The scalar-normalized levels converge smoothly to $S^2\times\mathbb R$.
If the levels have unbounded radius, there is one fixed normal trajectory
with $u^2/c\to2$ and $c\mathrm{R}\to1$.
For $m> -\infty$, the same asymptotic conclusions hold with
$c=f-m$, since $(f-m)/f\to1$; thus they also apply to the
finite-infimum level parameter used below.
\end{lem}
\begin{proof}
We first suppose that $m$ is finite. The positive lower bounds for scalar curvature and the lapse in
the preceding proof, together with $(\Delta+X)\mathrm{R}=-2|\Ric|^2$ and
\eqref{reflcrit:derivatives}, imply
\[
 0\le-X\mathrm{R}\le C\mathrm{R}^2,\qquad \mathrm{R}\ge b/(1+f-m).
\]
Indeed, $\mathrm{d}(1/\mathrm{R})/\mathrm{d}c\le C/|X|^2$ on a high superlevel, and integration starts
from a level with a positive scalar infimum. On the thick region,
\eqref{reflcrit:Ru} consequently gives
$b\le f\mathrm{R}\le C_d$ and $b_d\le \mathrm{R}u^2\le C$. A hypothetical Bryant-product
limit has a tip line point at bounded scalar distance; by
Lemma~\ref{reflcrit:equivariant-models}, it gives actual axis
points with $|f(y_i)-f(x_i)|=O(\sqrt{c_i})$ and $\mathrm{R}(y_i)\asymp \mathrm{R}(x_i)$.
They are outgoing, with bounded $f(y_i)\mathrm{R}(y_i)$. This is impossible:
if $L_+>0$ their product diverges trivially, and if $L_+=0$ it diverges by
Lemma~\ref{reflcrit:decay}.

Now suppose $m=-\infty$ and $L_+=0$. Every level has unbounded radius by
Lemma~\ref{reflcrit:envelope}, and uniform high-potential decay is available.
For a thick sequence put $\lambda_i:=\mathrm{R}(x_i)$ and pass so that
$\lambda_i u_i^2\to\ell\in[0,C]$. If $\ell>0$ and the limit is Bryant,
the normalized potentials $\sqrt{\lambda_i}(f-f(x_i))$ converge to a unit
parallel function. Thus, the levels converge to Bryant slices, and their
radial slopes at the basepoints tend to a positive number $v$.
Concavity gives, for the actual tip distance $D_i$ on the same level,
$D_i {\partial_{r} u}(x_i)\le u_i$. Hence $D_i\sqrt{\lambda_i}$ is bounded, so $\mathrm{R}_{\rm tip}\le C\lambda_i$. This contradicts $c_i\mathrm{R}_{\rm tip}\to\infty$
and $c_i\lambda_i\le C_d$. If $\ell=0$, the identities
\[
 0\le1-|\nabla_hu|_h^2\le \mathrm{R}u^2/2,\qquad
 \bigl(\partial_{r} u\bigr)^2=|\nabla_hu|_h^2-(Xu)^2/|X|^2,\qquad 0<Xu\le u\mathrm{R}
\]
give ${\partial_{r} u}(x_i)\to1$ and $D_i\sqrt{\lambda_i}\to0$.
Since \eqref{reflcrit:derivatives} gives $|\nabla (\mathrm{R}^{-1/2})|\le C$, this implies
$\mathrm{R}_{\rm tip}/\lambda_i\to1$, again contradicting the tip divergence.

Only cylindrical limits remain. By
Lemma~\ref{reflcrit:equivariant-models}, the surviving rotations are
standard on the sphere and trivial on the Euclidean factors, so their
radius is $\sqrt2$. This proves all the stated limits and convergence of
levels. Finally, start a normal trajectory on $\Sigma_{c_0}$ at $u^2=2c_0$. For fixed
$d^2<2<D^2$, the inequalities $u^2=d^2c$ and $u^2=D^2c$ cannot be crossed
outward at a first exit, since $\mathrm{d}W/\mathrm{d}c=XW/(1-\mathrm{R})\to2$ on either boundary.
The trajectory remains thick, and integration gives $u^2/c\to2$ and
$c\mathrm{R}\to1$.
The assertion follows.
\end{proof}

In the positive-end-limit case, we need a global bound for the
gradient of the squared orbit radius. We obtain it from the level geometry
before using the heat equation comparison.

\begin{lem}\label{lem:decorated-sector}
Suppose $|\nabla W|\le C$ on $M$.
Suppose the quotient blow-down is a sector of positive angle and its
two boundary rays are the limits of the axis ends, with unit
arclength speed. Orient axis arclength increasingly and put
\[
 a_-:=\lim_{s\to-\infty}\frac{\mathrm d}{\mathrm ds}(f\circ\gamma)(s),\qquad
 a_+:=\lim_{s\to\infty}\frac{\mathrm d}{\mathrm ds}(f\circ\gamma)(s).
\]
Assume $-1<a_-<a_+<1$. Then there are Euclidean coordinates on
the sector for which the vertex is $(0,0)$ and the limiting potential
has unit gradient $\partial_x$. Choose $y$ so that the incoming boundary
ray lies in $y>0$. Fix a base point $o$ on the axis, also denoting its
quotient image by $o$, and let $T\to\infty$ be the distance rescaling
parameter, so $h_T:=T^{-2}h$. Then, locally uniformly up to the boundary, as $T\to\infty$, we have
\[
 \frac{f-f(o)}T\longrightarrow x,\qquad
 \frac WT\longrightarrow
 \begin{cases}2x+s_-y,&y\ge0,\\2x+s_+y,&y\le0,\end{cases}
 \qquad\mathrm{where}\qquad s_\pm=\frac{2a_\pm}{\sqrt{1-a_\pm^2}}.
\]
Potential gradients converge strongly in local $L^2$ in the interior,
in the measured Gromov--Hausdorff sense. Explicitly, put
$\varphi_T=(f-f(o))/T$. For each relatively compact interior domain
$K$ with boundary of measure zero and its corresponding domains $K_T$
in the rescaled quotients, the weak convergence
$\nabla_{h_T}\varphi_T\rightharpoonup\partial_x$ is accompanied by
\[
 \int_{K_T}|\nabla_{h_T}\varphi_T|_{h_T}^2\,\mathrm d\mu_{h_T}
 \longrightarrow\int_K|\partial_x|^2\,\mathrm dx\,\mathrm dy
 =\operatorname{Area}(K).
\]
Moreover $XW\to2$ uniformly as $u\to\infty$.
\end{lem}
\begin{proof}
Rescale the quotient and retain the potential and squared radius.
Set $h_T:=T^{-2}h$, $\varphi_T:=(f-f(o))/T$,
$W_T:=W/T$, and $\rho_T:=u/\sqrt T$. The assumed gradient bound makes $W_T$
uniformly Lipschitz, with zero boundary trace, and makes $\rho_T$ concave
and locally Lipschitz in the sector interior. Any limiting nonnegative
concave $\rho$ is either strictly positive or identically zero. The latter
is impossible: local concavity would give $\nabla W_T\to0$, whereas testing
\[
 T^{-1}\Delta_{h_T}W_T+\langle \nabla_{h_T}\varphi_T,\nabla_{h_T}W_T\rangle=2
\]
on an interior ball gives a contradiction. Thus $W_T$ has a positive
lower bound on every compact interior set.

We set $b:=\log u$. The horizontal equation reads
$\Hess_h(f+2b)=Kh-2\mathrm{d}b\otimes \mathrm{d}b$. Its trace gives
\[
 \Delta_{h_T}\varphi_T=\frac2T
           (K_{h_T}-\Delta_{h_T}b-|\nabla_{h_T}b|_{h_T}^2)\ge0.
\]
Fix an interior ball $\mathcal B_T$ in the rescaled quotient and a
cutoff $\chi\in C_c^\infty(\mathcal B_T^+)$ equal to one on $\mathcal B_T$,
where $\mathcal B_T^+$ is a fixed larger interior ball.
Interior bounds for $\mathrm{d}b$, this cutoff, and finite total curvature imply
$\int_{\mathcal B_T^+}\chi\Delta_{h_T}\varphi_T\,\mathrm d\mu_{h_T}\le C/T$. Harmonic replacement on an interior
ball then has difference of Dirichlet integrals $O(T^{-1})$: if $\eta_T$ is harmonic with
boundary value $\varphi_T$, then we obtain
\[
 \int_{\mathcal B_T}|\nabla_{h_T}(\varphi_T-\eta_T)|_{h_T}^2\,\mathrm d\mu_{h_T}
 =-\int_{\mathcal B_T}(\varphi_T-\eta_T)\Delta_{h_T}\varphi_T\,\mathrm d\mu_{h_T}\le C/T.
\]
The uniform local Dirichlet Poincar\'e inequality and harmonic convergence
on noncollapsed $\mathrm{RCD}(0,2)$ surfaces give a harmonic limit and strong local
gradient convergence; the precise replacement and convergence results are
\cite[\S4, Lemma 4.7, Theorem 4.4 and Corollary 4.5]{AmbrosioHonda}.
The limit is also convex, hence affine. Since $\mathrm{R}W\le C$ and $W_T$ is
bounded below, $|\nabla \varphi_T|^2=1-\mathrm{R}\to1$; write the affine limit as $x$.
The same Dirichlet integral argument and integration by parts pass the radius equation
to $\partial_xW_\infty=2$, with zero traces on both sector edges.

The axis directions $e_\pm$ have scalar products $-a_-$ and $a_+$ with
$\partial_x$. This leaves two possible placements. If $\partial_x$ lay
outside the sector, some parallel segment would join one boundary edge to
the other. Integrating $\partial_xW_\infty=2$ on these segments contradicts
the zero traces. Thus choose the interior placement
\[
 e_-=(-a_-,\sqrt{1-a_-^2}),\qquad
 e_+=(a_+,-\sqrt{1-a_+^2}).
\]
Integration from each edge gives the stated tent. Its two slopes differ
because $a/\sqrt{1-a^2}$ is strictly increasing.

Finally, if $u(x_i)\to\infty$, the bound $\mathrm RW\le C$ gives
$\mathrm R(x_i)\to0$. The estimate
$|\nabla(\mathrm R^{-1/2})|\le C$ makes this decay uniform on
every fixed ball, and the curvature derivative estimates give an
unscaled flat pointed limit. Noncollapsing
on all scales makes it $\mathbb R^4$: a flat quotient with nontrivial
infinite deck group has smaller volume growth order, and a nontrivial
finite Euclidean group cannot act freely. The isotropy-fixed quotient
planes converge to $\mathbb R^2$. The centered $W$ have bounded gradients
and solve $\Delta_hW+XW=2$, while centered potentials converge to a unit
linear function. Elliptic interior estimates give a limit $w$ on
$\mathbb R^2$ satisfying $\Delta w+\partial_xw=2$, $|\nabla w|\le C$.
Each component of $dw$ is a bounded drift-harmonic function. Translating
it by $-t\partial_x$ gives a bounded ancient heat solution; the heat kernel
gradient bound at arbitrarily long times makes it constant. Thus $w$ is
affine and $\partial_xw=2$. Arbitrary-sequence compactness proves $XW\to2$.
\end{proof}

To apply Lemma~\ref{lem:decorated-sector}, we next derive its
radius-gradient bound from the hypothetical positive outgoing scalar
limit. This connects the axis geometry with the sector profile.
\begin{lem}
\label{reflcrit:tent}
Assume the setting of Lemma~\ref{reflcrit:foliation}, with
$L_+:=1-a_+^2>0$. Fix $o\in M$. Then $|\nabla W|\le C$ on $M$.
As $T\to\infty$, on the quotient sector limit of
$(M/O(3),T^{-2}h)$ there are coordinates $(x,y)$ such that
\[
 \frac{f-f(o)}T\longrightarrow x,\qquad
 \frac WT\longrightarrow
 \begin{cases}2x+s_-y,&y\ge0,\\2x+s_+y,&y\le0,\end{cases}
 \qquad\mathrm{where}\qquad s_\pm:=\frac{2a_\pm}{\sqrt{1-a_\pm^2}},\quad s_-<s_+.
\]
Convergence is locally uniform up to the sector boundary; potential
gradients converge strongly in local $L^2$ in the interior.
Moreover, $XW\to2$ uniformly as $u\to\infty$.
\end{lem}
\begin{proof}
We first compare axis arclength with ambient distance. We then bound
$|\nabla W|$ and identify the limiting potential and radius functions
on the quotient sector.

We first determine how the axis lengths behave at large scale.
The complete meridian plane $S$ in
Theorem~\ref{lem:metric-blowdown} has nonnegative curvature and finite total
curvature. The quotient cone cannot be a ray: by
Lemma~\ref{lem:ray-direct-converse}, a ray blow-down would imply uniform
scalar decay, contradicting the outgoing axis limit $L_+>0$.
Thus, the plane cone has positive angle. This argument precedes, and does
not assume, any comparison of axis arclength with ambient distance. A unit-speed axis geodesic passes to a
unit-speed quasigeodesic under this noncollapsed two-dimensional
convergence; this is \cite[\S5.1, properties \textup{(1)} and \textup{(6)}]{PetruninQC}.
In the punctured flat cone it is a straight geodesic, so both halves run
radially at unit speed. Thus $d_S(o,\gamma(s))/|s|\to1$.

We need the uniform version on shorter chords as well. If
$s_i\to\pm\infty$ and $r_i^\pm\in[(1-\epsilon)s_i,(1+\epsilon)s_i]$
(with the interval endpoints reordered when necessary), put
$\ell_i:=|r_i^+-r_i^-|$. When $\ell_i/|s_i|$ stays positive the cone limit
already gives $d_S(\gamma(r_i^-),\gamma(r_i^+))/\ell_i\to1$.
Otherwise center at their axis midpoint. A ball of radius $\rho|s_i|$,
with fixed small $\rho$, converges to a Euclidean ball away from the cone
vertex. Its area ratio tends to one by Hausdorff measure convergence under
noncollapsed Alexandrov convergence \cite[Theorem 10.8]{BGP}.
Bishop--Gromov makes the area ratios on all smaller balls, including every
fixed multiple of $\ell_i$, tend to one. Their rescaled limits are Euclidean
by volume rigidity. Quasigeodesic stability now makes the rescaled axis a
straight unit-speed line, proving the same chord conclusion. Distances
between axis points in $M$ and $S$ agree: rotate an ambient minimizing
segment's initial normal direction into $S$; the totally geodesic fixed
surface contains the rotated segment. Consequently this chord statement
also holds in the quotient and ambient distances.

We next extend the axis estimates across the collars to obtain the required global bound.
At either axis end $\frac{\mathrm d^2 F}{\mathrm ds^2}\to0$, since
it is integrable and uniformly continuous. Because
$L_-:=1-a_-^2>L_+>0$, isotropy gives a uniformly positive transverse Ricci
curvature at both late ends. On a normal minimizing segment from the axis,
therefore, we conclude that $\mathrm df(n)\ge b>0$ after one fixed length $D$.
Let $x$ have closest axis foot $\gamma(s)$ and $r:=d_{\mathcal A}(x)$,
with $r\le\epsilon_0|s|$. We choose $\ell:=r/\epsilon$ and the two comparison
points $\gamma(s\pm\ell)$. The uniform chord result and squared distance
Hessian comparison give
\[
 d(\gamma(s),\gamma(s\pm\ell))=(1+o(1))\ell,
 \quad d(\gamma(s-\ell),\gamma(s+\ell))=(2+o(1))\ell,
 \quad d(x,\gamma(s\pm\ell))^2\le r^2+\ell^2.
\]
We choose lateness so the relative chord errors are at most $\epsilon^2$.
The cosine law and Toponogov show that the two directions from $x$ to
these axis points are $\pm t+O(\epsilon)$, where $t\perp n$, and that
their inverse lengths are $O(\epsilon/r)$. Concavity of $u$ and its zero
values at their endpoints then imply
\begin{equation}\label{reflcrit:collar-tangent}
                   |\mathrm du(t)|\le C\epsilon u/r+C\epsilon|\mathrm du(n)|.
\end{equation}
These minimizing segments are horizontal by conservation of angular
Killing momenta. They cannot meet the axis before their endpoint: an angular
Jacobi field vanishing at two such contacts, extended by zero on the
remaining part, would have zero index form on a minimizing segment and a
nonzero derivative jump, contradicting regularity of the index form kernel.
Thus quotient concavity applies throughout them. The same reasoning applies
at each intermediate point of the original normal segment.

From \eqref{reflcrit:quotient} and \eqref{reflcrit:Ru},
\[
 0\le uXu=1-|\nabla u|^2-u\Delta_hu\le1+\mathrm{R}u^2/4\le C.
\]
For $\mathrm du(n)\ge0$, combine this with $\mathrm df(n)\ge b$ and
\eqref{reflcrit:collar-tangent}; for small fixed $\epsilon$ it gives
$\frac{\mathrm d}{\mathrm dr}(u^2)\le C+\alpha u^2/r$, where $\alpha<1$. When $\mathrm du(n)<0$ the inequality
holds automatically. Integrating from $D$ and using $u\le r$ before $D$
gives $u^2\le Cr$ in both proportional collars.

We next show that the estimate on one outgoing collar gives $u(x)^2\le C(1+d(o,x))$
everywhere. We set $T:=1+d(o,x)$ and join $\gamma(T)$ to $x$ by a minimizing
segment of length at most $2T$. A point $y$ at a small fixed fraction
$\delta T$ of this segment has a closest axis foot with parameter in
$[T/2,2T]$: otherwise the potential difference from $\gamma(T)$, using
$\frac{\mathrm d F}{\mathrm ds}\to a_+>0$, exceeds its distance, at most $2\delta T$.
Choosing $\delta$ small puts $y$ in the collar, so $u(y)^2\le C\delta T$.
Concavity from the zero initial radius gives
$u(x)\le2u(y)/\delta$, proving the coarse bound. If the segment is shorter
than $\delta T$, the endpoint itself is in the collar. Now if a closest
axis foot of $x$ is $\gamma(s)$, either $d_{\mathcal A}(x)$ is a fixed
fraction of $1+|s|$, when the coarse bound applies, or one of the two
collars applies. Thus $u^2\le C d_{\mathcal A}$ globally, after using
$u\le d_{\mathcal A}$ for distances less than one. Concavity along quotient
segments of length $d_{\mathcal A}/2$ in both directions gives
$|\nabla u|\le2u/d_{\mathcal A}$. Hence $|\nabla W|\le C$.

The axis comparison above gives unit-speed boundary rays, and
$0\le a_-<a_+<1$. Lemma~\ref{lem:decorated-sector} now proves
the stated decorated convergence and the assertion about $XW$.
This proves the stated conclusion.
\end{proof}

We isolate the parabolic comparison that will exclude the
hypothetical no-critical point geometry. Its hypotheses are
verified from the geometric estimates below.

\begin{lem}
\label{reflcrit:heat-ridge}
Let $(M^4,g,f)$ be a complete steady soliton with
$\Ric=\Hess f\ge0$, $\mathrm R+|\nabla f|^2=1$, and a
potential-preserving isometric $O(3)$ action with $S^2$ principal
orbits. Write $g=h+u^2g_{S^2}$, $W:=u^2$, $X:=\nabla f$, and
$c:=f$. Assume the rotational potential levels have positive
intrinsic radial curvature: if
$\gamma_c=\mathrm dr^2+u(r,c)^2g_{S^2}$ and $e$ is a unit orbit
tangent vector, then this means
\[
 K_{\gamma_c}(\partial_r,e)=-\frac{\partial_r^2u}{u}>0
 \quad\text{where }u>0.
\]
Assume also that the quotient $\mathcal B=M/O(3)$
satisfies
\[
 K_h\ge0,\qquad \int_{\mathcal B}K_h\,\mathrm d\mu_h<\infty.
\]
There is no sequence $p_i\in\Sigma_{c_i},\quad c_i:=f(p_i)$, $c_i\to\infty$, with
$W(p_i)\asymp c_i$ and the following properties. In signed level
arclength $z$, centered at $p_i$, there are samples
$z_i^-<0<z_i^+$ such that
\begin{gather*}
 |z_i^-|+|z_i^+|=o(c_i),\\
 (\partial_z W)(z_i^-,c_i)\ge s_H-o(1),\qquad
 (\partial_z W)(z_i^+,c_i)\le s_L+o(1),\qquad s_H>s_L,
\end{gather*}
and either of the following conditions holds. Here a \textbf{level meridian}
is the one-dimensional quotient of a rotational potential level by its
$S^2$ orbits, with the induced arclength metric; a \textbf{pole} is an endpoint
where the orbits collapse.
\begin{enumerate}[label=\textup{(\roman*)},ref=\textup{(\roman*)}]
\item There is $c_0$ such that on the union of the complete
potential levels $\{c\ge c_0\}$ under consideration,
$|\nabla W|\le C$, $\mathrm RW\le C$, and
\[
 \lim_{A\to\infty}\sup_{\{c\ge c_0,\ u\ge A\}}|XW-2|=0.
\]
These bounds include the future normal trajectories and level intervals
used below. Each finite meridian endpoint
is a smooth pole with $W=0$; other meridian ends have infinite length.
\item There are $\lambda_i\to\infty$ with
$\lambda_i/c_i\to0$ and $\lambda_i\ge\max\{|z_i^-|,|z_i^+|\}$.
Let $p_i(c)$ solve
\[
 \frac{\mathrm d p_i}{\mathrm dc}(c)=\frac{X}{|X|^2}(p_i(c)),
 \qquad p_i(c_i)=p_i,
\]
and on each level use signed arclength $z$ with $z(p_i(c))=0$.
For every fixed $R_0,T>0$, these trajectories and the tubes
\[
 c_i\le c\le c_i+T\lambda_i^2,\qquad |z|\le R_0\lambda_i,
\]
exist and satisfy $W\asymp c$, $|X|\ge q_0>0$, $\mathrm R\le C/c$,
and $|\nabla W|\le C$, with constants independent of $i$.
\end{enumerate}
\end{lem}
\begin{proof}
Under the first alternative, follow $p_i$ by $Z=X/q^2$, $q:=|X|$, and center signed level arclength $z$
on this trajectory. These trajectories exist for every future potential:
initially $q\ge1/2$ for large $i$, $q$ increases along $X$, and therefore, we conclude that
$|Z|=q^{-1}$ stays bounded on the complete manifold.
Since $W$ increases and $\mathrm{R}W\le C$,
$ZW=XW/(1-\mathrm{R})\to2$ uniformly along all future centers, so $W(p_i(c))\asymp c$.
We set $\lambda_i:=\max\{\sqrt{c_i},\text{absolute sample distances}\}$.
Then $\lambda_i\to\infty$ and $\lambda_i/c_i\to0$.
Each finite endpoint is at distance at least $W(p_i(c))/C$ from the center
by $|\nabla W|\le C$. Consequently for fixed $R_0,T$ the entire tubes
\[
 c_i\le c\le c_i+T\lambda_i^2,\qquad |z|\le R_0\lambda_i
\]
exist, and on them $W\asymp c$, $\mathrm{R}\le C/c$, $q\ge1/2$.
Thus, the first alternative supplies the tubes and estimates explicitly
assumed in the second. We now use only those tube estimates.
We write
\begin{gather*}
 h=q^{-2}\mathrm{d}c^2+(\mathrm{d}z-\beta \mathrm{d}c)^2,\qquad
 Z=\partial_c+\beta\partial_z,\qquad\beta(c,0)=0,\\
 \mathfrak v=ZW,\qquad P=q^2\mathfrak v=XW,\qquad
 \sigma=q^{-2}\Ric(\partial_z,\partial_z).
\end{gather*}
Normal variation of level arclength, the radius equation, and the horizontal
soliton equation give respectively
\begin{align}
 {\partial_{z} \beta}&=\sigma\ge0,\nonumber\\
 q^2Z\mathfrak v+{\partial_{z}^2 W}-({\partial_{z} q}/q){\partial_{z} W}+\mathfrak v-q^2\mathfrak v^2/W&=2,\label{reflcrit:heat-identities}\\
 \Ric(\partial_z,\partial_z)(1+\mathfrak v/W)
   &=K-{\partial_{z}^2 W}/W+\bigl(\partial_{z} W\bigr)^2/(2W^2).\nonumber
\end{align}
For example the two horizontal Hessian components of $W$ are
$q^2Z\mathfrak v+q(Zq)\mathfrak v-({\partial_{z} q}/q){\partial_{z} W}$ and
${\partial_{z}^2 W}+\Ric(\partial_z,\partial_z)\mathfrak v$; their sum in
$\Delta_hW+XW=2$ proves the middle identity, using the traced soliton
identity. Level concavity gives ${\partial_{z}^2 W}\le \bigl(\partial_{z} W\bigr)^2/(2W)\le C/c$.
We spell out the integrated estimates on the comparison tubes. Fix
$R,T<\infty$, and at each level take
$I_i(c):=[-(R+1)\lambda_i,(R+1)\lambda_i]$ for
$c_i\le c\le c_i+T\lambda_i^2$. Constants below may depend on
$R,T$ and the tube bounds, but not on $i$. Since
$|\partial_zW|\le C$ and $(\partial_z^2W)_+\le C/c$, the
fundamental theorem of calculus gives
\[
 \int_{I_i(c)}(-\partial_z^2W)_+\,\mathrm dz
 \le \int_{I_i(c)}(\partial_z^2W)_+\,\mathrm dz+2C
 \le C+C\lambda_i/c.
\]
The subscript $+$ denotes the positive part. Also
$\mathfrak v=2W\Ric_{\rm orb}/q^2\ge0$, so the last identity in
\eqref{reflcrit:heat-identities}, $W\asymp c$, and $q\ge q_0$
imply
\[
 0\le\sigma\le C K+\frac Cc(-\partial_z^2W)_++\frac C{c^2},
 \qquad
 \int_{I_i(c)}\sigma\,\mathrm dz
 \le C\int_{I_i(c)}K\,\mathrm dz+C/c+C\lambda_i/c^2.
\]
Because $\beta(c,0)=0$ and $\partial_z\beta=\sigma\ge0$, we have
\[
 \int_{I_i(c)}|\beta|\,\mathrm dz
       \le C\lambda_i\int_{I_i(c)}\sigma\,\mathrm dz.
\]
Let $\mathcal R_i$ be this enlarged tube in the $(c,z)$ chart.
The change of variables
$\tau=-(c-c_i)/\lambda_i^2$, $y=z/\lambda_i$ has Jacobian
$\mathrm d\tau\,\mathrm dy=\lambda_i^{-3}\mathrm dc\,\mathrm dz$.
Since $\mathrm d\mu_h=q^{-1}\mathrm dc\,\mathrm dz$ and
$q_0\le q\le1$, the preceding inequalities give, on
$[-R,R]\times[-T,0]$,
\begin{align*}
 \|\lambda_i\beta\|_{L^1}
 &\le C\lambda_i^{-1}\int_{\mathcal R_i}K\,\mathrm d\mu_h
       +CT\lambda_i/c_i+CT\lambda_i^2/c_i^2\longrightarrow0,\\
 \|\lambda_i\sigma\|_{L^1}
 &\le C\lambda_i^{-2}\int_{\mathcal R_i}K\,\mathrm d\mu_h
       +CT/c_i+CT\lambda_i/c_i^2\longrightarrow0.
\end{align*}
The $L^1$ norms here and below use $\mathrm d\tau\,\mathrm dy$.
Finite total curvature bounds the curvature integrals uniformly;
the remaining terms tend to zero because $\lambda_i\to\infty$
and $\lambda_i/c_i\to0$. No assumption that
$\lambda_i^2/c_i$ is small has been used.

We define $d_i^0:=W(p_i)-2c_i$ and
\[
 V_i=(W-2c-d_i^0)/\lambda_i,\qquad
 \widehat V_i=V_i+P/\lambda_i.
\]
Here $V_i(0,0)=0$, $|\partial_{y}\bigl(V_i\bigr)|\le C$, and
$\|\widehat V_i-V_i\|_\infty\le C/\lambda_i$. We retain the normal second
derivative in \eqref{reflcrit:heat-identities}; direct substitution gives
\begin{align*}
 \partial_{\tau}\bigl(\widehat V_i\bigr)-\partial_{y}^2\bigl(V_i\bigr)
 ={}&\lambda_i\beta\partial_{y}\bigl(V_i\bigr)+\beta {\partial_{y} P}
 -\lambda_i Z(q^2)\mathfrak v-\lambda_i({\partial_{z} q}/q)\partial_{y}\bigl(V_i\bigr)
 -\lambda_i q^2\mathfrak v^2/W.
\end{align*}
The correction $P/\lambda_i$ cancels the normal second
derivative term $Z\mathfrak v$. To estimate every remaining term,
use
\[
 |P|+|\mathfrak v|\le C,\qquad
 |Z(q^2)|+|\partial_zq|\le C/c,\qquad W^{-1}\le C/c.
\]
The first bounds follow from $P=2W\Ric_{\rm orb}$, $W\asymp c$,
$\mathrm R\le C/c$, and $q\ge q_0$; the gradient bounds follow
from $\Hess f=\Ric$. We keep the derivative of $P$ in divergence
form:
\[
 \beta\partial_yP=\partial_y(\beta P)-\lambda_i\sigma P.
\]
Thus, with
\begin{align*}
 E_i:={}&\beta P,\\
 e_i:={}&\lambda_i\beta\partial_yV_i-\lambda_i\sigma P
       -\lambda_i Z(q^2)\mathfrak v
       -\lambda_i(\partial_zq/q)\partial_yV_i
       -\lambda_i q^2\mathfrak v^2/W,
\end{align*}
we have exactly
\[
 \partial_\tau\widehat V_i-\partial_y^2V_i=e_i+\partial_yE_i.
\]
On each fixed rectangle, its finite measure and
$|\partial_yV_i|\le C$ give the term-by-term bounds
\begin{align*}
 \|E_i\|_{L^1}
 &\le C\lambda_i^{-1}\|\lambda_i\beta\|_{L^1},\\
 \|e_i\|_{L^1}
 &\le C\bigl(\|\lambda_i\beta\|_{L^1}
             +\|\lambda_i\sigma\|_{L^1}+\lambda_i/c_i\bigr)
       \longrightarrow0.
\end{align*}
In particular, no bound for $\partial_yP$ is needed.

We next give compactness up to the terminal time. Work on a fixed
rectangle and reserve a spatial margin of one. Let
$\varrho_\varepsilon$ be a nonnegative smooth spatial mollifier of
integral one, supported in $[-\varepsilon,\varepsilon]$, where
$0<\varepsilon<1$. Convolution of the last equation yields
\[
 \partial_\tau(\varrho_\varepsilon*\widehat V_i)
 =(\partial_y\varrho_\varepsilon)*\partial_yV_i
       +\varrho_\varepsilon*e_i+(\partial_y\varrho_\varepsilon)*E_i.
\]
Consequently, uniformly in the spatial center and for
$-T\le\tau_1\le\tau_2\le0$,
\begin{align*}
 &|(\varrho_\varepsilon*\widehat V_i)(y,\tau_2)
             -(\varrho_\varepsilon*\widehat V_i)(y,\tau_1)|\\
 &\hspace{1cm}\le A_\varepsilon|\tau_2-\tau_1|
       +A_\varepsilon(\|e_i\|_{L^1}+\|E_i\|_{L^1}),
\end{align*}
where $A_\varepsilon$ is independent of $i$. The first term uses
$\|(\partial_y\varrho_\varepsilon)\|_{L^1}|\partial_yV_i|\le A_\varepsilon$;
the other two use the supremum norms of the fixed mollifier and
its derivative. The spatial Lipschitz bound and the correction give
\[
 |V_i-\varrho_\varepsilon*\widehat V_i|
       \le C\varepsilon+C/\lambda_i.
\]
Taking one fixed mollifier, $V_i(0,0)=0$ bounds its terminal
average. The time-increment estimate then bounds this average
throughout $[-T,0]$, and the spatial Lipschitz estimate bounds
$V_i$ on the whole smaller rectangle. To obtain equicontinuity,
first choose $\varepsilon$ small, then $i$ large so that the
residual norms and $\lambda_i^{-1}$ are small, and finally choose
$|\tau_2-\tau_1|$ small. This order proves asymptotic uniform
equicontinuity, including at $\tau=0$, rather than just compactness
on open negative-time rectangles. Arzel\`a--Ascoli and a diagonal
subsequence therefore give local uniform convergence on
$\mathbb R\times(-\infty,0]$. Since
$\widehat V_i-V_i\to0$ uniformly and $e_i,E_i\to0$ in local
$L^1$, testing against compactly supported smooth functions gives
\[
 \partial_\tau V=\partial_y^2V\quad(\tau<0),\qquad
 |\partial_yV|\le C.
\]
The limit is smooth for $\tau<0$ and continuous through $\tau=0$.
The bounded ancient caloric function ${\partial_{y} V}$ is constant: represent it by the
heat kernel over an arbitrary elapsed time $S$ and differentiate to get
$|\partial_y{\partial_{y} V}|\le C\|{\partial_{y} V}\|_\infty/\sqrt S$, then let $S\to\infty$.
Thus $V$ has a single affine spatial slope.

At $\tau=0$, however,
$\partial_y^2V_i\le C\lambda_i/c_i\to0$. The rescaled sample
locations $y_i^\pm:=z_i^\pm/\lambda_i$ lie in $[-1,0)$ and
$(0,1]$, respectively. Integrating this one-sided second derivative
bound from a point in $[-3,-2]$ to $y_i^-$ gives
$\partial_yV_i\ge s_H-o(1)$ throughout $[-3,-2]$.
Likewise $\partial_yV_i\le s_L+o(1)$ throughout $[2,3]$.
Integration over these fixed intervals and local uniform convergence
at the terminal time give
\[
 V(-2,0)-V(-3,0)\ge s_H,\qquad
 V(3,0)-V(2,0)\le s_L.
\]
A single affine spatial slope would have to be both at least $s_H$
and at most $s_L$, contradicting $s_H>s_L$.
This finishes the proof.
\end{proof}

\begin{cor}\label{reflcrit:Lzero}
A soliton with no critical point as in Lemma~\ref{reflcrit:foliation}
satisfies $L_+:=1-a_+^2=0$, where
$a_+=\lim_{s\to\infty}\frac{\mathrm d}{\mathrm ds}(f\circ\gamma)(s)$ is defined in Lemma~\ref{reflcrit:foliation}.
\end{cor}
\begin{proof}
If $L_+>0$, Lemma~\ref{reflcrit:tent} gives the two distinct slopes at its
ridge. We justify the passage to exact levels without assuming smooth
convergence at the distance scale. We fix a sufficiently small $\delta>0$.
At scale $T_i\to\infty$, choose quotient points approximating the five
interior cone points $(1,j\delta)$, $j=-2,-1,0,1,2$. We set
$c:=f$ and project them onto the exact level $c=T_i$ along
$Z=X/q^2$, $q:=|X|$; call the projected points $p_{j,i}$.
The decorated convergence was stated for $(f-f(o))/T_i$, but
$f(o)/T_i\to0$, so it also gives $f/T_i\to x$.
Thus the initial potential discrepancy is $o(T_i)$. In a fixed strip with compact closure in the comparison domain
about these cone points, $W\ge bT_i$ by the decorated convergence, so $\mathrm{R}\le C/T_i$ and $q\to1$. Also, we have
\[
                 0<ZW=2W\Ric_{\rm orb}/q^2\le C
\]
there. A continuation argument shows that each projection trajectory,
whose potential length is $o(T_i)$, retains $W\ge bT_i/2$, has metric
length $o(T_i)$, and stays in this strip. Thus, the projected points have
the same cone positions and limiting values of $W/T_i$.

We require the intrinsic level distances as well. A quotient minimizing
chord between any two of the initial points converges to their Euclidean
segment: the limiting straight segment is the unique minimizer and lies
strictly inside the sector. Uniform convergence of the decorated functions
therefore gives $W\ge bT_i$ and $|c-T_i|=o(T_i)$ along the whole chord.
We project this chord to $c=T_i$. If $\varphi_t$ denotes the $Z$-flow, the
projection is $\Pi(x)=\varphi_{T_i-f(x)}(x)$ and its differential is
\[
 \mathrm{d}\Pi_x(v)=\mathrm{d}\varphi_{T_i-f(x)}
             \bigl(v-\mathrm{d}f(v)Z(x)\bigr).
\]
The vector in parentheses is precisely the orthogonal projection to the
level tangent space, so its norm is at most $|v|$. Along the short
projection trajectories, the preceding continuation estimate gives
$\mathrm{R}\le C/T_i$ and $q\ge1/2$. For tangent vectors normal-potential variation
is $\partial_c\gamma=2q^{-2}\Ric|_{T\{c=\mathrm{constant}\}}$.
Consequently $\mathrm{d}\varphi_t$ on these vectors has norm at most
$\exp(C|t|/T_i)=1+o(1)$; for negative $t$ it is in fact nonexpanding.
The projected chord thus has length at most $1+o(1)$ times its original
length. The quotient distance gives the reverse bound. We have proved,
for every pair among the five points,
\[
 d_{\{c=T_i\}/O(3)}(p_{j,i},p_{k,i})
                   =(|j-k|\delta+o(1))T_i.
\]

The level quotient is a ray by Lemma~\ref{reflcrit:foliation}. These
pairwise distances place the five points in the stated order or its
reverse for large $i$. The increase in $W/T_i$ from $j=-2$ to $j=-1$
is $s_+\delta>0$; since the radius strictly increases toward the
noncompact end, this selects the stated order in that orientation.
Center signed level arclength at $p_{0,i}$. The mean-value theorem on the
intervals $[p_{-2,i},p_{-1,i}]$ and $[p_{1,i},p_{2,i}]$ now supplies
a left sample with ${\partial_{z} W}=s_++o(1)$ and a right sample with
${\partial_{z} W}=s_-+o(1)$. This also covers $s_-=0$; no strictly positive lower
bound for the right slope is needed. We choose $\delta_i\downarrow0$
slowly and $T_i\to\infty$ sufficiently fast so that all distance and
radius errors are $o(\delta_iT_i)$. The resulting slope errors tend
to zero while both sample distances are $o(T_i)$. The central values
satisfy $W(p_{0,i})/T_i\to2$.
The future-tube hypotheses in Lemma~\ref{reflcrit:heat-ridge} hold by
Lemma~\ref{reflcrit:foliation}, and the quotient has finite total curvature
by its smooth meridian double. That lemma gives the contradiction.
The assertion follows.
\end{proof}

We next construct the scalar barrier needed to control the radius
profile. The conclusion will be used to restrict its possible
growth at infinity.

\begin{lem}
\label{lem:ambient-barrier-comparison}
Work on a region of regular potential levels in an $O(3)$-invariant
normalized steady soliton. Let $c=f$ in cases \textup{(ii)}--\textup{(iii)}
below, and let $c=f-m$, $m=\inf_M f$, in case \textup{(i)}.
Restrict to $c\ge c_0>0$. Put
$V:=|\nabla_hu|_h^2$ and assume the quotient equations
\eqref{reflcrit:quotient}, $0<|X|\le1$, $V\le1$, $0<Xu\le C/u$,
and the envelope $u^2-2c\le E_0u$. Write
$U(c):=\sup u$ over the level.
For every fixed $d>0$, assume $V\to0$ uniformly on the
region $u\ge d\sqrt c$ as $c\to\infty$.
Suppose the region $\mathcal U$ has one of the following three
geometries:
\begin{enumerate}[label=\textup{(\roman*)}]
\item $m:=\inf_Mf> -\infty$, $c=f-m$, and
$\mathcal U=\bigcup_{c\ge c_0}\Sigma_{m+c}$ consists of complete
one-tip levels, each parametrized by $r\in[0,\infty)$.
The only finite meridian endpoint is the smooth pole $u(0,c)=0$;
as $r\to\infty$ the limits \eqref{eq:finite-infimum-end} hold,
uniformly for $c$ in compact intervals.
\item The region is swept out by level segments
$0\le r\le\bar r(c)$, where $r$ is arclength from the smooth pole,
$\bar r(c)$ varies smoothly, and for every $c_0\le a<b<\infty$
the swept set $\{a\le c\le b,\ 0\le r\le\bar r(c)\}$ is compact.
Its outer boundary is $r=\bar r(c)$, where
\[
 \bar u(c):=u(\bar r(c),c)\longrightarrow\infty,
 \qquad \bar u(c)^2V(\bar r(c),c)\le C_0.
\]
\item The region consists of entire compact two-tip levels
$\Sigma_c$, whose meridians are $[0,D(c)]$, with smooth poles
$u(0,c)=u(D(c),c)=0$. Every retained slab
$\{a\le c\le b\}\cap\mathcal U$ is compact for finite $a,b$.
There is no noncompact spatial end or additional outer boundary.
\end{enumerate}
Then, for one fixed $K_0$, throughout this region where
$u>r_*\sqrt{K_0}$, we have
\[
                         V\le\phi_B(u/\sqrt{K_0}).
\]
Here $r_*>0$ is the fixed inner endpoint of the singular Bryant
profile used in the proof, and $\phi_B:(0,\infty)\to(0,\infty)$
is that squared slope profile, with
$\phi_B(t)=t^{-2}+2t^{-4}+O(t^{-6})$ as $t\to\infty$,
as in \cite[Theorem 2.1]{BrendleAncient}. We choose $r_*$ by
$\phi_B(r_*)=2$, as in the paragraph following
\cite[Remark 2.2]{BrendleAncient}.
\end{lem}
\begin{proof}
We specify the one-variable barrier functions. The construction in
\cite[Theorem 2.1, Lemmas 2.3--2.5, Definition 2.6 and Proposition 2.7]
{BrendleAncient}, with the quantitative residual in
\cite[Proposition 2.19]{MMS}, gives positive $C^1$, piecewise $C^2$ functions
$\psi_a$ on $[r_*/a,1+a^{-2}/100]$, for all sufficiently large $a$, such that
\begin{align*}
 \mathcal Q(\psi_a)-x\frac{\mathrm d \psi_a}{\mathrm dx}&\le-d_*a^{-4}x^{-5}
                    &&(x\ge N/a),\\
 \mathcal Q(\psi_a)-x\frac{\mathrm d \psi_a}{\mathrm dx}&\le-d_*a
                    &&(r_*/a\le x\le N/a),\\
 \mathcal Q(\psi)&=\psi\frac{\mathrm d^2 \psi}{\mathrm dx^2}-\tfrac12(\frac{\mathrm d \psi}{\mathrm dx})^2
                  +(1-\psi)\frac{\mathrm d \psi}{\mathrm dx}/x+2\psi(1-\psi)/x^2.
\end{align*}
They satisfy $\psi_a\ge a^{-4}/32$, $\psi_a(r_*/a)\ge3/2$.
On the inner piece $\psi_a\ge b$, $|\frac{\mathrm d \psi_a}{\mathrm dx}|\le Ca$,
$|\frac{\mathrm d^2 \psi_a}{\mathrm dx^2}|\le Ca^2$. On the middle piece, separated from the final endpoint,
$\psi_a\ge ka^{-2}x^{-2}$; on the entire outer piece
$\left|\frac{\mathrm d^j\psi_a}{\mathrm dx^j}\right|\le Ca^{-2}x^{-2-j}$ for $j\le2$.
Near the final endpoint the lower bound instead is
\[
          \psi_a(x)\ge a^{-2}(x^{-2}-1)+a^{-4}/16.
\]
Finally, $\psi_a(t/a)\to\phi_B(t)$ on each fixed compact subset of
$t>r_*$. These statements follow from the explicit two-piece formulas in
the cited construction.

We set $M:=2c+Ka^2$, $x:=u/\sqrt M$, and $B(u,c):=\psi_a(x)$.
The assumed envelope gives, for all sufficiently large $K$ independently
of $a$,
\begin{equation}\label{reflcrit:barrier-denominator}
 x<1+a^{-2}/100,\qquad u^2B\ge k_0K
 \quad\hbox{whenever }x\ge r_*/a.
\end{equation}
For the first assertion minimize
$d_au^2-(1+d_a)E_0u+(1+d_a)Ka^2$, with
$d_a=(1+a^{-2}/100)^2-1$. The inner and middle lower bounds give the second
assertion on those pieces. At the endpoint use specifically
\[
 u^2B\ge(M-u^2)/a^2+u^2/(16a^4)
       \ge K-E_0u/a^2+(u/a^2)^2/16\ge K-4E_0^2.
\]

We now compute the ambient differential inequality at a contact point.
In all three cases $\mathrm dc=\mathrm df$, so the differential
identities below are the same for $c=f$ and $c=f-m$.
We write $\mathcal L:=\Delta_h+X$, $q:=|X|$, and $w:=Xu$.
The equations \eqref{reflcrit:quotient} give
\[
 \mathcal Lu=(1-V)/u,\quad \mathcal Lc=1-2w/u,\quad
 \mathcal LV=2|\Hess_hu|^2-2V(1-V)/u^2,\quad 0<w\le C/u.
\]
At a zero local maximum of $V-B$, choose $e_1=\nabla u/\sqrt B$ and
put $b_\perp:=\Hess_hu(e_2,e_2)$. Resolving $\mathrm{d}V={\partial_{u} B}\,\mathrm{d}u+{\partial_{c} B}\,\mathrm{d}f$ in this
frame gives
\[
 2|\Hess_hu|^2=\tfrac12\bigl(\partial_{u} B\bigr)^2+{\partial_{u} B}{\partial_{c} B}w/B+\bigl(\partial_{c} B\bigr)^2q^2/B
                        -\bigl(\partial_{c} B\bigr)^2w^2/(2B^2)+2b_\perp^2.
\]
Consequently, with $\mathcal Q_u$ denoting the same expression as
$\mathcal Q$ in variable $u$,
\begin{align}
 \mathcal L(V-B)={}&-\mathcal Q_u(B)-{\partial_{c} B}
 +w({\partial_{u} B}{\partial_{c} B}/B-2{\partial_u\partial_c B}+2{\partial_{c} B}/u)\nonumber\\
 &+q^2(\bigl(\partial_{c} B\bigr)^2/B-{\partial_{c}^2 B})-\bigl(\partial_{c} B\bigr)^2w^2/(2B^2)+2b_\perp^2.
 \label{reflcrit:contact}
\end{align}
The derivatives are
\begin{align*}
 \partial_u B&=M^{-1/2}\frac{\mathrm d\psi_a}{\mathrm dx},&
 \partial_c B&=-\frac{x}{M}\frac{\mathrm d\psi_a}{\mathrm dx},\\
 \partial_u\partial_c B
 &=-\frac{\frac{\mathrm d\psi_a}{\mathrm dx}
                   +x\frac{\mathrm d^2\psi_a}{\mathrm dx^2}}{M^{3/2}},&
 \partial_c^2 B
 &=\frac{3x\frac{\mathrm d\psi_a}{\mathrm dx}
                   +x^2\frac{\mathrm d^2\psi_a}{\mathrm dx^2}}{M^2}.
\end{align*}
In particular, $-\mathcal Q_u(B)-{\partial_{c} B}
=-[\mathcal Q(\psi_a)-x\frac{\mathrm d \psi_a}{\mathrm dx}]/M$.
We now bound each possibly negative contribution explicitly. Set
\[
 \begin{aligned}
 \mathcal E_1&:=w|\partial_uB\,\partial_cB|/B,\\
 \mathcal E_2&:=2w\bigl(|\partial_u\partial_cB|
                                  +|\partial_cB|/u\bigr),\\
 \mathcal E_3&:=q^2|\partial_c^2B|,\qquad
 \mathcal E_4:=(\partial_cB)^2w^2/(2B^2).
 \end{aligned}
\]
Since $q\le1$, \eqref{reflcrit:contact} implies
\[
 \mathcal L(V-B)\ge
 -[\mathcal Q(\psi_a)-x\,\mathrm d\psi_a/\mathrm dx]/M
             -\sum_{j=1}^4\mathcal E_j.
\]
Here we discarded only the nonnegative terms
$q^2(\partial_cB)^2/B$ and $2b_\perp^2$.

On the outer piece $x\ge N/a$, use
$|\mathrm d^j\psi_a/\mathrm dx^j|\le Ca^{-2}x^{-2-j}$ for
$j=1,2$, $w\le C/(\sqrt Mx)$, and
$B^{-1}\le Mx^2/(k_0K)$ from
\eqref{reflcrit:barrier-denominator}. Substituting into the four
terms, without requiring a positive lower bound for $\psi_a$
independent of $a$, gives
\[
 \begin{aligned}
 \mathcal E_1&\le\frac{Ca^{-4}x^{-4}}{KM},&
 \mathcal E_2&\le\frac{Ca^{-2}x^{-4}}{M^2},\\
 \mathcal E_3&\le\frac{Ca^{-2}x^{-2}}{M^2},&
 \mathcal E_4&\le\frac{Ca^{-4}x^{-2}}{K^2M}.
 \end{aligned}
\]
For example, the two denominator factors in $\mathcal E_4$ give
\[
 \mathcal E_4\le
 C\frac{a^{-4}x^{-4}}{M^2}\frac1{Mx^2}
                 \frac{M^2x^4}{k_0^2K^2}.
\]
Relative to the strict residual scale $D:=a^{-4}x^{-5}/M$,
the four estimates become
\[
 \frac{\mathcal E_1}{D}\le\frac{Cx}{K},\quad
 \frac{\mathcal E_2}{D}\le\frac{Ca^2x}{M},\quad
 \frac{\mathcal E_3}{D}\le\frac{Ca^2x^3}{M},\quad
 \frac{\mathcal E_4}{D}\le\frac{Cx^3}{K^2}.
\]
Since $x\le1+a^{-2}/100\le2$ and $M\ge Ka^2$, their sum is
at most $C(K^{-1}+K^{-2})D$. This also treats the region near the
outer endpoint, where the uniform floor for $\psi_a$ alone would
not give the required denominator estimate.

On the inner piece $r_*/a\le x\le N/a$, the bounds
$\psi_a\ge b>0$, $|\mathrm d\psi_a/\mathrm dx|\le Ca$, and
$|\mathrm d^2\psi_a/\mathrm dx^2|\le Ca^2$ instead give
\[
 \begin{gathered}
 |\partial_uB|\le Ca/\sqrt M,\qquad
 |\partial_cB|\le C/M,\qquad
 |\partial_u\partial_cB|\le Ca/M^{3/2},\\
 |\partial_c^2B|\le C/M^2,\qquad
 w\le Ca/\sqrt M,\qquad B^{-1}\le C.
 \end{gathered}
\]
It follows that
\[
 \mathcal E_1+\mathcal E_2\le Ca^2/M^2,\qquad
 \mathcal E_3\le C/M^2,\qquad
 \mathcal E_4\le Ca^2/M^3.
\]
Dividing their sum by the inner strict residual scale $a/M$
gives
\[
 C\left(\frac aM+\frac1{aM}+\frac a{M^2}\right)
                  \le\frac C{Ka}\qquad(K,a\ge1).
\]
All constants here depend only on the fixed barrier and geometric
bounds. Choose $K_0$ large enough that
$C(K_0^{-1}+K_0^{-2})<d_*/2$, and then choose the lower threshold
for $a$ so that $C/(K_0a)<d_*/2$. The same choices work for every
$K\ge K_0$. On each piece the negative error is therefore smaller
than half the corresponding positive residual. Thus
$\mathcal L(V-B)>0$, contrary to a zero maximum. At the $C^1$ matching
interface the same conclusion follows from either one-sided Hessian test;
if $\mathrm{d}x=0$, the sole discontinuous chain-rule term $\frac{\mathrm d^2 \psi_a}{\mathrm dx^2}\mathrm{d}x\otimes \mathrm{d}x$
vanishes.

We compare on $\Omega_{a,K}=\{c\ge c_0, u\ge r_*\sqrt{2c+Ka^2}/a\}$
within the chosen region. We enlarge $K_0$ to make its initial $c=c_0$ face
empty. Its inner boundary has $B\ge3/2>V$. In the reference case,
\eqref{reflcrit:barrier-denominator} and $\bar u^2V\le C_0$ strictly
protect the outer boundary after increasing $K_0$.
In case (i) the noncompact end is also strictly protected.
Indeed, $C_\Phi:=\sup_{v\ge U(c_0)}v^2\Phi_B(v)<\infty$, by the regular
Bryant asymptotics, so $k_0K_0>C_\Phi$ makes
\[
 \psi_a(U(c)/\sqrt{2c+Ka^2})>\Phi_B(U(c)).
\]
Lemma~\ref{reflcrit:envelope} makes this a uniform limiting margin on
compact $c,K$ intervals. There is no contact through the upper endpoint
of the barrier's domain by \eqref{reflcrit:barrier-denominator}.

Now fix $a$. On all these domains $u\ge r_*\sqrt{2c}/a$, so thick
cylindrical convergence gives $V\to0$ as $c\to\infty$, uniformly in
$K\ge K_0$. The floor $B\ge a^{-4}/32$ therefore protects a common tail
$c\ge C(a)$. On the remaining finite $c$ interval, the radius is bounded:
by the envelope in case (i), by compactness of the swept segments
in case (ii), and by compactness of potential slabs in case (iii).
In case (iii) the poles have $u=0$ and lie outside the comparison
domain, so there is no additional outer boundary. Thus $\Omega_{a,K}$
is empty there for sufficiently large $K$. Starting from these values,
lower $K$ continuously. A first loss of comparison cannot occur at the
protected faces or ends, and hence would give an interior zero maximum,
already excluded. Therefore $V\le B_{a,K_0}$ for this fixed $a$.
The same $K_0$ works for every large $a$. At a fixed point with
$u>r_*\sqrt{K_0}$, let $a\to\infty$ to obtain the asserted singular-profile
bound. This order does not require a tail threshold uniform in $a$.
This proves the stated conclusion.
\end{proof}

We now apply the preceding comparison to the two level geometries
that arise when the potential has no critical point: full one-tip
levels and the swept tip-to-reference segments.
\begin{lem}
\label{reflcrit:barrier}
Assume that the potential has no critical point, and use the notation
of Lemma~\ref{reflcrit:foliation}. Use $c=f-m$ in the
finite-infimum full-level case below and $c=f$ in the
tip-to-reference case. Fix constants $c_0,E_0,C_0>0$,
and put $V=|\nabla_hu|_h^2$.
Suppose that $u^2-2c\le E_0u$ on the region under consideration,
where $c\ge c_0$. Suppose either that this region consists of the
full one-tip levels with $c=f-m$, $m> -\infty$, and the
finite-infimum end limits \eqref{eq:finite-infimum-end}, or that it
consists of smoothly varying
segments $0\le r\le\bar r(c)$, where $r$ is level arclength from the
unique tip. In the latter case assume that the swept segments over
bounded $c$ intervals are compact and that their outer endpoints satisfy
\[
 \bar u(c):=u(\bar r(c),c)\longrightarrow\infty,
 \qquad \bar u(c)^2V\le C_0.
\]
Assume the quotient identities \eqref{reflcrit:quotient},
$0<|X|\le1$, $|\nabla_hu|^2\le1$, and $0<Xu\le C/u$.
Assume also the thick-region limit of Lemma~\ref{reflcrit:thick}, namely
\[
 \sup_{\mathcal U\cap\{c\ge A,\ u\ge d\sqrt c\}}
             |\nabla_hu|^2\longrightarrow0
 \quad(A\to\infty),\qquad\text{for every fixed }d>0.
\]
The complete levels in the first case and the swept segments in the
second case have precisely the boundary and compactness properties
of cases \textup{(i)} and \textup{(ii)} of
Lemma~\ref{lem:ambient-barrier-comparison}, respectively. Then there is a fixed $K_0$
such that
\[
 V\le\phi_B(u/\sqrt{K_0})\qquad\text{where }u>r_*\sqrt{K_0}.
\]
\end{lem}

\begin{proof}
The quotient identities and Lemma~\ref{reflcrit:thick} supply the
geometric bounds and the uniform thick-band decay required by
Lemma~\ref{lem:ambient-barrier-comparison}. Its cases \textup{(i)}
and \textup{(ii)}, respectively, give the assertion.
\end{proof}

\begin{cor}\label{reflcrit:infinite-infimum}
In the no-critical point setting of Lemma~\ref{reflcrit:foliation},
$\inf_Mf=-\infty$.
\end{cor}
\begin{proof}
Suppose that $m$ is finite. Lemma~\ref{reflcrit:envelope} gives
$u^2-2c\le E_0u$ for all large $c=f-m$. We apply
Lemma~\ref{reflcrit:barrier} to the full levels, using
Lemma~\ref{reflcrit:thick}. By Corollary~\ref{reflcrit:Lzero} and
Lemma~\ref{reflcrit:decay}, scalar curvature tends uniformly to zero on
these high levels. We fix $u_0>r_*\sqrt{K_0}$ with
$\phi_B(u_0/\sqrt{K_0})<1$. Since $U(c)\to\infty$, each late level
reaches this radius. If $M(c)=\sup_{\{f=m+c\}}\mathrm{R}$, Gauss and $q^2\ge1/2$
give $0\le K_{\rm rad}\le C M(c)$. The meridian Jacobi equation,
$u(0)=0,{\partial_{r} u}(0)=1$, gives on every fixed bounded interval
\[
 {\partial_{r} u}(r)\ge1-CM(c)r^2/2,\qquad
 u(r)\ge r-CM(c)r^3/6.
\]
At $r=2u_0$ the second bound exceeds $u_0$ for large $c$, so the first
radius-$u_0$ point lies before it and has ${\partial_{r} u}\to1$. But there
$V\ge \bigl(\partial_{r} u\bigr)^2\to1$, contradicting the strict singular-profile bound.
This finishes the proof.
\end{proof}

The no-critical point branch has only one relevant cap direction.
We now analyze that branch using a fixed normal reference
trajectory and the associated profile equation.

\begin{lem}
\label{reflcrit:onecap}
Assume the setting of Lemma~\ref{reflcrit:foliation}. Fix the
normal trajectory of Lemma~\ref{reflcrit:thick}, write
$c:=f$ as the level parameter, and let $F(z,c)$ be the orbit radius in signed level
arclength from that trajectory, increasing toward the noncompact end.
There are $c_0,C>0$ such that for every $c\ge c_0$, we have
\[
 |F(0,c)-\sqrt{2c}|\le C,\qquad
 0<{\partial_{z} F}(0,c)\le Cc^{-1/2}.
\]
\end{lem}

\begin{proof}
We first derive the profile equation along the fixed reference
trajectory. We then estimate its residual and apply the one-variable
barrier to obtain the radius and slope bounds.

By Corollaries~\ref{reflcrit:Lzero} and \ref{reflcrit:infinite-infimum},
$m=-\infty$ and $L_+=0$. Levels have unbounded increasing radius, and
Lemma~\ref{reflcrit:thick} gives
$F(0,c)^2/c\to2$, $c\mathrm{R}\to1$ at the fixed reference. Every fixed rescaled
meridian interval around it lies in the level: otherwise the tip is at
bounded scalar distance, contradicting its $c\mathrm{R}_{\rm tip}\to\infty$.

We set $q:=|X|$ and write $\gamma_c$ for the metric induced on
the potential level $\Sigma_c:=\{f=c\}$, indexed by the potential value $c$. Its normal
variation is
\[
 \partial_c\gamma_c=2\Ric_{\gamma_c}+2E,\qquad
 E=q^{-2}\Ric|_{T\Sigma_c}-\Ric_{\gamma_c}.
\]
In a level orthonormal frame, Gauss gives
\[
 E_{ij}=\Rm(n_f,e_i,n_f,e_j)
       +q^{-2}\Ric_{ij}\Ric(n_f,n_f)
       +q^{-2}\sum_a\Ric_{ia}\Ric_{aj},\qquad n_f=X/q.
\]
In particular, $E\ge0$. The full-index soliton commutation identity is
\[
 \Rm(X,\cdot,X,\cdot)
 =-\tfrac12\Hess \mathrm{R}-\Ric^2+\Delta\Ric+2\Rm(\Ric).
\]
It follows by differentiating $\mathrm d\mathrm{R}=-2\Ric(X,\cdot)$, commuting the
covariant derivatives of Ricci using the soliton equation, and applying
$\Delta\Ric+\nabla_X\Ric+2\Rm(\Ric)=0$.
Together with \eqref{reflcrit:derivatives} and $q\ge1/2$ on high levels, it gives
$|\nabla_{\gamma_c}^jE|\le C_j\mathrm{R}^{2+j/2}$. Differentiating the tangent
projections adds second fundamental form factors $O(\mathrm{R})$ and preserves
these orders. If $E=E_{\rm rad}\mathrm{d}z^2+E_{\rm orb}g_{S^2}$, then we obtain
\begin{equation}\label{reflcrit:source-size}
 0\le E_{\rm rad}\le CF^{-4},\qquad
 0\le E_{\rm orb}\le CF^{-2},\qquad
 |\partial_{z}\bigl(E_{\rm orb}\bigr)|\le CF^{-3}.
\end{equation}
Here we used $\mathrm{R}F^2\le C$ and $|{\partial_{z} F}|\le1$.

Restoring arclength after normal variation, with zero arclength shift at
$z=0$, gives the exact equation
\begin{align*}
 -{\partial_{c} F}={}&{\partial_{z}^2 F}-(1+\bigl(\partial_{z} F\bigr)^2)/F
 +2{\partial_{z} F}\left[({\partial_{z} F}/F)(0,c)-\int_0^z\bigl(\partial_{z} F\bigr)^2/F^2\,\mathrm{d}\zeta\right]\\
 &-E_{\rm orb}/F+{\partial_{z} F}\int_0^zE_{\rm rad}\,\mathrm{d}\zeta.
\end{align*}
In particular the nonzero slope at the reference has not been discarded.
We set
\[
 \tau:=-\log c,\quad \xi=z/\sqrt c,\quad
 r=\sqrt2+G=c^{-1/2}F(\sqrt c\xi,c),\quad d={\partial_{\xi} G}(0,\tau).
\]
Then we have
\begin{align*}
	{\partial_{\tau} G}&=\mathscr LG+\mathcal N+\mathcal S
\end{align*}
where
\begin{align*}
 \mathscr L&=\partial_\xi^2-\tfrac\xi2\partial_\xi+1,
 \\
 \mathcal N&=-G^2/(2r)-\bigl(\partial_{\xi} G\bigr)^2/r+2d{\partial_{\xi} G}/r(0)
       -2{\partial_{\xi} G}\int_0^\xi \bigl(\partial_{\xi} G\bigr)^2/r^2\,d\eta,\nonumber\\
 \mathcal S&=-E_{\rm orb}/r+{\partial_{\xi} G}\int_0^\xi cE_{\rm rad}\,d\eta.
 \nonumber
\end{align*}
We next give the cutoff and spectral estimates, including the reference trace,
using monotonicity of the one-cap profile.

For large integers $k$, define
\[
 \delta_k:=e^{-k/2}+\sup_{\tau\le-k}
                  (|G(0,\tau)|+{\partial_{\xi} G}(0,\tau)).
\]
Then $\delta_k\downarrow0$. We choose $a_k^{-2}=B\delta_k$ with a fixed large
$B$. On the tip-to-reference interval use radius $v=F$ as coordinate and
let $b(v,c):=\bigl(\partial_{z} F\bigr)^2$. Its equation, obtained from normal metric variation, is
\begin{align*}
 -{\partial_{c} b}={}&b {\partial_{v}^2 b}-\tfrac12\bigl(\partial_{v} b\bigr)^2
       +v^{-2}(1-b)(v{\partial_{v} b}+2b)\\
       &+2bE_{\rm rad}+v^{-1}(\partial_v b)E_{\rm orb}
                         -2b\partial_{v}\bigl(E_{\rm orb}/v\bigr).
\end{align*}
Since ${\partial_{v} b}=2{\partial_{z}^2 F}<0$, \eqref{reflcrit:source-size} bounds the last line
above by $C\sqrt b\,v^{-4}$. We compare with the one-variable barrier
$\psi_{a_k}(v/\sqrt{2c})$ from Lemma~\ref{reflcrit:barrier}. Its strict
residual dominates this error once $c\ge Da_k^3$: on the outer piece the
ratio of the source to the strict defect is at most $Ca_k^3/c$, and the
inner estimates give the same threshold. This holds on $\tau\le-k$
because $a_k^3\le Ce^{3k/4}$. At the inner radius boundary the barrier is
at least $3/2$. At the reference, choice of $B$ gives
\[
 F(0,c)/\sqrt{2c}\le1+a_k^{-2}/100,\qquad
 b(F(0,c),c)=d^2\le\delta_k^2<a_k^{-4}/32.
\]
For fixed $k$, cylindrical convergence initializes comparison as
$c\to\infty$ on this fixed radius-fraction interval. Running toward smaller
$c$ gives no first contact by the strict residual and the two boundaries.
It follows that $\bigl(\partial_{z} F\bigr)^2\le C\delta_k$ where $F\ge\sqrt{c/2}$ on the inner
side. Integration from the reference, with a first-crossing argument,
extends the neck to $|\xi|\le2\delta_k^{-1/100}$; on the outer side
concavity and ${\partial_{z} F}\le d$ give the same extension. Uniform thick-region
curvature estimates and interpolation then yield on $[-k-1,-k]$
\begin{equation}\label{reflcrit:rough-window}
 |G|+|{\partial_{\xi} G}|\le C\delta_k^{1/4},\qquad
 |{\partial_{\xi}^2 G}|\le C\delta_k^{1/8},\qquad
 |\partial_\xi^jG|\le C_j\ (j\ge3).
\end{equation}
For example the variation of $G$ across this window is bounded by
$C\delta_k^{1/2-1/100}$; the displayed weaker exponent leaves room for
uniform interpolation. All estimates are on a fixed positive radius-fraction
band, so their constants are independent of $k$.

Let $\mathrm{d}\nu:=(4\pi)^{-1/2}e^{-\xi^2/4}\mathrm{d}\xi$, fix a cutoff equal to one on
$[-1,1]$ and zero outside $[-2,2]$, and freeze
$\chi_k(\xi)=\chi(\delta_k^{1/100}\xi)$ on each block.
We set $Y_k:=\chi_kG$ and $\gamma_k:=\|Y_k\|_{L^2(\mathrm{d}\nu)}^2$.
Increasing $k$ if necessary, monotonicity and concavity give
\begin{equation}\label{reflcrit:trace}
                  |G(0)|^2+d^2\le C\gamma_k.
\end{equation}
For the value, use $[0,1]$ if $G(0)\ge0$ and $[-1,0]$ otherwise.
If $|G(0)|<d/2$, concavity gives $G(-1)<-d/2$, and monotonicity preserves
this bound on $[-3,-2]$. Otherwise the value estimate controls $d$.
These intervals all lie in the cutoff plateau for large $k$.

We estimate every term in the residual and apply the barrier comparison. Integrating the derivative of
$\chi_k^4G \bigl(\partial_{\xi} G\bigr)^3e^{-\xi^2/4}$ and using
\eqref{reflcrit:rough-window} gives
\[
 \int_{\mathbb R}\chi_k^4\bigl(\partial_{\xi} G\bigr)^4\,\mathrm{d}\nu
 \le C\delta_k^{1/100}\gamma_k
       +C\exp(-b\delta_k^{-1/50}).
\]
Indeed, the two interior terms are bounded by
$C\delta_k^{1/8}\gamma_k^{1/2}
 (\int_{\mathbb R}\chi_k^4\bigl(\partial_{\xi} G\bigr)^4\,\mathrm{d}\nu)^{1/2}$, since
$|\xi {\partial_{\xi} G}|\le C\delta_k^{1/4-1/100}$; Young's inequality proves the
claim. Changing cutoff powers only creates the same Gaussian tail.
This controls $\bigl(\partial_{\xi} G\bigr)^2/r$. The term $G^2/r$ is bounded by
$C\delta_k^{1/4}\gamma_k^{1/2}$. The trace estimate and the preceding
quartic bound control $d{\partial_{\xi} G}$ by a positive power of $\delta_k$ times
$\gamma_k^{1/2}$. For the Volterra term, positivity and decrease of
${\partial_{\xi} G}$ give
\[
 0\le\int_0^\xi \bigl(\partial_\eta G(\eta,\tau)\bigr)^2/r(\eta,\tau)^2\,\mathrm d\eta
       \le d\,[r(0)^{-1}-r(\xi)^{-1}]\quad(\xi>0),
\]
whereas for $\xi<0$ its absolute value is at most
${\partial_{\xi} G}(\xi)[r(\xi)^{-1}-r(0)^{-1}]$.
Thus its product with ${\partial_{\xi} G}$ is bounded by
$C\delta_k^{1/4}(\bigl(\partial_{\xi} G\bigr)^2+d{\partial_{\xi} G})$.
Finally \eqref{reflcrit:source-size} gives
$|E_{\rm orb}/r|\le Ce^\tau$ and
$|\int_0^\xi cE_{\rm rad}(\eta,c)\,\mathrm d\eta|\le Ce^\tau|\xi|$.
The cutoff commutator is supported in the Gaussian tail and the cutoff
has no time derivative inside a block. We conclude, for fixed
$\epsilon>0$, for example any sufficiently small $\epsilon<1/500$, that
\begin{equation}\label{reflcrit:residual}
 \|\partial_\tau Y_k-\mathscr LY_k\|_{L^2(\mathrm{d}\nu)}
 \le C\delta_k^\epsilon\gamma_k^{1/2}
      +Ce^\tau+C\exp(-b\delta_k^{-1/50}).
\end{equation}

We normalize the first Hermite functions as
$h_0=1$, $h_1=\xi/\sqrt2$, $h_2=(\xi^2-2)/(2\sqrt2)$.
The positive eigenspace of $\mathscr L$ is $\operatorname{span}\{h_0,h_1\}$,
the zero-eigenvalue space is $\operatorname{span}\{h_2\}$, and the remaining
spectrum is at most $-1/2$. We adjoin the scalar coordinate
$j(\tau)=e^{\tau/2}$, which satisfies $\frac{\mathrm d j}{\mathrm d\tau}=j/2$, and define blockwise
\[
 p=\|P_+Y_k\|^2+j^2,\qquad n=\|P_0Y_k\|^2,\qquad
 b=\|P_-Y_k\|^2.
\]
Assign the cutoff $\chi_k$ to the half-open block
$I_k:=(-k-1,-k]$, and write $Y(\tau):=Y_k(\tau)$ there.
Thus each endpoint uses the cutoff from its left. This convention
fixes the value of $Y$ at every time before any past supremum is
taken. Let $P_k,N_k,B_k$ be the respective suprema of $p,n,b$ over
$\tau\le-k$, and put $\Gamma_k:=P_k+N_k+B_k$. Their defining sets
are nested, so all four sequences are nonincreasing. Changing $k$
does not retroactively change the cutoff at an already assigned
time. The trace estimate gives
\[
 \delta_k^2\le C\Gamma_k,\qquad \Gamma_k\ge e^{-k},\qquad
                         \Gamma_k\le C\delta_k^{1/2}.
\]
We check both the jumps and the differential inequalities before
taking the past suprema. At the boundary $\tau=-\ell-1$, the
left value uses $\chi_{\ell+1}$ and the right limit uses $\chi_\ell$.
Because $\delta_{\ell+1}\le\delta_\ell$, their difference is
supported where $|\xi|\ge\delta_\ell^{-1/100}$. Both profiles
satisfy the rough-window estimates at this common time. Integration
against the Gaussian weight therefore gives
\[
 \|Y(-\ell-1+)-Y(-\ell-1)\|_{L^2(\mathrm d\nu)}
 \le T_\ell:=C\exp(-b\delta_\ell^{-1/50}).
\]
For each fixed $m>0$, the elementary bound
$\exp(-b\delta^{-1/50})\le C_m\delta^{m+1}$ and
$\delta_\ell^2\le C\Gamma_\ell$ imply
\[
 T_\ell\le C_m\delta_\ell^m\Gamma_\ell^{1/2}.
\]
If $\ell\ge k$, both this jump and its adjacent traces are within
the past range controlled by $\Gamma_k$. Orthogonal projection
cannot increase the $L^2$ norm, and
$|\|u\|^2-\|v\|^2|\le(\|u\|+\|v\|)\|u-v\|$.
The continuous coordinate $j$ contributes no jump. Consequently,
\[
 |p^+-p^-|+|n^+-n^-|+|b^+-b^-|
 \le C\Gamma_k^{1/2}T_\ell+CT_\ell^2
 \le C_m\delta_k^m\Gamma_k.
\]
The same bounds hold for one-sided traces by limits from the
neighboring open blocks.

Inside a block $I_\ell$ with $\ell\ge k$, set
$\mathcal R:=\partial_\tau Y-\mathscr LY$. Equation
\eqref{reflcrit:residual} and $\|Y\|^2\le\Gamma_k$ give
\[
 \|\mathcal R\|\le
 C\delta_k^\epsilon\Gamma_k^{1/2}+Ce^{-k}+CT_\ell
 \le C\eta_k\Gamma_k^{1/2},\qquad
 \eta_k\le C\delta_k^\epsilon\longrightarrow0.
\]
For the source term we used
$e^{-k}/\Gamma_k^{1/2}\le e^{-k/2}\le\delta_k$;
choose $m>\epsilon$ to absorb the Gaussian tails and the jump
errors into the same $\eta_k$. Self-adjointness of $\mathscr L$
and its eigenvalue bounds yield, within every such block,
\[
 \frac{\mathrm dp}{\mathrm d\tau}\ge p-C\eta_k\Gamma_k,
 \qquad
 \left|\frac{\mathrm dn}{\mathrm d\tau}\right|
                    \le C\eta_k\Gamma_k,
 \qquad
 \frac{\mathrm db}{\mathrm d\tau}\le-b+C\eta_k\Gamma_k.
\]
Indeed, the least positive eigenvalue is $1/2$, the negative
spectrum is at most $-1/2$, and differentiating squared norms
supplies the factor two. The additional term $j^2$ satisfies
$\mathrm d(j^2)/\mathrm d\tau=j^2$. Each residual pairing is
bounded by $C\Gamma_k^{1/2}\|\mathcal R\|$.

For any $\tau\le-k$, integrate these inequalities from
$\tau-1$ to $\tau$, splitting the interval at cutoff boundaries.
A unit interval has at most two boundary contributions, including
one-sided endpoint limits. Their squared-norm jumps have the
bound just proved. The exponential integrating factors are bounded
by $e$ on this interval, so, after absorbing fixed constants,
\[
 \begin{aligned}
 p(\tau-1)&\le e^{-1}p(\tau)+C\eta_k\Gamma_k,\\
 n(\tau-1)&\ge n(\tau)-C\eta_k\Gamma_k,\\
 b(\tau-1)&\ge e\,b(\tau)-C\eta_k\Gamma_k.
 \end{aligned}
\]
Take the supremum in the first inequality, and choose approximate
maximizers for $N_k$ and $B_k$ in the second and third. The
translated times lie in $\{\tau\le-k-1\}$; conversely, every time
in that set is a translate by $-1$ of a time in
$\{\tau\le-k\}$. Together with the nested-set inequality for
$N_k$, this proves, on enlarging $\eta_k$ by a fixed factor,
\begin{align}\label{reflcrit:cutoff-recurrence}
 P_{k+1}&\le e^{-1}P_k+\eta_k\Gamma_k,\nonumber\\
 N_k-\eta_k\Gamma_k&\le N_{k+1}\le N_k,\\
 B_{k+1}&\ge eB_k-\eta_k\Gamma_k.\nonumber
\end{align}
The argument uses only finitely many jumps on each unit interval;
it does not sum jump bounds over the entire infinite past.

We spell out the elementary dichotomy. Monotonicity of $B_k$ and its last
recurrence give $B_k=o(\Gamma_k)$. Replacing $\Gamma_k$ by a fixed multiple
of $P_k+N_k$, the ratio $r_k=N_k/P_k$ obeys
\[
 r_{k+1}\ge\frac{r_k-\eta_k(1+r_k)}{e^{-1}+\eta_k(1+r_k)}.
\]
For any fixed $0<a<R_0$, at sufficiently late indices the right side is
at least $\lambda r_k$, $\lambda>1$, on $[a,R_0]$, and is at least $R_0$
when $r_k\ge R_0$. Hence either $r_k\to0$ or $r_k\to\infty$.
In the second case, choose approximate record times for $N_k$; the original
profiles divided by $N_k^{1/2}$ converge in Gaussian $L^2$ to $h_2$ or
$-h_2$. This is incompatible with increasing profiles: their normalized
averages on $[-3,-1],[-1,1],[1,3]$ must increase, whereas those of $h_2$
are positive, negative, positive, and the negatives have the reversed
violation. Thus $N_k+B_k=o(P_k)$.

The first inequality in \eqref{reflcrit:cutoff-recurrence} now gives geometric decay $\Gamma_k\le C\beta^k$
for some $\beta<1$, so $\delta_k\le C\beta^{k/2}$. Returning to the
residual makes $\sum\eta_k<\infty$. The sharper iteration
$P_{k+1}\le e^{-1}(1+C\eta_k)P_k$ gives $\Gamma_k\le Ce^{-k}$.
Finally \eqref{reflcrit:trace} implies
$|G(0,\tau)|+{\partial_{\xi} G}(0,\tau)\le Ce^{\tau/2}$.
Returning to $F$ proves the two asserted estimates.
This finishes the proof.
\end{proof}

\begin{lem}\label{lem:critical-properness-II}
Let $(M,g)$ be complete and connected and let $\nabla^2f>0$ everywhere.
If $\nabla f(p)=0$, then $p$ is the unique critical point and there is
a constant $C>1$ such that
\[
 f(x)\ge f(p)+C^{-1}d(p,x)-C.
\]
In particular $f$ is proper and bounded below.
\end{lem}
\begin{proof}
Choose $r_0>0$ and $\lambda_0>0$ so that
$\nabla^2f\ge\lambda_0g$ on $B(p,r_0)$. Along a unit-speed minimizing
geodesic $\gamma$ starting at $p$, put $v:=f\circ\gamma$.
Then $\frac{\mathrm d v}{\mathrm dt}(0)=0$, $\frac{\mathrm d^2 v}{\mathrm dt^2}\ge\lambda_0$ for $0\le t\le r_0$, and
$\frac{\mathrm d^2 v}{\mathrm dt^2}>0$ thereafter. If $d(p,x)\ge r_0$, integration gives
\[
 f(x)-f(p)\ge\lambda_0r_0d(p,x)-\tfrac12\lambda_0r_0^2.
\]
Changing the constant covers the remaining ball. Strict convexity
along a minimizing segment excludes a second critical point.
Hopf--Rinow and the lower bound make every sublevel compact.
This proves the lemma.
\end{proof}

We can now combine the radius, blow-down, and parabolic comparisons
to rule out the no-critical point branch.

\begin{prop}

\label{reflcrit:critical}
Let $(M^4,g,f)$ satisfy the hypotheses of Theorem~\ref{thm:main},
with the potential-preserving $O(3)$ action constructed above.
Then $f$ has a unique critical point and is proper. After an
additive normalization, we may assume $\min_Mf=0$.
\end{prop}
\begin{proof}
The four-dimensional Bryant soliton has a critical point. In the
other branch, we suppose for contradiction that $f$ has no critical
point. The preceding lemmas give $L_+=0$, $m=-\infty$, uniform scalar decay
on high-potential levels, and the reference estimate in
Lemma~\ref{reflcrit:onecap}. Write $c:=f$ and restrict to the swept tip-to-reference segments.
They have
\[
 u^2-2c\le E_0u.
\]
Indeed it holds at the reference by $F(0,c)=\sqrt{2c}+O(1)$, and propagates
to all smaller radii by monotonicity of $t^2-E_0t$ for $t\ge E_0/2$; below
that threshold it is automatic. At the reference the normal/tangential
splitting gives
\[
 V=\bigl(\partial_{z} F\bigr)^2+(Xu)^2/q^2,\qquad 0<Xu=u\Ric_{\rm orb}\le C/u.
\]
Since $q\to1$ and $F(0,c)^2\asymp c$, the one-cap estimate gives
$F(0,c)^2V\le C$. All hypotheses of Lemma~\ref{reflcrit:barrier} hold.
We choose again fixed $u_0$ with $\phi_B(u_0/\sqrt{K_0})<1$ and
$u_0>r_*\sqrt{K_0}$. The tip-to-reference interval reaches $u_0$ for all
large $c$. The fixed-interval Jacobi calculation in the proof of
Corollary~\ref{reflcrit:infinite-infimum} makes its level slope tend to one,
contradicting the ambient singular-profile bound. This excludes the
no-critical case.

Since $\Hess f=\Ric>0$, Lemma~\ref{lem:critical-properness-II}
now makes the critical point unique and gives the linear lower bound
for $f$, hence properness.
This proves the stated conclusion.
\end{proof}

\begin{proof}[Completion of Theorem~\ref{thm:main}]
The four-dimensional Bryant soliton has the required angular action
as a subgroup of its full rotation group. In the other branch, Corollary \ref{cor:globalfieldcriterion} implies the $O(3)$-action; Proposition~\ref{reflcrit:critical}
gives the unique critical point and proper potential. This finishes the proof.
\end{proof}

\bibliographystyle{alpha}

\vskip10pt

\noindent
Hanbing Fang, Department of Mathematics, The Chinese University of Hong
Kong, Shatin, N.T., Hong Kong; E-mail:
\href{mailto:hanbingfang@cuhk.edu.hk}{hanbingfang@cuhk.edu.hk}\par
\end{document}